\documentclass[11pt]{article}
\usepackage[english]{babel}
\usepackage{amsmath, amsfonts, amssymb, amsthm}
\usepackage{graphicx}
\usepackage{array}
\usepackage{xcolor}

\newcommand \bM {\overline{M}}

\newcommand \R {\mathbb{R}}

\newcommand \N {\mathbb{N}}

\newcommand \Z {\mathbb{Z}}
\newcommand \la {\langle}
\newcommand \ra {\rangle}
\newcommand \ind {\mathbf{1}}
\newcommand \E {\mathbb{E}}
\def\P{\mathbb{P}}
\newcommand \e {\varepsilon}
\newcommand \cO {\mathcal{O}}

\def\cM {\mathcal M}
\def\bw {\overline{w}}
\def\cG {\mathcal G}
\def\cL {\mathcal{L}}

\def \txi {\tilde{\xi}}
\def \Da {\mathcal{D}^\alpha}

\newcommand \rmd {\mathrm{d}}
\def\bfP {\mathbf{P}}
\def \bfE {\mathbf{E}}

\newtheorem{theorem}{Theorem}[section]

\newtheorem{lemma}[theorem]{Lemma}
\newtheorem{corollary}[theorem]{Corollary}
\newtheorem{proposition}[theorem]{Proposition}
\newtheorem{remark}[theorem]{Remark}
\newtheorem{definition}[theorem]{Definition}

\title{Rescaling limits of the spatial Lambda-Fleming-Viot process
with selection}

\author{A.M. Etheridge\thanks{AME supported in part by EPSRC Grant EP/I01361X/1}\\
Department of Statistics\\University of Oxford\\24-29 St Giles\\Oxford OX1 3LB, UK\\~\\
A. V\'eber\thanks{AV supported in part by the {\em chaire Mod\'elisation Math\'ematique et Biodiversit\'e} of Veolia
Environnement-\'Ecole Polytechnique-Museum National d'Histoire Naturelle-Fondation X.}\\
CMAP -- \'Ecole Polytechnique\\CNRS -- I.P. Paris\\ Route de Saclay\\91128 Palaiseau Cedex\\France \\~\\
F. Yu \thanks{FY supported in part by EPSRC Grant EP/I028498/1.}
\\ Department of Mathematics\\ University of Bristol\\ University Walk\\Bristol BS8 1TW, UK.}

\date{\today}

\begin{document}
\maketitle

\newpage
\begin{abstract}
We consider the spatial $\Lambda$-Fleming-Viot process model for
frequencies of genetic types in
a population living in $\R^d$, with two types of individuals ($0$ and $1$) and natural
selection favouring individuals of type $1$. We first prove that the model is well-defined and provide a measure-valued dual process encoding the locations of the `potential ancestors' of a sample taken from such a population, in the same spirit as the dual process for the SLFV without natural selection \cite{BEV2010}. We then consider two cases, one in which the dynamics of the process are driven by purely `local' events (that is, reproduction events of bounded radii) and
one incorporating large-scale extinction-recolonisation events whose radii have a polynomial tail distribution.
In both cases, we consider a sequence of spatial $\Lambda$-Fleming-Viot processes indexed by $n$, and we assume that the fraction of individuals replaced during a
reproduction event and the relative
frequency of events during which natural selection acts tend to $0$ as $n$ tends to infinity.
We choose the decay of these parameters in such a way that when reproduction is only local, the measure-valued process describing the local frequencies of the less favoured type
converges in distribution to a (measure-valued) solution to the stochastic Fisher-KPP equation in one dimension, and to a (measure-valued) solution to the deterministic Fisher-KPP equation in more than one dimensions. When large-scale extinction-recolonisation events occur, the sequence of processes converges instead to the solution to the analogous equation in which the Laplacian
is replaced by a fractional Laplacian  (again, noise can be retained in the limit only in one spatial dimension).
We also consider the process of `potential ancestors' of a sample of individuals taken from these
populations, which we see as (the empirical distribution of) a system of branching and coalescing symmetric jump processes. We show their convergence in distribution towards a system of Brownian or stable motions
which branch at some finite rate. In one dimension, in the limit, pairs of particles also coalesce at a
rate proportional to their collision local time. In contrast to previous proofs of scaling limits for the
spatial $\Lambda$-Fleming-Viot process, here the convergence of the more complex \emph{forwards in time} processes is used to prove the convergence of the dual process of potential ancestries.


\end{abstract}

\tableofcontents

\section{Introduction}
\label{introduction}

The principal aim of mathematical population genetics is to understand the
influence of the different forces of evolution that act on a population, and the
interactions between them, in shaping the patterns of genetic diversity that we
see in the present-day population. One important aspect of this is the interplay
between spatial structure of the population and the intrinsic randomness due to reproduction
in a finite population (known as genetic drift).
This is particularly mathematically challenging in one of the most biologically
important situations, when the
population is distributed across a two-dimensional spatial continuum. The obstructions
to producing a mathematically consistent and analytically tractable model in
this setting were highlighted
in \cite{felsenstein:1975} and dubbed ``the pain in the torus''.
The \emph{spatial $\Lambda$-Fleming-Viot process} (SLFV), introduced in
\cite{BEV2010,etheridge:2008}, provides one route to overcoming those obstructions,
and its relatively simple mathematical structure makes it a powerful tool for
investigating genetic diversity in spatially structured populations.
In fact, it is not so
much a process as a general framework for modelling frequencies of
different genetic types in populations which evolve in a spatial continuum.
For example, it is readily adapted to include things like the large-scale
extinction/recolonisation events which have dominated the demographic
history of many species. In this paper, we shall be interested in an extension of
this measure-valued process in which some individuals have higher reproductive success
than others, modelling the evolution of a spatially structured population subject
to \emph{natural selection}.

Variants of the SLFV that incorporate forms of natural selection
already appear in a number of studies
\cite{biswas/etheridge/klimek:2018,EFP2016,EFPS2016,EFS2017,FP2017}, but without a detailed discussion of
the construction of the stochastic processes, or whether they are well-defined when the
geographic space in which the population evolves is infinite. Our first contribution is to
formulate and construct an SLFV with natural selection. The methods that we employ
can be readily adapted to capture all of the forms of selection considered to date, and
indeed the form of selection considered here contains many
of them as special cases.

We shall then turn to using our model to study the interaction between natural
selection, spatial structure, and genetic drift. In particular, we are interested in
identifying the spatial and temporal scales over which one can expect to see a non-trivial signature
of the interaction between these forces.
More precisely, we investigate rescaling limits of the model which capture the resultant
patterns of genetic diversity over large spatial and temporal scales.
In particular, our second contribution is to
find suitable scalings of time, space and of the strength of selection for which,
in the limit as the scaling parameter $n$ tends to infinity, we recover the Fisher-KPP
equation~\cite{FIS1937, KPP1937} and, in one spatial dimension, its stochastic counterpart.
In the presence of large-scale demographic
events, the appropriate rescalings are different and lead to analogous equations with
the Laplacian replaced by the fractional Laplacian, but, intriguingly, no other trace
of the large-scale events survives. The limits obtained here assume that the local
population densities are high, 
thus
complementing results of \cite{EFPS2016,EFS2017} which address the interaction of
natural selection and genetic drift when local population densities are small.

The Fisher-KPP equation
\begin{equation}\label{FKPP}
\partial_t p = \frac{\sigma^2}{2}\Delta p +sp(1-p)
\end{equation}
was introduced independently by Fisher \cite{FIS1937},
specifically to model the spread of an advantageous gene through a spatially distributed
population, and Kolomogorov, Petrovsky \& Piskunov \cite{KPP1937}, who also
highlighted the applications to biology.
Fisher considered a population living in a one-dimensional space, whereas
Kolmogorov et al.~worked in two dimensions (although they then assumed that the distribution
of types was independent of the second coordinate, thus reducing it to the one-dimensional case).
The equation has been extensively studied (and extended in many ways), and is now a standard
model of invasion in biology. A major focus of work has been on the
travelling wave solutions.
When the motion of individuals or genes is not local but has a heavy-tailed distribution,
one replaces the Laplacian in~(\ref{FKPP}) by a fractional Laplacian
$-(-\Delta)^{\alpha}$. This, notably, modifies the speed of the travelling wave
solutions, which is constant in the diffusive case
and increases exponentially
in the fractional case; see \cite{CR2013} and references therein.

To take into account the stochasticity inherent in reproduction
in a finite population, in one dimension one can add a noise term of the form
\[
\e \sqrt{p(1-p)} \dot{\mathcal{W}},
\]
to the right hand side of~(\ref{FKPP}), where $\dot{\mathcal{W}}$ is a space-time white noise.
This yields the natural continuous space analogue of the classical stepping-stone model
of population genetics, introduced without selection in \cite{kimura:1953}, and
studied in more generality in, for example, \cite{shiga/shimizu:1980}.
The (continuous space) stochastic Fisher-KPP equation can be obtained from the discrete
space counterpart through rescaling (\emph{c.f.}~\cite{barton/depaulis/etheridge:2002}, where the
case without selection is treated) and
was also
obtained as the limit (over appropriate
large spatial and temporal scales)
of a family of long-range contact processes in \cite{MT1995}.
It has been the object of intensive study, with
the perturbations of solutions due to the noise when $\e$ is very small receiving
particular attention,
\emph{e.g.}~\cite{conlon/doering:2005, mueller/mytnik/quastel:2011,mueller/sowers:1995}
and a huge body of closely related work inspired by work of Brunet, Derrida and
coworkers, \emph{e.g.}~\cite{brunet/derrida:2001}. Our results here
provide the parameter regimes under which the SLFV with selection can be
thought of as a noisy perturbation of the Fisher-KPP equation. Crucially, they apply
in two or more spatial dimensions, where the stochastic PDE has no solution. In particular, the rescaled process $\bM^n$ introduced in Section~\ref{ss: main results} of this work provides a tractable analogue in dimension $d\geq 2$ to the one-dimensional stochastic Fisher-KPP equation with small noise when $n$ is large.

\subsection{The spatial $\Lambda$-Fleming-Viot process with selection}\label{ss: def SLFVS}

The main innovation in the SLFV is that reproduction in the population is
based on a Poisson point process of \emph{events}, rather than on individuals.
It is this which overcomes the pain in the torus.
This is discussed in detail in \cite{BEV2010} and so we do not repeat the motivation here.
Each event determines the region of space in which reproduction (or
extinction/recolonisation) will take place and an \emph{impact} $u$.  As a
result of the event, a proportion $u$ of the individuals living in the
region is replaced by offspring of a parent chosen from the population
immediately before the event (a precise definition of the process is given below).
The Poisson structure renders the process particularly amenable to analytic study.
In the neutral setting, which has been studied rather extensively
(see \cite{BEV2013} for a somewhat out of date review), the parent is chosen uniformly
at random from the affected region, irrespective of type.
There are many possible ways to incorporate
natural selection. Here we shall focus on one of the simplest, but also
most important, in which in the selection of the parent, individuals
are weighted according to their genetic type.

To motivate our definition of the process with (fecundity) selection,
suppose that there are two possible types in the population, which we shall
denote by $0$ and $1$.  In order to give a slight selective advantage to
type~$1$, we fix a \emph{selection coefficient} $s>0$
and suppose that, when an event falls, if the proportion of type~$0$
individuals in the affected region immediately before the event is $\bar{w}$,
then the probability of picking a type~$0$ parent is
$p(\bar{w},s)=\bar{w}/(1+s(1-\bar{w}))$. In other words, in the choice of the parent we give a weight $1$ to type $0$ individuals, and a weight $1+s>1$ to type $1$ individuals, so that the probability of picking a parent of type $0$ is $\bar{w}/(\bar{w}+(1+s)(1-\bar{w}))=p(\bar{w},s)$. Typically one is
interested in weak selection, so
that $s\ll 1$ and, in this case, we can estimate this
probability by $(1-s)\bar{w}+s\bar{w}^2$.  Here again we reap the
benefit of the Poisson structure of events: we can think of events as
being of one of two types. A proportion $(1-s)$ of events are ``neutral'':
the parent is selected exactly as in the neutral setting and has
probability $\bar{w}$ of being of type $0$.  On the other hand, a
proportion $s$ of events are ``selective'' and then the probability of
a type $0$ parent is $\bar{w}^2$.  One way to achieve this is to dictate
that at selective events we choose \emph{two} potential parents, independently, and only if
both are type $0$ will the offspring be type $0$. The Poisson structure
allows us to view neutral and selective events as being driven by
independent Poisson processes. This approach exactly
parallels that usually adopted to incorporate genic selection into the
classical Moran model of population genetics
(see, \emph{e.g.}, Definition 5.6 in \cite{ETH2011}).
Of course there are many ways to modify the selection mechanism. For
example, as in Definition~\ref{defn of process} below, we can allow both
the distribution of
the size of the region affected and of the impact to differ between selective
and neutral events, or we can consider \emph{density dependent} selection, in which
the fitness of an individual depends on the local distribution of genetic types,
\emph{e.g.}~\cite{EFP2016}.

Let us turn to a precise definition. All the random objects in this section are defined on some probability space $(\Omega,\mathcal{F},\P)$.

First we describe the state space of the process, borrowing some results
from \cite{VW2012} in the special case in which the compact space of
possible genetic types is $K=\{0,1\}$. We suppose that the population
evolves in $\R^d$ (although the space of geographical locations could equally,
for example, be taken to be some subset of $\R^d$, or a $d$-dimensional torus).
At each time $t$, the population is represented by a measure $M_t$ on $\R^d\times K$
whose first marginal is Lebesgue measure on $\R^d$. As in the neutral setting,
this corresponds to assuming that individuals are uniformly distributed over
$\R^d$ and for any measurable subset $E$ of $\R^d$ and
$\kappa\in \{0,1\}$, $\mathrm{Vol}(E)^{-1}M_t(E\times\{\kappa\})$ gives the
proportion of individuals of type $\kappa$ in $E$. The space
\begin{equation}\label{def M-lambda}
\cM_\lambda\!:=\! \Big\{M \hbox{ measure on }\R^d\times \{0,1\}\!: \forall f\in C_c(\R^d),\int_{\R^d\times \{0,1\}}f(x)M(\rmd x,\rmd \kappa)= \int_{\R^d}f(x)\rmd x\Big\}
\end{equation}
of such measures is equipped with the topology of vague convergence, which makes it a compact set
(\emph{c.f.}~Lemma~1.1 in \cite{VW2012}). Here $C_c(\R^d)$ denotes the space of all compactly supported continuous functions on $\R^d$. A standard decomposition theorem (see \emph{e.g.}~\cite{KAL2002},
p.561) gives us the existence of a measurable mapping $w_t:\R^d \rightarrow [0,1]$ such that
\begin{equation}\label{relation M-w}
M_t(\rmd x,\rmd \kappa)= \big(w_t(x)\delta_{0}(\rmd \kappa)+(1-w_t(x))\delta_1(\rmd\kappa)\big)\, \rmd x.
\end{equation}
Morally, $w_t(x)$ represents the local fraction of individuals of type $0$ at site $x\in \R^d$ at time $t$, and we abuse notation and call it the ``density'' of $M$. Note that $w_t$ is defined up to a Lebesgue null set, that is two mappings $w_t$ and $\tilde{w}_t$ will be equivalent if and only if
\[
\mathrm{Vol}\big(\big\{ x\in \R^d\,:\, w_t(x)\neq \tilde{w}_t(x)\big\}\big)=0.
\]
In what follows, $w_t$ will denote any representative of the equivalence class of densities for $M_t$. We shall thus equally speak of $M_t$ or $w_t$, depending on what makes the notation more fluid. However, it should be understood that the object of interest in all our results is the measure-valued evolution $(M_t)_{t\geq 0}$.

For every $f\in C_c(\R^d)$ and every $F\in C^1(\R)$ (the space of all continuously differentiable functions on $\R$), let us set
\begin{equation}\label{def crochets}
\la w,f\ra := \int_{\R^d}w(x)f(x)\rmd x
\end{equation}
and let us define the function $\Psi_{F,f}$ on $\cM_\lambda$ by
\begin{align}\label{test function 1}
\Psi_{F,f}(M):= F(\la w,f\ra)=F\bigg(\int_{\R^d\times \{0,1\}}f(x)\ind_{\{0\}}(\kappa)M(\rmd x,\rmd \kappa)\bigg),
\end{align}
where $w$ is any representative of the density of $M$. These functions will prove particularly useful for the following reason.
\begin{lemma}\label{lem: dense set}
The set of functions of the form $\Psi_{F,f}$, $F\in C^1(\R)$ and $f\in C_c(\R^d)$, is dense in $C(\cM_\lambda)$ for the supremum norm topology.
\end{lemma}

\begin{proof}[Proof of Lemma~\ref{lem: dense set}]
Since we endow $\cM_\lambda$ with the topology of vague convergence, the set of all functions of the form
\begin{equation}\label{sup density}
M\mapsto G\left(\int_{\R^d \times \{0,1\}}\varphi(x,\kappa)M(\rmd x,\rmd \kappa)\right),
\end{equation}
with $G\in C^1(\R)$ and $\varphi\in C_c(\R^d\times \{0,1\})$ is dense in $C(\cM_\lambda)$. But if $w$ is a representative of the density of $M$, we can write
\begin{align*}
\int_{\R^d \times \{0,1\}}\varphi(x,\kappa)M(\rmd x,\rmd \kappa)& = \int_{\R^d}\varphi(x,0)w(x)\rmd x + \int_{\R^d}\varphi(x,1)(1-w(x))\rmd x \\
& = \int_{\R^d}(\varphi(x,0)-\varphi(x,1))w(x)\rmd x + \int_{\R^d}\varphi(x,1)\rmd x,
\end{align*}
and so the mapping \eqref{sup density} can be rewritten in the form $F(\la w,f\ra)$, with
\[
F(y) = G\left(y + \int_{\R^d}\varphi(x,1)\rmd x\right) \qquad \hbox{and} \qquad f(x)= \varphi(x,0)-\varphi(x,1).
\]
By construction, in the above we have $F\in C^1(\R)$ and $f \in C_c(\R^d)$. The set of functions of the form \eqref{sup density} is thus included in the set of functions of the form \eqref{test function 1} and the conclusion follows.
\end{proof}

In order to gain a feeling for the process, let us first give a non-rigorous description based on the two independent Poisson point processes of ``neutral'' and ``selective'' events mentioned above. This intuitive idea of how the SLFV with fecundity selection should evolve suggests a natural choice of operator ${\mathcal L}$ on functions of the form \eqref{test function 1}, see \eqref{other def L}, and we shall show in Theorem~\ref{th:existence} that for any probability measure $P$ on $\cM_\lambda$ describing the law of the initial condition, the martingale problem for $(\mathcal{L},P)$ has a unique solution on the space of all measurable $\cM_\lambda$-valued paths. Furthermore, this solution is a Markov process with a.s. c\`adl\`ag paths, and it has the Feller property. The SLFV with selection, with initial distribution $P$, can then be defined as the unique solution to this well-posed martingale problem (see Definition~\ref{defn of process}).

So first, the idea. Let $\mu,\mu'$ be two $\sigma$-finite measures on $(0,\infty)$, and let $\nu=\{\nu_r,\, r>0\}$, $\nu'=\{\nu_r',\, r>0\}$ be two collections of probability measures on $[0,1]$ such that
\begin{equation} \label{cond existence}
\int_0^\infty r^d \int_0^1 u\nu_r(\rmd u)\mu(\rmd r) <\infty,\qquad \hbox{and}\qquad \int_0^\infty r^d \int_0^1 u\nu'_r(\rmd u)\mu'(\rmd r) <\infty.
\end{equation}
Further, let $\Pi^N$ and $\Pi^S$ be two independent Poisson point processes on $\R\times \R^d\times (0,\infty)\times [0,1]$ with respective intensity measures $\rmd t\otimes \rmd x\otimes \mu(\rmd r)\nu_r(\rmd u)$ and
$\rmd t\otimes \rmd x\otimes \mu'(\rmd r)\nu'_r(\rmd u)$. Let $M_0\in \cM_\lambda$ be the (for now, deterministic) initial value of the process. The dynamics of $(M_t)_{t\geq 0}$ are as follows. If $(t,x,r,u)\in \Pi^N$, a neutral event occurs at time $t$, within the closed ball $B(x,r)$:
\begin{enumerate}
\item Sample a type $\kappa$ according to the type distribution within $B(x,r)$ just before the event. That is, $\kappa=0$ with probability $V_r^{-1}M_{t-}(B(x,r)\times\{0\})$, where $V_r$ is the volume of a $d$-dimensional ball of radius $r$; otherwise, $\kappa=1$.
\item Update the value of $M_t$ (only) within $B(x,r)$ by setting
\[
M_t\Big|_{B(x,r)\times\{0,1\}}:= (1-u)M_{t-}\Big|_{B(x,r)\times\{0,1\}} + u\,\rmd x\Big|_{B(x,r)}\otimes \delta_\kappa.
\]
In words, at every site $y\in B(x,r)$ we keep a fraction $(1-u)$ of the population as it was
just before the event, and we replace the remaining fraction $u$ by descendants of the
individual with type $\kappa$ chosen during the first step. These offspring all inherit the
type $\kappa$ of their parent. Thus, a representative of the density of $M_t$ can be taken
to be $w_t(y)=w_{t-}(y)$ if $y\notin B(x,r)$, and
\[
w_t(y) = (1-u)w_{t-}(y)+u \ind_{\{\kappa=0\}} \qquad \hbox{if }y\in B(x,r).
\]
\end{enumerate}
Similarly, if $(t,x,r,u)\in \Pi^S$, a selective event occurs at time $t$, within the closed ball $B(x,r)$:
\begin{enumerate}
\item Sample two types $\kappa$ and $\kappa'$ independently, according to the type distribution
within $B(x,r)$ just before the event. We interpret them as the types of two ``potential'' parents.
\item Update the value of $M_t$ (only) within $B(x,r)$ by setting
\[
M_t\Big|_{B(x,r)\times\{0,1\}}:= (1-u)M_{t-}\Big|_{B(x,r)\times\{0,1\}} + u\,\rmd x\Big|_{B(x,r)}\otimes \delta_{\max\{\kappa,\kappa'\}}.
\]
That is, the offspring are of type $0$ if and only if both potential parents are of type $0$.
This time, a representative of the density of $M_t$ can be taken to be $w_t(y)=w_{t-}(y)$ if
$y\notin B(x,r)$, and
\[
w_t(y) = (1-u)w_{t-}(y)+u \ind_{\{\kappa=\kappa'=0\}} \qquad \hbox{if }y\in B(x,r).
\]
\end{enumerate}


Let us now introduce the operator that will encode this dynamics. For every potential density $w:\R^d\rightarrow [0,1]$, $x\in \R^d$, $r>0$ and $u\in [0,1]$, let us define
\begin{align}
\Theta^+_{x,r,u}(w)& := \ind_{B(x,r)^c}w + \ind_{B(x,r)}((1-u)w+u),\qquad \hbox{and}\nonumber\\
\Theta^-_{x,r,u}(w)& := \ind_{B(x,r)^c}w + \ind_{B(x,r)}(1-u)w. \label{notation Delta}
\end{align}
These quantities will correspond to the value of the density immediately after an event $(t,x,r,u)$ if the offspring are of type $0$ or type $1$ respectively.

Assuming that the above description corresponds to a well-posed martingale problem, we would expect the corresponding operator $\mathcal{L}$ to act on functions of the form \eqref{test function 1} as follows (recall that $V_r$ stands for the volume of a ball of radius $r$): for every $M\in \cM_\lambda$,
\begin{align}
{\mathcal L}\Psi_{F,f}(M) & =
\int_{\R^d} \int_0^\infty\int_0^1 \int_{B(x,r)}\frac{1}{V_r}\, \Big[w(y) F(\la \Theta^+_{x,r,u}(w), f \ra)\label{other def L} \\
& \qquad \quad  + (1-w(y))F(\la \Theta^-_{x,r,u}(w),f\ra) - F(\la w,f \ra)\Big] \rmd y\,\nu_r(\rmd u)\,\mu(\rmd r)\,\rmd x\nonumber \\
& \quad  + \int_{\R^d} \int_0^\infty\int_0^1  \int_{B(x,r)^2}
\frac{1}{V_r^2}\, \Big[w(y)w(z)F\big(\la \Theta^+_{x,r,u}(w), f \ra\big)\nonumber \\
& \qquad \quad   + (1-w(y)w(z))F(\la \Theta^-_{x,r,u}(w), f\ra) - F(\la w,f \ra)\Big] \rmd y\,\rmd z\,\nu'_r(\rmd u)\,\mu'(\rmd r)\,\rmd x. \nonumber
\end{align}
Note that $\cL\Psi_{F,f}(M)$ can also be expressed without referring to the density of $M$:
\begin{align}
& {\mathcal L}\Psi_{F,f}(M)\label{other expression LPsi}\\
& = \int_{\R^d} \int_0^\infty\int_0^1 \int_{B(x,r)\times \{0,1\}}\frac{1}{V_r}\, \Big[\ind_{\{0\}}(\kappa) F\Big(\big\la M,\ind_{B(x,r)^c\times \{0\}}f\big\ra + (1-u)\big\la M,\ind_{B(x,r)\times \{0\}}f\big\ra \nonumber\\
& \qquad  + u\int_{B(x,r)}f(x')\rmd x' \Big) + \ind_{\{1\}}(\kappa) F\Big(\big\la M,\ind_{B(x,r)^c\times \{0\}}f\big\ra + (1-u)\big\la M,\ind_{B(x,r)\times \{0\}}f\big\ra \Big) \nonumber\\
& \qquad  - F\Big(\big\la M, \ind_{\R^d\times \{0\}}f\big\ra\Big) \Big]M(\rmd y,\rmd \kappa)\,\nu_r(\rmd u)\,\mu(\rmd r)\,\rmd x\nonumber\\
& \quad + \int_{\R^d} \int_0^\infty\int_0^1 \int_{(B(x,r)\times \{0,1\})^2}\frac{1}{V_r^2}\, \Big[\ind_{\{0\}}(\kappa)\ind_{\{0\}}(\kappa') F\Big(\big\la M,\ind_{B(x,r)^c\times \{0\}}f\big\ra \nonumber\\
& \qquad + (1-u)\big\la M,\ind_{B(x,r)\times \{0\}}f\big\ra + u \int_{B(x,r)}f(x')\rmd x'\Big)\nonumber \nonumber\\
& \qquad  + \big(1-\ind_{\{0\}}(\kappa)\ind_{\{0\}}(\kappa')\big) F\Big(\big\la M,\ind_{B(x,r)^c\times \{0\}}f\big\ra + (1-u)\big\la M,\ind_{B(x,r)\times \{0\}}f\big\ra\Big)\nonumber\\
&\qquad - F\Big(\big\la M, \ind_{\R^d\times \{0\}}f\big\ra\Big)\Big]M(\rmd y,\rmd \kappa)M(\rmd y',\rmd \kappa')\,\nu'_r(\rmd u)\,\mu'(\rmd r)\,\rmd x,\nonumber
\end{align}
where we have used the bracket notation for some of the integrals to ease the notation.

Let $B_{\cM_\lambda}[0,\infty)$ (\emph{resp.}, $D_{\cM_\lambda}[0,\infty)$) denote the space of all paths (\emph{resp.}, c\`adl\`ag paths) with values in $\cM_\lambda$. When needed, $D_{\cM_\lambda}[0,\infty)$ is endowed with the standard Skorokhod topology and the associated Borel $\sigma$-field. Our first main result is the following.
\begin{theorem}\label{th:existence}Suppose that Condition~\eqref{cond existence} holds. Then, for every probability measure $P$ on $\cM_\lambda$ we have:

\noindent $(i)$ The $B_{\cM_\lambda}[0,\infty)$-martingale problem for $(\mathcal{L},P)$ is well-posed. That is, there exists a unique measurable process $(M_t)_{t\geq 0}$ with values in $\cM_\lambda$ such that $M_0$ has law $P$ and for every function $\Psi_{F,f}$ of the form \eqref{test function 1},
\begin{equation}\label{MP}
\left(\Psi_{F,f}(M_t) - \Psi_{F,f}(M_0) -\int_0^t \cL \Psi_{F,f}(M_s)\rmd s \right)_{t\geq 0}
\end{equation}
is a martingale.

\noindent $(ii)$ The process $(M_t)_{t\geq 0}$ in $(i)$ is a Markov process and its semigroup is Feller. Moreover, it has c\`adl\`ag paths almost surely.
\end{theorem}
We can finally define the SLFV with fecundity selection in a rigourous way, assuming that Condition~\eqref{cond existence} is satisfied.
\begin{definition}[SLFV with fecundity selection (SLFVS)]
\label{defn of process}
Let $P$ be a probability measure on $\cM_\lambda$. We call \emph{spatial $\Lambda$-Fleming-Viot process with fecundity selection, with initial distribution $P$}, the unique solution $(M_t)_{t\geq 0}$ to the martingale problem for $(\mathcal{L},P)$ obtained in Theorem~\ref{th:existence}$(i)$. In particular, by Theorem~\ref{th:existence}$(ii)$, the SLFVS is a strong Markov process with c\`adl\`ag paths a.s.
\end{definition}

The proof of Theorem~\ref{th:existence}, given in Section~\ref{s: proof existence} to ease the exposition, proceeds as follows. First, the result would be an obvious consequence of the Poisson point process formulation
if we had chosen a compact set $E$ in place of $\R^d$ for the geographical space in which the
population evolves, and if the intensities of the Poisson point processes $\Pi^N$ and $\Pi^S$
were finite, as then the global rate at which events fall and $M_t$ is updated would be finite.
We thus start from this simple case and take a sequence of Poisson point processes whose intensities converge to the (possibly infinite) intensities of $\Pi^N$ and $\Pi^S$ on $\R\times E\times (0,\infty) \times [0,1]$. We then take a sequence of hypercubes growing to $\R^d$, and construct the process $(M_t)_{t\geq 0}$ of
Theorem~\ref{th:existence} as a potential limit for the corresponding processes. Uniqueness of such a limit is proved via a duality relation between any solution to the martingale problem \eqref{MP} and a given family of solutions to the martingale problem satisfied by the particle system $(\Xi_t)_{t\geq 0}$ introduced in Section~\ref{s: existence}. This duality argument is a natural analogue of the argument guaranteeing  uniqueness of the neutral SLFV \cite{BEV2010}, for which the dual process is a system of coalescing random
walks interpreted as tracing the locations of the ancestors of individuals in a sample from the population. In the case with selection, we shall see the dual process as a system of branching and coalescing random walks that describes the locations of all \emph{potential} ancestors of individuals in a sample from the population modelled by $(M_t)_{t\geq 0}$. The technical Condition~(\ref{cond existence}) corresponds to Assumption~2.4 in \cite{BEV2010} and expresses the fact that each ``ancestral lineage'' is affected by an event at a finite rate.

Observe that the reproduction events encoded by the Poisson point process $\Pi^S$
favour the subpopulation of individuals of type $1$, since during an event determined by
$\Pi^S$, offspring are of type $0$ only if both the potential parents sampled are of type $0$.
Since we only consider this particular form of selection in this paper,
there should be no ambiguity in simply calling this process the SLFV \emph{with selection},
but we emphasise that, although this is certainly one of the most natural, there are many
alternative models.
For example, one could modify the construction so that one first selects a parental type and
then an impact depending on that type, or one could ``kill'' with differential weights
(\emph{c.f.}~\cite{bah/pardoux:2013, foucart:2013, MIL2015} in the non-spatial setting).

We note that
\cite{etheridge/kurtz:2014} describes two constructions of the SLFV.
The first gives the building blocks for the existence of an SLFV with type-dependent killing,
under somewhat weaker conditions than~(\ref{cond existence}). The proof of existence is
given (only) in the neutral case, but uniqueness remains open.
The second construction, which requires Condition~(\ref{cond existence}),
allows for the sort of selection considered here, although, again, the actual proof
of existence is only provided in the neutral case.

\subsection{A measure-valued dual process of ``potential ancestors''}
\label{ss: general dual}
In this section, we first introduce a process $(\Xi_t)_{t\geq 0}$ with values in the set of all finite point measures on $\R^d$, whose evolution is driven by an independent copy of the Poisson point processes $\Pi^N$ and $\Pi^S$. In Section~\ref{sss: duality}, we state a duality relation between any solution to the martingale problem \eqref{MP} and the process $\Xi$ starting from suitable initial distributions. This duality is the analogue of the relation between the neutral SLFV and its ``genealogical process'' (see Theorem~4.2 in \cite{BEV2010} for a general version of this relation, and Equation~(8) in \cite{BEV2012} for the particular case of two types of individuals). This is the content of Proposition~\ref{prop: dual}, whose proof is deferred to Section~\ref{s: existence} to ease the exposition. Although the duality presented here is very reminiscent of the standard notion of duality between two martingale problems (see \cite{EK1986}, pp.188--189, with $\alpha=\beta=0$ for us), it differs in that the natural duality function
\[
f(M,\Xi) = \prod_{i=1}^k w(x_i)
\]
(for every $\Xi=\sum_{i=1}^k \delta_{x_i}$ and $M\in \cM_\lambda$ with ``density'' $w$) suggested by classical population genetics is not well defined (see the discussion at the beginning of Section~\ref{sss: duality}). Indeed, another representative $w'$ of the density of $M$ may differ from $w$ at some of the $x_i$, yielding a different value for $f(M,\Xi)$. Consequently, we must modify the Ethier \& Kurtz approach to duality, but Relation \eqref{dual bancale} stated in Proposition~\ref{prop: dual} will still take the same form as Relation~(4.35) in \cite{EK1986}.

\subsubsection{Definition of the dual process}\label{sss: def dual}
In contrast with the strategy adopted in Section~\ref{ss: def SLFVS} to construct the SLFVS, here we do not base the definition of the dual process on a martingale problem but, instead, we provide an explicit construction of this finite rate jump process in Definition~\ref{defn of dual}. In Proposition~\ref{prop:well-defined dual}, we show that this definition gives rise to a well-defined Markov process which also solves a martingale problem. This will be sufficient to obtain the duality relation stated in Proposition~\ref{prop: dual} and which is required to prove uniqueness of the solution to the martingale problem for $(\mathcal{L},P)$ stated in \eqref{MP}. 

Let us start with some heuristics on the form and dynamics of the dual process before formulating Definition~\ref{defn of dual}. Recall that during a neutral event $(t,x,r,u)\in \Pi^N$, a single parental type is
chosen according to the type distribution
\[
\frac{1}{V_r}\int_{B(x,r)}M_{t-}(\rmd z,\rmd \kappa)= \frac{1}{V_r}\int_{B(x,r)}\big(w_{t-}(z)\delta_0(\rmd \kappa)+(1-w_{t-}(z))\delta_1(\rmd \kappa)\big)\rmd z
\]
in $B(x,r)$ at time $t-$. Although, strictly speaking, the density $w_{t-}$ is only defined up to a
Lebesgue null set (and so for a given $z$ the value of $w_{t-}(z)$ may differ
between two representatives of the density of $M_{t-}$), this sampling can informally be seen as
picking a spatial location $z$ uniformly at random within $B(x,r)$, and then choosing a parent
from the population at $z$ immediately before the event. Thus the parent is of
type $0$ with probability $w_{t-}(z)$, or $1$ with probability $1-w_{t-}(z)$.
Similarly, the independent sampling of two types within $B(x,r)$ during a selective
event can be interpreted as choosing two locations $z$ and $z'$ independently and
uniformly at random within $B(x,r)$, and then potential parental types
according to the type distributions at $z$ and $z'$ just before the event.

Suppose now that we sample $k\in \N$ individuals at some locations $x_1,\ldots,x_k\in \R^d$
at time $0$, ``the present'', assuming that the population has been evolving for some very large time (that we do not specify). We want to trace back the locations of the ``ancestors'' of the individuals in the sample: that is, we want to go back into the past and describe at every earlier time~$t$ the set of locations in $\R^d$ from which the collection of types seen in our sample may have originated. To motivate the introduction of the process $(\Xi_t)_{t\geq 0}$ below, let us first analyse from a genealogical perspective what happens during each reproduction event. If the event is neutral (\emph{i.e.}, belongs to $\Pi^N$), when an ancestor finds itself in the region affected by the event just after the latter has occurred, the probability that it belongs to the fraction $u$ of the local population replaced during the event is precisely $u$. In this case, the ``parent'' of this ancestor was the ``individual'' whose type was chosen to be the one reproduced during the event, and as expounded above, the location of this ``parent'' is uniformly distributed over the affected area. Consequently, precisely at the time of this event in the past, the ancestral lineage corresponding to the ancestor found in this area jumps onto the location of the ``parent''. On the other hand, if (with probability $1-u$) the ancestor does not belong to the fraction replaced, it is not an offspring of the ``parent'' and its ancestral lineage is not affected by the event (\emph{i.e.}, it remains at the same spatial location). Finally, if there is more than one ancestor in this area, each of them belongs to the fraction of the population just replaced with probability $u$ independently of each other, and the ancestral lineages of all those (and only those) who lie in this ``offspring'' population merge into a single ancestral lineage located at the position of the ``parent''. Note that this procedure is independent of the type of the ``parent''. During a selective event (\emph{i.e.}, belonging to $\Pi^S$),
this can no longer be the case; since we only follow the spatial locations from
which the sampled types originate, and not their types, we are unable to decide which of the two ``potential'' parents
is the true parent of the event. Instead we follow the
locations of all ``potential'' ancestors. More precisely, as in a neutral event, every ancestor present in the area of the event just after it occurred belongs to the fraction of the local population just replaced with probability $u$, independently of each other. At the time of the event in the past, the ancestral lineages corresponding to the ancestors who belong to the ``offspring'' population merge, since they all have the same ``parent''. However, we do not know \emph{a priori} from which of the two potential ``parents'' they inherit their types and so the new ancestral lineage instantly splits into two potential lineages, starting from the positions of the two potential ``parents'', independently and uniformly distributed over the area covered by the event. This parallels the construction of the \emph{ancestral selection graph} and its duality relation with the Wright-Fisher diffusion with selection in the case of a panmictic population \cite{krone/neuhauser:1997, neuhauser/krone:1997}. We shall sometimes use this informal description to see our dual process as a system of branching and coalescing jump processes, although this interpretation will appear much clearer when we describe the limiting ``ancestral'' processes that arise in the regimes of parameters on which we focus in Theorems~\ref{th: conv duals} and \ref{th: conv duals alpha}.

We now give a formal definition of the process $(\Xi_t)_{t\geq 0}$ which will keep track of the locations of the potential ancestors of a sample taken from the current state of the population. To this end, observe that the time-reversed point processes
\begin{equation}\label{time-reversed PPP}
\overleftarrow{\Pi}^i:= \big\{(-t,x,r,u)\, :\, (t,x,r,u)\in \Pi^i\big\},\qquad i\in \{N,S\},
\end{equation}
also form two independent Poisson point processes on $\R\times \R^d\times (0,\infty)\times [0,1]$ with the same intensity measures as the corresponding forwards in time processes. The way in which events happen in both directions of time is thus the same in distribution. Hence, let $\widetilde{\Pi}^N$ and $\widetilde{\Pi}^S$ be independent copies of $\Pi^N$ and $\Pi^S$ respectively, defined on another probability space $(\mathbf{\Omega},\mathcal{F}',\bfP)$ (and so is the process $\Xi$ introduced below).

Let $\cM_p(\R^d)$ denote the set of all finite point measures on $\R^d$, which we endow with the topology of weak convergence. The process $(\Xi_t)_{t\geq 0}$ will take its values in $\cM_p(\R^d)$: each atom of $\Xi_t$ will represent the location of a potential ancestor $t$ units of time in the past.
\begin{definition}
\label{defn of dual}
Let $\Xi^0$ be an $\cM_p(\R^d)$-valued random variable, and let us define the process $(\Xi_t)_{t\geq 0}$ with initial value $\Xi^0$ as follows. We set $\Xi_0=\Xi^0$ and, for convenience, at every time $t\geq 0$ we write
\[
\Xi_t = \sum_{i=1}^{N_t}\delta_{\xi_t^i}
\]
where $N_t= \Xi_t(\R^d)$ and some of the $\xi_t^i$ may be identical (by Lemma~2.3 in \cite{KAL1976}, the elements of this decomposition are measurable functions of $\Xi_t$). Note that the ordering by $1,\ldots,N_t$ of the atoms is arbitrary and will play no role in the updating of $\Xi_t$.

Then:

\noindent For every $(t,x,r,u)\in \widetilde{\Pi}^N$:
\begin{enumerate}
\item To each $\xi_{t-}^i \in B(x,r)$, we independently give a mark with probability $u$, or not with probability $1-u$;
\item If at least one atom $\xi_{t-}^i$ is marked, to form $\Xi_t$ we remove all the marked atoms from $\Xi_{t-}$ and we add a Dirac mass at a location which is drawn uniformly at random from within $B(x,r)$.
\end{enumerate}
For every $(t,x,r,u)\in \widetilde{\Pi}^S$:
\begin{enumerate}
\item To each $\xi_{t-}^i \in B(x,r)$, we independently give a mark with probability $u$, or not with probability $1-u$;
\item If at least one atom $\xi_{t-}^i$ is marked, to form $\Xi_t$ we remove all the marked atoms from $\Xi_{t-}$ and we add two Dirac masses at locations which are drawn independently and uniformly from within $B(x,r)$.
\end{enumerate}
In both cases, if no particles in $\Xi_{t-}$ are marked, then nothing happens.
\end{definition}
Note that the point measure $\Xi_t$ always has at least one atom (unless $\Xi_0=0$), since any removal is accompanied by the insertion of at least one new atom.

Before stating the result showing that this definition gives rise to a well-defined Markov process, let us introduce the operator $\cG$ which will turn out to be the extended generator of $(\Xi_t)_{t\geq 0}$ (\emph{i.e.}, the operator on which the martingale problem satisfied by $\Xi$ is based). Let $C_b^1(\R)$ denote the set of all functions on $\R$ which are bounded, of class $C^1$ and whose first derivatives are bounded. Let also $\mathcal{B}_b(\R^d)$ denote the set of all bounded measurable functions on $\R^d$. For every $F\in C^1_b(\R)$ and $f\in \mathcal{B}_b(\R^d)$, we define the function $\Phi_{F,f}$ by
\begin{equation}\label{test functions for dual}
\Phi_{F,f}(\Xi) := F(\la \Xi,f\ra), \qquad \quad \forall\Xi\in \cM_p(\R^d),
\end{equation}
where $\la \Xi,f\ra = \int f(x)\Xi(\rmd x)$. Finally, we define the function $\cG \Phi_{F,f}$ as follows. For every $\Xi = \sum_{i=1}^l\delta_{x_i}\in \cM_p(\R^d)$,
\begin{align}
\cG \Phi_{F,f}(\Xi):= & \int_{\R^d}\int_0^\infty \int_0^1 \int_{B(x,r)}\frac{1}{V_r} \Bigg[\sum_{\substack{D\subseteq I_{x,r}(\Xi) \\ |D|\geq 1}} u^{|D|}(1-u)^{|I_{x,r}(\Xi)\setminus D|} \label{def G}\\
& \qquad \qquad\times \bigg(F\Big(\la \Xi,f\ra - \sum_{i\in D}f(x_i)+f(z)\Big)-F\big(\la \Xi ,f\ra\big)\bigg)\Bigg]\rmd z \nu_r(\rmd u)\mu(\rmd r)\rmd x \nonumber\\
&+ \int_{\R^d}\int_0^\infty \int_0^1 \int_{B(x,r)^2}\frac{1}{V_r^2} \Bigg[\sum_{\substack{D\subseteq I_{x,r}(\Xi) \\ |D|\geq 1}} u^{|D|}(1-u)^{|I_{x,r}(\Xi)\setminus D|}\nonumber \\
&  \times \bigg(F\Big(\la \Xi,f\ra - \sum_{i\in D}f(x_i)+f(z)+f(z')\Big)-F\big(\la \Xi ,f\ra\big)\bigg)\Bigg]\rmd z\rmd z' \nu'_r(\rmd u)\mu'(\rmd r)\rmd x, \nonumber
\end{align}
where
\begin{equation}\label{def Ixr}
I_{x,r}(\Xi)= \{i\in \{1,\ldots,l\}\,:\, x_i \in B(x,r)\}
\end{equation}
is the set of atoms of $\Xi$ sitting in the closed ball $B(x,r)$ and by convention, the sum over $D\subset I_{x,r}(\Xi), |D|\geq 1$ is set to $0$ if $I_{x,r}(\Xi)$ is empty. Note again that by Lemma~2.3 in \cite{KAL1976}, the elements $l,x_1,\ldots,x_l$ of the decomposition of $\Xi$ are measurable functions of $\Xi$, and so the mapping $\cG\Phi_{F,f}$ is a well-defined measurable function on $\cM_p(\R^d)$.

\begin{proposition}\label{prop:well-defined dual}
The process $(\Xi_t)_{t\geq 0}$ of Definition~\ref{defn of dual} is a well-defined Markov jump process with values in $\cM_p(\R^d)$. In addition, if there exists $K>0$ such that $\bfP[\Xi^0(\R^d)\leq K]=1$, then for every $F\in C^1_b(\R)$ and $f\in \mathcal{B}_b(\R^d)$, the process
\begin{equation}\label{MP dual}
\bigg(\Phi_{F,f}(\Xi_t) - \Phi_{F,f}(\Xi_0) -\int_0^t \cG \Phi_{F,f}(\Xi_s)\,\rmd s\bigg)_{t\geq 0}
\end{equation}
is a martingale.
\end{proposition}

\begin{proof}[Proof of Proposition~\ref{prop:well-defined dual}]
Let us first argue that the process $(\Xi_t)_{t\geq 0}$ of Definition~\ref{defn of dual} is well defined for all time $t\geq 0$. Let us focus on a given atom in $\Xi_t$ (for some $t\geq 0$), say at $z\in \R^d$. Since it is affected by a reproduction event only if it lies in the area of the event and if it is marked (which happens with a prescribed probability $u$), by construction the rate at which this atom is impacted by an event is given by
\begin{align}
& \int_{\R^d}\int_0^\infty \int_0^1 \ind_{\{|z-x|\leq r\}}u\, \nu_r(\rmd u)\mu(\rmd r)\rmd x + \int_{\R^d}\int_0^\infty \int_0^1 \ind_{\{|z-x|\leq r\}}u\, \nu'_r(\rmd u)\mu'(\rmd r)\rmd x \nonumber\\
& = \int_0^\infty \int_0^1 V_ru\, \nu_r(\rmd u)\mu(\rmd r) + \int_0^\infty \int_0^1 V_ru\, \nu'_r(\rmd u)\mu'(\rmd r) := C_0 <\infty \label{jump rate}
\end{align}
(where the finiteness of $C_0$ comes from Condition~(\ref{cond existence})), and so the total rate at which any of the atoms of $\Xi_t$ is affected, and hence $\Xi$ jumps, is bounded from above by $C_0\Xi_t(\R^d)$. Furthermore, the number of atoms in $\Xi_t$ can increase only during an event of $\widetilde{\Pi}^S$, and by at most one (if only one atom is erased and two atoms are created during a selective event). Consequently, the total number of atoms in $\Xi_t$ is stochastically bounded by the number of particles in a Yule process starting with $\Xi_0(\R^d)$ particles, each of which splits into two at constant rate
\begin{equation}\label{branching rate}
\int_0^\infty \int_0^1 V_ru\, \nu'_r(\rmd u)\mu'(\rmd r),
\end{equation}
independently of each other. Combining the above with the fact that $\Xi_0(\R^d)$ is finite a.s., we obtain that with probability one the total mass of $\Xi_t$ is finite for every $t\geq 0$ and there is no accumulation of jumps in finite time. That is, $(\Xi_t)_{t\geq 0}$ is a finite rate jump process defined for all time $t\geq 0$. The fact that $\Xi$ is Markovian then comes from the Poisson point process structure of its evolution.

Let us now give a bound on $\cG \Phi_{F,f}(\Xi)$, defined in \eqref{def G}, to prove first that the operator $\cG$ is well-defined on the set of test functions considered, and second that $(\Xi_t)_{t\geq 0}$ is indeed solution to the martingale problem \eqref{MP dual}. To this end, let $F\in C^1_b(\R)$, $f\in \mathcal{B}_b(\R^d)$ and $\Xi \in \cM_p(\R^d)$. Denoting the sup norm by $\|\cdot\|$ and applying Taylor's theorem to the function $F$, we can write
\begin{align}
\big|\cG& \Phi_{F,f}(\Xi)\big| \label{bound on G}\\
&\leq \|F'\|\Xi(\R^d)\|f\|\int_{\R^d}\int_0^\infty\int_0^1 \Bigg[\sum_{\substack{D\subseteq I_{x,r}(\Xi) \\ |D|\geq 1}} u^{|D|}(1-u)^{|I_{x,r}(\Xi)\setminus D|}\Bigg] \nu_r(\rmd u)\mu(\rmd r)\rmd x \nonumber\\
& \quad + \|F'\|\big(\Xi(\R^d)+1\big)\|f\|\int_{\R^d}\int_0^\infty\int_0^1 \Bigg[\sum_{\substack{D\subseteq I_{x,r}(\Xi) \\ |D|\geq 1}} u^{|D|}(1-u)^{|I_{x,r}(\Xi)\setminus D|}\Bigg] \nu'_r(\rmd u)\mu'(\rmd r)\rmd x.\nonumber
\end{align}
Next, using the bounds $|I_{x,r}(\Xi)|=\Xi(B(x,r))\leq \Xi(\R^d)$,
\[
\sum_{\substack{D\subseteq I_{x,r}(\Xi) \\ |D|\geq 1}} u^{|D|}(1-u)^{|I_{x,r}(\Xi)\setminus D|} = 1-(1-u)^{|I_{x,r}(\Xi)|} \leq u\,|I_{x,r}(\Xi)| \leq u\,\Xi(\R^d)\ind_{\{\Xi(B(x,r))>0\}},
\]
and
\begin{align*}
\int_{\R^d}\int_0^\infty\int_0^1& u\,\Xi(\R^d)\ind_{\{\Xi(B(x,r))>0\}} \nu_r(\rmd u)\mu(\rmd r)\rmd x \\
& \leq \Xi(\R^d) \int_0^\infty \int_0^1 u\, \mathrm{Vol}\big(\mathrm{Supp}(\Xi)+B(0,r)\big)\nu_r(\rmd u)\mu(\rmd r)\\
& \leq C_d \, \Xi(\R^d)^2 \int_0^\infty \int_0^1 ur^d \nu_r(\rmd u)\mu(\rmd r),
\end{align*}
where $\mathrm{Supp}(\Xi)$ denotes the (discrete) support of $\Xi$ and $C_d$ is the volume of a $d$-dimensional ball of radius $1$, we obtain that
\begin{align}
\big|\cG \Phi_{F,f}(\Xi)\big| \leq & \ C_d \|F'\|\|f\|\, \Xi(\R^d)^2\big(\Xi(\R^d)+1\big) \nonumber \\
& \qquad \qquad \times \left(\int_0^\infty\int_0^1 ur^d \nu_r(\rmd u)\mu(\rmd r)+\int_0^\infty\int_0^1 ur^d \nu'_r(\rmd u)\mu'(\rmd r)\right). \label{bound on G ii}
\end{align}
From this we can first conclude that the operator $\cG$ is indeed well-defined on the set of functions of the form $\Phi_{F,f}$, with $F\in C^1_b(\R)$ and $f\in \mathcal{B}_b(\R^d)$. It is then straightforward to see that for every such test function and every $\Xi\in \cM_p(\R^d)$,
\begin{equation}\label{extended generator}
\frac{\rmd}{\rmd t}\bfE_{\Xi}\big[\Phi_{F,f}(\Xi_t)\big]\Big|_{t=0} = \cG \Phi_{F,f}(\Xi).
\end{equation}
Since $\Xi^0(\R^d)\leq K$ a.s., the Yule process with branching rate given in \eqref{branching rate} that dominates the number of particles in $\Xi$ has finite moments at any time $t\geq 0$ (see Equation~(5) in \cite{YUL1925} for the original derivation of the distribution of the number of individuals at any time $t$ in a Yule process, which is negative binomial for any initial number of individuals), and so the expression on the r.h.s. of \eqref{bound on G ii} applied to $\Xi_t$ is integrable for any $t\geq 0$. Combined with the boundedness of $F$ and Fubini's theorem, this yields that
\[
\Phi_{F,f}(\Xi_t) - \Phi_{F,f}(\Xi_0) -\int_0^t \cG \Phi_{F,f}(\Xi_s)\,\rmd s
\]
is integrable for every $t\geq 0$. Together with \eqref{extended generator}, this allows us to conclude that $\Xi$ is indeed a solution to the martingale problem \eqref{MP dual} with initial distribution the law of $\Xi^0$.
\end{proof}

\subsubsection{Duality relation between $(M_t)_{t\geq 0}$ and $(\Xi_t)_{t\geq 0}$}\label{sss: duality}
A key feature of our model, that we shall use repeatedly, is the fact that the processes $(M_t)_{t\geq 0}$ and $(\Xi_t)_{t\geq 0}$ are dual to each other if we restrict our attention to initial distributions on $\cM_p(\R^d)$ of a particular form (in essence, the atoms of $\Xi_0$ should be random and have a distribution absolutely continuous with respect to Lebesgue measure - see below). As in the neutral case \cite{BEV2010,BEV2012}, this will allow us to transfer the information we obtain on $(M_t)_{t\geq 0}$ onto $(\Xi_t)_{t\geq 0}$, and vice versa. Because we want to use this property in the proof of existence of $(M_t)_{t\geq 0}$ in Section~\ref{s: proof existence} (more precisely, to show that there is at most one solution to the martingale problem for $(\cL,\delta_{M^0})$), Proposition~\ref{prop: dual} is phrased in a more general way and relates $(\Xi_t)_{t\geq 0}$ to any solution to the martingale problem for $\cL$.

The difficulty that we face is that the density of any element of $\cM_\lambda$ is only defined Lebesgue a.e.~and so the usual test functions used to establish such dualities in population genetics when the underlying geographical space is discrete, which take the form
\begin{equation}\label{dual function}
D(M,\Xi) := \exp\left(\int_{\R^d}\ln w(x)\Xi(\rmd x)\right)=\prod_{i=1}^k w(x_i)
\end{equation}
for $M\in\cM_\lambda$ with density $w$ and $\Xi=\sum_{i=1}^k \delta_{x_i}$, will not make sense. However, if, instead of taking deterministic points $x_1,\ldots, x_k$, we take random points, with a distribution which has a density $\psi$ with respect to Lebesgue measure on $(\R^d)^k$, then writing $\mu_\psi$ for the law of the random measure constructed in this way, we have for any $M\in \cM_\lambda$
\begin{align}
\int_{\cM_p(\R^d)} &D(M,X)\mu_\psi(dX) =\int_{(\R^d)^k}\psi(x_1,\ldots ,x_k)\bigg\{\prod_{j=1}^kw(x_j)\bigg\}
\,\rmd x_1\cdots \rmd x_k \label{average test function} \\
& = \int_{(\R^d\times \{0,1\})^k}\psi(x_1,\ldots,x_k)\bigg\{\prod_{j=1}^k\ind_{\{0\}}(\kappa_j)\bigg\}\, M(\rmd x_1,\rmd \kappa_1)\cdots M(\rmd x_k,\rmd \kappa_k), \nonumber
\end{align}
which \emph{is} well-defined (and independent of the representative $w$ of the density of $M$). The following property will therefore be very useful for the main result of this section, Proposition~\ref{prop: dual}.
\begin{lemma}\label{lem:absolutely cont}
Suppose that the distribution of $\Xi_0$ has the form $\mu_\psi$, for some $k\geq 1$ and some density function $\psi$ on $(\R^d)^k$. Then for every $t\geq 0$ and every $j\in \{1,2,\ldots\}$, conditionally on $N_t=j$, the law of $(\xi_t^1,\ldots,\xi_t^j)$ is absolutely continuous with respect to Lebesgue measure on $(\R^d)^j$.
\end{lemma}
\begin{proof}[Proof of Lemma~\ref{lem:absolutely cont}]
The desired property follows from the fact that during every event of $\widetilde{\Pi}^N$ or $\widetilde{\Pi}^S$, the distribution of each ``potential parent'' is uniformly distributed over the area of the event, independently of the current locations of the atoms of $\Xi_s$. Hence, each time a point from $\Xi_s$ is removed, the one or two atoms that are added have a location whose law is absolutely continuous with respect to Lebesgue measure on $\R^d$, while the thinning procedure used to remove the points already in $\Xi_s$ preserves the property that the distribution of the locations of the remaining atoms has a density with respect to Lebesgue measure.
\end{proof}

Setting for every vector of $k$ locations $(x_1,\ldots,x_k)\in (\R^d)^k$
\begin{equation}\label{def Xi[x]}
\Xi[x_1,\ldots,x_k] := \sum_{i=1}^k \delta_{x_i}\in \cM_p(\R^d),
\end{equation}
writing $\Xi_0 \sim \mu_\psi$ to denote the fact that the random variable $\Xi_0$ has law $\mu_\psi$, and recalling that $\P$ (\emph{resp.}, $\bfP$) is the probability measure on the space on which $(M_t)_{t\geq 0}$ (\emph{resp.}, $(\Xi_t)_{t\geq 0}$) is defined, we can now state the following result, whose proof is given in Section~\ref{s: existence}.
\begin{proposition}\label{prop: dual}
Let $M^0\in \cM_\lambda$, $k\in \{1,2,\ldots\}$ and let $\psi$ be a density function on $(\R^d)^k$. Then any solution $(M_t)_{t\geq 0}$ to the $B_{\cM_\lambda}[0,\infty)$-martingale problem for $(\cL,\delta_{M^0})$ satisfies: for every $t\geq 0$
\begin{equation}\label{dual bancale}
\int_{\cM_p(\R^d)} \E\big[D(M_t,X)\,|\, M_0=M^0\big]\mu_\psi(dX) = \bfE\big[D(M^0,\Xi_t)\, |\, \Xi_0\sim \mu_\psi\big].
\end{equation}
Equivalently, by Fubini's theorem and \eqref{average test function},
\begin{align}
\E_{M^0}\bigg[\int_{(\R^d)^k} & \psi(x_1,\ldots,x_k)\bigg\{\prod_{j=1}^k w_t(x_j)\bigg\}\, \rmd x_1\cdots \rmd x_k\bigg] \nonumber\\
& = \int_{(\R^d)^k} \psi(x_1,\ldots,x_k)\bfE_{\Xi[x_1,\ldots,x_k]}\bigg[\prod_{j=1}^{N_t} w^0\big( \xi_t^j\big)\bigg]\, \rmd x_1 \cdots \rmd x_k. \label{dual formula}
\end{align}
\end{proposition}
\begin{remark}\label{rk: complements}
By linearity, \eqref{dual formula} also holds for every $\psi\in \mathbb{L}^1((\R^d)^k)$. In addition, to keep the notation simple we have restricted our attention to deterministic initial values $M^0$, but the proof of Proposition~\ref{prop: dual} shows that a similar duality formula holds when $M^0$ is any $\cM_\lambda$-valued random variable. See in particular \eqref{lemme EK}, which only needs to be integrated with respect to the law of $M^0$ to yield the result.
\end{remark}
\begin{remark}\label{rk: dual not enough}
One may try to use Proposition~\ref{prop: dual} to prove uniqueness of the solution to the martingale problem for $(\mathcal{G},\mu_\psi)$, which is satisfied by $(\Xi_t)_{t\geq 0}$. To this end, in the statement of Proposition~\ref{prop: dual}, we would like to replace the process $(\Xi_t)_{t\geq 0}$ of Definition~\ref{defn of dual} by any process $(\widetilde{\Xi}_t)_{t\geq 0}$ solving the same martingale problem \eqref{MP dual}. However, in contrast with the explicit construction of $\Xi$ which immediately yields Lemma~\ref{lem:absolutely cont}, one cannot see from the martingale problem formulation that at any time $t\geq 0$, conditionally on $\widetilde{\Xi}_t(\R^d)$, the law of the locations in $\R^d$ of the atoms of $\widetilde{\Xi}_t$ is absolutely continuous with respect to Lebesgue measure on $(\R^d)^{\widetilde{\Xi}_t(\R^d)}$. But this property is crucial to the proof of Proposition~\ref{prop: dual}, and therefore we cannot prove that \eqref{dual formula} holds more generally than for the process $\Xi$ of Definition~\ref{defn of dual}.
\end{remark}

\subsection{Convergence of the rescaled SLFVS to Fisher-KPP processes}\label{ss: main results}

Now that we have introduced the spatial $\Lambda$-Fleming-Viot process with selection and its dual process of ``potential ancestors'', we turn to the main questions of this work: can we recover the solution to the deterministic or the stochastic Fisher-KPP equation as a scaling limit of the SLFVS, and how does the introduction of (a particular form of) rare but geographically extended extinction-recolonisation events impact the law of the limiting process under analogous scaling assumptions? Recall that the Fisher-KPP equation is a classical model for the wave of advance of a slightly favourable allele in a very dense population, in which individuals reproduce locally, so that changes in local allele frequencies are continuous in time and space. In our framework, this corresponds to focusing on a regime of parameters in which selective events are rare compared to neutral events, the impact of every event (\emph{i.e.}, the fraction of the local population actually affected by the event) is very small, and the event radii have a bounded variance. Therefore, writing $n$ for a parameter that we shall let tend to infinity, in what follows we shall assume that there exist $\delta,\gamma>0$ such that the relative frequency of selective events to neutral events scales like $n^{-\delta}$, and the impact of every event scales like $n^{-\gamma}$. Furthermore, in the first case that we consider below, all events will have the same radius (but the assumption of bounded radii would lead to the same type of results), and this assumption will be relaxed in the second case we consider. Since we are interested in the patterns of variation that we see under this model if we look over large spatial and temporal scales, we shall need a third parameter $\beta\geq 0$ to describe the relevant spatial scale to be considered: time will be scaled by a factor $n$ when space will be scaled by a factor $n^\beta$.

Let us be more precise about our assumptions. First, we concentrate on the particular case in which the intensity measures of the Poisson point processes of reproduction events (see Definition~\ref{defn of process}) satisfy
\begin{equation}\label{cond intensity}
\mu'(\rmd r)\nu_r'(\rmd u)=s_n \mu(\rmd r)\nu_r(\rmd u)
\end{equation}
for a parameter $s_n$ of the form $\sigma n^{-\delta}$, with $\sigma>0$ independent of $n$. That is, the distribution of radii and impacts are the same for neutral and selective events, but neutral events happen $n^\delta/\sigma$ times faster than events during which type $1$ individuals are favoured.  We also choose very
special forms for the measures $\mu(\rmd r)$ and $\nu_r(\rmd u)$.  Our results will certainly hold under much more general conditions, but the proofs become obscured by notation. More precisely, we assume that all events (neutral and selective) have impact $u_n =un^{-\gamma}$, where $u>0$ is independent of $n$. In formulae:
\begin{equation}\label{cond intensity 2}
\nu_r(\rmd u)=\nu'_r(\rmd u)=\delta_{u_n}(\rmd u) \qquad \mathrm{for\ every}\ r>0,
\end{equation}
implying in particular that $\mu' = s_n\mu$ by \eqref{cond intensity}. The assumptions that
\begin{equation}\label{assumptions par}
u_n = \frac{u}{n^\gamma}, \qquad \hbox{and} \qquad s_n=\frac{\sigma}{n^\delta}
\end{equation}
mirror the usual assumptions in the classical Moran
and Wright-Fisher models, in the absence of spatial
structure, in which one is interested in the scaling limits that
are obtained as population size $N$ tends to infinity while $Ns_N$ remains
${\mathcal O}(1)$ (see, \emph{e.g.}, Chapter~5 in \cite{ETH2011}).

We shall consider the following two cases:
\begin{itemize}
\item \textbf{Fixed radius:} $\mu(\rmd r)= \delta_R(\rmd r)$, for some fixed $R>0$.
In this case, we choose $\gamma = 1/3$, $\delta=2/3$, $\beta=1/3$ and set (in any
dimension)
\begin{equation}\label{def wn1}
\bw^n_t(x):= \frac{1}{V_R}M_{nt}\big(B(n^{\beta}x,R)\times \{0\}\big)=\frac{1}{V_R}\, \int_{B(n^{1/3}x,R)}w_{nt}(y)\, \rmd y,
\end{equation}
where we recall that $V_R$ stands for the volume of a $d$-dimensional ball of radius $R$. Writing $w_t^n(\cdot)=w_{nt}(n^{1/3}\cdot)$ and $B_n(x)=B(x,n^{-1/3}R)$, we see that
\[\bw^n_t(x)=\frac{n^{d/3}}{V_R}\int_{B_n(x)}w_t^n(x),\]
and so this scaling corresponds to scaling down the spatial coordinate by
$n^\beta=n^{1/3}$ (so that distance one in the new units corresponds to distance $n^{1/3}$
in the original units), and to considering the timescale $(nt,\, t\geq 0)$.
The random variable $\bw^n_t(x)$ gives the local proportion of
individuals of the unfavoured type $0$ in a small neighbourhood
(of radius $n^{-1/3}R$) of the point $x$ and at time $t$ in these new units.
\item \textbf{Stable radii:} For some $\alpha\in (1,2)$, we set
\begin{equation}\label{form of mu}
\mu(\rmd r)=\frac{\ind_{\{r\geq 1\}}}{r^{d+\alpha+1}}\, \rmd r,
\end{equation}
and
\begin{equation}\label{def wn2}
\bw^n_t(x):=\frac{1}{V_1}M_{nt}\big(B(n^{\beta}x,1)\times \{0\}\big) = \frac{1}{V_1}\, \int_{B(n^{\beta}x,1)}w_{nt}(y)\, \rmd y,
\end{equation}
with
\begin{equation}\label{def constants}
\beta = \frac{1}{2\alpha -1},\quad \gamma= \frac{\alpha-1}{2\alpha-1}\quad \hbox{and}\quad \delta = \frac{\alpha}{2\alpha-1}.
\end{equation}
\end{itemize}
In both cases, we write $\bM^n_t$ for the random measure (taking its values in $\cM_\lambda$) with density $\bw^n_t$. It is straightforward to check that the integrability conditions (\ref{cond existence}) are satisfied; in particular, the indicator function $\ind_{\{r\geq 1\}}$ in the definition of $\mu$ in the stable case prevents microscopic events from accumulating at a rate which would violate these conditions. Consequently, the unscaled $\cM_\lambda$-valued process corresponding to each $n$ is well-defined, and so is its scaled and locally averaged version $(\bM^n_t)_{t\geq 0}$. Note however that the process $\bM^n=(\bM_t^n)_{t\geq 0}$ is not Markovian. Indeed, it is not simply obtained by a change in space and time coordinates of the measures $(M_t)_{t\geq 0}$ (with parameters $s_n$, $u_n$, ...) but its density $\bw^n_t$ at any time is defined as an average over a ball of fixed radius $R$ (or $1$ in the stable case) of the density of $M_{nt}$. Therefore, the law of the ``parental'' type(s) picked during an event cannot be expressed in terms of a sampling from the current value of $\bM^n$ and, additionally, the change in the value of each $\bw^n_t(y)$ due to an event centered in $B(x,r)$ will depend on the geometry of the intersection $B(n^\beta y,R)\cap B(x,r)$. Hence, the evolution of quantities of the form $\la \bw^n_t,f\ra$, with $f\in C_c(\R^d)$, cannot be fully described in terms of $\bM^n_t$.
\begin{remark}
We recover the parameters for the fixed radius case from those for stable radii on setting $\alpha =2$, and so there is some sort of continuity between the two regimes. In the fixed radius case, we are able to provide an informal argument which explains why our choice for the parameters $\beta$, $\gamma$, $\delta$ is appropriate (\emph{c.f.}~Section~\ref{heuristics}). These heuristics also partly explain the choice of the parameter values in the stable case. The missing condition on $\beta,\gamma,\delta$ in this case is less intuitive and arises from a generator calculation, see also Section~\ref{heuristics}.
\end{remark}

Recall that the space $\cM_\lambda$ is equipped with the topology of vague convergence. Let $C_c^{\infty}(\R^d)$ denote the set of all smooth compactly supported functions on $\R^d$ and recall the notation $\la w,f\ra$ from (\ref{def crochets}). Our main results are as follows, starting with the case of ``local'' reproduction.
\begin{theorem}[Fixed radius]\label{th:fixed}
Suppose that $(\bM^n_0)_{n\geq 1}$ converges in distribution to some $M_0\in \cM_\lambda$.  Then, as $n\rightarrow \infty$, the process
$(\bM_t^n)_{t\geq 0}$
converges weakly in $D_{\cM_\lambda}[0,\infty)$ towards a Markov process
$(M_t^\infty)_{t\geq 0}$ with continuous sample paths,
starting at $M^{\infty}_0=M_0$. The limiting process is characterised as follows. Let
\begin{equation}\label{def Gamma}
\Gamma_R= \frac{1}{V_R}\int_{B(0,R)}\int_{B(x,R)}(z_1)^2\rmd z\rmd x
\end{equation}
(where $z_1$ denotes the first coordinate of $z$).

$(i)$ When $d=1$, $(M_t^\infty)_{t\geq 0}$ is the unique process for which,
for every choice of the representative $w^\infty_s$ of the density of $M^\infty_s$ at every time $s$, and for every $f,g\in C^{\infty}_c(\R)$,
\[
\mathcal{Z}^f:=\bigg(\la w^\infty_t,f\ra - \la w_0^\infty,f\ra -\int_0^t
\bigg\{\frac{u\Gamma_R}{2}\, \la w_s^\infty ,\Delta f\ra - 2Ru\sigma\,
\la w_s^\infty (1-w_s^\infty),f\ra\bigg\}\, \rmd s\bigg)_{t\geq 0}
\]
is a continuous zero-mean martingale with quadratic variation at time $t$ equal to
\[
4R^2 u^2 \int_0^t \la w_s^\infty (1-w_s^\infty),f^2\ra\, \rmd s.
\]
Furthermore, the bracket process between $\mathcal{Z}^f$ and $\mathcal{Z}^g$ is given by
\[
\big[\mathcal{Z}^f,\mathcal{Z}^g\big]_t = 4R^2u^2 \int_0^t \la w_s^\infty (1-w_s^\infty),fg\ra\, \rmd s.
\]

$(ii)$ When $d\geq 2$, $(M_t^\infty)_{t\geq 0}$
is the unique (deterministic) process
for which, for every choice of the representative $w^\infty_s$ of the density of $M^\infty_s$ at every time $s$, and for every $f\in C^{\infty}_c(\R^d)$ and $t\geq 0$,
\[
\la w^\infty_t,f\ra = \la w_0^\infty,f\ra +
\int_0^t \bigg\{\frac{u\Gamma_R}{2}\, \la w_s^\infty ,\Delta f\ra - u\sigma
V_R\, \la w_s^\infty (1-w_s^\infty),f\ra\bigg\}\, \rmd s.
\]
\end{theorem}

Informally, in one space dimension, one can see the time-indexed family of densities of the limiting process $(M_t^\infty)_{t\geq 0}$ as a weak solution to the
stochastic partial differential equation
\[
\frac{\partial w}{\partial t} =
\frac{u\Gamma_R}{2}\, \Delta w - 2Ru\sigma w(1-w)
+ 2Ru\sqrt{w(1-w)}\, \dot{\mathcal{W}}
\]
(independently of the representative chosen at every time $t$), where $\dot{\mathcal{W}}$ a space-time white noise.
In dimension $d\geq 2$, on the other hand, the noise term disappears
in the limit and the time-indexed family of densities of $(M_t^\infty)_{t\geq 0}$ can be seen as a weak solution to the
deterministic Fisher-KPP equation
\[
\frac{\partial w}{\partial t} = \frac{u\Gamma_R}{2}\, \Delta w - u\sigma
V_R\, w(1-w).
\]

\begin{remark}\label{rk:1d deterministic}
As we shall explain in Section~\ref{heuristics}, our choice of $\beta=1/3=\gamma$ and $\delta=2/3$ is obtained by solving
\[
1-\gamma=2\beta,\qquad 1-\delta-\gamma=0,\qquad \hbox{and }\beta=\gamma.
\]
This set of three equations guarantees that in one dimension, the limiting process $M^\infty$ is solution to the stochastic Fisher-KPP equation. If we replace the last condition by the inequality $0<\beta<\gamma$, then the sequence of processes $(\bM^n)_{n\geq 1}$ still converges, to a limit which is solution to the deterministic Fisher-KPP equation in any dimension (including $d=1$).
\end{remark}

Theorem~\ref{th:fixed} has a counterpart for the correspondingly rescaled dual process. For every $n\in \N$, let $(\Xi_t)_{t\geq 0}$ be the
process of Definition~\ref{defn of dual} with parameters $\mu=\delta_R$,
$\mu'= s_n\delta_R$, $\nu_R=\nu_R'=\delta_{u_n}$, where $s_n=\sigma n^{-2/3}$
and $u_n=u n^{-1/3}$. $(\Xi_t)_{t\geq 0}$ is thus dual to the unscaled process $(M_t)_{t\geq 0}$ with the same parameters, in the sense of Proposition~\ref{prop: dual} (to ease the notation, the dependence on $n$ of these processes is not reported). Now, define the rescaled process $(\Xi_t^n)_{t\geq 0}$ so that for every $t\geq 0$,
\begin{equation}\label{rescaled dual}
\Xi^n_t = \sum_{i=1}^{N_t^n}\delta_{\xi_t^{n,i}}:= \sum_{i=1}^{N_{nt}}\delta_{n^{-1/3}\xi_{nt}^i}.
\end{equation}

Recall that the space $\cM_p(\R^d)$ of finite point measures on $\R^d$ is endowed with the topology of \emph{weak} convergence, and recall also the definition of the law $\mu_\psi$ on $\cM_p(\R^d)$ given in the paragraph below \eqref{dual function}.
\begin{theorem}[Fixed radius - Dual]\label{th: conv duals}
Let $k\in \{1,2,\ldots\}$, $\psi$ be a probability density on $(\R^d)^k$ and suppose that for any $n\geq 1$, $\Xi^n_0$ has law $\mu_\psi$. Then, as $n\rightarrow \infty$, $(\Xi_t^n)_{t\geq 0}$
converges in distribution in $D_{\cM_p(\R^d)}[0,\infty)$ to a limiting Markov process $(\Xi_t^\infty)_{t\geq 0}$ characterised as follows: $\Xi_0^\infty$ has law $\mu_\psi$ and

$(i)$ When $d=1$, $(\Xi_t^\infty)_{t\geq 0}$ is a system of branching and coalescing Brownian motions,
in which particles follow independent Brownian motions with variance
parameter $u\Gamma_R$, and branch at rate $u\sigma V_R$ into two new
particles, started at the location of the parent. In addition to branching and diffusing,
each pair of particles, independently, also coalesces at
rate $4R^2u^2$ times their collision local time.

$(ii)$ When $d\geq 2$, $(\Xi_t^\infty)_{t\geq 0}$ is a branching Brownian motion (with no coalescence), in which particles follow independent Brownian motions with variance parameter $u\Gamma_R$, and branch at rate $u\sigma V_R$ into two new particles, started at the location of the parent.
\end{theorem}

To state the corresponding results for stable radii, we need some more notation. We write $V_r(x,y)$ for the volume of $B(x,r)\cap B(y,r)$ and define
\[
\Phi(|z-y|):= \int_{\frac{|z-y|}{2}}^\infty \frac{1}{r^{d+1+\alpha}}\frac{V_r(y,z)}{V_r}\, \rmd r.
\]
Now, for every $f\in C_c^\infty(\R^d)$ we set
\begin{equation}\label{generator of stable motion}
\Da f(y) = u \int_{\R^d} \Phi(|z-y|) (f(z)-f(y)) \rmd z.
\end{equation}
We shall check in Lemma~\ref{lemma generator alpha-stable}
that this defines the infinitesimal generator of a symmetric stable process (that is, it is a constant
multiple of the fractional Laplacian).  Our results for stable radii are
then as follows.
\begin{theorem}[Stable radii]\label{th:stable}
Suppose that $\bM^n_0$ converges in distribution to some $M_0\in \cM_\lambda$. Then,
as $n\rightarrow \infty$, the process $(\bM_t^n)_{t\geq 0}$
converges weakly in
$D_{\cM_\lambda}[0,\infty)$ towards a Markov process $(M_t^\infty)_{t\geq 0}$ starting at $M_0$.
Furthermore, if $\Da$ denotes the generator of the symmetric $\alpha$-stable
process defined in (\ref{generator of stable motion}), then

$(i)$ When $d=1$, $(M_t^\infty)_{t\geq 0}$ is the unique process for which,
for every choice of the representative $w^\infty_s$ of the density of $M^\infty_s$ at every time $s$, and for every $f,g\in C_c^{\infty}(\R)$,
\[
\mathcal{Z}^f:=\bigg(\la w_t^\infty , f\ra - \la w_0^\infty , f\ra -\int_0^t
\bigg\{\la w_s^\infty, \Da f\ra - \frac{2u\sigma}{\alpha}\,
\la w_s^\infty(1-w_s^\infty),f\ra \bigg\}\, \rmd s\bigg)_{t\geq 0}
\]
is a continuous zero-mean martingale with quadratic variation at time $t$ equal to
\[
\frac{4u^2}{\alpha-1}\int_0^t  \la w_s^\infty(1-w_s^\infty), f^2\ra\, \rmd s.
\]
Furthermore, the bracket process between $\mathcal{Z}^f$ and $\mathcal{Z}^g$ is given by
\[
\big[\mathcal{Z}^f,\mathcal{Z}^g\big]_t = \frac{4u^2}{\alpha-1}\int_0^t  \la w_s^\infty(1-w_s^\infty), fg\ra\, \rmd s.
\]

$(ii)$ When $d\geq 2$, $(M_t^\infty)_{t\geq 0}$ is the unique (deterministic) process for which,
for every choice of the representative $w^\infty_s$ of the density of $M^\infty_s$ at every time $s$, and for every $f\in C_c^{\infty}(\R^d)$ and $t\geq 0$,
\[
\la w_t^\infty , f\ra = \la w_0^\infty , f\ra +\int_0^t
\bigg\{\la w_s^\infty, \Da f\ra - \frac{u\sigma V_1}{\alpha}\,
\la w_s^\infty(1-w_s^\infty),f\ra \bigg\}\, \rmd s.
\]
\end{theorem}
\begin{remark}\label{rk:1d deterministic stable}
Again, our choice of values for $\beta$, $\gamma$ and $\delta$ is obtained by solving
\[
1-\gamma=\alpha\beta,\qquad 1-\delta-\gamma=0,\qquad \hbox{and }(\alpha-1)\beta=\gamma.
\]
(see Section~\ref{heuristics}) in order to obtain a limiting process $M^\infty$ which is stochastic in one dimension. If we replace the last condition by the inequality $0<(\alpha-1)\beta<\gamma$, then (in any dimension) $(\bM^n)_{n\geq 0}$ converges to a deterministic limit which is characterised as in the statement of Theorem~\ref{th:stable}$(ii)$.
\end{remark}

Likewise, letting $(\Xi_t)_{t\geq 0}$ be the $\cM_p(\R^d)$-valued process which is dual to the
unscaled process $(M_t)_{t\geq 0}$ corresponding to the case of stable radii with parameters $u_n=u/n^{-\gamma}$ and
$s_n=\sigma/n^{-\delta}$, and defining the rescaled process $(\Xi_t^n)_{t\geq 0}$ in such a way that for every $t\geq 0$,
\begin{equation} \label{rescaled dual stable}
\Xi^n_t = \sum_{i=1}^{N^n_t} \delta_{\xi_t^{n,i}} := \sum_{i=1}^{N_{nt}} \delta_{n^{-\beta}\xi_{nt}^i},
\end{equation}
(with the values of $\beta,\gamma,\delta$ given in \eqref{def constants}), we have the following convergence result.
\begin{theorem}[Stable radii - Dual] \label{th: conv duals alpha}
Let $k\in \{1,2,\ldots\}$, $\psi$ be a probability density on $(\R^d)^k$ and suppose that for any $n\geq 1$, $\Xi^n_0$ has law $\mu_\psi$. Then, as $n\rightarrow \infty$, $(\Xi_t^n)_{t\geq 0}$
converges in distribution in $D_{\cM_p(\R^d)}[0,\infty)$ to a limiting Markov process $(\Xi_t^\infty)_{t\geq 0}$ characterised as follows: $\Xi_0^\infty$ has law $\mu_\psi$ and

$(i)$ When $d=1$, $(\Xi_t^\infty)_{t\geq 0}$ is a branching and coalescing stable process, in which particles follow independent symmetric $\alpha$-stable processes which branch at rate $u\sigma V_1/\alpha$
into two particles starting at the location of their parent.
The motion of a single particle is fully described by the generator $\Da$ defined
in (\ref{generator of stable motion}). In addition, each pair of particles, independently, coalesces at rate $4u^2/(\alpha-1)$ times their collision local time.

$(ii)$ When $d\geq 2$, $(\Xi_t^\infty)_{t\geq 0}$ is a branching stable process (with no coalescence), in which particles follow independent symmetric $\alpha$-stable processes with generator $\Da$, and branch at rate $u\sigma V_1/\alpha$ into two new particles, started at the location of the parent.
\end{theorem}

In fact, we shall use knowledge of the limiting ``population model'' $(M_t^\infty)_{t\geq 0}$ to
recover the corresponding limiting results for our rescaled duals. The difficulty with proving Theorems~\ref{th: conv duals} and \ref{th: conv duals alpha} directly
stems from problems with identifying the limiting coalescence mechanism in one
dimension.  This contrasts with the situation of uniformly bounded local population densities (\emph{i.e.}, the impact $u$ not tending to zero) considered in \cite{BEV2012} in the neutral case and in
\cite{EFPS2016,EFS2017} in the selective case,
where it is the ability to identify the
limiting behaviour of the (analytically tractable)
coalescent dual that allows us to prove results about
the large scale evolution of the spatial pattern of allele frequencies.

We close this section with a few remarks. First, one may observe from the expression of $\Da$ given in (\ref{generator of stable motion}) that, as in the fixed radius case, the
drift component of the limiting process is proportional to $u$
and the quadratic variation is proportional to $u^2$, so that $u$ can be thought
of as scaling time.  Moreover, the limiting process
that we obtain in the stable radius case can be seen as a weak solution
to a (stochastic) PDE which only differs from that obtained in the fixed radius
case in that the Laplacian has been replaced by the generator of a
symmetric stable process. This is, perhaps, at first sight rather surprising.
The only effect of the large scale events is on the spatial motion of
individuals in the population, and we see no trace of the correlations in
their movement, or of the selection or genetic drift acting over large scales,
that we have in the prelimiting model.  Notice also that the scaling of $s_n$
(relative to $u_n$) that leads to a nontrivial limit is independent of spatial
dimension. In contrast, in \cite{FP2017}, the authors consider a different scaling for the parameters and prove a similar convergence result and a central limit theorem, in which the order of magnitude and the limit of the fluctuations around the deterministic limiting process are dimension-dependent (despite the fact that the impact $u_n$ tends to zero, while the dependence on dimension mostly occurs when the impact remains fixed.).

As remarked above, we would obtain the same results under much more general conditions. For example, in selecting the regions to be affected by events, not only could one take more general measures $\mu$ (it is the tail behaviour of
$\mu(\rmd r)$ that we see in our limits), but also reproduction events do not need to be based on balls. We anticipate
that this robustness will also be maintained if one considers more general selection mechanisms, in which the strength and direction of selection depends on the local frequencies of different types in the population, and it should be clear how to modify our proofs in such cases.

\subsection{Structure of the paper}

The rest of the paper is laid out as follows. In Section~\ref{s: proof existence}, we prove Theorem~\ref{th:existence}. In Section~\ref{s: existence}, we prove the duality relation stated in Proposition~\ref{prop: dual}. In Section~\ref{heuristics}, we provide heuristic arguments to explain our rescalings.
In Section~\ref{s: convergence 1}, we turn to proving Theorem~\ref{th:fixed}, the scaling limit in the case of fixed radii, and Theorem~\ref{th: conv duals}
which provides the
corresponding result for the rescaled duals. In Section~\ref{s: convergence 2}, we
prove Theorems~\ref{th:stable} and~\ref{th: conv duals alpha},
the analogous results for stable radii. In Appendices~\ref{continuity fixed} and \ref{continuity stable}, we obtain
continuity estimates for the rescaled SLFVS of Sections~\ref{s: convergence 1} and \ref{s: convergence 2}.
In particular, these rather technical estimates are key ingredients in
(and nice complements to) the proofs of Theorems~\ref{th:fixed} and \ref{th:stable}.

\section{Proof of Theorem~\ref{th:existence} (Existence of the SLFVS)}
\label{s: proof existence}
The strategy of the proof is the following. We start with a version of the process in which $\R^d$ is replaced by a hypercube $E$ of finite sidelength and the measures $\mu$ and $\mu'$ are assumed to be finite. In this case, the total rate at which events happen is finite and the corresponding process is a well-defined measure-valued Markov jump process with a.s. c\`adl\`ag trajectories. We then proceed in two steps:
\begin{itemize}
\item[$(i)$] We show existence when $E$ has finite sidelength but $\mu$ and $\mu'$ are only $\sigma$-finite, by taking sequences of finite measures $(\mu^n)_{n\geq 1}$ and $(\mu'^{n})_{n\geq 1}$ such that $\mu^n(\rmd r)$ converges to $\mu(\rmd r)$ (and the same with primes), and proving that the corresponding sequence of processes converges to a well-defined limit.
\item[$(ii)$] Given $(i)$, we extend to $\R^d$ by considering a sequence of processes obtained by restricting to
an increasing family of hypercubes $(E_n)_{n\geq 1}$ which exhaust the space, and proving that this sequence converges to the process $(M_t)_{t\geq 0}$ that we are seeking.
\end{itemize}
Both steps rely on Theorem~4.8.10 in \cite{EK1986}, which states that provided we can show that
\begin{itemize}
\item[$(a)$] The operator $\bf L$ on which the limiting martingale problem is based is included in the set $C_b(\cM_\lambda)\times C_b(\cM_\lambda)$, where $C_b(\cM_\lambda)$ is the set of all bounded continuous functions on $\cM_\lambda$;
\item[$(b)$] The limiting $D_{\cM_\lambda}[0,\infty)$-martingale problem for $(\mathbf{L},P^E)$ (where $P^E$ is the distribution of the limit of $(M^{(n)}_0)_{n\geq 1}$, in particular $P^{\R^d}=P$) has at most one solution;
\item[$(c)$] For every $n$, $M^{(n)}$ is a process with sample paths in $D_{\cM_\lambda}[0,\infty)$  (here we follow Ethier and Kurtz in taking $({\cal G}^n_t)_{t\geq 0}$ to be the natural filtration associated to $M^{(n)}$) and the sequence $(M^{(n)})_{n\geq 1}$ is relatively compact;
\item[$(d)$] There exists a countable set $\Gamma\subset [0,\infty)$ such that for every $(F,G)\in {\bf L}$,
\begin{equation}\label{conv generator}
\lim_{n\rightarrow \infty} \E\left[\left(F\big(M^{(n)}_{t+s}\big)- F\big(M^{(n)}_t\big)-\int_t^{t+s}G\big(M^{(n)}_u\big)\rmd u\right)\bigg(\prod_{i=1}^k h_i\big(M^{(n)}_{t_i}\big)\bigg)\right] = 0
\end{equation}
for all $k\geq 0$, $0\leq t_1<t_2< \cdots < t_k\leq t < t+s$ with $t_i,\, t,\, t+s\notin \Gamma$, and $h_i\in C_b(\cM_\lambda)$;
\end{itemize}
then there exists a solution $(M_t)_{t\geq 0}$ to the $D_{\cM_\lambda}[0,\infty)$-martingale problem for $({\bf L},P^E)$ and $(M^{(n)})_{t\geq 0}\Rightarrow (M_t)_{t\geq 0}$ as $n\rightarrow \infty$. Item~$(b)$ will be a consequence of Proposition~\ref{prop: dual} and Remark~\ref{rk: complements}, whose proofs are postponed until Section~\ref{s: existence} for the sake of clarity. The other items will be checked one by one below. Once this is done, existence of a solution to the martingale problem in $D_{\cM_\lambda}[0,\infty)$ will imply the existence of a solution in the larger space $B_{\cM_\lambda}[0,\infty)$, and since uniqueness holds in $B_{\cM_\lambda}[0,\infty)$ too (see the proof of item $(b)$ below), this will show that the $B_{\cM_\lambda}[0,\infty)$-martingale problem for $(\mathbf{L},P^E)$ is well-posed. That is, Theorem~\ref{th:existence}$(i)$ and the property that the trajectories of $(M_t)_{t\geq 0}$ are c\`adl\`ag a.s. will be proved. Furthermore, since $(b)$ is satisfied for any distribution $P^E$ on $\cM_\lambda$, we shall be able to deduce from $(a)$, $(b)$ and Theorem~4.4.2$(a)$ in \cite{EK1986} that the limiting process $M$ is a Markov process with respect to its natural filtration. The last step will consist in showing that its semigroup is Feller.

Recall the definitions of $\Theta^+_{x,r,u}(w)$ and $\Theta^-_{x,r,u}(w)$ given in (\ref{notation Delta}), and the notation $\|f\|$ (\emph{resp.}, $\|f\|_1$) for the supremum (\emph{resp.}, $\mathbb{L}^1$) norm of the function $f$. To simplify notation, we shall restrict our attention to initial distributions $P$ of the form $\delta_{M^0}$ for some $M^0\in \cM_\lambda$. Indeed, the extension of $(b)$ to a general $P$ is covered by Remark~\ref{rk: complements}, the bounds on the elements of the semi-martingale decomposition on which the proof of $(c)$ rely are independent of the choice of the initial values for the processes of interest, and the proof of $(d)$ can easily be generalised by using the linearity of the expectation and integrating all key equations with respect to $P(dM^0)$. From now on, we thus fix $M^0\in \cM_\lambda$.

\subsubsection*{Proof of $(i)$.}
Let $E$ be some hypercube with sidelength $\ell$, and let $\mu$, $\mu'$ be the
$\sigma$-finite measures on $(0,\infty)$ of Theorem~\ref{th:existence}. Let $(\mu^n)_{n\geq 1}$ and
$(\mu'^{n})_{n\geq 1}$ be two sequences of finite measures on $(0,\infty)$ such that, as $n\rightarrow \infty$,
\begin{equation}\label{approx}
\int_0^\infty\varphi(r)\mu^n(\rmd r) \nearrow \int_0^\infty \varphi(r)\mu(\rmd r) \qquad \hbox{and}\qquad \int_0^\infty\varphi(r)\mu'^{n}(\rmd r) \nearrow \int_0^\infty \varphi(r)\mu'(\rmd r),
\end{equation}
for all measurable $\varphi\geq 0$. Let $\cM_\lambda(E)$ be the analogue of $\cM_\lambda$ (see \eqref{def M-lambda}) when the ``geographical'' space $\R^d$ is replaced by $E$. That is, $\cM_\lambda(E)$ is the set of all measures on $E\times \{0,1\}$ whose first marginal distribution is Lebesgue measure on $E$. It is also a compact space when endowed with the topology of vague convergence (note that since $E\times \{0,1\}$ is compact, the topology of vague convergence is the same as the topology of weak convergence). Let
\[
M^0_E:=M^0\bigg|_{E\times \{0,1\}}\in \cM_\lambda(E)
\]
be the measure induced by $M^0$ (the initial value fixed above) on $E\times \{0,1\}$, and for every $n\geq 1$, let $\Pi^{N,n}_E$ and $\Pi^{S,n}_E$ be independent Poisson point processes on $\R\times E\times (0,\infty)\times [0,1]$ with respective intensity measures $\rmd t\otimes \rmd x \otimes \mu^n(\rmd r)\nu_r(\rmd u)$ and  $\rmd t\otimes \rmd x \otimes \mu'^{n}(\rmd r)\nu'_r(\rmd u)$. Finally, let $(M_t^{(n)})_{t\geq 0}$ be defined as in Definition~\ref{defn of process}, with $\Pi^N$ replaced by $\Pi^{N,n}_E$, $\Pi^S$ replaced by $\Pi^{S,n}_E$, and with $M^{(n)}_0=M^0_E$ (what we call the ball $B(x,r)$ in this case is $B_E(x,r):=B(x,r)\cap E$). For a given $n\geq 1$, the total rate at which an event of $\Pi^{N,n}_E$ or $\Pi^{S,n}_E$ happens is
\begin{align}
 \int_E\int_0^\infty \int_0^1 \nu_r(\rmd u)\mu^n(\rmd r)\rmd x\ + &\int_E\int_0^\infty \int_0^1 \nu'_r(\rmd u)\mu'^{n}(\rmd r)\rmd x \nonumber\\
&\qquad = \mathrm{Vol}(E)\big[\mu^n((0,\infty))+\mu'^{n}((0,\infty))\big]<\infty, \label{bounded jump rate}
\end{align}
and so $M^{(n)}$ is a Markov jump process with jump rates uniformly bounded by the quantity in \eqref{bounded jump rate} and with c\`adl\`ag paths, solution to the martingale problem: $M^{(n)}_0=M^0_E$ and for every $F\in C^1(\R)$ and $f\in C(E)$,
\[
\left(\Psi_{F,f}\big(M^{(n)}_t\big) - \Psi_{F,f}\big(M^{(n)}_0\big)-\int_0^t {\cal L}^{(n)}\Psi_{F,f}\big(M^{(n)}_s\big)\rmd s \right)_{t\geq 0}
\]
is a martingale (for the natural filtration associated to $M^{(n)}$), where $\Psi_{F,f}$ is defined as in \eqref{test function 1} and the bounded continuous function ${\mathcal L}^{(n)}\Psi_{F,f}$ is defined by
\begin{align}
{\mathcal L}^{(n)}\Psi_{F,f}(M) & =
\int_E \int_0^\infty\int_0^1 \int_{B_E(x,r)}\frac{1}{\mathrm{Vol}(B_E(x,r))}\, \Big[w(y) F(\la \Theta^+_{x,r,u}(w), f \ra)\label{def Ln} \\
& \qquad \quad  + (1-w(y))F(\la \Theta^-_{x,r,u}(w),f\ra) - F(\la w,f \ra)\Big] \rmd y\,\nu_r(\rmd u)\,\mu^n(\rmd r)\,\rmd x\nonumber \\
& \quad  + \int_E \int_0^\infty\int_0^1  \int_{B_E(x,r)^2}
\frac{1}{\mathrm{Vol}(B_E(x,r))^2}\, \Big[w(y)w(z)F\big(\la \Theta^+_{x,r,u}(w), f \ra\big)\nonumber \\
& \qquad \quad   + (1-w(y)w(z))F(\la \Theta^-_{x,r,u}(w), f\ra) - F(\la w,f \ra)\Big] \rmd y\,\rmd z\,\nu'_r(\rmd u)\,\mu'^{n}(\rmd r)\,\rmd x. \nonumber
\end{align}
for every $M\in \cM_\lambda(E)$. (As earlier, here $w$ is any representative of the density of $M$ and we have kept the notation $\Theta^{\pm}_{x,r,u}$ for the change in $w$ during an event even though $w$ is now defined on $E$ only.)

Let us show that as $n\rightarrow \infty$, $M^{(n)}$ converges in distribution in $D_{\cM_\lambda(E)}[0,\infty)$ to the unique solution $M^{(\infty)}$ to the $D_{\cM_\lambda(E)}[0,\infty)$-martingale problem for $({\cal L}^{(\infty)},\delta_{M_E^0})$, where ${\cal L}^{(\infty)}$ is defined as in \eqref{def Ln} with $\mu^n$ and $\mu'^{n}$
respectively replaced by $\mu$ and $\mu'$. We check items~$(a)-(d)$ one by one, with $\mathbf{L}= {\cal L}^{(\infty)}$ whose domain ${\cal D}({\cal L}^{(\infty)})$ is taken to be the set of all functions of the form $\Psi_{F,f}$ with $f\in C(E)$ and $F\in C^1(\R)$.

For item~$(a)$, observe that since every mapping $w$ that we consider takes its values in $[0,1]$ (and so does its image by any $\Theta^{\pm}_{x,r,u}$) and $E$ is compact, for every $f\in C(E)$ and every $x\in E$, $r>0$ and $u\in [0,1]$, we have
\begin{equation}\label{ineq theta}
\big|\la \Theta^{\pm}_{x,r,u}(w), f \ra - \la w,f\ra \big| \leq u \|f\| \mathrm{Vol}(B_E(x,r)).
\end{equation}
Consequently, for any $F\in C^1(\R)$, by Taylor's theorem we have
\begin{equation}\label{increments}
\big| F(\la \Theta^{\pm}_{x,r,u}(w), f \ra) - F(\la w,f\ra) \big| \leq  \left(\sup_{|z|\leq \|f\|\mathrm{Vol}(E)}|F'(z)|\right)u\|f\|C_d r^d,
\end{equation}
where $C_d$ is a constant that depends only on the dimension $d$. Writing
\begin{align*}
&\big| w(y) F(\la \Theta^+_{x,r,u}(w), f \ra) + (1-w(y))F(\la \Theta^-_{x,r,u}(w),f\ra) - F(\la w,f \ra) \big| \\
& \leq w(y) \big|F(\la \Theta^+_{x,r,u}(w), f \ra)-F(\la w,f \ra)\big| + (1-w(y)) \big|F(\la \Theta^-_{x,r,u}(w), f \ra)-F(\la w,f \ra)\big| \\
& \leq \left(\sup_{|z|\leq \|f\|\mathrm{Vol}(E)}|F'(z)|\right)u\|f\|C_d r^d,
\end{align*}
we obtain that
\begin{align}
\big|{\cal L}^{(\infty)}\Psi_{F,f}(M)\big| \leq & \left(\sup_{|z|\leq \|f\|\mathrm{Vol}(E)}|F'(z)|\right)\|f\|C_d \mathrm{Vol}(E) \bigg(\int_0^\infty \int_0^1 ur^d\nu_r(\rmd u)\mu(\rmd r) \nonumber\\
& \qquad \qquad + \int_0^\infty \int_0^1 ur^d\nu'_r(\rmd u)\mu'(\rmd r)\bigg), \label{Linfty}
\end{align}
and the quantity on the r.h.s. is a finite constant independent of $M$ by Condition~\eqref{cond existence}. This result proves that the operator ${\cal L}^{(\infty)}$ is indeed well defined on $\mathcal{D}({\cal L}^{(\infty)})$. Since $\cM_\lambda(E)$ is a compact subset of the set of all measures on $E\times \{0,1\}$ and since functions of the form $\Psi_{F,f}$ are continuous on the latter, each $\Psi_{F,f}$ belongs to $C_b(\cM_\lambda(E))$. Recalling \eqref{test function 1}, we can rewrite ${\cal L}^{(\infty)}\Psi_{F,f}(M)$ in terms of the measure $M$ as follows:
\begin{align}
{\cal L}^{(\infty)}&\Psi_{F,f}(M)\label{TF with measures}\\
 = & \int_E\int_0^\infty\int_0^1 \int_{B_E(x,r)\times \{0,1\}} \bigg\{F\bigg(\int_{E\times \{0,1\}}f(z) \ind_{\{0\}}(k)M(\rmd z,\rmd k) \nonumber\\
&\qquad \qquad - u\int_{B_E(x,r)\times \{0,1\}}f(z)\ind_{\{0\}}(k)M(\rmd z,\rmd k)+ u\ind_{\{0\}}(\kappa)\int_{B_E(x,r)}f(z)\rmd z\bigg) \nonumber\\
& \qquad - F\bigg(\int_{E\times \{0,1\}} f(z)\ind_{\{0\}}(k)M(\rmd z,\rmd k)\bigg)  \bigg\} \frac{M(\rmd y,\rmd \kappa)}{\mathrm{Vol}(B_E(x,r))}\nu_r(\rmd u)\mu(\rmd r)\rmd x \nonumber\\
& + \int_E\int_0^\infty\int_0^1 \int_{(B_E(x,r)\times \{0,1\})^2}  \bigg\{F\bigg(\int_{E\times \{0,1\}}f(z) \ind_{\{0\}}(k)M(\rmd z,\rmd k) \nonumber\\
& \qquad\qquad - u\int_{B_E(x,r)\times \{0,1\}}f(z)\ind_{\{0\}}(k)M(\rmd z,\rmd k)+ u\ind_{\{0\}}(\kappa\vee \kappa')\int_{B_E(x,r)}f(z)\rmd z\bigg) \nonumber\\
& \qquad - F\bigg(\int_{E\times \{0,1\}} f(z)\ind_{\{0\}}(k)M(\rmd z,\rmd k)\bigg)  \bigg\}\frac{M^{\otimes 2}(\rmd y,\rmd \kappa,\rmd y',\rmd \kappa')}{\mathrm{Vol}(B_E(x,r))^2}\nu'_r(\rmd u)\mu'(\rmd r)\rmd x. \nonumber
\end{align}
Using the Dominated Convergence Theorem and \eqref{Linfty}, it is then straightforward to check that if $(M_l)_{l\geq 0}$ is a sequence in $\cM_\lambda(E)$ converging to $M$, then ${\cal L}^{(\infty)}\Psi_{F,f}(M_l)$ converges to ${\cal L}^{(\infty)}\Psi_{F,f}(M)$ as $l\rightarrow \infty$ and the function ${\cal L}^{(\infty)}\Psi_{F,f}$ is (sequentially) continuous on $\cM_\lambda(E)$. Together with \eqref{Linfty}, this implies that ${\cal L}^{(\infty)}\Psi_{F,f}\in C_b(\cM_\lambda(E))$ and item~$(a)$ is proved.

Item~$(b)$ is a consequence of the exact analogue of Proposition~\ref{prop: dual} in which the ``geographical'' space $\R^d$ is replaced by $E$ (and $\cM_\lambda$ by $\cM_\lambda(E)$). Indeed, by Lemma~2.1(c) in \cite{VW2012}, the linear span of the set of constant functions and of functions of the form
\begin{equation}\label{functions}
M\mapsto \int_{(\R^d)^k}\psi(x_1,\ldots,x_k) \bigg\{ \prod_{j=1}^k w(x_j)\bigg\}\, \rmd x_1\cdots \rmd x_k,
\end{equation}
(where $M$ has density $w$) for $k\geq 1$ and $\psi\in \mathbb{L}^1((\R^d)^k)\cap C((\R^d)^k)$, is dense in the set of all continuous functions on the compact space $\cM_\lambda(E)$ (and the same holds with $E=\R^d$). This set of functions is therefore separating on the space of all probability distributions on $\cM_\lambda(E)$. We can then proceed exactly as in the proof of Proposition~4.4.7 in \cite{EK1986} to conclude that for every $M^0\in \cM_\lambda(E)$, uniqueness holds for the $B_{\cM_\lambda}[0,\infty)$-martingale problem for $(\mathcal{L}^{(\infty)},\delta_{M^0})$ (or more generally for any distribution for the initial value $M_0$). Indeed, in short \eqref{dual formula} allows us to conclude that any two solutions to the martingale problem have the same one-dimensional distributions, and then Theorem~4.4.2 in \cite{EK1986} gives us that these two solutions necessarily have the same finite dimensional distributions and thus uniqueness in $B_{\cM_\lambda}[0,\infty)$ holds. Item~$(b)$ is proved.

We now turn to item~$(c)$, the relative compactness of $(M^{(n)})_{n\geq 1}$. Since $\cM_{\lambda}(E)$ (equipped with the topology of vague convergence) is a compact space and since by Lemma~\ref{lem: dense set} the set $\mathcal{D}({\cal L}^{(\infty)})$ is dense in $C(\cM_\lambda(E))$, by Theorem~3.9.1 in \cite{EK1986} the relative compactness of $(M^{(n)})_{n\geq 1}$ is equivalent to the relative compactness of the sequence of real-valued processes $(\Psi_{F,f}(M^{(n)}))_{n\geq 1}$ for all $\Psi_{F,f}\in \mathcal{D}({\cal L}^{(\infty)})$. Thus let $F\in C^1(\R)$ and $f\in C(E)$. Using the standard Aldous-Rebolledo criterion \cite{ALD1978,REB1980} and writing $(\Phi^n_t)_{t\geq 0}$ for the predictable finite variation part of $\Psi_{F,f}(M^{(n)})$ and $(Q^n_t)_{t\geq 0}$ for the predictable quadratic variation of its martingale part, we only have to show that
\begin{itemize}
\item[(1)] For every $t\geq 0$, the sequence $(\Psi_{F,f}(M_t^{(n)}))_{n\geq 1}$ is tight.
\item[(2)] For every $T>0$, given a sequence of stopping times $(\tau_n)_{n\geq 1}$ bounded by $T$, for every $\e>0$ there exists $\delta>0$ such that
\begin{equation}\label{Phi}
\limsup_{n\rightarrow \infty} \sup_{\theta \in [0,\delta]} \P\left[\big|\Phi^n_{\tau_n+\theta} - \Phi^n_{\tau_n}\big| >\e \right]\leq \e,
\end{equation}
and
\begin{equation}\label{Qn}
\limsup_{n\rightarrow \infty} \sup_{\theta \in [0,\delta]} \P\left[\big|Q^n_{\tau_n+\theta} - Q^n_{\tau_n}\big| >\e \right]\leq \e.
\end{equation}
\end{itemize}
(1) is straightforward, since for any potential density $w$ we have $|\la w,f\ra| \leq \|f\|\mathrm{Vol}(E)$ and $F$ is continuous, which implies that $\Psi_{F,f}(M^{(n)}_t)$ is bounded uniformly in $n$ and $t$. To deal with (2), for every time $t$ we fix a representative $w^{(n)}_t$ of the density of $M^{(n)}_t$. Since each $(M_t^{(n)})_{t\geq 0}$ is a Markov jump process with bounded jump rates, we have
\[
\Phi^n_t=\int_0^t  \mathcal{L}^{(n)}\Psi_{F,f}(M_s^{(n)}) \rmd s
\]
and
\begin{align*}
& Q^n_t
 = \int_0^t \int_E \int_0^{\infty} \int_0^1 \int_{B_E(x,r)}\frac{1}{V(B_E(x,r))}\,
\Big\{w^{(n)}_{s}(y)\big[F(\la \Theta^+_{x,r,u}(w^{(n)}_{s}),f\ra)-
F(\la w^{(n)}_{s},f \ra)\big]^2\\
& \qquad \qquad   + (1-w^{(n)}_{s}(y))\big[F(\la
\Theta^-_{x,r,u}(w^{(n)}_{s}),f\ra) - F(\la w^{(n)}_{s},f \ra)\big]^2\Big\}\,
 \rmd y\,\nu_r(\rmd u)\,\mu^n(\rmd r)\,\rmd x\, \rmd s \\
& \qquad + \int_0^t \int_E \int_0^{\infty} \int_0^1 \int_{B_E(x,r)^2}
\frac{1}{V(B_E(x,r))^2}\, \Big\{w^{(n)}_{s}(y)w^{(n)}_{s}(z)
\big[F(\la \Theta^+_{x,r,u}(w^{(n)}_{s}),f\ra) \\
& \qquad \qquad \qquad \qquad\qquad \qquad \qquad \qquad \qquad \qquad\qquad \qquad \qquad \qquad \qquad \qquad \qquad- F(\la w^{(n)}_{s},f \ra)
\big]^2  \\ & \ + (1-w^{(n)}_{s}(y)w^{(n)}_{s}(z))
\big[F(\la \Theta^-_{x,r,u}(w^{(n)}_{s}),f\ra) - F(\la w^{(n)}_{s},f \ra)
\big]^2\Big\}\, \rmd y\,\rmd z\,\nu'_r(\rmd u)\,\mu'^{n}(\rmd r)\,\rmd x\,\rmd s.
\end{align*}
Using the expression
for $\mathcal{L}^{(n)}$ given in (\ref{def Ln}) and the bound \eqref{increments}, as in \eqref{Linfty} we obtain that for every $M\in \cM_{\lambda}(E)$,
\begin{align}
\big|\mathcal{L}^{(n)}\Psi_{F,f}(M)\big| &  \leq C \left(\sup_{|z|\leq \|f\|\mathrm{Vol}(E)}|F'(z)|\right)\|f\| \bigg(\int_E \int_0^\infty \int_0^1 \, ur^d \nu_r(\rmd u) \mu^n(\rmd r)\rmd x \nonumber\\
& \qquad \qquad \qquad \qquad \qquad \qquad \qquad \quad +
\int_E \int_0^\infty \int_0^1 \, ur^d \nu'_r(\rmd u)\mu'^{n}(\rmd r)\rmd x\bigg)\nonumber \\
& \leq C' \bigg( \int_0^\infty \int_0^1 \, ur^d \nu_r(\rmd u)\mu(\rmd r) + \int_0^\infty \int_0^1 \, ur^d\nu'_r(\rmd u)\mu'(\rmd r)\bigg), \label{bound on LN}
\end{align}
where we have used that $E$ has finite volume and, by assumption,
$\mu^n(\rmd r)\nearrow \mu(\rmd r)$ (and the corresponding statement with primes).
By Condition (\ref{cond existence}), the expression on the r.h.s. is finite (and independent of $M$), and so with probability $1$ we have for every $\theta >0$
\[
\big|\Phi^n_{\tau_n+\theta} - \Phi^n_{\tau_n}\big| \leq C''\theta,
\]
where $C''$ is independent of $n$ (and even of $T$). Hence, it suffices to choose $\delta>0$ small enough for \eqref{Phi} to hold.

Similarly, the integrands in the expression of $Q^n_t$ are bounded by
\begin{equation}\label{increment Q}
C''\left(\sup_{|z|\leq \|f\|\mathrm{Vol}(E)}|F'(z)|\right)^2 \|f\|^2 u^2 \min(r^{2d},\ell^{2d}),
\end{equation}
where $\ell$ is the sidelength of $E$. But $u^2\leq u$ and there exists a constant $\mathcal{C}_E<\infty$ such that $\min(r^{2d},\ell^{2d})\leq \mathcal{C}_E r^d$, and  so the same reasoning shows that \eqref{Qn} holds too for $\delta>0$ small enough. The relative compactness of every sequence $(\Psi_{F,f}(M^{(n)}))$ is proved, and the relative compactness of $(M^{(n)})_{n\geq 1}$ in $D_{\cM_{\lambda}(E)}[0,\infty)$ thus follows.

Finally, we prove item~$(d)$. Let $F\in C^1(\R)$ and $f\in C(E)$. Let also $k\geq 0$, $0\leq t_1<t_2<\cdots<t_k \leq t<t+s$ and $h_i\in C_b(\cM_\lambda(E))$ for all $i\in \{1,\ldots,k\}$. Since $M^{(n)}$ satisfies the martingale problem for $({\cal L}^{(n)},\delta_{M_E^0})$, we have
\begin{equation}\label{martingale conv}
\E\left[\left(\Psi_{F,f}\big(M^{(n)}_{t+s}\big)- \Psi_{F,f}\big(M^{(n)}_t\big)-\int_t^{t+s}\mathcal{L}^{(n)}\Psi_{F,f}\big(M^{(n)}_u\big)\rmd u\right)\bigg(\prod_{i=1}^k h_i\big(M^{(n)}_{t_i}\big)\bigg)\right]=0.
\end{equation}
We can thus write
\begin{align}
&\E\left[\left(\Psi_{F,f}\big(M^{(n)}_{t+s}\big)- \Psi_{F,f}\big(M^{(n)}_t\big)-\int_t^{t+s}\mathcal{L}^{(\infty)}\Psi_{F,f}\big(M^{(n)}_u\big)\rmd u\right)\bigg(\prod_{i=1}^k h_i\big(M^{(n)}_{t_i}\big)\bigg)\right] \nonumber\\
&\qquad  = 0 + \E\left[\left(\int_t^{t+s}\left(\mathcal{L}^{(n)}\Psi_{F,f}\big(M^{(n)}_u\big)-\mathcal{L}^{(\infty)}\Psi_{F,f}\big(M^{(n)}_u\big)\right)\rmd u\right)\bigg(\prod_{i=1}^k h_i\big(M^{(n)}_{t_i}\big)\bigg)\right].\label{difference}
\end{align}
Using (\ref{approx}) and applying the estimate~(\ref{bound on LN}) with the positive measures $\mu(\rmd r)\nu_r(\rmd u)-\mu^n(\rmd r)\nu_r(\rmd u)$ and $\mu'(\rmd r)\nu'_r(\rmd u)-\mu'^n(\rmd r)\nu'_r(\rmd u)$, we see that we can use the Dominated Convergence Theorem to argue that the second term on the r.h.s. of \eqref{difference} converges to $0$ as $n\rightarrow \infty$. Hence, \eqref{conv generator} holds true with $\Gamma=\emptyset$ and item~$(d)$ is proved.

As detailed at the beginning of this section, we can therefore conclude that there exists a unique solution $M^{(\infty)}$ to the $B_{\cM_\lambda(E)}[0,\infty)$-martingale problem for $(\cL^{(\infty)},\delta_{M_E^0})$ and this process is Markov with respect to its natural filtration and has c\`adl\`ag paths a.s.

Concerning the Feller property of the semigroup of $M^{(\infty)}$, the fact that for every $t\geq 0$ and every $\varphi\in C(\cM_\lambda(E))$, $M\mapsto \E_M[\varphi(M_t^{(\infty)})]$ is a continuous function is a consequence of the continuity in $M^0$ of the quantity on the r.h.s. of \eqref{dual formula} (which is more easily seen in \eqref{average test function} when we replace $\psi$ by the density at time $t$ -- conditional on $N_t$ -- of the locations of the atoms $\xi_t^1,\ldots,\xi_t^{N_t}$ ordered in some arbitrary way) and the property already mentioned in item~$(b)$ that the linear span of the set of constant functions and of functions of the form \eqref{functions} is dense in $C(\cM_\lambda(E))$. The strong continuity of the semigroup is a consequence of the fact, proved in item~$(d)$, that for every $\Psi_{F,f}\in \mathcal{D}(\cL^{(\infty)})$, we have for every $t\geq 0$ and every $M\in \cM_\lambda(E)$
\[
\E_{M}\big[\Psi_{F,f}\big(M_t^{(\infty)}\big)\big] - \Psi_{F,f}(M) = \E_M\bigg[\int_0^t \cL^{(\infty)}\Psi_{F,f}\big(M_s^{(\infty)}\big)\rmd s\bigg].
\]
Together with the uniform bound \eqref{Linfty}, it shows that there exists $C_{F,f}>0$ such that for every $t\geq 0$
\[
\sup_{M\in \cM_\lambda(E)}\big|\E_{M}\big[\Psi_{F,f}\big(M_t^{(\infty)}\big)\big] - \Psi_{F,f}(M)\big|\leq C_{F,f}t,
\]
and the quantity on the l.h.s. indeed converges to $0$ as $t\rightarrow 0$. This property can then be extended to any $\varphi\in C(\cM_\lambda(E))$ by Lemma~\ref{lem: dense set}.

The proof of $(i)$ is thus complete.

\subsubsection*{Proof of $(ii)$.}

The proof of $(ii)$ follows exactly the same pattern, but now the task of bounding the integrals
defining $\Phi^n$ and $Q^n$ becomes more delicate. The resolution is to exploit the fact
that $f$ has compact support $S_f$.

Let $\mu,\mu', \nu$ and $\nu'$ satisfy (\ref{cond existence}) and
let $\{E_n\}_{n\geq 1}$ be a sequence of hypercubes increasing to $\R^d$.
We embed each $\cM_{\lambda}(E_n)$ into
$\cM_\lambda=\cM_{\lambda}(\R^d)$ by setting $w(x)\equiv 0$ outside $E_n$. For every $n\geq 1$, let $M^{[n]}$ denote the $\cM_\lambda$-valued Markov process obtained by imposing that $M^{[n]}\big|_{E_n\times\{0,1\}}$ should evolve like the SLFVS on $E_n\times \{0,1\}$ obtained in $(i)$, and $M^{[n]}\big|_{E_n^c\times\{0,1\}}=\rmd x\big|_{E_n^c}\otimes \delta_1(d\kappa)$ (\emph{i.e.}, $w^{[n]}\equiv 0$ on $E_n^c$). For each $n\in \N$, $M^{[n]}$ is an a.s. c\`adl\`ag process and we assume that it starts from the measure $M_{E_n}^0$ obtained by restricting $M^0$ to $E_n$ (as in $(i)$) and by assuming that its ``density'' $w^0_{E_n}$ is $0$ outside $E_n$ (obviously, $M_{E_n}^0$ converges vaguely to $M^0$ as $n\rightarrow \infty$). According to the previous paragraph, it satisfies the property that for every $F\in C^1(\R)$ and every $f\in C_c(\R^d)$,
\begin{equation}\label{prelimiting MP}
\left(\Psi_{F,f}\big(M^{[n]}_t\big) - \Psi_{F,f}\big(M^{[n]}_0\big) -\int_0^t \cL^{[n]}\Psi_{F,f}\big(M^{[n]}_s\big)\rmd s\right)_{t\geq 0}
\end{equation}
is a martingale, where
\begin{align}
\mathcal{L}^{[n]}&\Psi_{F,f}(M) \nonumber\\
=& \int_0^\infty\int_{(S_f+B(0,r))\cap E_n}
\int_0^1  \int_{B_{E_n}(x,r)}\frac{1}{\mathrm{Vol}(B_{E_n}(x,r))}
\Big[w(y)F(\la \Theta^+_{n,x,r,u}(w),f\ra) \nonumber\\
& \qquad \qquad \qquad     + (1-w(y))F(\la \Theta^-_{n,x,r,u}(w),f\ra) - F(\la w,f\ra)\Big]\, \rmd y\,\nu_r(\rmd u)\,\rmd x\,\mu(\rmd r) \nonumber\\
& +
\int_0^\infty\int_{(S_f+B(0,r))\cap E_n} \int_0^1  \int_{B_{E_n}(x,r)^2}
\frac{1}{\mathrm{Vol}(B_{E_n}(x,r))^2}\Big[w(y)w(z)F\big(\la \Theta^+_{n,x,r,u}(w),f\ra\big)
\nonumber \\
&\quad  + (1-w(y)w(z))F\big(\la \Theta^-_{n,x,r,u}(w),f\ra\big) -
F(\la w,f\ra)\Big]\, \rmd y\,\rmd z\,\nu'_r(\rmd u)\,\rmd x\,\mu'(\rmd r). \label{generator}
\end{align}
Here we have written $S_f+B(0,r):=\{x+y:\, x\in S_f, y\in B(0,r)\}$ (motivated by the fact that if the centre of an event of radius $r$ does not belong to this set, then the event does not intersect the support of $f$ and therefore it does not affect the value of $\Psi_{F,f}(M)$) and we have chosen to report the dependence of the operations $\Theta^{\pm}_{n,x,r,u}$ on $n$ since they modify the value of $w$ only within $E_n$.
The key observation is that
\begin{equation}
\label{bound for increments}
|\la \ind_{B(x,r)}w,f\ra|\leq \|f\| \mathrm{Vol}(S_f\cap B(x,r))\leq C_1\|f\| (r^d\wedge 1),
\end{equation}
and
\begin{equation}
\label{bound for integrating increments}
\mathrm{Vol}(S_f + B(0,r))\leq C_2 (r^d \vee 1),
\end{equation}
where $C_1$ and $C_2$ are independent of $r$ and depend only on the support of $f$.
Moreover, the estimate~(\ref{bound for increments})
is uniform in $w$ and, in particular, the same bound holds if we replace
$w$ by $1-w$.

To see how to apply this, consider the part of~(\ref{generator}) corresponding to neutral events. We split
the integral over $(0,\infty)$ at some radius $R_0>1$. We have that
\begin{align}
\bigg|&\int_{R_0}^\infty \int_{(S_f+B(0,r))\cap E_n}
\int_0^1\int_{B_{E_n}(x,r)}\frac{1}{\mathrm{Vol}(B_{E_n}(x,r))}\Big[w(y)F(\la \Theta^+_{n,x,r,u}(w),f\ra) \nonumber \\
& \qquad \qquad\qquad\qquad+ (1-w(y))F(\la \Theta^-_{n,x,r,u}(w),f\ra) -
F(\la w,f\ra)\Big]\, \rmd y\,\nu_r(\rmd u)\,\rmd x\,\mu(\rmd r)\bigg|\nonumber \\
& \leq C_3(F',f) \int_{R_0}^\infty \int_{(S_f+B(0,r))\cap E_n}
\int_0^1 \, u \mathrm{Vol}(B(x,r)\cap S_f) \, \nu_r(\rmd u)\,\rmd x\,\mu(\rmd r)\nonumber\\
& \leq C_4(F',f) \mathrm{Vol}(S_f)\int_{R_0}^\infty \int_0^1\, ur^d\, \nu_r(\rmd u)\mu(\rmd r), \label{vanishing bound}
\end{align}
where $C_3(F',f)$ and $C_4(F',f)$ depend only on $F'$ and $f$ and the last line uses \eqref{bound for integrating increments} and the fact that $\mathrm{Vol}(B(x,r)\cap S_f)\leq \mathrm{Vol}(S_f)$. To control the second part of the integral corresponding to the neutral part,
notice that a simple estimate using the fact that the corresponding events
have radius bounded above by $R_0$, yields
\begin{align}
\bigg|&\int_0^{R_0} \int_{(S_f+B(0,r))\cap E_n}
\int_0^1\int_{B_{E_n}(x,r)}\frac{1}{\mathrm{Vol}(B_{E_n}(x,r))}\Big[w(y)F(\la \Theta^+_{n,x,r,u}(w),f\ra) \nonumber \\
& \qquad \qquad\qquad\qquad+ (1-w(y))F(\la \Theta^-_{n,x,r,u}(w),f\ra) -
F(\la w,f\ra)\Big]\, \rmd y\,\nu_r(\rmd u)\,\rmd x\,\mu(\rmd r)\bigg|\nonumber \\
& \leq C_3(F',f) \int_0^{R_0} \int_{(S_f+B(0,r))\cap E_n}
\int_0^1 \, u \mathrm{Vol}(B(x,r)\cap S_f) \, \nu_r(\rmd u)\,\rmd x\,\mu(\rmd r)\nonumber\\
& \leq C_5(F',f) \mathrm{Vol}(S_f+B(0,R_0))\int_0^{R_0} \int_0^1 ur^d\, \nu_r(\rmd u)\mu(\rmd r),\label{second bound}
\end{align}
where we have used \eqref{bound for increments} to bound $\mathrm{Vol}(B(x,r)\cap S_f)$ by $Cr^d$ (independently of $x$), and then we have bounded the remaining integral of $\rmd x$ over $(S_f+B(0,r))\cap E_n$ by $\mathrm{Vol}(S_f+B(0,R_0))$. Observe that both bounds \eqref{vanishing bound} and \eqref{second bound} are finite by Condition~\eqref{cond existence}, and they are independent of $M$ (or $w$). Exactly the same arguments control the selection part of the generator $\cL^{[n]}$. Furthermore, the same bounds apply if we replace the operator $\cL^{[n]}$ by $\cL$ defined in \eqref{other def L}. Consequently, we can proceed as in $(i)$ to prove that $\cL$ is included in $C_b(\cM_\lambda)\times C_b(\cM_\lambda)$, which was item~$(a)$ to check. Item~$(b)$ is a direct consequence of Proposition~\ref{prop: dual} and of the same arguments as in the analogous part of the proof of $(i)$. The fact that each $M^{[n]}$ is a process with sample paths in $D_{\cM_\lambda}[0,\infty)$ is part of the conclusion of $(i)$. Next, the estimates \eqref{vanishing bound} and \eqref{second bound} enable us to proceed as in the proof of item~$(c)$ for $(i)$ to show that the sequence $(M^{[n]})_{n\geq 1}$ is relatively compact, which proves item~$(c)$ for $(ii)$.

To check item~$(d)$, notice that
by Condition~(\ref{cond existence}), by taking $R_0$ sufficiently large, the right
hand side of~(\ref{vanishing bound}) can be made arbitrarily small, uniformly in $M$. This is enough to ensure that the missing contribution of the events centered \emph{outside} $E_n$ is negligible, that is that
\begin{align}
\bigg|&\int_{0}^\infty \int_{(S_f+B(0,r))\cap E_n^c}
\int_0^1\int_{B(x,r)}\frac{1}{\mathrm{Vol}(B(x,r))}\Big[w(y)F(\la \Theta^+_{x,r,u}(w),f\ra) \nonumber \\
& \qquad\qquad+ (1-w(y))F(\la \Theta^-_{x,r,u}(w),f\ra) -
F(\la w,f\ra)\Big]\, \rmd y\,\nu_r(\rmd u)\,\rmd x\,\mu(\rmd r)\bigg| \nonumber \\
& \leq C\mathrm{Vol}(S_f)\int_{d(S_f,E_n^c)}^\infty \int_0^1 ur^d \nu_r(\rmd u)\mu(\rmd r) \rightarrow 0
\label{extra events}
\end{align}
uniformly in $M$ (or $w$) as $n\rightarrow\infty$, where $d(S_f,E_n^c)$ is the minimal distance between a point of $S_f$ and a point of $E_n^c$ (which tends to infinity as $n$ tends to infinity). The same estimates hold for the selection term, and can also be used to control the error due to the vanishing difference between $\la \Theta^{\pm}_{x,r,u}(w),f\ra$ and $\la \Theta^{\pm}_{n,x,r,u}(w),f\ra$ (the latter modifying $w$ on $B(x,r)\cap E_n$ only). We can then argue as in \eqref{difference} to conclude.

Since items~$(a)-(d)$ are now checked, we can conclude that there exists a unique measurable process $(M_t)_{t\geq 0}$ such that $M_0=M^0$ and satisfying the martingale problem \eqref{MP}, and this process is Markov and has c\`adl\`ag paths a.s. The Feller property of its semigroup can then be derived using the same arguments as in the proof of $(i)$, recalling the uniform bound on $\cL$ obtained in \eqref{vanishing bound} and \eqref{second bound} (see the paragraph following \eqref{second bound}).

\section{Proof of Proposition~\ref{prop: dual} (Duality)}
\label{s: existence}
To prove Proposition~\ref{prop: dual}, we first show that we can extend the operator $\cL$ to a larger class of functions on $\cM_\lambda$ and that any solution to the $B_{\cM_\lambda}[0,\infty)$-martingale problem for $(\cL,\delta_{M^0})$ stated in \eqref{MP} satisfies the corresponding extended martingale problem. We then use this new set of test functions to complete the proof of Proposition~\ref{prop: dual}.

For every $k\in \{1,2,\ldots\}$ and $\psi\in \mathbb{L}^1((\R^d)^k)$, let us define the function $D_\psi$ by: for every $M\in \cM_\lambda$ with density $w$,
\begin{align}
D_\psi (M)&:= \int_{(\R^d\times \{0,1\})^k}\psi(x_1,\ldots,x_k)\bigg\{\prod_{j=1}^k\ind_{\{0\}}(\kappa_j)\bigg\}M(\rmd x_1,\rmd \kappa_1)\cdots M(\rmd x_k,\rmd \kappa_k) \nonumber\\
&=\int_{(\R^d)^k}\psi(x_1,\ldots,x_k) \, \bigg\{\prod_{j=1}^k w(x_j)\bigg\}\, \rmd x_1\cdots \rmd x_k.  \label{test function 2}
\end{align}
Let us also set
\begin{align}
&\mathcal{L}D_\psi(M)\nonumber\\
&:= \int_{\R^d} \int_0^\infty \int_0^1 \int_{B(x,r)}\frac{1}{V_r}\int_{(\R^d)^k} \psi(x_1,\ldots,x_k) \bigg\{\prod_{j\in I^c}w(x_j)\bigg\} \bigg[w(y)\prod_{j\in I}\big((1-u)w(x_j)+u\big)\nonumber\\
& \qquad  + (1-w(y))\prod_{j\in I}\big((1-u)w(x_j)\big) - \prod_{j\in I}w(x_j)\bigg]\rmd x_1\ldots \rmd x_k \rmd y\nu_r(\rmd u)\mu(\rmd r)\rmd x\nonumber\\
&\ \ + \int_{\R^d} \int_0^\infty \int_0^1 \int_{B(x,r)^2}\frac{1}{V_r^2}\int_{(\R^d)^k} \psi(x_1,\ldots,x_k) \bigg\{\prod_{j\in I^c}w(x_j)\bigg\} \nonumber\\
& \qquad \times \bigg[w(y)w(z)\prod_{j\in I}\big((1-u)w(x_j)+u\big) + (1-w(y)w(z))\prod_{j\in I}\big((1-u)w(x_j)\big) - \prod_{j\in I}w(x_j)\bigg]\nonumber\\
& \qquad \qquad \qquad \qquad\qquad \qquad\qquad \qquad \rmd x_1\ldots \rmd x_k \rmd y \rmd z\nu'_r(\rmd u)\mu'(\rmd r) \rmd x, \label{def L}
\end{align}
where we use the notation $I:=\{i\, :\, x_i\in B(x,r)\}$ and the convention that a product over the empty set is equal to $1$. Note that for every $\psi\in C_c(\R^d)$, the function $D_\psi$ coincides with the function $\Psi_{Id,\psi}$ defined in \eqref{test function 1} and, likewise, the function $\cL D_\psi$ defined in \eqref{def L} coincides with the function $\cL\Psi_{Id,\psi}$ defined in \eqref{other def L}. To see this, let us observe that the first part of \eqref{def L} can be rewritten (using the convention that the product over an empty set is equal to $1$, and then Fubini's theorem to pass from the first line to the next)
\begin{align*}
& \int_{\R^d} \int_0^\infty \int_0^1 \int_{B(x,r)}\frac{1}{V_r}\int_{\R^d}\psi(x_1)\ind_{\{x_1\in B(x,r)\}}\big[w(y)\big((1-u)w(x_1)+u\big)\\
& \qquad \qquad \qquad \qquad \qquad \qquad \qquad \qquad +(1-w(y))(1-u)w(x_1) - w(x_1)\big]\, \rmd x_1\rmd y\nu_r(\rmd u)\mu(\rmd r)\rmd x \\
& = \int_{\R^d}\! \int_0^\infty\! \int_0^1\! \int_{B(x,r)}\frac{1}{V_r} \Big\{w(y) \big\langle \ind_{B(x,r)}\big((1-u)w+u\big),\psi\big\rangle + (1-w(y))\big\langle \ind_{B(x,r)}(1-u)w,\psi\big\rangle\\
& \qquad \qquad \qquad \qquad \qquad \qquad \qquad \qquad \qquad \qquad \qquad \qquad -\big\langle \ind_{B(x,r)}w,\psi\big\rangle\Big\} \, \rmd y\nu_r(\rmd u)\mu(\rmd r)\rmd x,
\end{align*}
which is equal to the first integral on the r.h.s. of \eqref{other def L} when $F=Id$ and $f=\psi$. The same reasoning can be made on the second part of \eqref{def L}.

The following result shows that the expression in \eqref{def L} in fact extends the operator $\cL$ defined in \eqref{other def L} to all functions of the form $D_\psi$.

\begin{lemma}\label{property 2}
Let $(M_t)_{t\geq 0}$ be a solution to the martingale problem \eqref{MP}. Then for every $k\geq 1$ and $\psi\in \mathbb{L}^1((\R^d)^k)$,
\begin{equation}\label{alt MP}
\left(D_\psi(M_t)-D_\psi(M_0)-\int_0^t \cL D_\psi(M_s)\rmd s\right)_{t\geq 0}
\end{equation}
is a martingale.
\end{lemma}

\begin{proof}[Proof of Lemma~\ref{property 2}]
We have already checked that the desired property held true for $k=1$ and $\psi\in C_c(\R^d)$, since then $D_\psi=\Psi_{Id,\psi}$ and $\cL D_\psi = \cL \Psi_{Id,\psi}$. We first extend the result to $\psi\in \mathbb{L}^1(\R^d)$ by a density argument, and then consider the case $k\geq 2$.

To complete both parts of the programme, we need a general bound on functions of the form $\cL D_\psi$ that we derive now. Let $m\geq 1$ and $\widetilde{\psi}\in \mathbb{L}^1((\R^d)^m)$. For every $M\in \cM_\lambda$ (with density $w$), the expression for $\cL D_{\widetilde{\psi}}(M)$ given in (\ref{def L}) can be rewritten as
\begin{align*}
&\int_{\R^d} \int_0^\infty \int_0^1 \int_{B(x,r)}\frac{1}{V_r}\int_{(\R^d)^m} \widetilde{\psi}(x_1,\ldots,x_m) \bigg\{\prod_{j\in I^c}w(x_j)\bigg\}\\
&\qquad \times \bigg[\bigg(\sum_{J\subset I,J\neq I} (1-u)^{|J|} u^{|I\setminus J|}\prod_{j\in J}w(x_j) \bigg) w(y) +\big((1-u)^{|I|}-1\big)\prod_{j\in I}w(x_j)\bigg] \\
& \qquad \qquad \qquad \qquad\qquad \qquad\qquad \qquad\qquad \qquad \qquad \qquad \qquad \rmd x_1\ldots \rmd x_m \rmd y\nu_r(\rmd u)\mu(\rmd r)\rmd x \\
& + \int_{\R^d} \int_0^\infty \int_0^1 \int_{B(x,r)^2}\frac{1}{V_r^2}\int_{(\R^d)^m} \widetilde{\psi}(x_1,\ldots,x_m) \bigg\{\prod_{j\in I^c}w(x_j)\bigg\} \\
 & \qquad \times \bigg[w(y)w(z)\bigg(\sum_{J\subset I,J\neq I} (1-u)^{|J|} u^{|I\setminus J|}\prod_{j\in J}w(x_j) \bigg) + \big((1-u)^{|I|}-1\big)\prod_{j\in I}w(x_j)\bigg]\\
&\qquad \qquad \qquad \qquad\qquad \qquad\qquad \qquad\qquad \qquad \qquad \qquad \qquad \rmd x_1\ldots \rmd x_k \rmd y\rmd z\nu'_r(\rmd u)\mu'(\rmd r)\rmd x.
\end{align*}
Bounding $w$ by $1$ and using the facts that
\[
\sum_{J\subsetneq I} (1-u)^{|J|}u^{|I\setminus J|} = 1- (1-u)^{|I|}\leq u|I|\ind_{\{|I|\geq 1\}}\leq  mu \ind_{\{B(x,r)\cap \{x_1,\ldots,x_m\}\neq \emptyset\}}
\]
and
\[
\mathrm{Vol}\big(\{x\in \R^d:\, B(x,r)\cap \{x_1,\ldots,x_m\}\neq \emptyset\}\big) \leq m C_dr^d
\]
for a constant $C_d$ depending only on $d$, we obtain that the first term in $|\cL D_{\widetilde{\psi}}(M)|$ is bounded by
\begin{align}
& \int_{\R^d}\int_0^\infty \int_0^1 \int_{B(x,r)}\frac{1}{V_r}\int_{(\R^d)^m}\big|\widetilde{\psi}\big|(x_1,\ldots,x_m)\ind_{\{B(x,r)\cap \{x_1,\ldots,x_m\}\neq \emptyset\}}2mu \nonumber \\
& \qquad \qquad \qquad \qquad\qquad \qquad \qquad \qquad \qquad \qquad \qquad \qquad \rmd x_1\cdots \rmd x_m\rmd y\nu_r(\rmd u)\mu(\rmd r)\rmd x \nonumber \\
& \leq \int_0^\infty \int_0^1 \int_{(\R^d)^m} \big|\widetilde{\psi}\big|(x_1,\ldots,x_m) 2m^2 C_d r^du \, \rmd x_1\cdots \rmd x_m\nu_r(\rmd u)\mu(\rmd r) \nonumber \\
& = 2m^2C_d \|\widetilde{\psi}\|_1 \int_0^\infty \int_0^1 r^d u\, \nu_r(\rmd u)\mu(\rmd r). \label{first bound for L}
\end{align}
We can then bound the second part of $|\cL D_{\widetilde{\psi}}(M)|$ in a similar way and obtain that
\begin{equation} \label{bound for L}
|\cL D_{\widetilde{\psi}}(M)| \leq 2 m^2C_d \|\widetilde{\psi}\|_1\bigg(\int_0^\infty \int_0^1 r^d u\, \nu_r(\rmd u)\mu(\rmd r)+ \int_0^\infty \int_0^1 r^d u\, \nu'_r(\rmd u)\mu'(\rmd r)\bigg),
\end{equation}
and the expression on the r.h.s. is finite by Condition~(\ref{cond existence}).

Using the fact that $C_c(\R^d)$ is dense in $\mathbb{L}^1(\R^d)$ (for the $\mathbb{L}^1$ norm), the bound \eqref{bound for L} and dominated convergence, we can then use the same approach as in the proof of items~$(d)$ in Section~\ref{s: proof existence} (see Equations~\eqref{martingale conv} and \eqref{difference}) to conclude that the process in \eqref{alt MP} is indeed a martingale when $\psi\in \mathbb{L}^1(\R^d)$.

Let us now consider $k\geq 2$. Any integrable function $\psi$ on $(\R^d)^k$ can be approximated (in $\mathbb{L}^1$ norm) by linear combinations of functions of the product form $\psi_1(x_1)\cdots \psi_k(x_k)$ with $\psi_i\in C_c(\R^d)$ for every $i$. Furthermore, by polarisation, the test function
\begin{equation}\label{D tensor}
D_{\otimes \psi_i}(M)=\prod_{i=1}^k\bigg(\int_{\R^d\times \{0,1\}}\psi_i(x_i)\ind_{\{0\}}(\kappa_i)M(\rmd x_i,\rmd \kappa_i)\bigg) = \prod_{i=1}^k \la w,\psi_i\ra
\end{equation}
can in turn be written as a linear combination of functions of the
form $\la w, f\ra^m$, with $m\in \N$ and $f\in C_c(\R^d)$, for which we can use (\ref{other def L}) to obtain
\begin{align}
{\mathcal L}\Psi_{(\cdot)^m,f}(M) &=
\int_{\R^d} \int_0^\infty\int_0^1 \int_{B(x,r)}\frac{1}{V_r}\, \Big[w(y) \big\la \ind_{B(x,r)^c}w + \ind_{B(x,r)}((1-u)w+u), f \big\ra^m \nonumber \\
&\  + (1-w(y))\big\la \ind_{B(x,r)^c}w + \ind_{B(x,r)}(1-u)w,f\big\ra^m - \la w,f \ra^m\Big] \rmd y\nu_r(\rmd u)\mu(\rmd r)\rmd x \nonumber\\
&\ + \int_{\R^d} \int_0^\infty\int_0^1  \int_{B(x,r)^2}
\frac{1}{V_r^2}\, \Big[w(y)w(z)\big\la \ind_{B(x,r)^c}w + \ind_{B(x,r)}((1-u)w+u), f \big\ra^m \nonumber \\
&\  + (1-w(y)w(z))\big\la \ind_{B(x,r)^c}w + \ind_{B(x,r)}(1-u)w, f\ra^m - \la w,f \big\ra^m \Big] \nonumber \\
& \qquad \qquad  \qquad \qquad \qquad \qquad \qquad \qquad \qquad \qquad  \qquad\hfill \rmd y\,\rmd z\nu'_r(\rmd u)\mu'(\rmd r)\rmd x. \label{coincidence}
\end{align}
On the other hand, taking $\psi(x_1,\ldots,x_m)= \prod_{i=1}^m f(x_i)$ in (\ref{test function 2}), we obtain that
\[
D_\psi (M)= \la w,f\ra^m.
\]
Let us thus show that, in this case, the expression on the r.h.s. of (\ref{def L}) coincides with (\ref{coincidence}). We focus on the first term on the r.h.s. of (\ref{def L}), since the computations are the same for the other terms. For fixed $x,r,u,y$, and writing $B$ for $B(x,r)$ to simplify the notation, we have
\begin{align*}
&\int_{(\R^d)^m}  f(x_1)\cdots f(x_m) \bigg[\prod_{j:x_j\notin B}w(x_i)\bigg]\bigg[\prod_{j:x_j\in B}\big((1-u)w(x_j)+u\big)\bigg]\rmd x_1 \ldots \rmd x_m \\
& = \sum_{J\subseteq \{1,\ldots,m\}}\int_{(\R^d)^m}f(x_1)\cdots f(x_m)\bigg[\prod_{j\in J^c}\ind_{\{x_j\notin B\}}w(x_j)\bigg]\bigg[\prod_{j\in J}\ind_{\{x_j\in B\}}((1-u)w(x_j)+u)\bigg] \\
& \qquad \qquad \qquad \qquad\qquad \qquad\qquad \qquad\qquad \qquad\qquad \qquad\qquad \qquad\qquad \qquad \qquad\rmd x_1\ldots \rmd x_m\\
& = \sum_{J\subseteq \{1,\ldots,m\}} \big\la \ind_{B^c} w,f\big\ra^{m-|J|} \big\la \ind_B ((1-u)w+u),f\big\ra^{|J|}\\
& = \sum_{j=0}^m \binom{m}{j} \big\la \ind_{B^c} w,f\big\ra^{m-j} \big\la \ind_B ((1-u)w+u),f\big\ra^{j}
= \big\la \ind_{B^c} w + \ind_B\big((1-u)w+u\big),f \big\ra^m,
\end{align*}
which coincides with the integrand in the first part of (\ref{coincidence}). Checking that the same holds for the three other parts of (\ref{coincidence}), we can conclude that the two expressions for the action of $\cL$ on functions of the form $\la w,f\ra^m$ coincide. Consequently, the process in \eqref{alt MP} is a martingale for $\psi$ form $f\otimes \cdots \otimes f$ with $f\in C_c(\R^d)$, and by linearity for $\psi$ of the form $\psi_1\otimes \cdots \otimes \psi_k$ with $\psi_i\in C_c(\R^d)$ for every $i$. Using the same density argument as in the case $k=1$, together with the bound \eqref{bound for L}, we can finally conclude that the process in \eqref{alt MP} is a martingale for every $k\geq 1$ and every $\psi\in \mathbb{L}^1((\R^d)^k)$, and Lemma~\ref{property 2} is proved.
\end{proof}

\begin{remark}\label{rk:test functions}
Note that either of these two sets of test functions, (\ref{test function 1}) or (\ref{test function 2}), is sufficient to characterise the law of the SLFVS (see Lemma~\ref{lem: dense set} for the first set, and Lemma~2.1(c) in \cite{VW2012} for the second), and so we can use them interchangeably.
In particular, the family (\ref{test function 1}) will be more convenient in proving
the convergence of our rescaled $\cM_\lambda$-valued processes, whereas the duality relation that will give us the uniqueness of the limit is based on the family (\ref{test function 2}).
\end{remark}

Armed with Lemma~\ref{property 2}, we can now prove Proposition~\ref{prop: dual}.

\begin{proof}[Proof of Proposition~\ref{prop: dual}]
Despite the fact that Theorem~4.4.11 in \cite{EK1986} does not directly apply, we follow its proof closely (with $\alpha\equiv 0\equiv \beta$). Let $(M_t)_{t\geq 0}$ be solution to the martingale problem~\eqref{MP}. For every $s,t\geq 0$, let
\[
F(s,t):= \E_{M^0}\Big[\bfE\big[ D(M_s,\Xi_t)\, |\, \Xi_0 \sim \mu_\psi \big] \Big].
\]
By Lemma~\ref{lem:absolutely cont}, since $\Xi_0$ has law $\mu_\psi$, at every time $t\geq 0$ the locations of the atoms of $\Xi_t$ have a joint distribution which is absolutely continuous with respect to Lebesgue measure. Let us write $\psi_t^{(n)}$ for the density of these points conditionally on the event $\{N_t=n\}$. We thus have, by \eqref{average test function},
\begin{align}
F(s,t) & = \E_{M^0}\bigg[\bfE\bigg[\int_{(\R^d)^{N_t}} \psi_t^{(N_t)}\big(x_1,\ldots,x_{N_t}\big)\bigg(\prod_{j=1}^{N_t}w_s(x_j)\bigg)\rmd x_1 \cdots \rmd x_{N_t} \, \Big|\, \Xi_0 \sim \mu_\psi \bigg] \bigg] \label{def Fst} \\
&= \bfE\left[\E_{M^0}\big[D_{\psi_t^{(N_t)}}(M_s)\big]\, \big|\, \Xi_0\sim \mu_\psi\right],\nonumber
\end{align}
where the last line uses Fubini's theorem. Since for every $n\in \{1,2,\ldots\}$ we have $\psi_t^{(n)}\in \mathbb{L}^1((\R^d)^n)$, we can use Lemma~\ref{property 2} and write that
\begin{equation}\label{increment F}
F(s,t)-F(0,t) = \bfE\bigg[\E_{M^0}\bigg[\int_0^s \cL D_{\psi_t^{(N_t)}}(M_\tau)\rmd\tau\bigg]\, \bigg|\, \Xi_0\sim \mu_\psi\bigg].
\end{equation}
On the other hand, for any fixed $s\geq 0$ and any $t\geq 0$, we can rewrite the expression on the r.h.s. of \eqref{def Fst} as
\begin{align}
\E_{M^0}\bigg[ \bfE\bigg[ \prod_{j=1}^{N_t}w_s\big(\xi_t^j\big)\bigg|\, \Xi_0\sim \mu_\psi \bigg]\bigg]& = \E_{M^0}\Big[ \bfE\Big[\exp \big(\la \Xi_t,\ln w_s\ra \big)\, \big|\, \Xi_0\sim \mu_\psi \Big]\Big] \nonumber \\
& = \E_{M^0}\Big[\bfE\Big[\Phi_{\exp,\ln w_s}(\Xi_t)\, \big|\, \Xi_0\sim \mu_\psi \Big]\Big], \label{other expr}
\end{align}
where we have fixed a representative $w_s$ of the density of $M_s$ (since by \eqref{average test function}, the r.h.s. of \eqref{other expr} is independent of the choice of this representative) and $\Phi_{\exp,\ln w_s}$ is defined as in \eqref{test functions for dual} with $F=\exp$ and $f=\ln w_s$. Here we use the convention that $\Phi_{\exp,\ln w_s}(\Xi)=0$ whenever at least one of the atoms $x$ of $\Xi$ is such that $w_s(x)=0$. With this convention, the definition of $\cG \Phi_{\exp,\ln w_s}(\Xi)$ given in \eqref{def G} still makes sense, the function $\Phi_{\exp,\ln w_s}$ takes its values in $[0,1]$ and the same bound as in \eqref{bound on G} controls the expectation of $\cG \Phi_{\exp,\ln w_s}(\Xi_\tau)$ for every $\tau\in [0,t]$ (since in this case $F'(\la \Xi_\tau,f\ra)=\prod_i w_s(\xi^i_\tau)\in [0,1]$). We may therefore use
Proposition~\ref{prop:well-defined dual} to write that
\begin{equation}\label{increment F2}
F(s,t)-F(s,0) = \E_{M^0}\bigg[\bfE\bigg[\int_0^t \cG \Phi_{\exp,\ln w_s}(\Xi_\tau)\rmd\tau\, \bigg|\, \Xi_0\sim\mu_\psi \bigg]\bigg].
\end{equation}
It remains to show that for every $s,t\geq 0$,
\begin{equation}\label{lemme EK}
\bfE\Big[\E_{M^0}\big[\cL D_{\psi_{t-s}^{(N_{t-s})}}(M_s)\big]\, \big|\, \Xi_0\sim\mu_\psi \Big] = \E_{M^0}\Big[\bfE\big[\cG \Phi_{\exp,\ln w_s}(\Xi_{t-s})\, \big|\, \Xi_0\sim\mu_\psi \big]\Big],
\end{equation}
so that we may use Lemma~4.4.10 in \cite{EK1986} to conclude that
\[
F(t,0)=F(0,t),
\]
which is equivalent to \eqref{dual bancale}.

Using \eqref{def G} and shortening the notation $I_{x,r}(\Xi_{t-s})$ into $I$, we have,
\begin{align}
& \cG \Phi_{\exp,\ln w_s}(\Xi_{t-s})\nonumber\\
& = \int_{\R^d} \int_0^\infty\int_0^1 \int_{B(x,r)}\frac{1}{V_r}\bigg(\prod_{j\in I^c} w_s(\xi_{t-s}^j)\bigg) \nonumber \\
& \qquad \  \times \Bigg[\sum_{D\subseteq I;|D|\geq 1}u^{|D|}(1-u)^{|I\setminus D|}\bigg(w_s(y)\prod_{i\in D^c}w_s(\xi_{t-s}^i) - \prod_{i\in I}w_s(\xi_{t-s}^i)\bigg)\Bigg]\rmd y \nu_r(\rmd u)\mu(\rmd r)\rmd x \nonumber\\
&\quad + \int_{\R^d} \int_0^\infty\int_0^1 \int_{B(x,r)^2}\frac{1}{V_r^2}\bigg(\prod_{j\in I^c} w_s(\xi_{t-s}^j)\bigg) \Bigg[\sum_{D\subseteq I;|D|\geq 1}u^{|D|}(1-u)^{|I\setminus D|} \nonumber\\
& \qquad \qquad  \times \bigg(w_s(y)w_s(z)\prod_{i\in D^c}w_s(\xi_{t-s}^i) - \prod_{i\in I}w_s(\xi_{t-s}^i)\bigg)\Bigg]\rmd y\rmd z \nu'_r(\rmd u)\mu'(\rmd r)\rmd x \nonumber \\
&= \int_{\R^d}\! \int_0^\infty\!\int_0^1\! \int_{B(x,r)}\frac{1}{V_r}\bigg(\prod_{j\in
I^c} w_s(\xi_{t-s}^j)\bigg) \Bigg[w_s(y)\Bigg(\sum_{D\subseteq I;|D|\geq 1}u^{|D|}(1-u)^{|I\setminus D|}\prod_{i\in D^c}w_s(\xi_{t-s}^i)\Bigg)\nonumber\\
& \qquad \qquad - \big(1-(1-u)^{|I|}\big)\prod_{i\in I}w_s(\xi_{t-s}^i)\Bigg] \rmd y\nu_r(\rmd u)\mu(\rmd r)\rmd x \nonumber\\
& \qquad  + \int_{\R^d} \int_0^\infty\int_0^1 \int_{B(x,r)^2}\frac{1}{V_r^2}\bigg(\prod_{j\in I^c} w_s(\xi_{t-s}^j)\bigg)\Bigg[w_s(y)w_s(z)\nonumber\\
& \qquad \qquad \times \Bigg(\sum_{D\subseteq I;|D|\geq 1}u^{|D|}(1-u)^{|I\setminus D|}\prod_{i\in D^c}w_s(\xi_{t-s}^i)\Bigg)- \big(1-(1-u)^{|I|}\big)\prod_{i\in I}w_s(\xi_{t-s}^i)\Bigg]\nonumber \\
& \qquad \qquad \qquad \qquad\qquad \qquad\qquad \qquad \qquad \qquad\qquad \qquad\qquad \qquad  \rmd y\rmd z\nu'_r(\rmd u)\mu'(\rmd r)\rmd x. \label{identity of generators1}
\end{align}

On the other hand, using \eqref{def L} and writing $\psi_{t-s}$ for $\psi_{t-s}^{(N_{t-s})}$ to ease the notation, we obtain
\begin{align}
 \cL &D_{\psi_{t-s}}(M_s) \nonumber\\
=&\int_{\R^d} \int_0^\infty \int_0^1 \int_{B(x,r)}\frac{1}{V_r}\int_{(\R^d)^{N_{t-s}}} \psi_{t-s}(x_1,\ldots,x_{N_{t-s}}) \bigg(\prod_{j\in I^c}w_s(x_j)\bigg)\nonumber\\
& \qquad \times \bigg[w_s(y)\prod_{j\in I}\big((1-u)w_s(x_j)+u\big) + (1-w_s(y))\prod_{j\in I}\big((1-u)w_s(x_j)\big) - \prod_{j\in I}w_s(x_j)\bigg] \nonumber\\
& \qquad \qquad\qquad \qquad\rmd x_1\ldots \rmd x_{N_{t-s}} \rmd y\nu_r(\rmd u)\mu(\rmd r)\rmd x\nonumber\\
&+ \int_{\R^d} \int_0^\infty \int_0^1 \int_{B(x,r)^2}\frac{1}{V_r^2}\int_{(\R^d)^{N_{t-s}}} \psi_{t-s}(x_1,\ldots,x_{N_{t-s}}) \bigg(\prod_{j\in I^c}w_s(x_j)\bigg) \nonumber \\
& \qquad \times \bigg[w_s(y)w_s(z)\prod_{j\in I}\big((1-u)w_s(x_j)+u\big) + (1-w_s(y)w_s(z))\prod_{j\in I}\big((1-u)w_s(x_j)\big) \nonumber \\
&\qquad \qquad\qquad \qquad - \prod_{j\in I}w_s(x_j)\bigg] \rmd x_1\ldots \rmd x_{N_{t-s}} \rmd y \rmd z\nu'_r(\rmd u)\mu'(\rmd r) \rmd x. \label{def L2}
\end{align}
Now, the first integral in the above is equal to
\begin{align}
&\int_{\R^d} \int_0^\infty \int_0^1 \int_{B(x,r)}\frac{1}{V_r}\int_{(\R^d)^{N_{t-s}}}\, \psi_{t-s}(x_1,\ldots,x_{N_{t-s}}) \prod_{j\in I^c}w_s(x_j) \nonumber \\
& \qquad \times \Bigg[w_s(y) \sum_{D\subseteq I;|D|\geq 1}u^{|D|}(1-u)^{|I\setminus D|}\prod_{j\in D^c}w_s(x_j)
 - \big(1-(1-u)^{|I|}\big)\prod_{j\in I}w_s(x_j)\Bigg] \nonumber\\
& \qquad \qquad \qquad \qquad\qquad \qquad\qquad \qquad \qquad\qquad  \rmd x_1\ldots \rmd x_{N_{t-s}} \rmd y \nu_r(\rmd u)\mu(\rmd r) \rmd x.\label{identity of generators}
\end{align}
Taking the expectation of \eqref{identity of generators} with respect to $\bfE[\,\cdot\,| \Xi_0\sim \mu_\psi]$, we obtain that it is equal to the expectation of the first term on the r.h.s. of \eqref{identity of generators1}. The same holds for the expectation of the second part of \eqref{def L2} and that of the
second part of \eqref{identity of generators1}. Taking then the expectation with respect to $\E_{M_0}$ (and using Fubini's theorem), we arrive at the desired equality \eqref{lemme EK}. Hence, we can use Lemma~4.4.10 in \cite{EK1986} to conclude that
\[
F(t,0)=F(0,t),
\]
and \eqref{dual bancale} is proved.
\end{proof}

\section{Heuristics for the large-scale behaviour of the SLFVS and its dual process}
\label{heuristics}
In this section, we provide an informal
justification of our choices for the parameters $\beta$, $\gamma$ and
$\delta$ in our scalings. Recall from the introduction that we should like to establish
scalings of the selection and impact parameters, $s_n$ and $u_n$, for which selection
will leave a trace in the long-term evolution of the population, without leading to an instantaneous invasion by the favoured allele. In particular, we wish to complement the work of \cite{EFPS2016,EFS2017}, in
which the impact (or fraction of the local population replaced during an event) is kept of order $\cO(1)$ while the selection coefficient goes to $0$ as the scaling parameter $n$ tends to infinity, modelling a population in which the local densities of individuals are low and therefore any event leads to the replacement of a macroscopic fraction of the local population. In our work, we always assume that local population densities are very high and reproduction impacts only a very small fraction of the individuals present in the affected area (\emph{i.e.}, we assume that $u \ll 1$). We consider long timescales $(nt,\,t\geq 0)$ and identify the orders of
magnitude of $u_n$ and $s_n$, and the corresponding spatial scale
$(n^\beta x,\, x\in\R^d)$ over which the dynamics of the population will
converge as $n\to\infty$ to a Fisher-KPP type evolution. In passing, we shall also see why there is no way to obtain a stochastic limit in more than one dimension when we assume that both the impact $u_n$ of each event and the relative frequency $s_n$ of selective to neutral events tend to $0$, and we shall justify the claims made in Remarks~\ref{rk:1d deterministic} and \ref{rk:1d deterministic stable} about the range of parameters leading to a deterministic Fisher-KPP process in one dimension.

Let us start with the case of local reproduction. As usual in the spatial Lambda-Fleming-Viot framework, it is easier to first think about the corresponding scaled dual processes. Chapter~7 in \cite{LIA2009} suggests that if there is a regime of parameters in which the SLFVS converges to the solution to the Fisher-KPP equation, and in one dimension to its stochastic counterpart, then the sequence of corresponding dual processes should converge to a branching Brownian motion in which, in dimension $1$, pairs of particles coalesce at a rate proportional to their collision local time (we explain the correspondence between our frameworks in the paragraph on uniqueness of the limit in the proof of Theorem 1.11, see Section 5.1). Let us thus analyse the different types of event which can affect the dual process when time is sped up by $n$ and space is scaled down by $n^\beta$, focusing first on what happens to a single particle. During a neutral event, if this particle is marked (with probability $u_n=\cO(n^{-\gamma})$), then it is removed and replaced by another particle whose distribution is uniformly distributed in the region of the event (of radius $Rn^{-\beta}$ in our new units). We see this as a \emph{jump} of the particle. Because the component of the intensity of $\Pi^N$ corresponding to the centres of events is Lebesgue measure on $\R^d$,
when the particle jumps, the location of the centre of the corresponding
event is uniformly distributed in the ball of radius $Rn^{-\beta}$ around the current location $x$ of the particle. Consequently, the position of the particle ``after the jump'' belongs to $B(x,2Rn^{-\beta})$ a.s. and has a radially symmetric distribution around $x$. Summing up the above, in our new units a single particle jumps at rate $\cO(n^{1-\gamma})$ and makes mean zero, finite variance, jumps of size bounded by a constant times $1/n^\beta$. For this jump process to converge to some non-trivial process (Brownian motion, in fact) as $n\rightarrow \infty$, we thus have to assume that
\begin{equation}\label{eq BM}
1-\gamma=2\beta.
\end{equation}
Now consider what happens at a selective event. Again our particle is marked with probability $\cO(n^{-\gamma})$, and in this case, the two particles which replace it are created at a separation of order $\cO(n^{-\beta})$. Consequently, they may be overlapped by a new event very quickly (after a time of order $\cO(n^{-1})$), and
then with probability $u_n^2=\cO(n^{-2\gamma})$ they are both erased and
replaced by a single ``parental'' particle (we see this type of event as a
``coalescence''). In the limit as $n\to\infty$, we will only ``see'' the
branching event before it is erased by such a coalescence if the two
particles have positive probability of moving apart to a distance of order
one before (perhaps) coalescing. Let us find conditions under which we can
expect this to hold.

In our new timescale, each particle is overlapped by an event at rate $\cO(n)$. The probability that only one of the two particles is marked during such an event, and therefore ``jumps'' to a location at distance $\cO(n^{-\beta})$ while the other stands still, is $un^{-\gamma}(1-\cO(n^{-\gamma}))$. Furthermore, as soon as the two particles are at distance larger than $2Rn^{-\beta}$, they cannot be overlapped by the same event and so they jump independently of each other according to a continuous time random walk. Hence, what we actually have to understand is how many times the two particles come back to a separation less than $2Rn^{-\beta}$ before they manage to move apart to a separation of $\cO(1)$. Indeed, the same type of analysis as the one carried out in the proof of Lemma~6.6 in \cite{BEV2010} (see Lemma~\ref{lem: incursions} and below in Section~\ref{s: convergence 1} of the present work) shows that when they come together, the two particles remain at distance less than $2Rn^{-\beta}$ during a number of events \emph{affecting them} of order $\cO(1)$ (which translates into a number of events simply \emph{overlapping them} of order $\cO(n^\gamma)$). Hence, the probability that they are both affected by an event and coalesce before separating again to a distance more than $2Rn^{-\beta}$ is of the same order as the probability of coalescence during a single event conditionally on at least one of the particles being marked, which is $\cO(n^{-2\gamma}/n^{-\gamma})=\cO(n^{-\gamma})$. From this, we can conclude in particular that the two particles will need to come back
``together'' $\cO(n^\gamma)$ times before they have a positive probability of
both being affected by the same event and therefore
coalescing into common ancestral particle(s).

When they are more than $2Rn^{-\beta}$ apart, the two particles jump
independently according to a continuous time symmetric random walk with
step sizes of order $\cO(n^{-\beta})$, and so the separation between them is
also a symmetric random walk with step size of this order. We are interested in the probability that during an excursion away from $B(0,2Rn^{-\beta})$, the difference walk reaches a distance of order $1$, \emph{i.e.} $n^{\beta}$ times larger than its initial value. It is convenient to work in our original space units. For a symmetric continuous-time random walk with step size of order $\cO(1)$, starting at distance slightly larger than $2R$ from $0$, the probability of reaching distance $n^\beta$ from $0$ before reentering $B(0,2R)$ has the same order as the probability that the number of steps to come back within $B(0,2R)$ is
larger than $n^{2\beta}$. This in turn will have the same order as the
corresponding quantity for simple symmetric random walk.
Using Proposition~5.1.1 in \cite{LL2010} when $d=1$, Theorem~1 in \cite{RR1966} when $d=2$ and the transience of simple symmetric random walk when $d\geq 3$, we obtain that in one dimension, this probability is of order $\cO(n^{-\beta})$; in two dimensions, this probability is of order $\cO(1/\ln n)$; in dimension $d\geq 3$, this probability tends to $p\in (0,1)$ as $n\rightarrow \infty$. Consequently, when $d\geq 2$ the probability that the two particles come back together $\cO(n^\gamma)$ times before they separate to a distance of $\cO(1)$
tends to $0$ and, in the limit, a given branching event is never followed by
instantaneous coalescence. Since branching events happen at a rate $\cO(ns_nu_n)=\cO(n^{1-\delta-\gamma})$, we need to impose that
\begin{equation}\label{rate selection}
1-\delta-\gamma = 0
\end{equation}
if they are to occur at rate $\cO(1)$ in the limit. On the other hand, in
one dimension, we see that if $\beta>\gamma$, with probability tending to one,
the two particles will coalesce back together before they can separate to a
distance of $\cO(1)$ and in the limit, all branching events are
cancelled (\emph{i.e.}, there is no branching in the limiting dual).
If $\beta<\gamma$, the two particles become separated at distance $\cO(1)$
before they have any chance to coalesce, and all branching events are
conserved in the limit; in contrast coalescence will never be seen in the limit
for particles starting at any separation.
Finally, when
\begin{equation}\label{branch and coalesce}
\beta=\gamma
\end{equation}
the particles have positive probability of separating to a distance
of $\cO(1)$ before coalescing, but coalescence of particles
happens in finite time a.s. in the limit. In the last two cases, we also
need to impose Condition~\eqref{rate selection} for branching events to occur
at rate $\cO(1)$ in the limit.

Solving the system given by Conditions~\eqref{eq BM}, \eqref{rate selection} and \eqref{branch and coalesce}, we obtain $\beta=1/3=\gamma$ and $\delta=2/3$ as specified in Section~\ref{ss: main results}. Observe that if we keep Conditions~\eqref{eq BM} and \eqref{rate selection} and replace Condition~\eqref{branch and coalesce} by
\begin{equation}\label{branch only}
\beta<\gamma,
\end{equation}
even in one dimension two particles can branch but not coalesce in the limit and, as explained in the paragraph on uniqueness of the limit in the proof of Theorem~\ref{th:fixed} (see Section~\ref{ss:proof1}), the limit of the SLFVS is the (measure-valued) solution to the deterministic Fisher-KPP equation.

\begin{remark}
In two dimensions, our heuristics suggest that if instead of scaling space by $n^\beta$, we were to scale it by a more general factor $k_n$, then the relation between the parameters allowing the number of excursions necessary for the particles to separate to distance $\cO(k_n)$ to be of the same order as the number
of times the particles need to come ``together'' before they coalesce is
\begin{equation}\label{alternate cond}
\ln k_n \approx 1/u_n.
\end{equation}
(Here we use $\approx$ to mean that the two quantities are of the same order of magnitude.) Together with the relations $nu_n \approx k_n^2$ ensuring that the limiting motion is Brownian motion, and $nu_ns_n \approx 1$ guaranteeing that branching occurs at rate $\cO(1)$ in the limit, this gives us that
\begin{equation}
\ln k_n \approx \frac{n}{k_n^2}, \qquad \hbox{and so }k_n \approx \sqrt{\frac{n}{\ln n}}.
\end{equation}
But with this choice of parameters, the bound \eqref{control Bn} we shall establish for the predictable quadratic variation of $F(\la\bw^n_\cdot,f\ra)$ reads
\begin{equation}
u_n^2nk_n^d\times k_n^{-2d} \approx \frac{k_n^2}{n}\approx \frac{1}{\ln k_n}\rightarrow 0 \quad \hbox{as }n\rightarrow \infty.
\end{equation}
Therefore, something more subtle happens here and even with a more general form of scaling of space, the limiting process of allele frequencies is still deterministic.
\end{remark}

We now turn to the stable case. As before, we first consider the ``jumps'' of a single particle. Again, such a jump occurs at rate $\cO(n^{1-\gamma})$, and its new position is chosen uniformly over a ball whose radius is
given by the intensity measure $\mu$ with polynomial decay described in (\ref{def constants}). Consequently, if
we choose $nu_n\propto n^{\alpha\beta}$, \emph{i.e.}
\begin{equation}\label{cond stable}
1-\gamma=\alpha\beta,
\end{equation}
then in the limit as $n\rightarrow\infty$ the jump process will converge to a
symmetric $\alpha$-stable process with index $\alpha$.
Second, in order to see any branching of particles due to selective events
at all, as before
we need $ns_nu_n$ to be order one,
that is Condition~\eqref{rate selection} to be fulfilled. Finally, let us
consider the simultaneous removal of two particles, to be
replaced by ``parental particle(s)'' (what we called a coalescence earlier). Since $u_n\rightarrow 0$ as $n\rightarrow\infty$,
although it is now the case that two particles can always be affected
by the same event (the radii of the events are not bounded), ``most of the time'' they will not and their jumps are
almost independent. Consequently, the difference of their positions is also approximately described by an
$\alpha$-stable process. Now, because events of radius $\cO(1)$ (in our original units) are much more
frequent than events of large radii
$\cO(n^a)$ for any $a>0$, and the probability that both particles belong to the fraction of the local population
replaced during an event is tiny ($u_n^2=u^2n^{-2\gamma}$), if coalescence is to happen in the limit, then
we expect it to be driven by the smaller events. Note that because there is no bound on the event radii, we can no longer perform the same decomposition into excursions away from some ball as in the fixed radius case. The following argument only gives the intuition behind Lemma~\ref{lem: limit coal} in Section~\ref{s: convergence 2}, which allows us to control the coalescence rate of the two particles.

In more than one dimension, the
rotation-invariant $\alpha$-stable processes with $\alpha\in (1,2)$ are transient
(see Example~37.19$(ii)$ in \cite{SAT1999}), and so as in the fixed radius case, this tells us
that the two particles do not spend enough time close together for coalescence to occur, whatever our choice of
$\beta,\gamma,\alpha$ consistent with the previous conditions.
In one dimension, we have not found a simple heuristic explanation for the last
condition on the parameters (which one would expect to be analogous to the comparison
between the number and lengths of visits
in a neighbourhood of zero for the difference process, and the coalescence rate of the two particles,
carried out in the fixed radius case).
Instead, the condition
\begin{equation}\label{cond coal stable}
\gamma=(\alpha -1)\beta
\end{equation}
will emerge when we control the second term on the r.h.s. of
(\ref{finite variation}), which corresponds to the variance term in the limiting process (and thus to the coalescence term in the dual process). See also Equations (\ref{vq1}) and (\ref{vq2}) and the surrounding paragraphs.
In the end, we have three equations in
three unknowns (in one dimension) and solving gives the values in Equation~(\ref{def constants}). As in the fixed radius case, we may replace Condition~\eqref{cond coal stable} by
\begin{equation}\label{cond coal stable 2}
(\alpha -1)\beta<\gamma
\end{equation}
and obtain a limiting dual in which, in any dimension, branching occurs at rate $\cO(1)$ but coalescence never occurs.

As a final comment, notice that in this work we have chosen a particular form for the parameters $u_n$ and $s_n$, given in \eqref{assumptions par}, and for the scaling of space (by $n^\beta$). We have argued that, within this particular framework (which nonetheless covers a wide range of scenarios), in more than one dimension it was not possible to find values for $\beta,\gamma,\delta$ such that the martingale problem characterising the limiting process $(M_t^\infty)_{t\geq 0}$ contains a Laplacian (or fractional Laplacian) term, a drift term due to the slight advantage of type~1 individuals \emph{and} a martingale term corresponding to the noise in the Fisher-KPP-like equation satisfied by the process. Of course this does not prove that other forms of parameters and spatial scalings would not yield a stochastic limit even under the assumption that $(u_n)_{n\geq 1}$ and $(s_n)_{n\geq 1}$ tend to $0$ as $n$ tends to infinity that we have imposed. One way to investigate this question is to use the correspondence between the fact that the limiting process is stochastic and the property that the events that we loosely call ``coalescence of particles'' have positive probability to happen in the limiting dual process (a general property of continuous-site stepping-stone models, see Section~5 in \cite{EVA1997}). We have not been able to find scalings for which, in the limiting dual, particles may ``move in space'', ``branch'' and ``coalesce'' (even in a non-local way), with positive probability when $d\geq 2$. In fact, due to the facts that rotation-invariant $\alpha$-stable processes are transient for $d>\alpha$ and that coalescence happens $u_n$ times more slowly than movement of particles before taking the limit, we conjecture that such scalings do not exist, even in the case of $\alpha$-stable radii. We leave this delicate question to the interested reader.

\section{Convergence of the rescaled SLFVS and its dual - the fixed radius case} \label{s: convergence 1}
In this section, we prove Theorem~\ref{th:fixed} and, from it, deduce Theorem~\ref{th: conv duals}.

\subsection{Proof of Theorem~\ref{th:fixed}.}\label{ss:proof1}
The proof proceeds in the usual way. First, we show that the sequence of non-Markovian processes $(\bM^n)_{n\geq 1}$ is tight in $D_{\cM_\lambda}[0,\infty)$ and that any limit point has a.s. continuous trajectories. Next, we prove that any limit point $M^\infty$ satisfies the martingale problem stated in Theorem~\ref{th:fixed}. Finally, we show that there is at most one solution in $D_{\cM_\lambda}[0,\infty)$ to this martingale problem (again thanks to a duality argument), which will allow us to conclude that indeed $(\bM^n)_{n\geq 1}$ converges to this solution.

\medskip
\textbf{1) Tightness and continuity of the limit.}

\medskip
\noindent First, let $n\geq 1$. Since the (unscaled) SLFVS $(M_t)_{t\geq 0}$ with parameters given in \eqref{cond intensity}, \eqref{cond intensity 2}, \eqref{assumptions par} and $\mu =\delta_R$ has sample paths in $D_{\cM_\lambda}[0,\infty)$ by Theorem~\ref{th:existence}, so has the locally averaged and scaled process $\bM^n$ (recall that its density is defined by Relation~\eqref{def wn1}).

Let us now show that the sequence $(\overline{M}^n)_{n\geq 1}$ is relatively compact. We proceed as in the proof of Theorem~\ref{th:existence} and refer to the paragraph Proof of $(i)$, item~$(c)$, in Section~\ref{s: proof existence} for a more detailed justification of the steps taken below. Again, due to the compactness of $\cM_\lambda$ endowed with the topology of vague convergence and the fact that the set of functions of the form $\Psi_{F,f}$ with $F\in C^3(\R)$ and $f\in C_c^\infty(\R^d)$ is dense (for the supremum norm) in the set of functions of the form $\Psi_{F,f}$ with $F\in C^1(\R)$ and $f\in C_c(\R^d)$, and is therefore dense in $C(\cM_\lambda)$ by Lemma~\ref{lem: dense set}, Theorem~3.9.1 in \cite{EK1986} tells us that it suffices to show the relative compactness of the sequence of real-valued processes $(\Psi_{F,f}(\overline{M}^n))_{n\geq 1}$ for every $F\in C^3(\R)$ and $f\in C_c^\infty(\R^d)$. Second, for each such function $\Psi_{F,f}$, we use the Aldous-Rebolledo criterion \cite{ALD1978,REB1980} to reduce the problem to tightness of the sequences of the predictable finite variation parts and of the predictable quadratic variation of the martingale parts
of $(\Psi_{F,f}(\bM^n_t))_{t\geq 0}$. More precisely, since $\Psi_{F,f}$ is a bounded function on $\cM_\lambda$, we directly have that for every $t\geq 0$, the sequence $(\Psi_{F,f}(\bM^n_t))_{n\geq 1}$ is tight. Writing $({\cal A}^n_t)_{t\geq 0}$ (\emph{resp.}, $({\cal Q}^n_t)_{t\geq 0}$) for the finite variation part (\emph{resp.}, the quadratic variation of the martingale part) of $(\Psi_{F,f}(\bM^n_t))_{t\geq 0}$, it remains to prove that for every $T>0$, every sequence of stopping times $(\tau_n)_{n\geq 1}$ bounded by $T$, and every $\e>0$, there exists $\eta>0$ such that
\begin{equation}\label{tightness cal An}
\limsup_{n\rightarrow \infty} \sup_{\theta\in [0,\eta]} \P\big[|{\cal A}^n_{\tau_n+\theta} -{\cal A}^n_{\tau_n}|>\e\big]\leq \e,
\end{equation}
and
\begin{equation}\label{tightness cal Qn}
\limsup_{n\rightarrow \infty} \sup_{\theta\in [0,\eta]} \P\big[|{\cal Q}^n_{\tau_n+\theta} -{\cal Q}^n_{\tau_n}|>\e\big]\leq \e.
\end{equation}
In what follows, we first establish an expression for the processes ${\cal A}^n$ and ${\cal Q}^n$ (see \eqref{n-finite variation} and \eqref{n-quad variation}), and then we prove \eqref{tightness cal An} and \eqref{tightness cal Qn} by decomposing these expressions into several parts that we control separately.

To find expressions for the predictable finite and quadratic variation parts of $\Psi_{F,f}(\bM^n)$ for any fixed $n\geq 1$, let us begin by considering the unscaled SLFVS $(M_t)_{t\geq 0}$ with reproduction events of fixed radius $R$, and parameters $u_n$, $s_n$ (notice that for simplicty we have suppressed the dependence of $(M_t)_{t\geq 0}$ on $n$ in the notation). Also, the function $F$ will be fixed but for a moment we replace the function $f$ by any function $\varphi\in C^{\infty}_c(\R^d)$. A judicious choice of $\varphi$ will then yield the desired function of the scaled and locally averaged process $\bM^n$.

Let $\varphi\in C^{\infty}_c(\R^d)$. By Theorem~\ref{th:existence}, we know that before scaling space and time, the extended generator of the SLFVS, acting on the function $\Psi_{F,\varphi}$, is given by
\begin{align}
{\cL} \Psi_{F,\varphi}(M)&
= \int_{\R^d}\int_{B(x,R)^2}\frac{1}{V_R^2}\, \Big\{w(y)(1+s_nw(z))\big[ F(\la \Theta^+_{x,R,u_n}(w),\varphi \ra) - F(\la w,\varphi\ra)\big] \\
 &\quad + (1-w(y)+s_n(1-w(y)w(z)))\big[ F(\la \Theta^-_{x,R,u_n}(w),\varphi \ra) - F(\la w,\varphi\ra)\big]\Big\} \rmd y \rmd z \rmd x, \nonumber
\end{align}
where $w$ is a representative of the density of $M$. This gives us that the
predictable finite variation part of $(\Psi_{F,\varphi}(M_t))_{t\geq 0}$ is
\begin{equation}\label{cal A}
{\cal A}_t=\int_0^t {\cL} \Psi_{F,\varphi}(M_s)\, \rmd s,\qquad t\geq 0.
\end{equation}
Furthermore, the martingale problem stated in Theorem~\ref{th:existence} applies to
$\Psi_{F^2,\varphi}=F(\la \cdot,\varphi\ra)^2$, which allows us to obtain (using It\^o's formula) that the predictable quadratic variation of $(\Psi_{F,\varphi}(M_t))_{t\geq 0}$ at any time $t\geq 0$ is given by
\begin{align}
{\cal Q}_t & = \int_0^t \int_{\R^d}\int_{B(x,R)^2}\frac{1}{V_R^2}\, \Big\{w_s(y)(1+s_nw_s(z))\big[ F(\la \Theta^+_{x,R,u_n}(w_s),\varphi \ra) - F(\la w_s,\varphi\ra)\big]^2 \label{cal Q}\\
 &\quad + (1-w_s(y)+s_n(1-w_s(y)w_s(z)))\big[ F(\la \Theta^-_{x,R,u_n}(w_s),\varphi \ra) - F(\la w_s,\varphi\ra)\big]^2\Big\} \rmd y \rmd z \rmd x \rmd s.\nonumber
\end{align}
Let us now consider the Markov process $(M^n_t)_{t\geq 0}$ whose density at time $t$ is $w^n_t(\cdot):=w_{nt}(n^{1/3}\,\cdot)$.
We set
\begin{equation}\label{def Bn}
B_n(x)=B(x,n^{-1/3}R)
\end{equation}
and write $\bw(x)=n^{d/3}V_R^{-1}\int_{B_n(x)}w(z)\rmd z$.
In particular, in the notation of Section~\ref{ss: main results} we have for every $t\geq 0$
\begin{equation}
\frac{n^{d/3}}{V_R}\int_{B_n(x)}w^n_t(z)\rmd z = \frac{1}{V_R}\int_{B(n^{1/3}x,R)}w_{nt}(y)\rmd y=\bw^n_t(x).
\end{equation}
From our expression for $\cL$,
accelerating time by a factor $n$ and performing several changes of the spatial variables, we obtain that the extended generator of $M^n$ is given by
\begin{align}
&{\cL}^n \Psi_{F,\varphi}(M) \nonumber\\
& = n\int_{\R^d}\int_{B(x,R)^2}\frac{1}{V_R^2}\, \Big\{w(n^{-1/3}y)\big(1+s_n w(n^{-1/3}z))\big[ F(\la \Theta^+_{n^{-1/3}x,n^{-1/3}R,u_n}(w),\varphi \ra) \nonumber\\
&\qquad \qquad \qquad \qquad \qquad \qquad \qquad \qquad \qquad \qquad \qquad \qquad \qquad \qquad \qquad \qquad \quad   - F(\la w,\varphi\ra)\big] \nonumber\\
& \quad + (1-w(n^{-1/3}y)+s_n(1-w(n^{-1/3}y)w(n^{-1/3}z)))\big[ F(\la \Theta^-_{n^{-1/3}x,n^{-1/3}R,u_n}(w),\varphi \ra) \nonumber\\
& \qquad \qquad \qquad \qquad \qquad \qquad \qquad \qquad \qquad \qquad \qquad \qquad \qquad \qquad \qquad   - F(\la w,\varphi\ra)\big]\Big\}\rmd y\rmd z \rmd x \nonumber \\
& = n^{1+\frac{d}{3}}\int_{\R^d} \Big\{\bw(x)(1+s_n\bw(x))\big[F(\la \Theta^+_{x,n^{-1/3}R,u_n}(w),\varphi\ra)-F(\la w,\varphi\ra)\big] \nonumber \\
& \qquad \qquad  +\big(1-\bw(x)+s_n(1-\bw(x)^2)\big)\big[F(\la \Theta^-_{x,n^{-1/3}R,u_n}(w),\varphi\ra)-F(\la w,\varphi\ra)\big]\Big\}\rmd x. \label{generator wn}
\end{align}
The predictable finite variation part of $(\Psi_{F,\varphi}(M^n_t))_{t\geq 0}$ is thus $(\int_0^t{\cL}^n\Psi_{F,\varphi}(M^n_s)\rmd s)_{t\geq 0}$. Likewise, its predictable quadratic variation at time $t$ is equal to
\begin{align}
& n^{1+\frac{d}{3}}\int_0^t \int_{\R^d} \Big\{\bw^n_s(x)(1+s_n\bw^n_s(x))\big[F(\la \Theta^+_{x,n^{-1/3}R,u_n}(w^n_s),\varphi\ra)-F(\la w^n_s,\varphi\ra)\big]^2 \nonumber \\
& \quad +\big(1-\bw^n_s(x)+s_n(1-\bw^n_s(x)^2)\big)\big[F(\la \Theta^-_{x,n^{-1/3}R,u_n}(w_s^n),\varphi\ra)-F(\la w^n_s,\varphi\ra)\big]^2\Big\}\rmd x\rmd s.
\end{align}

Finally, it remains to evaluate the above expressions with $\varphi$ of the form
\begin{equation}\label{def varphif}
\varphi_f(x)=\frac{n^{d/3}}{V_R}\int_{B_n(x)} f(y)\rmd y
\end{equation}
for our fixed $f\in C_c^{\infty}(\R^d)$ and to use the fact that, by Fubini's Theorem,
\begin{equation}
\la w^n, \varphi_f\ra = \int_{\R^d}\int_{\R^d}\, w^n(y)\frac{n^{d/3}}{V_R}f(z) \ind_{\{|z-y|\leq n^{-1/3}R\}} \rmd y\rmd z = \la \bw^n,f\ra,
\end{equation}
to obtain that the predictable finite variation part of $(\Psi_{F,f}(\bM^n_t))_{t\geq 0}$ is given by
\begin{equation}\label{n-finite variation}
\mathcal{A}^n_t = \int_0^t {\cL}^n \Psi_{F,\varphi_f}\big(M^n_s\big)\rmd s,
\end{equation}
with ${\cL}^n\Psi_{F,\varphi_f}$ as in (\ref{generator wn}), and its predictable quadratic variation is given by
\begin{align}
\mathcal{Q}^n_t&= n^{1+\frac{d}{3}}\int_0^t \int_{\R^d} \Big\{\bw^n_s(x)(1+s_n\bw^n_s(x))\big[F(\la \Theta^+_{x,n^{-1/3}R,u_n}(w^n_s),\varphi_f\ra)-F(\la w^n_s,\varphi_f\ra)\big]^2 \nonumber \\
& \quad +\big(1-\bw^n_s(x)+s_n(1-\bw^n_s(x)^2)\big)\big[F(\la \Theta^-_{x,n^{-1/3}R,u_n}(w_s^n),\varphi_f\ra)-F(\la w^n_s,\varphi_f\ra)\big]^2\Big\}\rmd x\rmd s.\label{n-quad variation}
\end{align}

Note that
\begin{align}
\la \Theta^+_{x,n^{-1/3}R,u_n}(w),\varphi_f\ra - \la w,\varphi_f\ra & = u_n\, \la \ind_{B_n(x)}(1-w),\varphi_f \ra \nonumber\\
\la \Theta^-_{x,n^{-1/3}R,u_n}(w),\varphi_f\ra - \la w,\varphi_f\ra & = -u_n\, \la \ind_{B_n(x)}w,\varphi_f \ra, \label{small increments}
\end{align}
so that both increments are of the order of $u_n n^{-d/3}$. Moreover, $f$ has compact support $S_f$ in $\R^d$ and thus so has $\varphi_f$. This will enable us to control the integrals over space of these increments.

Using this observation, we first show that
$|\mathcal{A}^n_t|$ is bounded by a constant independent of $n$. To this end, we
write it as the sum of a neutral term and a selective term and perform a Taylor expansion of
$F$ (truncating at second order in the neutral term and at first order in the selective term).
This yields, for any $t\geq 0$,
\begin{equation}\label{decomp cal A}
\mathcal{A}^n_t = \int_0^t \big(A_n(s)+B_n(s)+C_n(s)+D_n(s)+E_n(s)\big)\, \rmd s,
\end{equation}
where
\begin{align}
A_n(s) &= u_nn^{1+\frac{d}{3}}F'(\la \bw^n_s,f\ra)\int_{\R^d}  \Big[\bw^n_s(x)\la \ind_{B_n(x)}(1-w^n_s),\varphi_f\ra \nonumber\\
& \qquad \qquad \qquad \qquad\qquad \qquad \qquad \qquad\qquad \qquad - (1-\bw^n_s(x))\la \ind_{B_n(x)}w^n_s,\varphi_f\ra\Big] \rmd x, \nonumber\\
B_n(s) &= u_n^2 n^{1+\frac{d}{3}}\frac{F''(\la \bw^n_s,f\ra)}{2}\int_{\R^d} \Big[\bw^n_s(x)\la \ind_{B_n(x)}(1-w^n_s),\varphi_f\ra^2 \nonumber\\
& \qquad \qquad \qquad \qquad\qquad \qquad \qquad \qquad\qquad \qquad + (1-\bw^n_s(x)) \la \ind_{B_n(x)}w^n_s,\varphi_f\ra^2 \Big] \rmd x ,\nonumber \\
C_n(s) & \leq \mathcal{C}n^{1+\frac{d}{3}} \int_{\R^d}  \big(u_n \mathrm{Vol}(B_n(x))\big)^3 \ind_{\{B_n(x)\cap S_f\neq \emptyset\}} \rmd x , \nonumber\\
D_n(s) & = u_ns_nn^{1+\frac{d}{3}} F'(\la \bw^n_s,f\ra)\int_{\R^d} \Big[\bw^n_s(x)^2\la \ind_{B_n(x)}(1-w^n_s),\varphi_f\ra\nonumber \\
& \qquad \qquad \qquad \qquad\qquad \qquad \qquad \qquad \qquad- (1-\bw^n_s(x)^2)\la \ind_{B_n(x)}w^n_s,\varphi_f\ra\Big]\rmd x ,\nonumber\\
E_n(s) & \leq \mathcal{C'} n^{1+\frac{d}{3}}s_n u_n^2 \int_{\R^d}  \mathrm{Vol}(B_n(x))^2 \ind_{\{B_n(x)\cap S_f\neq \emptyset\}}\rmd x,
\end{align}
for some constant ${\cal C},{\cal C}'$ independent of $n$ and $s$. To control these expressions, we take a Taylor expansion of $\varphi_f$. We illustrate with
the term $A_n(s)$.
In fact, in identifying the limiting process we shall need a precise expression
for the limit of $A_n(s)$ and so we perform the expansion slightly more carefully than would
be required to simply conclude boundedness.

Let us write $D\varphi_f$ for the vector of first derivatives of $\varphi_f$ and $H\varphi_f$ for
the corresponding Hessian ($H\varphi_f=DD\varphi_f$). Recall that $S_f$
denotes the compact support of $f$. Then
\begin{align}
A_n(s)& = u_nn^{1+\frac{d}{3}}F'(\la \bw^n_s,f\ra)\int_{\R^d}  \Big[\bw^n_s(x)\la \ind_{B_n(x)},\varphi_f\ra - \la \ind_{B_n(x)}w^n_s,\varphi_f\ra\Big]\rmd x \label{approx An}\\
& = u_nn^{1+\frac{d}{3}}F'(\la \bw^n_s,f\ra) \int_{\R^d} \frac{n^{d/3}}{V_R}\int \int \ind_{\{|y-x|\leq n^{-1/3}R\}}\ind_{\{|z-x|\leq n^{-1/3}R\}} w^n_s(y)\nonumber\\
& \qquad \qquad \qquad \qquad\qquad \qquad\qquad \qquad\qquad \qquad \qquad \qquad\times (\varphi_f(z)-\varphi_f(y)) \rmd z \rmd y\rmd x \nonumber \\
& = u_nn^{1+\frac{d}{3}}F'(\la \bw^n_s,f\ra) \int_{\R^d} \frac{n^{d/3}}{V_R}\int_{\R^d} \int_{\R^d} \ind_{\{|y-x|\leq n^{-1/3}R\}}\ind_{\{|z-x|\leq n^{-1/3}R\}} w^n_s(y) \nonumber\\
& \qquad \times \big[D\varphi_f(y)(z-y) +
\frac{1}{2}(z-y)H\varphi_f(y)(z-y)+\mathcal{O}(|z-y|^3)\ind_{\{y\in S_f\}}\big]\rmd z  \rmd y \rmd x. \nonumber
\end{align}
Consider the first term on the right. Integrating first with respect to $x$ (using Fubini's Theorem)
this term is
\begin{equation}
\frac{u_nn^{1+\frac{2d}{3}}}{V_R} F'(\la \bw^n_s,f\ra) \int_{\R^d} w^n_s(y) \int_{\R^d}
\mathrm{Vol}(B_n(y)\cap B_n(z))D\varphi_f(y)(z-y) \rmd z \rmd y,
\end{equation}
and since $\mathrm{Vol}(B_n(y)\cap B_n(z))$ is a function of $|z-y|$ alone, the integrand is
antisymmetric as a function of $z-y$ and so the integral with respect to $z$ vanishes.

Similarly, the integrals corresponding to the off-diagonal terms in the Hessian will vanish, leaving
\begin{align}
u_nn^{1+\frac{d}{3}}F'(\la \bw^n_s,f\ra)& \int_{\R^d}
\frac{n^{d/3}}{V_R}\int_{\R^d} \int_{\R^d}
\ind_{\{|y-x|\leq n^{-1/3}R\}}\ind_{\{|z-x|\leq n^{-1/3}R\}} w^n_s(y)\nonumber\\
& \times \frac{1}{2}\sum_{i=1}^d(z_i-y_i)^2\frac{\partial^2}{\partial y_i^2}\varphi_f(y)\rmd y\rmd z\rmd x \label{approx An2}
\end{align}
plus a lower order term. Now observe that since $f\in C_c^{\infty}(\R^d)$, another Taylor expansion argument enables us to write that
\begin{equation}
\frac{\partial^2}{\partial y_i^2}\varphi_f(y) =\varphi_{\frac{\partial^2 f}{\partial y_i^2}}(y)
=\frac{\partial^2f}{\partial y_i^2}(y)+\cO(n^{-2/3})\ind_{\{B_n(y)\cap S_f\neq \emptyset\}}
\end{equation}
(where the term $\cO(n^{-2/3})$ is independent of $y$). This yields
\begin{align}
A_n(s) &=
u_nn^{1+\frac{d}{3}}F'(\la \bw^n_s,f\ra) \int_{\R^d}
\frac{n^{d/3}}{V_R}\int_{\R^d} \int_{\R^d}
\ind_{\{|y-x|\leq n^{-1/3}R\}}\ind_{\{|z-x|\leq n^{-1/3}R\}} w^n_s(y)\nonumber
\\ &\qquad\qquad\qquad\qquad\qquad\qquad
\times\frac{1}{2}\sum_{i=1}^d(z_i-y_i)^2\frac{\partial^2f}{\partial y_i^2}(y)\rmd y\rmd z\rmd x
\nonumber
\\ +& \cO(n^{-2/3})n^{\frac{2}{3}(1+d)}
\int_{\R^d}\int_{\R^d}\int_{\R^d}
\ind_{\{|y-x|\leq n^{-1/3}R\}}\ind_{\{|z-x|\leq n^{-1/3}R\}}
|z-y|^2\ind_{\{B_n(y)\cap S_f\neq \emptyset\}}\nonumber\\
&=
\frac{un^{\frac{2}{3}(1+d)}}{2V_R}F'(\la \bw^n_s,f\ra) \int_{\R^d}
\int_{B_n(x)^2}\, w^n_s(y) \sum_{i=1}^d (z_i-y_i)^2 \frac{\partial^2 f}{\partial y_i^2}(y)\rmd y\rmd z \rmd x + \mathcal{O}(n^{-2/3})\nonumber\\
& = \frac{u\Gamma_R}{2} F'(\la \bw^n_s,f\ra) \int_{\R^d}
\, w^n_s(y) \Delta f(y) \rmd y + \mathcal{O}(n^{-2/3})\nonumber\\
& = \frac{u\Gamma_R}{2} F'(\la \bw^n_s,f\ra)
\, \la \bw^n_s,\Delta f\ra + \mathcal{O}(n^{-2/3}),
\label{careful form of AN}
\end{align}
where
\begin{equation}
\Gamma_R= \frac{n^{\frac{2}{3}(1+d)}}{V_R} \int_{B_n(y)} \int_{B_n(x)} (z_1-y_1)^2 \rmd z\rmd x = \frac{1}{V_R}\int_{B(0,R)}\int_{B(x,R)}(z_1)^2 \rmd z \rmd x
\end{equation}
was defined in (\ref{def Gamma}), and the last equality uses another Taylor expansion to show that for any $s$,
\begin{equation}\label{comparison w}
\la w^n_s,\Delta f\ra = \la \bw^n_s,\Delta f\ra + \cO(n^{-2/3})
\end{equation}
with an error term uniformly bounded in $s$. In particular, since $|\la \bw^n_s,f\ra|\leq \|f\| \mathrm{Vol}(S_f)$, we can conclude that $|A_n(s)|\leq \mathcal{C}_A$ uniformly in $s$ and $n$.

Very similar arguments allow us to control the other terms:
\begin{equation}\label{control Bn}
|B_n(s)| \leq  \frac{u_n^2n^{1+\frac{d}{3}}}{2}\, |F''(\la \bw^n_s,f \ra)|\int  \, 2\mathrm{Vol}(B_n(x))^2 \ind_{\{x\in S_f\}} \|f\|^2 \rmd x \leq  \mathcal{C}_B\, n^{\frac{1-d}{3}},
\end{equation}
and, again by the same arguments,
\begin{equation}\label{control Cn}
|C_n(s)| \leq \mathcal{C}_C\, n^{-\frac{2d}{3}}, \quad |D_n(s)| \leq \mathcal{C}_D \quad \hbox{and}\quad |E_n(s)| \leq \mathcal{C}_E\, n^{-\frac{1+d}{3}},
\end{equation}
where the constants $\mathcal{C}_B,\mathcal{C}_C,\mathcal{C}_D,\mathcal{C}_E$ are all independent of $n$ and $s$.
Coming back to \eqref{decomp cal A} and combining all the estimates we just obtained, for every $s<t$ we have
\begin{equation}\label{FV An}
|\mathcal{A}^n_t-\mathcal{A}^n_s| \leq \big(\mathcal{C}_A + \mathcal{C}_B\, n^{\frac{1-d}{3}} + \mathcal{C}_C\, n^{-\frac{2d}{3}} + \mathcal{C}_D + \mathcal{C}_E\, n^{-\frac{1+d}{3}} \big) (t-s).
\end{equation}
From there it is easy to deduce that for every $T>0$, given a sequence of stopping times $(\tau_n)_{n\geq 1}$ bounded by $T$, for every $\varepsilon>0$ there exists $\eta>0$ such that
\begin{equation}\label{tightness An}
\limsup_{n\rightarrow \infty} \sup_{\theta\in [0,\eta]} \P\big[\big|\mathcal{A}^n_{\tau_n+\theta}-\mathcal{A}^n_{\tau_n}\big|>\e \big] =0,
\end{equation}
which corresponds to \eqref{tightness cal An} and shows that the sequence of finite variation parts of $(\Psi_{F,f}(\bM^n_t))_{t\geq 0}$ is tight.

Similarly, we obtain that
\begin{equation}
\big[F(\la \Theta^{\pm}_{x,n^{-1/3}R,u_n}(w^n_s),\varphi_f\ra)-F(\la w^n_s,\varphi_f\ra)\big]^2 \leq \mathcal{C}_F''\|f\|^2 u_n^2 \hbox{Vol}(B_n(x))^2 \ind_{\{B_n(x)\cap S_f\neq \emptyset\}}.
\end{equation}
Notice that this bound is independent of the value of $w_s^n$. Substituting into the definition of $\mathcal{Q}^n_t$ given in (\ref{n-quad variation}), we obtain that for every $s<t$,
\begin{equation}\label{estimates qv}
|\mathcal{Q}^n_t-\mathcal{Q}^n_s|\leq \mathcal{C}_F\, n^{\frac{1-d}{3}}(t-s),
\end{equation}
for a constant $\mathcal{C}_F$ independent of $n$ (and $s,t$). Therefore, for every $T>0$, every sequence of stopping times $(\tau_n)_{n\geq 1}$ bounded by $T$ and every $\e>0$, there exists $\eta>0$ such that
\begin{equation}\label{tightness Qn}
\limsup_{n\rightarrow \infty} \sup_{\theta\in [0,\eta]} \P\big[\big|\mathcal{Q}^n_{\tau_n+\theta}-\mathcal{Q}^n_{\tau_n}\big|>\e \big] =0
\end{equation}
and the sequence of predictable quadratic variations of the martingale part of the process $(\Psi_{F,f}(\bM^n_t))_{t\geq 0}$ is not only tight, but also when $d\geq 2$ it tends to $0$ uniformly over compact time intervals. By the Aldous-Rebolledo criterion (see again $(i)$-item~$(c)$ in Section~\ref{s: proof existence}), we conclude that $(\bM^n)_{n\geq 1}$ is tight in $D_{\cM_\lambda}[0,\infty)$, as required.

Finally, coming back to \eqref{small increments}, we see that every increment of $\la \bM^n,f\ra = \la M^n,\varphi_f\ra$ is bounded by $\|f\|V_R u_n n^{-d/3}$. Consequently, for every $T>0$ we have
\begin{equation}
\sup_{t\in [0,T]}\sup_{f\in C_c^\infty(\R^d):\|f\|\leq 1} \big|\big\la \bM^n_t,f\big\ra - \big\la \bM^n_{t-},f\big\ra\big| \leq V_Run^{-(1+d)/3},
\end{equation}
and thus any potential limit for $(\bM^n)_{n\geq 1}$ has continuous paths in $\cM_\lambda$.

\medskip
\textbf{2) Identifying the limit.}

\medskip
\noindent In what follows, we suppose that $(M_t^\infty)_{t\geq 0}\in D_{\cM_\lambda}[0,\infty)$ is the weak limit of a subsequence $(\bM^{n_k})_{k\geq 1}$ and for any $t\geq 0$, we write $w_t^\infty$ for (some representative of) the density of $M_t^\infty$.

In order to show that $M^\infty$ satisfies the martingale problem stated in Theorem~\ref{th:fixed}, we use the fact (established in the previous paragraph) that for every $f\in C_c^\infty(\R^d)$ and every $n\geq 1$,
\begin{equation}\label{approximate MP}
\bigg( \Psi_{\mathrm{Id},f}\big(\bM^n_t\big)- \Psi_{\mathrm{Id},f}\big(\bM^n_0\big) - \int_0^t \cL^n \Psi_{\mathrm{Id},\varphi_f}(M^n_s)\rmd s\bigg)_{t\geq 0}
\end{equation}
is a martingale with predictable quadratic variation \eqref{n-quad variation} (with $F=\mathrm{Id}$), where $\cL^n$ was defined in \eqref{generator wn} and $\varphi_f$ in \eqref{def varphif}. We first show that for every $t\geq 0$,
\begin{equation}
\lim_{k\rightarrow \infty}\E\bigg[\bigg| \int_0^t \cL^{n_k} \Psi_{\mathrm{Id},\varphi_f}(M^{n_k}_s)\rmd s - \int_0^t \bigg\{\frac{u\Gamma_R}{2}\, \la w_s^\infty ,\Delta f\ra - 2Ru\sigma\,
\la w_s^\infty (1-w_s^\infty),f\ra\bigg\}\, \rmd s\bigg| \bigg] = 0, \label{L1 limit}
\end{equation}
so that we can then use the fact that the quantity in \eqref{approximate MP} is a martingale, the fact that $\Psi_{\mathrm{Id},f}$ is a bounded continuous function and the Dominated Convergence Theorem to conclude that for every $0\leq t<t'$, $m\in \N$, $0\leq t_1 < \cdots < t_m\leq t$ and $h_1,\ldots,h_m\in C_b(\cM_\lambda)$,
\begin{align}
\E\bigg[\bigg(\la w^\infty_{t'},f\ra -\la w^\infty_t,f\ra -\int_t^{t'}\bigg\{\frac{u\Gamma_R}{2}\, \la w_s^\infty ,\Delta f\ra -& 2Ru\sigma\,
\la w_s^\infty (1-w_s^\infty),f\ra\bigg\}\,\rmd s\bigg)\nonumber\\
& \qquad \qquad \times \bigg(\prod_{i=1}^m h_i\big(M^\infty_{t_i}\big)\bigg)\bigg]=0 \label{MP limit}
\end{align}
and consequently that ${\cal Z}^f$ is a martingale (with respect to the natural filtration of $M^\infty$). In the case $d\geq 2$ this property will be sufficient to conclude, since we showed in \eqref{estimates qv} that the quadratic variation of the martingale \eqref{approximate MP} tended to $0$ as $n\rightarrow \infty$, and therefore the limit ${\cal Z}^f$ is the constant process equal to $0$. In one dimension, we shall still have to prove that the quadratic variation of ${\cal Z}^f$ is non-trivial and has the announced form. This is what we do in the last part of this point $\mathbf{2)}$.

Let us prove \eqref{L1 limit}. Specialising the computation of $\mathcal{A}^n$ in \eqref{generator wn} to the case $F=\mathrm{Id}$, we have by \eqref{careful form of AN}
\begin{align}
A_{n_k}(s)
& = \frac{u\Gamma_R}{2} \, \la \bw^{n_k}_s,\Delta f\ra + \mathcal{O}(n_k^{-2/3})\nonumber \\
& = \frac{u\Gamma_R}{2}\int_{\R^d\times\{0,1\}}\Delta f(x)\ind_{\{0\}}(\kappa)\bM^{n_k}_s(\rmd x,\rmd \kappa)+ \mathcal{O}(n_k^{-2/3})  \nonumber\\
& \rightarrow \frac{u\Gamma_R}{2} \, \la w_s^\infty,\Delta f\ra \qquad \hbox{as }k\rightarrow \infty. \label{first term}
\end{align}
These quantities being bounded by $(u\Gamma_R/2)\|\Delta f\| \mathrm{Vol}(S_f)+\mathcal{O}(n_k^{-2/3})$, independently of $s$, the convergence also happens in $\mathbb{L}^1$ norm. Next, Taylor-expanding $f$ to write that for every $y\in B(x,Rn_k^{-1/3})$,
\begin{equation}
\varphi_f(y)=f(x)+\mathcal{O}(|y-x|)= f(x)+ \mathcal{O}\big(n_k^{-1/3}\big),
\end{equation}
we obtain that
\begin{align}
D_{n_k}(s) & = \sigma u n_k^{d/3} \int_{\R^d}  \Big\{\bw_s^{n_k}(x)^2 \big\la \ind_{B_{n_k}(x)}(1-w^{n_k}_s),f(x) + \mathcal{O}\big(n_k^{-1/3}\big)\big\ra \nonumber \\
&\qquad \qquad \qquad \qquad \qquad - (1-\bw_s^{n_k}(x)^2)\big\la \ind_{B_{n_k}(x)}w^{n_k}_s,f(x) + \mathcal{O}\big(n_k^{-1/3}\big)\big\ra \Big\}\rmd x\nonumber \\
& = \sigma uV_R \int_{\R^d}  \big\{\bw_s^{n_k}(x)^2(1-\bw_s^{n_k}(x)) - (1-\bw_s^{n_k}(x)^2)\bw^{n_k}_s(x)\big\}f(x) \rmd x + \mathcal{O}\big(n_k^{-1/3}\big) \nonumber\\
& = - \sigma uV_R\, \la \bw_s^{n_k}(1- \bw_s^{n_k}), f\ra + \mathcal{O}\big(n_k^{-1/3}\big).
\label{sigma in Dn}
\end{align}
As above, the part of $D_{n_k}(s)$ which is linear in $\bw^{n_k}_s$ converges (weakly and in $\mathbb{L}^1$) towards
\begin{equation}\label{second term}
-\sigma u V_R\, \la w_s^\infty,f\ra.
\end{equation}
We now would like to show that the ``quadratic'' part of $D_{n_k}(s)$ converges to
\begin{equation}
\sigma uV_R\, \la (w_s^\infty)^2,f\ra.
\end{equation}
Note that this is not a simple consequence of the weak convergence of $\bM^{n_k}$ to $M^\infty$, as $\la (\bw^{n_k}_s)^2,f\ra$ cannot be written as an integral with respect to $\bM^{n_k}_s$ or $(\bM^{n_k}_s)^{\otimes 2}$. Instead, we shall approximate this expression by an integral with respect to $(\bM^{n_k}_s)^{\otimes 2}$ and use the continuity estimates obtained in Proposition~\ref{prop: spatial continuity A} to bound the remaining terms. (The statement and proof of this proposition are postponed until Appendix~\ref{continuity fixed} to ease the reading).

Let $\e\in (0,1/2)$, and let $p_\e$ be a continuous probability density function on $\R^d$ supported in $B(0,\e)$. For every $k\geq 1$ and $s\geq 0$, we have
\begin{align}
& \big|\la (\bw^{n_k}_s)^2,f \ra - \la (w_s^\infty)^2,f\ra \big| \nonumber \\
&\leq \bigg|\int_{\R^d}f(x)\bw_s^{n_k}(x)^2\rmd x - \int_{\R^d}\int_{\R^d}f(x)\bw^{n_k}_s(x)\bw^{n_k}_s(y)p_\e(y-x)\rmd y\rmd x \bigg| \nonumber \\
&\quad  + \bigg|\int_{\R^d}\int_{\R^d}f(x)\bw^{n_k}_s(x)\bw^{n_k}_s(y)p_\e(y-x)\rmd y\rmd x - \int_{\R^d}\int_{\R^d}f(x)w_s^\infty(x)w_s^\infty(y)p_\e(y-x)\rmd y\rmd x \bigg| \nonumber \\
& \quad + \bigg|\int_{\R^d}\int_{\R^d}f(x)w_s^\infty(x)w_s^\infty(y)p_\e(y-x)\rmd y\rmd x - \int_{\R^d}f(x)w_s^\infty(x)^2\rmd x\bigg|.\label{conv quadratic part}
\end{align}
The second term on the r.h.s. can be rewritten as
\begin{align}
&\int_{(\R^d\times\{0,1\})^2}f(x)p_\e(y-x)\ind_{\{0\}}(\kappa)\ind_{\{0\}}(\kappa')\bM_s^{n_k}(\rmd y,\rmd \kappa')\bM_s^{n_k}(\rmd x,\rmd \kappa) \nonumber\\
& \rightarrow \int_{(\R^d\times\{0,1\})^2}f(x)p_\e(y-x)\ind_{\{0\}}(\kappa)\ind_{\{0\}}(\kappa')M_s^\infty(\rmd y,\rmd \kappa')M_s^\infty(\rmd x,\rmd \kappa)\nonumber\\
& = \int_{\R^d}\int_{\R^d}f(x)w_s^\infty(x)w_s^\infty(y)p_\e(y-x)\rmd y\rmd x
\end{align}
as $k$ tends to infinity (since the mapping $(x,y)\mapsto f(x)p_\e(y-x)$ belongs to $C_c((\R^d)^2)$, and since these terms are bounded uniformly in $k$ (and $\e,s$), this convergence also happens in $\mathbb{L}^1$ norm. That is, the expectation of the second term in (\ref{conv quadratic part}) tends to $0$ as $k\rightarrow \infty$.

Concerning the first term on the r.h.s. of (\ref{conv quadratic part}), because $\bw^{n_k}_s$ takes its values in $[0,1]$, we have, by Fubini's Theorem,
\begin{align}
&\E\bigg[\bigg|\int_{\R^d}f(x)\bw_s^{n_k}(x)^2\rmd x - \int_{\R^d}\int_{\R^d}f(x)\bw^{n_k}_s(x)\bw^{n_k}_s(y)p_\e(y-x)\rmd y\rmd x\bigg|\bigg]\nonumber\\
& \qquad \leq \|f\| \int_{S_f}\int_{B(x,\e)}\E\big[\big|\bw_s^{n_k}(x)-\bw_s^{n_k}(y)\big|\big]p_\e(y-x)\rmd y\rmd x.
\end{align}
By Proposition~\ref{prop: spatial continuity A} applied with $\e=Rn_k^{-1/3}$, there exists $a,v,\lambda,C>0$ independent of $k$ such that for every $x,y\in \R^d$ satisfying $|x-y|<1$ and every $s\in [0,t]$, we have
\begin{align}
\E\big[\big|\bw_s^{n_k}(x)-\bw_s^{n_k}(y)\big|\big]& \leq C\Big\{n_k^{-a} + \tau_{n_k}(x,y) + \Big(|x-y|^{1/4} + \tau_{n_k}(x,y)^{1/2}\Big)e^{\lambda(|x|+Rn_k^{-1/3})}\nonumber\\
& \qquad \qquad +n_k^{(1-d)/6}\tau_{n_k}(x,y)^{(2-d)/4}\Big\},
\end{align}
where
\begin{equation}
\tau_n(x,y) = n^{-v}\vee |x-y|^{2/(d+1)}.
\end{equation}
Thus, using the facts that the support $S_f$ of $f$ is compact, that $p_\e$ is a probability density supported in $B(0,\e)$, and that $\tau_n(x,y)\leq \e^{2/(d+1)}$ for $n$ large enough whenever $|x-y|\leq \e$, we can write that the first term on the r.h.s. of (\ref{conv quadratic part}) is bounded by
\begin{equation}\label{continuity bound}
C'\Big(n_k^{-a}+ \e^{2/(d+1)}+\e^{1/4}+\e^{1/(d+1)}+ n_k^{(1-d)/6}\e^{1/(d+1)}\Big).
\end{equation}
Likewise, by taking $n\rightarrow \infty$ in Proposition~\ref{prop: spatial continuity A} (along the converging subsequence), we obtain that the last term on the r.h.s. of (\ref{conv quadratic part}) is bounded by
\begin{equation}\label{continuity bound 2}
C'\Big(\e^{1/4} + \e^{2/(d+1)} +\e^{1/(d+1)}+ \e^{1/(d+1)}\mathbf{1}_{\{d=1\}}\Big).
\end{equation}
Combining the above, we have that for every $\e\in (0,1/2)$,
\begin{equation}
\limsup_{k\rightarrow \infty}\E\big[\big| \la (\bw^{n_k}_s)^2,f \ra - \la (w_s^\infty)^2,f\ra \big|\big]\leq C\big(\e^{1/4} + \e^{1/(d+1)}\big),
\end{equation}
and letting $\e$ tend to $0$ we can conclude that the part of the expression~(\ref{sigma in Dn})
for $D_{n_k}(s)$ which is quadratic in $\bw^{n_k}_s$ indeed converges in $\mathbb{L}^1$ towards
\begin{equation}\label{third term}
\sigma uV_R\, \la (w_s^\infty)^2,f\ra.
\end{equation}

Combining (\ref{first term}), (\ref{second term}) and (\ref{third term}), and using the facts
that $B_n(s)=0$ since $F=\mathrm{Id}$, and that $C_n(s)$ and $E_n(s)$ tend to zero uniformly in all possible values of $\bM^n$, we conclude that \eqref{L1 limit} is satisfied. As we explained above, this is sufficient to conclude in the case $d\geq 2$ since the quadratic variation of the martingale ${\cal Z}^f$ is then $0$.

We now turn to the case $d=1$. Defining
\begin{align}
W^n_t(f)& := \la \bw^n_t,f\ra - \la \bw^n_0,f\ra -\int_0^t \cL^n \Psi_{\mathrm{Id},\varphi_f}\big(M^n_s\big) \rmd s \\
& = \la \bw^n_t,f\ra - \la \bw^n_0,f\ra - \int_0^t \bigg\{\frac{u\Gamma_R}{2} \, \la \bw^n_s,\Delta f\ra - \sigma uV_R\, \la \bw_s^n(1- \bw_s^n), f\ra \bigg\} \, \rmd s+ \mathcal{O}(n^{-1/3}),\nonumber
\end{align}
we know that $W^n(f)$ is a zero-mean martingale with predictable quadratic variation
\begin{align}
u_n^2n^{4/3}&\int_0^t \int_{\R^d} \Big\{\bw^n_s(x)(1+s_n\bw^n_s(x))\la \ind_{B_n(x)}(1-w^n_s),f+\mathcal{O}(n^{-1/3})\ra^2 \nonumber\\
& \qquad \qquad \qquad  + (1-\bw^n_s(x)+s_n(1-\bw^n_s(x)^2))\la \ind_{B_n(x)}w^n_s,f+\mathcal{O}(n^{-1/3})\ra^2\Big\}\,  \rmd x\rmd s \nonumber\\
& = u^2V_R^2\int_0^t \la \bw^n_s(1-\bw^n_s),f^2\ra \rmd s  + \mathcal{O}(n^{-1/3}),
\end{align}
where, more precisely, the remainder term is bounded by a constant times $n^{-1/3}t$. As a consequence, for every $n\geq 1$, $0\leq t<t'$, $m\in \N$, $0\leq t_1<\cdots<t_m\leq t$ and $h_1,\ldots,h_m\in C_b(\cM_\lambda)$,
\begin{equation}\label{id qv}
\E\bigg[ \bigg(\big(W^n_{t'}(f)\big)^2- \big(W^n_t(f)\big)^2 - u^2V_R^2\int_t^{t'}\la \bw^n_s(1-\bw^n_s),f^2\ra \rmd s  + \mathcal{O}(n^{-1/3})\bigg)\bigg(\prod_{i=1}^m h_i\big(\bM^n_{t_i}\big)\bigg)\bigg]=0
\end{equation}
Observe that for every $n\geq 1$ and every $t\geq 0$,
\begin{equation}
|W^n_t(f)|\leq \mathrm{Vol}(S_f)\left[2\|f\| +t\left(\frac{u\Gamma_R}{2}\|\Delta f\|+\sigma uV_R\|f\|+ \mathcal{O}(n^{-1/3}) \right)\right],
\end{equation}
and so we can let $n\rightarrow \infty$ in \eqref{id qv} (along the converging subsequence) and use the Dominated Convergence Theorem, together with (\ref{first term}), (\ref{second term}) and (\ref{third term}), to conclude that
\begin{equation}\label{final qv}
\E\bigg[ \bigg(\big({\cal Z}^f_{t'}\big)^2- \big({\cal Z}^f_t\big)^2 - u^2V_R^2\int_t^{t'}\la w^\infty_s(1-w^\infty_s),f^2\ra \rmd s\bigg)\bigg(\prod_{i=1}^m h_i\big(M^\infty_{t_i}\big)\bigg)\bigg]=0.
\end{equation}
This allows us to identify the quadratic variation of the martingale ${\cal Z}^f$ as
\begin{equation}
[{\cal Z}^f]_t= u^2V_R^2\int_0^t\la w^\infty_s(1-w^\infty_s),f^2\ra \rmd s, \qquad t\geq 0.
\end{equation}
Since by \eqref{small increments} the jumps of $W^n(f)$ are all bounded by $Cn^{-2/3}$, ${\cal Z}^f$ is a continuous square-integrable martingale, starting at $0$. By the Dubins-Schwarz Theorem (see Remark~\ref{rk:DS} below), ${\cal Z}^f$ is therefore a time-changed Brownian motion, solution to the stochastic differential equation
\begin{equation}
dW_t = uV_R \sqrt{\la w_t^\infty(1-w_t^\infty), f^2\ra}\, dB^f_t,
\end{equation}
where $B^f$ denotes standard Brownian motion. The bracket process between ${\cal Z}^f$ and ${\cal Z}^g$ is then obtained by the same kind of calculations, writing first the bracket process for a fixed $n$ and then identifying the limit by letting $n_k\rightarrow \infty$.
\begin{remark}\label{rk:DS}
We cannot \emph{a priori} prove that $[{\cal Z}^f]_\infty=+\infty$ a.s., as required by the classical Dubins-Schwarz Theorem. Note however that this condition can be removed, at the expense of extending the probability space on which we work. Indeed, if we introduce a Brownian motion $(\beta^f_t)_{t\geq 0}$ independent of all other processes (possibly on some enlarged space) and set for every $t\geq 0$
\begin{equation}
B_t^f =\left\{\begin{array}{ll} {\cal Z}^f_{\tau_t} & \hbox{if }t< [{\cal Z}^f]_\infty,\vspace{0.1cm}\\
{\cal Z}^f_{\infty} + \beta^f_{t-[{\cal Z}^f]_\infty} & \hbox{if }t\geq  [{\cal Z}^f]_\infty,
\end{array}\right.
\end{equation}
where
\begin{equation}
\tau_t:= \inf\big\{s\geq 0:\, [{\cal Z}^f]_s>t\big\},
\end{equation}
then by Theorem~1.7 in Chapter~V of \cite{RY1999} we have that $(B_t^f)_{t\geq 0}$ is a standard Brownian motion and for every $t\geq 0$, ${\cal Z}^f_t=B^f_{[{\cal Z}^f]_t}$.
\end{remark}

To summarise, we have shown that any limit point $(M_t^\infty)_{t\geq 0}$ of
$(\bM^n)_{n\geq 1}$ satisfies the following system of
stochastic differential equations: for every $f\in C_c^\infty(\R^d)$,
\begin{align}
\rmd\la w_t^\infty,f\ra =& \bigg\{\frac{u\Gamma_R}{2} \la w_t^\infty,\Delta f\ra - \sigma
uV_R\, \la w_t^\infty(1- w_t^\infty), f\ra \bigg\} \rmd t \nonumber \\
& +  \ind_{\{d=1\}}uV_R \sqrt{\la w_t^\infty(1-w_t^\infty),
f^2\ra} \rmd B^f_t,\label{system sde}
\end{align}
with initial value $\la w_0,f\ra$, and in one dimension, by polarisation the covariation between $\la w_\cdot^\infty,f\ra$ and $\la w_\cdot^\infty,g\ra$ is as in the statement of Theorem~\ref{th:fixed}$(i)$.

\medskip
\textbf{3) Uniqueness of the limit.}

\medskip

Let us finally show that the system of equations (\ref{system sde}) has at most one solution. We start with the case $d\geq 2$. Any test function of the form
\begin{align}
M\mapsto & \int_{(\R^d\times \{0,1\})^k}\psi(x_1,\ldots,x_k)\bigg\{\prod_{j=1}^k\ind_{\{0\}}(\kappa_j)\bigg\}\, M(\rmd x_1,\rmd \kappa_1)\cdots M(\rmd x_k,\rmd \kappa_k) \nonumber\\
& = \int_{(\R^d)^k}\psi(x_1,\ldots,x_k)\bigg\{\prod_{i=1}^k w(x_i)\bigg\} \rmd x_1 \cdots \rmd x_k,\label{test function again}
\end{align}
with $\psi$ continuous and integrable on $(\R^d)^k$ (where as before $w$ is any representative of the density of $M$),  can be uniformly approximated by linear combinations of functions of the form $\prod_{i=1}^k \la \cdot,f_i\ra$
with $f_i\in C_c^\infty(\R^d)$ for every $i$. Thus, we can extend (\ref{system sde}) to this more general class of functions. Then in Chapter~7 of \cite{LIA2009}, it is proved that, when $\sigma=0$,
any solution to (\ref{system sde}) is dual, through the set of functional relations (\ref{dual formula}), to a system of independent Brownian motions with variance parameter $u\Gamma_R$,
in which particles never coalesce. This is easily modified to $\sigma>0$, in which case particles branch into two at rate $u\sigma V_R$, independently of each other. Since the set of all test functions of the form \eqref{test function again} is separating by Lemma~2.1(c) in \cite{VW2012}, we can proceed as in the proof of Proposition 4.4.7 in \cite{EK1986} to conclude that the system of equations~(\ref{system sde}) has at most one solution.
Hence, this solution exists and the full sequence $(\bM^n)_{n\geq 0}$ converges to it in distribution, as stated in Theorem~\ref{th:fixed}$(ii)$.

When $d=1$, we follow the same route and use It\^o's Formula to extend (\ref{system sde}) to functions of the product form $\prod_{i=1}^k\la \cdot,f_i\ra$ and then to the full class of functions \eqref{test function again} by the same density argument as before. Again in Chapter~7 of \cite{LIA2009}, it is proved that in one dimension and when $\sigma=0$, any solution to these equations is dual, through the set of relations (\ref{dual formula}), to a
system of independent Brownian motions with variance parameter $u\Gamma_R$,
in which, this time, particles coalesce pairwise at an instantaneous rate given by $u^2V_R^2$ times the local time at $0$ of their separation (independently of the other pairs). As earlier, this is easily modified to cover the case $\sigma>0$, by imposing that particles should
also branch into two at rate $u\sigma V_R$. By the same chain of arguments as in the case $d\geq 2$, we can therefore conclude that the system of equations (\ref{system sde}) has a unique solution, to which the full sequence $(\bM^n)_{n\geq 0}$ thus converges in distribution as $n$ tends to infinity. Theorem~\ref{th:fixed}$(i)$ is proved.

\begin{remark}
Liang's notation is very different from ours. To see that his process (with selection added and the coalescence rate multiplied by $u^2V_R^2$) and our limiting process do coincide, notice that $m(\rmd x)=\rmd x$
in our case and $\hat{X}_t(x)=w_t^\infty(x)\delta_0 + (1-w_t^\infty(x))\delta_1$. Hence, taking $\chi(\kappa)=\ind_{0}(\kappa)=\rho(\kappa)$ and $\psi(x)=f(x),\phi(x)=g(x)$ in Proposition~7.2 in \cite{LIA2009} indeed leads to
\begin{align}
\rmd [\mathcal{Z}^f ,\mathcal{Z}^g]_t &= u^2V_R^2\int_{\R^d} w_t^\infty(x)f(x)g(x)\rmd x - u^2V_R^2\int_{\R^d} \big(w_t^\infty(x)\big)^2f(x)g(x)\rmd x \nonumber\\
& = u^2V_R^2\int_{\R^d} w_t^\infty(x)\big(1-w_t^\infty(x)\big)f(x)g(x)\rmd x.
\end{align}
\end{remark}

\subsection{Proof of Theorem~\ref{th: conv duals}.}

We divide the proof into two parts.  The first, and simpler, shows that
the only possible limit for $(\Xi^n)_{n\geq 1}$
is the system of branching and coalescing
Brownian motions $\Xi^\infty$.
The second part, tightness of the sequence
$(\Xi_t^n)_{n\geq 1}$, is rather more
involved and will be broken into a number of smaller steps.

Recall that $\Xi^n$ is defined on the probability space $(\mathbf{\Omega},\mathcal{F}',\bfP)$ and takes its values in the set $\cM_p(\R^d)$ of all finite point measures on $\R^d$, which we have endowed with the topology of weak convergence. The linear hull of the set of test functions (recall \eqref{test functions for dual})
\begin{equation}\label{test functions dual}
\Phi_{\exp,\ln f}:\Xi \mapsto \prod_{i=1}^{|\Xi|} f(\xi^i) = \exp\bigg\{\int_{\R^d}(\ln f(x))\Xi(\rmd x)\bigg\},
\end{equation}
where $f\in C^1(\R^d)$ takes values in $[0,1]$, is dense in $C_0(\cM_p(\R^d))$ (the space of continuous functions on $\cM_p(\R^d)$ tending to $0$ at ``infinity'') for the topology of the uniform convergence (\emph{cf.} Lemma~0.2 in \cite{INW1968}, where the formalism is different but the result is equivalent to our claim). Consequently, the linear span of the set of functions \eqref{test functions dual} is dense in $C_b(\cM_p(\R^d))$ for the topology of uniform convergence over compact sets, and functions of the above form will thus be sufficient to characterise the law of an $\cM_p(\R^d)$-valued random variable. In this section, the atoms of the point measures considered will be viewed as particles evolving in $\R^d$.

We start with the following result.
\begin{lemma}\label{lem: conv dual}
The finite dimensional distributions of the system of scaled processes $\Xi^n$ converge as $n\rightarrow\infty$ to those of the system of branching and coalescing Brownian motions $\Xi^\infty$,
described in the statement of Theorem~\ref{th: conv duals}. In particular, the only possible
limit point for the sequence $(\Xi^n)_{n\geq 1}$ is $\Xi^\infty$.
\end{lemma}

\begin{proof}[Proof of Lemma~\ref{lem: conv dual}]
Suppose first that the density $\psi$ of the locations of the atoms of $\Xi^n_0$ can be factorised as $\psi(x_1,\ldots,x_k)=\psi_1(x_1)\cdots \psi_k(x_k)$, with $\psi_i\in C_c(\R^d)$ being a probability density function on $\R^d$ for every $i$.

Let us write $(M^{(n)}_t)_{t\geq 0}$ for the unscaled SLFVS with parameters $s_n$, $u_n$ (in the fixed radius case), and $w^{(n)}_t$ for a representative of the density of $M_t^{(n)}$, for every $t\geq 0$. Recall the notation $(M^n_t)_{t\geq 0}$ for the scaled process whose density at time $t$ is $w_{nt}^{(n)}(n^{1/3}\cdot)$. Let $w^0\in C^1(\R^d)$, and suppose that $M_0^{(n)}$ is such that $M_0^n$ has density $w^0$ for every $n\geq 1$. With this initial condition and Relation~\eqref{comparison w} (where $\Delta f$ can be replaced by any function $f\in C_c^2(\R^d))$, it is easy to check that $\bM^n_0$, as defined in Theorem~\ref{th:fixed}, converges to the measure $M^0\in \cM_\lambda$ with density $w^0$ as $n\rightarrow \infty$. Hence, by Theorem~\ref{th:fixed} the sequence of processes $(\bM^n)_{n\geq 1}$ converges weakly to $M^\infty$ starting at $M^0$. Using the approximation (\ref{comparison w}) to replace $\la \bw^n_t,\psi_i\ra$ by $\la w^{(n)}_{nt}(n^{1/3}\cdot),\psi_i\ra + \cO(n^{-2/3})$ on the third line,
together with Fubini's Theorem and the duality formula (\ref{dual formula}), we have
\begin{align}
&\E_{M_0^{(n)}}\bigg[\prod_{i=1}^k \bigg(\int_{\R^d\times \{0,1\}}\psi_i(x_i)\ind_{\{0\}}(\kappa_i)\bM^n_t(\rmd x_i,\rmd \kappa_i)\bigg)\bigg] \nonumber\\
& = \E_{M_0^{(n)}}\bigg[\int_{(\R^d)^k} \psi_1(x_1)\cdots \psi_k(x_k)\bigg\{\prod_{i=1}^k
\bw_t^n(x_i)\bigg\}\rmd x_1\ldots \rmd x_k\bigg]  \nonumber\\
& = \E_{M_0^{(n)}}\bigg[\int_{(\R^d)^k} \psi_1(x_1)\cdots \psi_k(x_k)\bigg\{\prod_{i=1}^k w^{(n)}_{nt}(n^{1/3}x_i)\bigg\}\rmd x_1\ldots \rmd x_k\bigg] + \mathcal{O}(n^{-2/3}) \nonumber\\
& = n^{-dk/3}\E_{M_0^{(n)}}\bigg[\int_{(\R^d)^k} \psi_1(n^{-1/3}x_1)\cdots \psi_k(n^{-1/3}x_k)\bigg\{\prod_{i=1}^k w^{(n)}_{nt}(x_i)\bigg\}\rmd x_1\ldots \rmd x_k\bigg]\nonumber\\
& \qquad \qquad \qquad \qquad\qquad \qquad\qquad \qquad\qquad \qquad \qquad \qquad\qquad \qquad \qquad \qquad+ \mathcal{O}(n^{-2/3}) \nonumber\\
& = n^{-dk/3} \int_{(\R^d)^k} \psi_1(n^{-1/3}x_1)\cdots \psi_k(n^{-1/3}x_k) \bfE_{\Xi[x_1,\ldots, x_k]}\bigg[\prod_{j=1}^{N_{nt}}w^{(n)}_0(\xi_{nt}^j)\bigg]\rmd x_1\ldots \rmd x_k \nonumber\\
& \qquad \qquad \qquad \qquad\qquad \qquad\qquad \qquad\qquad \qquad \qquad \qquad\qquad \qquad \qquad \qquad+ \mathcal{O}(n^{-2/3}) \nonumber\\
& = \int_{(\R^d)^k} \psi_1(x_1) \cdots \psi_k(x_k) \bfE_{\Xi[n^{1/3}x_1,\ldots, n^{1/3}x_k]}\bigg[\prod_{j=1}^{N_{nt}}w^{(n)}_0(n^{1/3}(n^{-1/3}\xi_{nt}^j))\bigg]\rmd x_1\ldots \rmd x_k \nonumber\\
& \qquad \qquad \qquad \qquad\qquad \qquad\qquad \qquad\qquad \qquad \qquad \qquad\qquad \qquad\qquad \qquad+ \mathcal{O}(n^{-2/3}) \nonumber\\
& = \int_{(\R^d)^k} \psi_1(x_1) \cdots \psi_k(x_k) \bfE_{\Xi[x_1,\ldots,x_k]} \bigg[\prod_{j=1}^{N^n_{t}}w_0^n(\xi_{t}^{n,j})\bigg]\rmd x_1\ldots \rmd x_k
+ \mathcal{O}(n^{-2/3})\nonumber\\
& = \bfE_{\Xi_0^n} \bigg[\prod_{j=1}^{N^n_{t}}w^0(\xi_{t}^{n,j})\bigg] + \mathcal{O}(n^{-2/3}). \label{scaled dual}
\end{align}
Now the expression on the l.h.s. of (\ref{scaled dual}) converges to the corresponding expression for $M^\infty$ as $n\rightarrow \infty$. Therefore, if $\Xi$ is the limit of a subsequence of $(\Xi^n)_{n\geq 1}$, then for every $t\geq 0$
\begin{equation}\label{limit duality}
\E_{M^0}\bigg[\prod_{i=1}^k \bigg(\int_{\R^d\times \{0,1\}}\psi_i(x_i)\ind_{\{0\}}(\kappa_i)M^\infty_t(\rmd x_i,\rmd \kappa_i)\bigg)\bigg]= \bfE_{\Xi_0} \bigg[\prod_{j=1}^{N^\infty_t}w^0(\xi_t^j)\bigg].
\end{equation}
On the other hand, as explained in Point {\bf 3)} of the proof of Theorem~\ref{th:fixed}, the same equality \eqref{limit duality} holds for any $w^0$ if we replace $\Xi_t$ in the r.h.s. by the empirical distribution at time $t$, $\Xi_t^\infty$, of the system of independent branching (and in dimension~1, coalescing) Brownian motions described in Theorem~\ref{th: conv duals}. As mentioned in the paragraph around (\ref{test functions dual}), test functions of the form used in the r.h.s. of (\ref{limit duality}) are separating. We can therefore conclude that the one-dimensional distributions of $(\Xi_t^n)_{t\geq 0}$ converge to those of $(\Xi_t^\infty)_{t\geq 0}$.
The generalisation to the finite-dimensional distributions is straightforward
since the duality formula (\ref{dual formula}) holds on any time interval
$[s,t]$ (if we replace $w_0$ by $w_s$ and $\xi_t^j$ by $\xi_{t-s}^j$).

Finally, since linear combinations of functions of the product form
\begin{equation}
x\mapsto\psi_1(x_1)\cdots \psi_k(x_k),
\end{equation}
with $\psi_i\in C_c(\R^d)$ a probability density function for every $i$, are dense (for the $\mathbb{L}^1$ norm) in the set of probability densities $\psi$ on $(\R^d)^k$, an analogue of Relation (\ref{limit duality}) can be established for this more general class of initial densities $\psi$. The same chain of arguments is then sufficient to conclude the proof of Lemma~\ref{lem: conv dual}.
\end{proof}

\medskip
\textbf{Tightness}

We now show tightness of the sequence $(\Xi^n)_{n\geq 1}$. To ease the notation, we write $\bfP_\psi$ for the probability measure on $D_{\cM_p(\R^d)}[0,\infty)$ under which the locations of the atoms of each $\Xi_0^n$ have density $\psi$. We first show that the compact containment condition holds if we see $(\Xi^n)_{n\geq 1}$ as a sequence of $\cM_p(\widehat{\R^d})$-valued Markov processes, where $\widehat{\R^d}$ is the one-point compactification of $\R^d$. We can then use Theorem~3.9.1 in \cite{EK1986}, together with the fact that the linear span of functions of the form \eqref{test functions dual} is dense in $C_b(\cM_p(\widehat{\R^d}))$ for the topology of uniform convergence on compact sets to reduce the tightness of $(\Xi^n)_{n\geq 1}$ to that of $(\Phi_{\exp,\ln f}(\Xi^n))_{n\geq 1}$ for every $f\in C^\infty(\widehat{\R^d})$ with values in $[0,1]$. More precisely, we show that for every such $f$, every $T>0$, every sequence of stopping times $(\tau_n)_{n\geq 1}$ bounded by $T$ and every $\e>0$, there exists $\delta=\delta(f,T,\psi,\e)$
such that
\begin{equation}\label{aim tightness}
\limsup_{n\rightarrow \infty}\, \bfP_{\psi}\Bigg[\sup_{0\leq t\leq \delta}\bigg|\prod_{i=1}^{N^n_{\tau_n+t}}f(\xi^{n,i}_{\tau_n+t}) - \prod_{i=1}^{N^n_{\tau_n}}f(\xi^{n,i}_{\tau_n})\bigg|>\e \Bigg] \leq \e.
\end{equation}
(This is actually stronger than the classical Aldous criterion based on stopping times \cite{ALD1978}, which considers the supremum over $t\in [0,\delta]$ of the probability that the increment between times $\tau_n$ and $\tau_n+t$ is larger than $\e$.) Finally, using Lemma~\ref{lem: conv dual} and Corollary~3.9.3 in \cite{EK1986}, we shall be able to conclude that $(\Xi^n)_{n\geq 1}$ is tight in $D_{\cM_p(\R^d)}[0,\infty)$, as desired (and furthermore that $\Xi^n$ converges weakly to $\Xi^\infty$ in $D_{\cM_p(\R^d)}[0,\infty)$).

We shall proceed in a number of steps.  First we control the maximum
number of particles in $\Xi_t^n$ up to time $T+1$. Not only does this
give us the compact containment condition, but conditional on this result, it is then easy to control the probability that there is a branch in an interval of
length $\delta$ (by branch, we mean that a particle is replaced by two ``parental'' particles during a selective event).  If we can also show that with high probability
there is no coalescence (\emph{i.e.}, no group of at least two particles is
ever removed during the same event and replaced by one or two ``parental'' particles), so that the number of particles in the system does
not change, then the problem is reduced to controlling
the jumps in a random walk. The most involved step, which is the substance
of Proposition~\ref{prop: controlling coalescence},
is showing that indeed there is no accumulation of coalescence events.

Let us replace $\R^d$ by its one-point compactification $\widehat{\R^d}$, so that the set of finite point measures with a total mass less than $K$ is compact for every $K>0$. Recall the notation $|\Xi|=\la \Xi,1\ra$ for the total mass of the measure $\Xi$. The following lemma thus implies the compact containment condition.
\begin{lemma}\label{lem: max number}
Let $T>0$. Given $\e>0$, there exists $K >0$ such that for every $n\geq 1$,
\begin{equation}
\bfP_\psi\bigg[\sup_{0\leq t\leq T+1} |\Xi^n_t| >K \bigg] \leq \frac{\e}{4}.
\end{equation}
\end{lemma}
\begin{proof}[Proof of Lemma~\ref{lem: max number}]
Recall that two particles are created when at least one of the extant particles is affected by a selective event. For a given particle of $\Xi^n$, this happens at rate $ns_n V_R u_n = u\sigma V_R$.
Furthermore, the presence of more than one particle in the area affected
by the event does not speed up the branching.  Consequently, the number of
particles in $(\Xi_t^n)_{t\geq 0}$ is stochastically bounded by the
number of particles
in a Yule process in which particles split
(independently of one another) into two offspring at rate $u\sigma V_R$.
Let $T>0$. Since the initial value, $\Xi_0^n$, has $k<\infty$ particles,
we conclude that there
exists $K\in \N$ such that for every $n\geq 1$,
\begin{equation}
\bfP_\psi\bigg[\sup_{0\leq t\leq T+1} |\Xi^n_t| >K \bigg] \leq \frac{\e}{4},
\end{equation}
as required.
\end{proof}

From now on, all our calculations proceed conditional on
the event
\begin{equation}
A_n= \bigg\{\sup_{t\in [0,T+1]} |\Xi^n_t| \leq K\bigg\}.
\end{equation}
From our reasoning above, we already see that for any $t\in [0,T]$,
conditional on $A_n$, the probability that at least one particle is
created during the time interval $(t,t+\delta]$ is bounded by
\begin{equation}
K\, \bfP_\psi\big[\hbox{a given particle branches in }(t,t+\delta]\big] \leq K
\big(1- e^{-u\sigma V_R \delta}\big) \leq u\sigma KV_R\, \delta.
\end{equation}
This bound is uniform in $n$ and so we see that there exists $\delta_1\in (0,1)$ such that for every $n\geq 1$,
\begin{equation}\label{proba branching}
\bfP_\psi\big[ \hbox{at least 1 particle created in }(\tau_n,\tau_n+\delta_1] \, ;\, A_n\big] \leq \frac{\e}{4}.
\end{equation}

We also want to control the probability of coalescence events.
Because of the calculation above, it is enough to do so in the absence of
branching.
\begin{proposition}
\label{prop: controlling coalescence}
Let $B_\delta^c$ denote the event that there is no branching event in
$(\tau_n,\tau_n+\delta]$.
There exists $\delta_2\in (0,\delta_1]$ such that
\begin{equation}
\bfP_\psi[\hbox{at least 1 coalescence in }(\tau_n,\tau_n+\delta_2]\, ;\, A_n, B_{\delta_2}^c] \leq  \frac{\e}{4}.
\end{equation}
\end{proposition}

Before proving Proposition~\ref{prop: controlling coalescence}, let us turn
to the final ingredient in the proof and
control the ``jumps'' of a single particle.

From the description in Section~\ref{ss: general dual}, after rescaling
of time and space, $\xi^{n,1}$ ``jumps'' (\emph{i.e.}, is removed and replaced by another particle seen as its parent) at rate
$nu_n V_R(1+s_n) = n^{2/3}uV_R (1+o(1))$, to a new location whose
distribution is symmetric about its current location. Furthermore, the locations of
the particle both before and after the jump belong to the same ball of
radius $Rn^{-1/3}$, and so the length of the jump is bounded by $2Rn^{-1/3}$.
Doob's Maximal Inequality and standard estimates for the variance of a compound Poisson process then imply that there exists $C_1>0$ such that for every $n$, any $s,\eta >0$, and every stopping time $T_n$,
\begin{equation}\label{max displacement}
\bfP_\psi\bigg[\sup_{t\in [0,s]}\big|\xi^{n,1}_{T_n+t} - \xi^{n,1}_{T_n}\big| > \eta \bigg] \leq \frac{C_1}{\eta^2}\, s,
\end{equation}
where we have used the strong Markov property of $\xi^{n,1}$ at time
$T_n$. From this, we can draw two conclusions. The first one, which is not necessary for the rest of the proof but gives some nice insight on our sequence of processes, is that taking $s=T$ and $T_n=0$, we can find a compact set $E\subset \R^d$ such that for every $n\geq 1$,
\begin{equation}\label{compact support}
\bfP_\psi\bigg[\sup_{t\in [0,T]}\Xi^n_t(E^c)>0;\, A_n \bigg] \leq \e.
\end{equation}
Indeed, since $\psi$ is integrable, there exists a compact set $\tilde{E}$ such that $\bfP_\psi[\Xi^n_0(\tilde{E}^c)>0]<\e /2$. Conditionally on all the initial particles belonging to $\tilde{E}$, by (\ref{max displacement}) we can then find a radius $\eta>0$ such that the probability that any of the (at most) $K$ particles leaves $E=\tilde{E}+B(0,\eta)$ is less than $\e/2$.

Second, conditional on the number of individuals not changing during a time interval
of length $\delta$, we can index the particles of $\Xi^n_{\tau_n}$ and
$\Xi^n_{\tau_n+\delta}$ by a common indexing set which we denote $I_n$, in such a way that a particle in $\Xi^n_{\tau_n+\delta}$ has the same label as a particle in $\Xi^n_{\tau_n}$ if and only if the position of the former can be seen as the result of a (potentially empty) series of jumps carried out by the latter during $(\tau_n,\tau_n+\delta]$. Under this assumption, for any $f\in C^\infty(\widehat{\R^d})$ with values in $[0,1]$, a Taylor expansion yields
\begin{equation}\label{diff products}
\bigg|\prod_{i\in I_n}f(\xi^{n,i}_{\tau_n+t}) - \prod_{i\in I_n}f(\xi^{n,i}_{\tau_n})\bigg| \leq C \big\|\nabla f\big\| \sum_{i\in I_n}\big|\xi^{n,i}_{\tau_n+t} - \xi^{n,i}_{\tau_n}\big|,
\end{equation}
for some $C>1$, where the sup norm of $\nabla f$ is finite since $\widehat{\R^d}$ is compact. Together with (\ref{max displacement}) and the choice $s=\delta$, $T_n=\tau_n$ and $\eta = \e/(K C\|\nabla f\|)$, this shows that there exists $\delta_3\in (0,\delta_2]$ such that
for $n$ large enough, writing $C_{\delta}^c$ for the event that there is no coalescence
in $(\tau_n, \tau_n+\delta]$,
\begin{equation}\label{proba moving}
\bfP_{\psi}\bigg[\sup_{t\in [0,\delta_3]}\bigg|\prod_{i\in I_n}
f(\xi^{n,i}_{\tau_n+t}) - \prod_{i\in I_n}f(\xi^{n,i}_{\tau_n})\bigg|>\e\, ;\,
A_n, B_{\delta_3}^c, C_{\delta_3}^c\bigg]\leq \frac{KC_1}{\eta^2}\, \delta_3 \leq \frac{\e}{4}.
\end{equation}
Combining Lemma~\ref{lem: max number}, (\ref{proba branching}), Proposition~\ref{prop: controlling coalescence} and (\ref{proba moving}), we obtain (\ref{aim tightness}) with $\delta=\delta_3$.

It remains to prove Proposition~\ref{prop: controlling coalescence}. Let us remark that it is not enough to consider particles at an
initial separation of order $\mathcal{O}(1)$ (or $\mathcal{O}(n^{1/3})$
before rescaling). In particular, when two particles are created through a
selective event, their (rescaled) initial distance is of order
$\mathcal{O}(n^{-1/3})$ and so we also need to control the
coalescence of particles starting from very small initial separations.

\begin{proof}[Proof of Proposition~\ref{prop: controlling coalescence}]
It suffices to consider just two particles and find $\delta_2>0$
such that the probability that they coalesce in a time interval
of length $\delta_2$ is bounded by $\e/(2K(K-1))$, irrespective of their
initial separation. Once this bound has been established, we can write
\begin{equation}
\bfP_\psi\big[\hbox{at least 1 coalescence in }(\tau_n,\tau_n+\delta_2]\, ;\, A_n, B_{\delta_2}^c\big] \leq \frac{K(K-1)}{2} \frac{\e}{2K(K-1)} = \frac{\e}{4},
\end{equation}
since, on the event $A_n$, there are at most $K(K-1)/2$ pairs of particles at
any time.

Recall that before scaling, each particle jumps at rate proportional to $u_n=un^{-1/3}$.
This makes it convenient to work in the timescale
$(n^{1/3}t,\, t\geq 0)$ and without rescaling space.
We shall write $\txi^{n,i}_{t} = \xi_{n^{1/3}t}^i$, $i\in \{1,2\}$.

When $\txi^{n,1}$ and $\txi^{n,2}$ are separated by more than $2R$,
they cannot be contained in the same reproduction event, and so they
evolve independently of one another.
The $i$th particle jumps at rate $n^{1/3}u_nV_R(1+s_n) = uV_R(1+o(1))$
to a new location, which is uniformly
distributed over the ball $B(Z,R)$, where $Z$ itself is chosen uniformly
at random from $B(\txi^{n,i},R)$. In what follows, we only need
that the jump made by each particle is an independent realisation
of a random variable $X$ taking values in $B(0,2R)$, whose distribution is
symmetric about the origin.

On the other hand, when $|\txi^{n,1}-\txi^{n,2}|<2R$, the two particles can
both lie in a region affected by a given reproduction event and their jumps
become correlated.
In particular, if they are both affected by this event, they merge together.
The infinitesimal generator of $((\txi_t^{n,1},\txi_t^{n,2}))_{t\geq 0}$, applied to any function $\phi\in C_0((\widehat{\R^d})^2)$ (the space of continuous functions on $(\widehat{\R^d})^2$ vanishing at infinity) takes the form
\begin{align}
&u (1+s_n) \int_{B(\txi^1,R)\setminus B(\txi^2,R)} \int_{B(x,R)} \frac{1}{V_R}\, \big(\phi(z,\txi^2)-\phi(\txi^1,\txi^2)\big)\rmd z\rmd x \nonumber\\
& +  u (1+s_n)\int_{B(\txi^2,R)\setminus B(\txi^1,R)} \int_{B(x,R)} \frac{1}{V_R}\, \big(\phi(\txi^1,z)-\phi(\txi^1,\txi^2)\big)\rmd z\rmd x \nonumber\\
& + u(1-un^{-1/3})(1+s_n)\int_{B(\txi^1,R)\cap B(\txi^2,R)} \int_{B(x,R)} \frac{1}{V_R}\, \big(\phi(z,\txi^2)+\phi(\txi^1,z)-2\phi(\txi^1,\txi^2)\big)\rmd z\rmd x \nonumber\\
& + u^2n^{-1/3}(1+s_n) \int_{B(\txi^1,R)\cap B(\txi^2,R)} \int_{B(x,R)} \frac{1}{V_R}\, \big(\phi(z,z)-\phi(\txi^1,\txi^2)\big)\rmd z\rmd x.
\end{align}
We can think of this as composed of two parts: the process
$((\hat{\xi}_t^{n,1},\hat{\xi}_t^{n,2}))_{t\geq 0}$ whose generator is
determined by the first three lines above,
on top of which a coalescence event occurs at instantaneous rate
$u^2n^{-1/3}(1+s_n)V_R(0,\hat{\xi}^{n,1}_t - \hat{\xi}^{n,2}_t)$
(recall that $V_R(0,a)$ is the volume of the intersection $B(0,R)\cap B(a,R)$).

With this description, the probability that the two particles have not
coalesced by time $\delta n^{2/3}$ (which corresponds to a time span of
$\delta$ on the timescale of $\xi^{n,i}$) is given by
\begin{equation}\label{coal proba}
\bfP_\psi\big[\tilde{T}> \delta n^{2/3}\big] =
\bfE_\psi\bigg[\exp\bigg\{-\frac{u^2(1+s_n)}{n^{1/3}}
\int_0^{\delta n^{2/3}}V_R(0,\hat{\xi}_s^{n,1}-\hat{\xi}_s^{n,2})\,
\rmd s\bigg\}\bigg],
\end{equation}
where we have written $\tilde{T}$ for the coalescence time of the two
particles.

Since $V_R(0,x)=0$ when $x\geq 2R$, it just remains to
establish how much time $\hat{\xi}^{n,1}-\hat{\xi}^{n,2}$ spends in the
ball $B(0,2R)$ by time $\delta n^{2/3}$. To do this, we define two
sequences of stopping times, $(\sigma^n_k)_{k\geq 1}$ and
$(\tau^n_k)_{k\geq 1}$ by
\begin{equation}
\sigma^n_1 = \inf\{t\geq 0\, :\, |\hat{\xi}_t^{n,1}-\hat{\xi}_t^{n,2}|\leq 2R\},
 \qquad \tau^n_1 = \inf\{t\geq \sigma^n_1\, :\, |\hat{\xi}_t^{n,1}-
\hat{\xi}_t^{n,2}|>2R\},
\end{equation}
and for every $k\geq 1$,
\begin{equation}
\sigma^n_k = \inf\{t\geq \tau^n_{k-1}\, :\, |\hat{\xi}_t^{n,1}-
\hat{\xi}_t^{n,2}|\leq 2R\}, \qquad
\tau^n_k = \inf\{t\geq \sigma^n_k\, :\, |\hat{\xi}_t^{n,1}-
\hat{\xi}_t^{n,2}|>2R\}.
\end{equation}
Now, we have the following result.
\begin{lemma}\label{lem: incursions}
There exists $\mathcal{C}>0$ such that for every $n,k\geq 1$,
\begin{equation}
\bfE_\psi\big[\tau^n_k - \sigma^n_k\big] \leq \mathcal{C}.
\end{equation}
\end{lemma}
In words, although the two particles are correlated when they are close
together, each ``incursion'' of $\hat{\xi}^{n,1}-\hat{\xi}^{n,2}$ inside
$B(0,2R)$ lasts only $\mathcal{O}(1)$ units of time, uniformly in $n$.
The proof of Lemma~\ref{lem: incursions} is similar to that of Lemma~6.6
in \cite{BEV2010} (based on the facts that the difference walk jumps at a
rate bounded from below by a positive constant, independent of its current
value, and that the probability that this jump leads to a sufficient
increase of their separation for $\hat{\xi}_t^{n,1}-\hat{\xi}_t^{n,2}$
to leave $B(0,2R)$ is
also bounded from below by a positive constant). Therefore, we omit it here.

Outside $B(0,2R)$, the difference $\hat{\xi}_t^{n,1}-\hat{\xi}_t^{n,2}$ has
the same law as a symmetric random walk, with jumps of size at most $2R$,
jumping at rate $2uV_R(1+s_n)$. Its behaviour will be determined by the
spatial dimension.

\medskip

\noindent
{${\mathbf d\geq 3}$:}
When $d\geq 3$, transience of the random walk guarantees that the number
of times $\hat{\xi}^{n,1}-\hat{\xi}^{n,2}$ returns to $B(0,2R)$ is a.s.~finite.
Since the parameter $n$ appears only in the jump rates and not in the embedded
chain of locations (during an excursion outside $B(0,2R)$), the probability
that the difference walk enters $B(0,2R)$ at least $k$ times
decays to $0$, uniformly in $n$, as $k\rightarrow \infty$. Together
with Lemma~\ref{lem: incursions} and the fact that $V_R(0,\cdot)$ is bounded,
this shows that for every $\eta>0$,
\begin{equation}
\lim_{n\rightarrow \infty}\bfP_\psi \bigg[\int_0^{\delta n^{2/3}}
V_R(0,\hat{\xi}_s^{n,1}-\hat{\xi}_s^{n,2})\, \rmd s > \eta\, \frac{n^{1/3}}{u^2(1+s_n)} \bigg] =0
\end{equation}
As a consequence, coming back to~(\ref{coal proba}), observing that
\begin{equation}
\bfP_\psi\big[\tilde{T}\leq  \delta n^{2/3}\big] = \bfP_\psi\bigg[\mathrm{Exp}(1)\leq \frac{u^2(1+s_n)}{n^{1/3}} \int_0^{\delta n^{2/3}} V_R(0,\hat{\xi}_s^{n,1}-\hat{\xi}_s^{n,2})\, \rmd s \bigg]
\end{equation}
(where $\mathrm{Exp}(1)$ denotes an exponential random variable with parameter $1$) and choosing $\eta$
small enough that $\bfP[\mathrm{Exp}(1)\leq \eta]\leq \e/(2K(K-1))$, we can conclude
that for any $\delta>0$,
\begin{equation}\label{proba of coal}
\limsup_{n\rightarrow \infty}\, \bfP_\psi\big[\tilde{T}\leq  \delta n^{2/3}\big] \leq \frac{\e}{2K(K-1)}.
\end{equation}

\medskip

\noindent
{$\mathbf d=2$:}
When $d=2$, we claim that there exists $\mathcal{C}'>0$, independent of $n$,
such that for every $x_1,x_2$ with $|x_1-x_2|> 2R$,
\begin{equation}
\bfP_{\{x_1,x_2\}}\big[\sigma^n_1> \delta n^{2/3}\big] \geq \frac{\mathcal{C}'}{\log(\delta n^{2/3})},
\end{equation}
where we have written $\bfP_{\{x_1,x_2\}}$ for the probability measure under which the two particles start at locations $x_1,x_2$. The proof of this claim is very similar to the beginning of the proof of
Lemma~4.2 in \cite{BEV2012}, and so we only sketch the main ideas. We can
a.s.~embed the trajectories of the difference process
$\hat{\xi}_t^{n,1}-\hat{\xi}_t^{n,2}$ into the trajectories of a
two-dimensional Brownian motion, in the same spirit as Skorokhod's embedding in
one dimension (see \emph{e.g.}~\cite{BIL1995}). Now, since the jumps of the difference
process (when outside $B(0,2R)$) are rotationally invariant, we have
\begin{equation}
\inf_{|x_1-x_2|>2R} \, \bfP_{\{x_1,x_2\}}\big[\hat{\xi}^{n,1}-\hat{\xi}^{n,2}
\hbox{ leaves }B(0,4R)\hbox{ before entering }B(0,2R)\big] >0,
\end{equation}
and the result then follows from that for Brownian motion, namely Theorem~2 in \cite{RR1966} applied with $a=2R$ and $r\geq 4R$.
As a consequence, the number $N_E^n$ of excursions outside $B(0,2R)$ that
the difference walk makes before starting an excursion of (time) length
at least $\delta n^{2/3}$ is stochastically bounded by a geometric random
variable with
success probability $C/\log(\delta n^{2/3})$. Now, once the difference walk
has started such a long excursion (say, the $k$th one), it is sure not to
come back within $B(0,2R)$ before time $\delta n^{2/3}$ and the number of
incursions in $B(0,2R)$ in the time interval $[0,\delta n^{2/3}]$ is
bounded by $k$. Thus, fixing $\eta>0$ as before and observing that $V_R(x,y)$ is bounded by the volume $V_R$ of a ball of radius $R$, we obtain that
\begin{align}
\bfP_\psi \bigg[\int_0^{\delta n^{2/3}}& V_R(0,\hat{\xi}_s^{n,1}-
\hat{\xi}_s^{n,2})\, \rmd s > \eta\, \frac{n^{1/3}}{u^2(1+s_n)} \bigg]\nonumber\\
 & \leq \bfP_\psi\big[N_E^n > C_E^n \log (\delta n^{2/3})\big] + \bfP_\psi\bigg[\sum_{k=1}^{\lceil C_E^n \log (\delta n^{2/3})\rceil} (\tau^n_k -\sigma^n_k) >\eta\, \frac{n^{1/3}}{u^2(1+s_n)V_R}\bigg] \nonumber\\
& \leq e^{-C_E^n \mathcal{C}'} + \frac{u^2(1+s_n)V_R}{\eta n^{1/3}}\, C_E^n\log(\delta n^{2/3}) \mathcal{C},
\end{align}
where the last inequality uses the stochastic bound of $N_E^n$ first,
and then Markov's inequality. Choosing $C_E^n = \log n$, for instance, we
deduce that for any $\delta >0$,
\begin{equation}
\lim_{n\rightarrow \infty} \bfP_\psi \bigg[\int_0^{\delta n^{2/3}}
V_R(0,\hat{\xi}_s^{n,1}-\hat{\xi}_s^{n,2})\, \rmd s > \eta\, \frac{n^{1/3}}{u^2(1+s_n)} \bigg] =0,
\end{equation}
and we conclude as in (\ref{proba of coal}).

\medskip

\noindent
{$\mathbf d=1$:}
Finally, when $d=1$ it is shown in \cite{PS1971}
that there exists $\mathcal{C}'>0$ such that for every $x_1,x_2$ such
that $|x_1-x_2|> 2R$,
\begin{equation}
\bfP_{\{x_1,x_2\}}[\sigma^n_1> \delta n^{2/3}] \geq \frac{\mathcal{C}'}{\sqrt{\delta}\, n^{1/3}}.
\end{equation}
Proceeding as before, and with the same notation, we therefore have
\begin{align}
\bfP_\psi \bigg[\int_0^{\delta n^{2/3}}&
V_R(0,\hat{\xi}_s^{n,1}-\hat{\xi}_s^{n,2})\, \rmd s >
\eta\, \frac{n^{1/3}}{u^2(1+s_n)} \bigg] \nonumber\\
 & \leq \bfP_\psi\big[N_E^n > C_E^n \sqrt{\delta} n^{1/3}\big] + \bfP_\psi\bigg[\sum_{k=1}^{\lceil C_E^n \sqrt{\delta} n^{1/3}\rceil} (\tau^n_k -\sigma^n_k) >\eta\, \frac{n^{1/3}}{u^2(1+s_n)\|V_R\|}\bigg] \nonumber\\
& \leq e^{-C_E^n \mathcal{C}'} + \frac{u^2(1+s_n)\|V_R\|}{\eta n^{1/3}}\, C_E^n \sqrt{\delta} n^{1/3} \mathcal{C}.
\end{align}
Choosing $C_E^n$ to be a constant large enough for the first term to be less
than $\e/(6K(K-1))$, and then $\delta_3>0$ small enough for the second term to be
less than $\e/(6(K(K-1))$, and finally taking $\eta$ small enough, we obtain that for any $\delta\leq \delta_3$,
\begin{align}
\bfP_\psi\big[\tilde{T}\leq  \delta n^{2/3}\big] & \leq
\bfP_\psi \bigg[\int_0^{\delta n^{2/3}}
V_R(0,\hat{\xi}_s^{n,1}-\hat{\xi}_s^{n,2})\, \rmd s > \eta\, \frac{n^{1/3}}{u^2(1+s_n)} \bigg] + \bfP[\mathrm{Exp}(1)\leq\eta]\nonumber \\
& \leq \frac{\e}{6K(K-1)}+\frac{\e}{6K(K-1)} + \frac{\e}{6K(K-1)} = \frac{\e}{2K(K-1)}.\label{proba of coal 2}
\end{align}
We have now proved the desired bound for the probability of a coalescence in
any dimension and the proof of Proposition~\ref{prop: controlling coalescence}
is complete.
\end{proof}

\section{Convergence of the rescaled SLFVS and its dual - the stable radius case} \label{s: convergence 2}
We proceed exactly as for the case of fixed radius.

\subsection{Proof of Theorem~\ref{th:stable}}
As in the proof of Theorem~\ref{th:fixed}, we first show that the sequence $(\bM^n)_{n\geq 1}$ is tight in $D_{\cM_\lambda}[0,\infty)$, then we show that any limit point of a subsequence satisfies the martingale problem stated in Theorem~\ref{th:stable}, and finally we prove that there exists at most one solution in $D_{\cM_\lambda}[0,\infty)$ to this martingale problem to conclude that $(\bM^n)_{n\geq 1}$ indeed converges to it.

\medskip
\textbf{1) Tightness.}

\medskip
We use the same method as in the proof of Theorem~\ref{th:fixed}, but the computations required are
different. Note that for every $n\geq 1$ the process $\bM^n$ has sample paths in $D_{\cM_\lambda}[0,\infty)$, since the unscaled process from which it is constructed has a.s. c\`adl\`ag paths by Theorem~\ref{th:existence}. Using again Theorem~3.9.1 in \cite{EK1986} and the compactness of $\cM_\lambda$, we reduce the proof of tightness of $(\bM^n)_{n\geq 1}$ to the proof of tightness of $(\Psi_{F,f}(\bM^n))_{n\geq 1}$ for every $F\in C^3(\R)$ and $f\in C_c^\infty (\R^d)$.

Hence, let us now fix $F$ and $f$ as above. Since $\Psi_{F,f}$ is a bounded function on $\cM_\lambda$ and consequently $(\Psi_{F,f}(\bM^n_t))_{n\geq 1}$ is a tight sequence for every $t\geq 0$, by the Aldous-Rebolledo criterion we only have to prove the equivalent of \eqref{tightness cal An} and \eqref{tightness cal Qn} after finding an expression for the predictable finite variation ${\cal A}^n$ of $(\Psi_{F,f}(\bM^n_t))_{t\geq 0}$ and for its predictable quadratic variation ${\cal Q}^n$. As earlier, we first replace $f$ by any function $\varphi\in C_c^\infty(\R^d)$ and then specialise the formulae we derive to a suitably chosen function $\varphi_f$ to conclude.

For any given $n\geq 1$, the extended generator of the unscaled process with parameters satisfying \eqref{cond intensity}, \eqref{cond intensity 2}, \eqref{assumptions par} and \eqref{form of mu}, acting on functions of the form $\Psi_{F,\varphi}$, is given by
\begin{align}
\cL \Psi_{F,\varphi}(M)& = \int_{\R^d}\int_1^{\infty}\int_{B(x,r)^2}\frac{1}{V_r^2}\, \Big\{w(y)(1+s_nw(z))\big[ F(\la \Theta^+_{x,r,u_n}(w),\varphi \ra) - F(\la w,\varphi\ra)\big] \label{unscaled generator} \\
 +& (1-w(y)+s_n(1-w(y)w(z)))\big[ F(\la \Theta^-_{x,r,u_n}(w),\varphi \ra) - F(\la w,\varphi\ra)\big]\Big\}\rmd y\rmd z\mu(\rmd r)\rmd x,\nonumber
\end{align}
where, as usual now, $w$ is a representative of the density of $M$. Arguing as in the part on tightness of the proof of Theorem~\ref{th:fixed} and using \eqref{unscaled generator} with $F$ and $F^2$, we obtain that the predictable finite variation part of $(\Psi_{F,\varphi}(M_t))_{t\geq 0}$ is given at any time $t\geq 0$ by
\begin{equation}
{\cal A}_t= \int_0^t \cL\Psi_{F,\varphi}(M_s)\, \rmd s,
\end{equation}
and its predictable quadratic variation by
\begin{align}
{\cal Q}_t = &\int_0^t \int_{\R^d}\int_1^{\infty}\int_{B(x,r)^2}\frac{1}{V_r^2}\, \Big\{w_s(y)(1+s_nw_s(z))\big[ F(\la \Theta^+_{x,r,u_n}(w_s),\varphi \ra) - F(\la w_s,\varphi\ra)\big]^2 \nonumber \\
 &\ + (1-w_s(y)+s_n(1-w_s(y)w_s(z)))\big[ F(\la \Theta^-_{x,r,u_n}(w_s),\varphi \ra) \nonumber\\
 & \qquad \qquad \qquad \qquad \qquad \qquad \qquad \qquad \qquad \qquad - F(\la w_s,\varphi\ra)\big]^2\Big\}\rmd y\rmd z\mu(\rmd r)\rmd x\rmd s.
\end{align}
To make the expressions easier to read, below we retain the notation $\beta$, $\gamma$ and $\delta$
from (\ref{def constants}). Let us now consider the process $(M^n_t)_{t\geq 0}$ whose density at time $t$ is $w_t^n(\cdot):=w_{nt}(n^\beta \cdot)$. Writing explicitly the martingale problem satisfied by $M^n$ and performing a change in the time and space variables, we obtain that the extended generator of this Markov process is given by
\begin{align}
&\cL^n \Psi_{F,\varphi}(M) \label{generator wn 2}\\
& = n\int_{\R^d}\int_1^{\infty}\int_{B(x,r)^2}\frac{1}{V_r^2}\, \Big\{w(n^{-\beta}y)\big(1+s_n w(n^{-\beta}z))\nonumber\\
& \qquad \qquad \qquad \qquad\qquad \qquad\qquad \qquad\times \big[ F(\la \Theta^+_{n^{-\beta}x,n^{-\beta}r,u_n}(w),\varphi \ra)  - F(\la w,\varphi\ra)\big]\big) \nonumber\\
& \quad + (1-w(n^{-\beta}y)+s_n(1-w(n^{-\beta}y)w(n^{-\beta}z)))\big[ F(\la \Theta^-_{n^{-\beta}x,n^{-\beta}r,u_n}(w),\varphi \ra) - F(\la w,\varphi\ra)\big]\Big\}\nonumber \\ & \qquad \qquad\qquad \qquad\qquad \qquad\qquad \qquad \qquad \qquad\qquad \qquad\qquad \qquad \qquad \qquad \rmd y\rmd z\mu(\rmd r) \rmd x  \nonumber \\
& = n^{1-\beta\alpha}\int_{\R^d}\int_{n^{-\beta}}^{\infty}\frac{1}{r^{d+1+\alpha}}\int_{B(x,r)^2}\frac{1}{V_r^2}\, \Big\{w(y)(1+s_nw(z))\big[F(\la \Theta^+_{x,r,u_n}(w),\varphi\ra)-F(\la w,\varphi\ra)\big] \nonumber \\
& \qquad  +\big(1-w(y)+s_n(1-w(y)w(z))\big)\big[F(\la \Theta^-_{x,r,u_n}(w),\varphi\ra)-F(\la w,\varphi\ra)\big]\Big\}\, \rmd y\rmd z\rmd r\rmd x, \nonumber
\end{align}
which allows us to write as in the fixed radius case that the predictable finite variation part of $\Psi_{F,\varphi}(M^n)$ is equal to  $(\int_0^t \cL^n \Psi_{F,\varphi}(M^n_s)\rmd s)_{t\geq 0}$, while its predictable quadratic variation is given by the integral with respect to time of the function in \eqref{generator wn 2} (applied to $M^n_s$) in which the increments $[F(\la \Theta^{\pm}_{x,r,u_n}(w^n_s),\varphi\ra)-F(\la w^n_s,\varphi\ra)]$ are squared. It remains to apply these results to $\varphi_f$ defined by
\begin{equation}
\label{stable varphi}
\varphi_f(x) = \frac{n^{d\beta}}{V_1}\int_{B(x,n^{-\beta})} f(y)\, \rmd y,
\end{equation}
and to use the fact that for every $t\geq 0$,
\begin{equation}
\Psi_{F,\varphi_f}(M^n_t)=F(\la w^n_t, \varphi_f\ra) = F(\la \bw^n_t,f\ra)=\Psi_{F,f}(\bM^n_t)
\end{equation}
(where $\bw^n_t$ is the density of $\bM^n_t$ defined in \eqref{def wn2}), to identify ${\cal A}^n$ as
\begin{equation}\label{expression An}
{\cal A}^n_t = \int_0^t \cL^n \Psi_{F,\varphi_f}(M^n_s)\, \rmd s,
\end{equation}
and ${\cal Q}^n$ as
\begin{align}
{\cal Q}^n_t=& n^{1-\beta\alpha}\int_0^t \int_{\R^d}\int_{n^{-\beta}}^{\infty}\frac{1}{r^{d+1+\alpha}}\int_{B(x,r)^2}\frac{1}{V_r^2}\, \Big\{w_s^n(y)(1+s_nw_s^n(z))\big[F(\la \Theta^+_{x,r,u_n}(w_s^n),\varphi_f\ra)\nonumber\\
& \qquad \qquad \qquad \qquad\qquad \qquad\qquad \qquad\qquad \qquad\qquad \qquad \qquad \qquad-F(\la \bw_s^n,f\ra)\big]^2 \nonumber \\
&  +\big(1-w_s^n(y)+s_n(1-w_s^n(y)w_s^n(z))\big)\big[F(\la \Theta^-_{x,r,u_n}(w_s^n),\varphi_f\ra)\nonumber\\
& \qquad \qquad \qquad \qquad\qquad \qquad\qquad \qquad\qquad \qquad-F(\la \bw^n_s,f\ra)\big]^2\Big\} \rmd y\rmd z\rmd r\rmd x \rmd s. \label{expression Qn}
\end{align}

Now that we have an expression for ${\cal A}^n$ and ${\cal Q}^n$, let us bound their increments to complete the proof of tightness of $(\Psi_{F,f}(\bM^n))_{n\geq 1}$. We start with ${\cal A}^n$.
As before, it is convenient to
split $\cL^n \Psi_{F,\varphi_f}(M_s^n)$ into its neutral and selective components.
Using a Taylor expansion of the function $F$, we obtain that the neutral part is equal to
\begin{align}
& n^{1-\beta\alpha}F'(\la \bw^n_s,f\ra)\int_{\R^d} \int_{n^{-\beta}}^{\infty}\frac{1}{r^{d+1+\alpha}}\int_{B(x,r)}\frac{1}{V_r}\Big[w^n_s(y)\big\la \Theta^+_{x,r,u_n}(w^n_s)-w^n_s,\varphi_f \big\ra \nonumber\\
& \qquad \qquad \qquad \qquad \qquad \qquad + (1-w^n_s(y))\big\la \Theta^-_{x,r,u_n}(w^n_s)-w^n_s,\varphi_f \big\ra \Big]\, \rmd y \rmd r \rmd x  \nonumber\\
& \quad + n^{1-\beta\alpha}\frac{F''(\la \bw^n_s,f\ra)}{2} \int_{\R^d} \int_{n^{-\beta}}^{\infty}\frac{1}{r^{d+1+\alpha}}\int_{B(x,r)}\frac{1}{V_r}\Big[w^n_s(y)\big\la \Theta^+_{x,r,u_n}(w^n_s)-w^n_s,\varphi_f \big\ra^2  \nonumber\\
& \qquad \qquad \qquad \qquad \qquad \qquad + (1-w^n_s(y))\big\la \Theta^-_{x,r,u_n}(w^n_s)-w^n_s,\varphi_f \big\ra^2 \Big]\, \rmd y \rmd r \rmd x  + \e_n \nonumber\\
& = n^{1-\beta\alpha-\gamma}u F'(\la \bw^n_s,f\ra)\int_{\R^d} \int_{n^{-\beta}}^{\infty}\frac{1}{r^{d+1+\alpha}}\int_{B(x,r)^2}\frac{1}{V_r}\, w^n_s(y)(\varphi_f(z)-\varphi_f(y))\, \rmd y\rmd z \rmd r \rmd x \nonumber\\
& \quad + n^{1-\beta\alpha-2\gamma}\frac{u^2}{2}\, F''(\la \bw^n_s,f\ra)\int_{\R^d} \int_{n^{-\beta}}^{\infty}\frac{1}{r^{d+1+\alpha}}\int_{B(x,r)}\frac{1}{V_r}\, \big\{w^n_s(y) \la \ind_{B(x,r)}(1-w^n_s),\varphi_f \ra^2 \nonumber\\
& \qquad \qquad \qquad \qquad\qquad + (1-w^n_s(y))\la \ind_{B(x,r)}w^n_s,\varphi_f\ra^2 \big\}\, \rmd y \rmd r \rmd x +\e_n, \label{finite variation}
\end{align}
with
\begin{equation}
|\e_n|\leq n^{1-\alpha\beta -3 \gamma}\frac{u^3C_F}{3!}
\int_{\R^d}\int_{n^{-\beta}}^\infty\frac{1}{r^{d+1+\alpha}}
\int_{B(x,r)}\frac{2}{V_r}\la\ind_{B(x,r)},\varphi_f\ra^3 \rmd y \rmd r \rmd x,
\end{equation}
where the constant $C_F$ is the supremum of $F^{(3)}$ over the bounded set in which its argument takes its values (recall that $\varphi_f\in C_c^\infty(\R^d)$).
Consider the first term on the right hand side of~(\ref{finite variation}).
Since $1-\alpha\beta -\gamma=0$, $n^{1-\beta\alpha-\gamma}=1$. We split the integral over
the radii into the sum of the integrals over $[n^{-\beta},1]$ and $[1,\infty)$.
By using a Taylor expansion of $\varphi_f$ and a symmetry argument to cancel the integral of $(z-y)\rmd z$,
we obtain that
\begin{align}
\bigg|u F'(\la \bw_s^n,f\ra)\int_{\R^d} & \int_{n^{-\beta}}^1\frac{1}{r^{d+1+\alpha}}\int_{B(x,r)^2}\frac{1}{V_r}\, w_s^n(y)(\varphi_f(z)-\varphi_f(y))\, \rmd y\rmd z \rmd r\rmd x\bigg|\label{taylor 1} \\
& \leq C \, \bigg|\int_{\R^d} \int_{n^{-\beta}}^1\frac{1}{r^{d+1+\alpha}}\int_{B(x,r)^2}\frac{1}{V_r}\, |z-y|^2 \ind_{\{B(x,r)\cap S_{\varphi_f}\neq \emptyset\}}\, \rmd y\rmd z\rmd r\rmd x\bigg| \nonumber\\
& \leq C' \hbox{Vol}(S_{\varphi_f}+B(0,1))\int_{n^{-\beta}}^1 \frac{1}{r^{d+1+\alpha}}\, r^{d+2} \rmd r = C''(1-n^{-\beta(2-\alpha)}) \nonumber
\end{align}
for some constants $C,C',C''>0$ (where $S_{\varphi_f}$ denotes the compact support of $\varphi_f$).
To control the integral over radii in $[1,\infty)$, the
cruder bound $|\varphi_f(y)-\varphi_f(z)|\leq 2\|f\|$ suffices and, using the fact that
\begin{equation}
\label{bound for integrating increments-bis}
\mathrm{Vol}\{x:S_{\varphi_f}\cap B(x,r)\neq\emptyset\}\leq C_2 (r^d \vee 1),
\end{equation}
we have
\begin{align}
\bigg|u F'(\la \bw^n_s,f\ra)\int_{\R^d} &\int_1^\infty\frac{1}{r^{d+1+\alpha}}\int_{B(x,r)^2}\frac{1}{V_r}\, w^n_s(y)(\varphi_f(z)-\varphi_f(y)) \rmd y\rmd z \rmd r\rmd x\bigg| \nonumber\\
& \leq C \int_{\R^d} \int_1^\infty \frac{1}{r^{d+1+\alpha}}\int_{B(x,r)^2}\frac{1}{V_r}\big(\ind_{\{y\in S_{\varphi_f}\}}+ \ind_{\{z\in S_{\varphi_f}\}}\big)  \rmd y\rmd z\rmd r \rmd x \nonumber\\
& \leq C' \int_1^\infty \frac{1}{r^{d+1+\alpha}}\int_{\R^d} \ind_{\{B(x,r)\cap S_{\varphi_f}\neq \emptyset\}}\hbox{Vol}(S_{\varphi_f}) \rmd x \rmd r\nonumber\\
& \leq C'' \int_1^\infty \frac{1}{r^{d+1+\alpha}} \ r^d \rmd r\leq C''', \label{taylor 2}
\end{align}
again for some constants $C,C',C''$ and $C'''$ which depend only on $d$, $F$ and $f$.

To control the second term on the right hand side of~(\ref{finite variation}), we use (\ref{bound for integrating increments-bis}) together with the inequality
\begin{equation}
\label{bound for increments-bis}
|\la \ind_{B(x,r)}w^n_s,\varphi_f\ra|\leq \|f\| \mathrm{Vol}(S_{\varphi_f}\cap B(x,r))\leq C_1\|f\| (r^d\wedge 1),
\end{equation}
to see that it is bounded by
\begin{align}
n^{-\gamma}& \frac{u^2C_F}{2}\times 2C_1^2\|f\|^2 \int_{\R^d}\int_{n^{-\beta}}^\infty \frac{1}{r^{d+1+\alpha}}(r^d \wedge 1)^2 \ind_{\{S_{\varphi_f}\cap B(x,r)\neq \emptyset\}} \rmd r \rmd x \nonumber\\
& = C_3 n^{-\gamma} \int_{n^{-\beta}}^\infty \frac{1}{r^{d+1+\alpha}}(r^d \wedge 1)^2 (r^d \vee 1) \rmd r \nonumber\\
& = C_3 n^{-\gamma} \int_{n^{-\beta}}^1\frac{r^{2d}}{r^{d+1+\alpha}} \rmd r + C_3 n^{-\gamma} \int_1^\infty \frac{r^d}{r^{d+1+\alpha}} \rmd r = C_4 n^{-\gamma}\big(1- n^{-\beta(d-\alpha)}\big).\label{vq1}
\end{align}
When $d\geq 2$, $d-\alpha>0$ and so this bound tends to $0$ as $n\rightarrow \infty$.
When $d=1$, $(\alpha -1)\beta -\gamma = 0$, and so this term is bounded by a constant
as $n\rightarrow \infty$. The same calculation shows that $\e_n\rightarrow 0$, uniformly in $w^n_s$,
as $n\rightarrow\infty$. As a consequence, in any dimension the absolute value of the neutral term of $\cL^n\Psi_{F,\varphi_f}(M)$ is bounded by a constant independent of $n$ and $M$.

Proceeding in the same way as for the second term above, we obtain that the ``selection'' term (\emph{i.e.}, that involving $s_n$) of $\cL^n\Psi_{F,\varphi_f}(M_s^n)$ is bounded by (recall that $1-\alpha\beta - \gamma=0$)
\begin{align}
2u\sigma n^{1-\beta\alpha -\gamma -\delta}&C_F\int_{\R^d}  \int_{n^{-\beta}}^\infty \frac{1}{r^{d+1+\alpha}}\int_{B(x,r)^3}\frac{1}{V_r^2}\, |\varphi_f(z')| \rmd y\,\rmd z\,\rmd z' \rmd r \rmd x \nonumber\\
& \leq Cn^{-\delta}\int_{n^{-\beta}}^\infty \frac{1}{r^{d+1+\alpha}}\int_{\R^d}  \big(1\wedge r^d\big)\, \ind_{\{B(x,r)\cap S_{\varphi_f}\neq \emptyset\}} \rmd x \rmd r \nonumber\\
& \leq  C'n^{-\delta}\int_{n^{-\beta}}^\infty \frac{1}{r^{d+1+\alpha}}\, \big(1\wedge r^d\big)\big(1\vee r^d\big) \rmd r \leq C'' n^{-\delta+\alpha\beta} = C'',\label{taylor 3}
\end{align}
since $\alpha\beta -\delta =0$. Combining \eqref{finite variation}, \eqref{taylor 1}, \eqref{taylor 2}, \eqref{vq1} and \eqref{taylor 3}, we obtain that there exists a constant $C$ independent of $n$ such that for every $0\leq t_1\leq t_2$,
\begin{equation}
\int_{t_1}^{t_2} \big|\cL^n\Psi_{F,\varphi_f}(M^n_s)\big|\rmd s\leq C(t_2-t_1),
\end{equation}
and therefore for every $T>0$, every sequence of stopping times $(\tau_n)_{n\geq 1}$ bounded by $T$, and every $\e>0$, we can choose $\eta>0$ small enough so that
\begin{equation}\label{target}
\limsup_{n\rightarrow \infty}\sup_{\theta\in [0,\eta]}\P\left[\big|{\cal A}^n_{\tau_n+\theta} - {\cal A}^n_{\tau_n}\big| >\e \right] = 0,
\end{equation}
which corresponds to the first part of the Aldous-Rebolledo criterion.

For the quadratic variation of the martingale part, a similar analysis yields that the integrand in
$\mathcal{Q}^n_t$ is bounded by
\begin{align}
Cn^{1-\beta\alpha-2\gamma} & \int_{\R^d}  \int_{n^{-\beta}}^\infty \frac{1}{r^{d+1+\alpha}} \int_{B(x,r)^3}\frac{1}{V_r}\, \varphi_f(z)\varphi_f(z')  \rmd y\, \rmd z\, \rmd z'\rmd r \rmd x \nonumber\\
& \leq C'n^{-\gamma}\int_{\R^d}  \int_{n^{-\beta}}^\infty \frac{1}{r^{d+1+\alpha}}\, \big(1\wedge r^d\big)^2 \ind_{\{B(x,r)\cap S_{\varphi_f}\neq \emptyset\}} \rmd r \rmd x \nonumber\\
& \leq C''n^{-\gamma} \int_{n^{-\beta}}^\infty \frac{1}{r^{d+1+\alpha}}\, \big(1\wedge r^d\big)^2 \big(1\vee r^d\big) \rmd r
\leq  C''' n^{-\gamma}\big(1+n^{-\beta(d-\alpha)}\big), \label{vq2}
\end{align}
which is bounded by a constant independent of $n$. As before, we conclude that the equivalent of \eqref{target} with ${\cal A}^n$ replaced by ${\cal Q}^n$ is satisfied for $\eta>0$ small enough. The Aldous-Rebolledo criterion allows us to conclude that the sequence of real-valued processes $(\Psi_{F,f}(\bM^n))_{n\geq 1}$ is tight, and since this is true for every $F\in C^3(\R)$ and $f\in C_c^\infty(\R^d)$, we obtain the tightness of $(\bM^n)_{n\geq 1}$ in $D_{\cM_\lambda}[0,\infty)$.

\medskip
\textbf{2) Identifying the limit.}

\medskip
Suppose $M^\infty \in D_{\cM_\lambda}[0,\infty)$ is the weak limit of a subsequence $(\bM^{n_k})_{k\geq 1}$, and for every $t\geq 0$, write $w^\infty_t$ for a representative of the density of $M^\infty_t$. We know from the previous paragraph that for every $f\in C_c^\infty(\R^d)$ and every $n\geq 1$,
\begin{equation}\label{approximate MP2}
\bigg( \Psi_{\mathrm{Id},f}\big(\bM^n_t\big)- \Psi_{\mathrm{Id},f}\big(\bM^n_0\big) - \int_0^t \cL^n \Psi_{\mathrm{Id},\varphi_f}(M^n_s)\rmd s\bigg)_{t\geq 0}
\end{equation}
is a martingale with predictable quadratic variation \eqref{expression Qn} (with $F=\mathrm{Id}$), where $\cL^n$ was defined in \eqref{generator wn 2} and $\varphi_f$ in \eqref{stable varphi}. As in the fixed radius case, we first show that for every $t\geq 0$,
\begin{equation}
\lim_{k\rightarrow \infty}\E\bigg[\bigg| \int_0^t \cL^{n_k} \Psi_{\mathrm{Id},\varphi_f}(M^{n_k}_s)\rmd s - \int_0^t \bigg\{\la w_s^\infty ,{\cal D}^\alpha f \ra - \frac{2u\sigma}{\alpha}\,
\la w_s^\infty (1-w_s^\infty),f\ra\bigg\}\, \rmd s\bigg| \bigg] = 0, \label{L1 limit2}
\end{equation}
so that we can then use the fact that the quantity in \eqref{approximate MP2} is a martingale, the fact that $\Psi_{\mathrm{Id},f}$ is a bounded continuous function and the Dominated Convergence Theorem to conclude that for every $0\leq t<t'$, $m\in \N$, $0\leq t_1 < \cdots < t_m\leq t$ and $h_1,\ldots,h_m\in C_b(\cM_\lambda)$,
\begin{align}
\E\bigg[\bigg(\la w^\infty_{t'},f\ra -\la w^\infty_t,f\ra -\int_t^{t'}\bigg\{\la w_s^\infty ,{\cal D}^\alpha f \ra -& \frac{2u\sigma}{\alpha}\,
\la w_s^\infty (1-w_s^\infty),f\ra\bigg\}\,\rmd s\bigg)\nonumber\\
& \qquad \qquad \times \bigg(\prod_{i=1}^m h_i\big(M^\infty_{t_i}\big)\bigg)\bigg]=0 \label{MP limit2}
\end{align}
and consequently that ${\cal Z}^f$ is a martingale (with respect to the natural filtration of $M^\infty$). In the case $d\geq 2$ this property is again sufficient to conclude, since we showed in \eqref{vq2} that the quadratic variation of the martingale \eqref{approximate MP2} tended to $0$ as $n\rightarrow \infty$, and therefore the limit ${\cal Z}^f$ is the constant process equal to $0$. We shall thus end this point $\mathbf{2)}$ by showing that in one dimension, the quadratic variation of ${\cal Z}^f$ and the bracket process between ${\cal Z}^f$ and ${\cal Z}^g$ have the required form.

Let us fix $f\in C_c^\infty(\R^d)$ and show \eqref{L1 limit2}. Let us first analyse the part of $\cL^n\Psi_{\mathrm{Id},\varphi_f}(M^n_s)$ corresponding to neutral events. By \eqref{finite variation} with $F=\mathrm{Id}$, since $1-\alpha\beta-\gamma =0$ this neutral part takes to form
\begin{align}
& un^{1-\alpha\beta -\gamma}\int_{\R^d} \int_{n^{-\beta}}^\infty \frac{1}{r^{d+1+\alpha}}\int_{B(x,r)^2}\frac{1}{V_r}\, w_s^n(y)(\varphi_f(z)-\varphi_f(y))\, \rmd y\rmd z \rmd r \rmd x \nonumber \\
&= u\int_{\R^d}  w_s^n(y)\int_{\R^d} \bigg(\int_{n^{-\beta}\vee \frac{|z-y|}{2}}^\infty \frac{1}{r^{d+1+\alpha}}\frac{V_r(y,z)}{V_r}\ \rmd r\bigg)(\varphi_f(z)-\varphi_f(y)) \, \rmd z \rmd y,
\end{align}
where $V_r(y,z)$ is again the volume of the intersection $B(y,r)\cap B(z,r)$. Now, a simple Taylor expansion to the second order gives us that
\begin{equation}
\varphi_f(z)-\varphi_f(y) = f(z)-f(y)+ \mathcal{O}(n^{-2\beta})
\big(\ind_{\{B_n(z)\cap S_f\neq\emptyset\}}+ \ind_{\{B_n(y)\cap S_f\neq\emptyset\}}\big),
\end{equation}
where $B_n(\cdot)=B(\cdot, n^{-\beta})$ and the error term is uniform in $y$ and $z$. Since
\begin{align}
& n^{-2\beta}\int_{\R^d}  w(y)\int_{\R^d}
\bigg(\int_{n^{-\beta}\vee \frac{|z-y|}{2}}^\infty
\frac{1}{r^{d+1+\alpha}}\frac{V_r(y,z)}{V_r}\, \rmd r\bigg)\,
\big(\ind_{\{B_n(z)\cap S_f\neq\emptyset\}}+
\ind_{\{B_n(y)\cap S_f\neq\emptyset\}}\big) \rmd z \rmd y \nonumber\\
& \qquad \qquad \leq C n^{-2\beta} \int_{S_f+B_n(0)}\int_{\R^d} \bigg(n^{-\beta} \vee \frac{|z-y|}{2}\bigg)^{-d-\alpha} \rmd z \rmd y \leq C' n^{-\beta(2-\alpha)} \rightarrow 0 \label{control error}
\end{align}
as $n\rightarrow \infty$, we can conclude that up to a vanishing error term,
the neutral part of $\cL^{n_k}\Psi_{F,\varphi_f}(M^{n_k}_s)$ is given by
\begin{equation}\label{first approx}
u\int_{\R^d}  w_s^n(y)\int_{\R^d} \bigg(\int_{n^{-\beta}\vee \frac{|z-y|}{2}}^\infty \frac{1}{r^{d+1+\alpha}}\frac{V_r(y,z)}{V_r}\ \rmd r\bigg)(f(z)-f(y)) \, \rmd z \rmd y.
\end{equation}
Now, our computations \eqref{taylor 1} and \eqref{taylor 2} in the proof of tightness imply that the function
\begin{equation}
a_n(y) : y\mapsto \int_{\R^d} \bigg(\int_{n^{-\beta}\vee \frac{|z-y|}{2}}^\infty \frac{1}{r^{d+1+\alpha}}\frac{V_r(y,r)}{V_r}\ \rmd r\bigg)(f(z)-f(y)) \, \rmd z
\end{equation}
is a continuous function, uniformly bounded in $y$ and $n$. Hence,
up to a vanishing error term we can first replace $w^n_s$ by $\bw^n_s$ in (\ref{first approx}) and, second,  use dominated convergence to pass to
the limit as $n\rightarrow \infty$ in (\ref{first approx}), along the converging subsequence. Doing so, and using the fact that all the error terms go to $0$ uniformly in $s$, we obtain that the limit in $\mathbb{L}^1$ norm of the neutral term in $\int_0^t\cL^{n_k}\Psi_{\mathrm{Id},\varphi_f}(M^{n_k}_s)\rmd s$ is equal to
\begin{equation}\label{limit motion}
u \int_0^t \int_{\R^d} w_s^\infty(y)\int_{\R^d}  \Phi(|z-y|)(f(z)-f(y))\rmd z\rmd y\rmd s,
\end{equation}
where, as in~(\ref{generator of stable motion}),
\begin{equation}
\Phi(|z-y|):= \int_{\frac{|z-y|}{2}}^\infty \frac{1}{r^{d+1+\alpha}}\frac{V_r(y,z)}{V_r}\, \rmd r.
\end{equation}

In passing, let us show the following property of the operator we obtain in the limit.
\begin{lemma}
\label{lemma generator alpha-stable}
For $f\in C_c^\infty(\R^d)$ write
\begin{equation}\label{generator stable motion}
\Da f(y) = u \int_{\R^d} \Phi(|z-y|) (f(z)-f(y)) \rmd z.
\end{equation}
Then $\Da$ is the infinitesimal generator of a symmetric $\alpha$-stable process
$(\zeta_t)_{t\geq 0}$.
\end{lemma}

\begin{proof}[Proof of Lemma~\ref{lemma generator alpha-stable}]
It is reassuring to first check that this is the generator of a well-defined L\'evy process:
\begin{align}
\int_{\R^d}(1\wedge |y|^2)& \int_0^\infty \frac{1}{r^{d+1+\alpha}}\, \frac{V_r(0,y)}{V_r}\, \rmd r\, \rmd y \nonumber\\
& \leq C \int_0^1\frac{1}{r^{d+1+\alpha}}\int_{B(0,2r)}|y|^2\, \rmd y \rmd r
+
C' \int_1^\infty\frac{1}{r^{d+1+\alpha}}\rmd r \nonumber\\
& \leq C''\int_0^1\frac{r^{d+2}}{r^{d+1+\alpha}}\rmd r +
C' \int_1^\infty\frac{1}{r^{d+1+\alpha}}\rmd r <\infty,
\end{align}
since $\alpha\in (1,2)$. Therefore the measure
\begin{equation}
\nu_\alpha(\rmd y)=\int_0^\infty \frac{1}{r^{d+1+\alpha}}\, \frac{V_r(0,y)}{V_r}\, \rmd r\, \rmd y
\end{equation}
is a L\'evy measure and there exists a unique L\'evy process with values in $\R^d$ whose L\'evy triplet is $(0,0,\nu_\alpha)$. By Theorem~6.8 in \cite{KS2016}, the operator $\Da$ is its infinitesimal generator.

To verify that the associated L\'evy process is a symmetric stable process,
we check the scaling property (the symmetry property is obvious from the form of $\nu_\alpha$). Let $b>0$.
The generator of $(b^{-1/\alpha}\zeta_{bt})_{t\geq 0}$ is given by
\begin{align}
\Da_bf(y) & = b u \int_{\R^d} \Phi(|z-b^{1/\alpha}y|) (f(b^{-1/\alpha}z)-f(y))\, \rmd z \nonumber\\
& = u b^{1+d/\alpha} \int_{\R^d} \Phi(|b^{1/\alpha}z-b^{1/\alpha}y|) (f(z)-f(y))\, \rmd z.
\end{align}
But a simple change of variables gives us that
\begin{align}
\Phi(|b^{1/\alpha}z-b^{1/\alpha}y|) &= \int_{\frac{b^{1/\alpha}|z-y|}{2}}^\infty \frac{1}{r^{d+1+\alpha}} \frac{V_r(b^{1/\alpha}y, b^{1/\alpha}z)}{V_r} \, \rmd r \nonumber\\
& = b^{-1-d/\alpha} \int_{\frac{|z-y|}{2}}^\infty \frac{1}{r^{d+1+\alpha}} \frac{V_r(y,z)}{V_r}\, \rmd r,
\end{align}
and so $\Da_b = \Da$ for all $b>0$. This shows the desired property of $\Da$.
\end{proof}

Having identified the neutral part of the limit, we now turn to the part of $\cL^n\Psi_{\mathrm{Id},\varphi_f}(M^n_s)$ corresponding to the selective events. It is given by
\begin{equation} \label{stable selection}
u\sigma n^{1-\beta\alpha-\gamma-\delta}\int_{\R^d} \int_{n^{-\beta}}^\infty \frac{1}{r^{d+1+\alpha}}\int_{B(x,r)^3}\frac{1}{V_r^2}\, (w^n_s(y)w^n_s(z)-w^n_s(z'))\varphi_f(z')\, \rmd y\rmd z\rmd z' \rmd r\rmd x.
\end{equation}
Now, the term which is linear in $w^n_s$ is easy to deal with: by Fubini's Theorem, it is equal to
\begin{align}
& u\sigma n^{-\delta} \int_{\R^d} \int_{n^{-\beta}}^\infty \frac{1}{r^{d+1+\alpha}}\int_{B(x,r)}w^n_s(z') \varphi_f(z')\, \rmd z'\rmd r\rmd x \nonumber \\
& \qquad= u\sigma n^{-\delta} \int_{\R^d} w^n_s(z')\varphi_f(z')
\bigg(\int_{n^{-\beta}}^\infty \frac{V_1r^d}{r^{d+1+\alpha}}\, \rmd r\bigg)\, \rmd z'
= \frac{u\sigma V_1}{\alpha}\, \la \bw^n_s,f\ra, \label{linear part selection}
\end{align}
where the last equality uses the fact that $\alpha\beta-\delta = 0$. It is then straightforward to obtain that
\begin{equation}\label{conv selection part}
\lim_{k\rightarrow \infty}\frac{u\sigma V_1}{\alpha}\, \E\bigg[\bigg|\int_0^t \la \bw^{n_k}_s,f\ra \rmd s - \int_0^t \la w^\infty_s,f\ra \rmd s\bigg|\bigg]=0.
\end{equation}

Similar calculations show that the ``quadratic'' term in (\ref{stable selection}) is equal to
\begin{equation}
u\sigma n^{-\delta}\int_{\R^d} \int_{n^{-\beta}}^{n^{-\beta}\log n} \frac{1}{r^{d+1+\alpha}}\bigg(\int_{B(x,r)}\frac{1}{V_r}\, w^n_s(y)\rmd y\bigg)^2 \int_{B(x,r)} \varphi_f(z)\rmd z \rmd r\rmd x+ \mathcal{O}\big( (\log n)^{-\alpha}\big).
\end{equation}
In contrast with the fixed radius case, here we first have to show that up to a vanishing error term, along the trajectories of the process $M^n$ we can replace the average of the density $w^n$ over a ball of radius at most $n^{-\beta} \log n$ by $\bw^n$, the average over a ball of radius $n^{-\beta}$ centered at the same point. In a second step, we use the same method as in the fixed radius case to prove that for every $t\geq 0$, $(\bw^{n_k}_t)^2$ converges to $(w^\infty_t)^2$ in the appropriate sense.

Concerning the first point, we have
\begin{align}
&\bigg(\int_{B(x,r)}\frac{1}{V_r} w^n_s(y)\rmd y\bigg)^2 \nonumber \\
& = \bigg(\int_{B(x,r)}\frac{1}{V_r} w^n_s(y)\rmd y + \bw^n_s(x)\bigg)\bigg(\int_{B(x,r)}\frac{1}{V_r} w^n_s(y)\rmd y - \bw^n_s(x)\bigg) + \bw^n_s(x)^2.
\end{align}
Suppose we have the following lemma (whose proof is quite technical and is given in Appendix~\ref{continuity stable}).
\begin{lemma}\label{lem: strong wish}
Under the conditions of Theorem~\ref{th:stable},
for every $r\in [n^{-\beta},n^{-\beta}\log n]$,
\begin{equation}
\lim_{n\rightarrow \infty} \E\bigg[\bigg|\int_{B(x,r)}\frac{w^n_s(y)}{V_r}\, \rmd y - \bw^n_s(x) \bigg|\bigg] =0
\end{equation}
uniformly in $x\in \R^d$ and uniformly in $s$ over compact time intervals $[0,t]$.
\end{lemma}

From this result, we can conclude from a dominated convergence argument
and a Taylor expansion of $\varphi_f$ that the
``quadratic'' part of (\ref{stable selection}) is equal to
\begin{equation}\label{quadratic part}
u\sigma n^{-\delta} \int_{\R^d} \int_{n^{-\beta}}^{n^{-\beta}\log n}
\frac{V_r}{r^{d+1+\alpha}}\,
\bw^n_s(x)^2 f(x)\, \rmd r\rmd x + \epsilon_n= \frac{u\sigma V_1}{\alpha}\, \la (\bw^n_s)^2,f\ra + \epsilon_n,
\end{equation}
where $\epsilon_n$ tends to zero  as $n\rightarrow \infty$ uniformly in $s$ over compact intervals of time.

As concerns the second point, we proceed as in (\ref{conv quadratic part}) and below. Using Proposition~\ref{prop: spatial continuity}$(ii)$ in Appendix~\ref{continuity stable}, the facts that the support of $f$ is bounded, and
that $p_\e$ is supported in $B(0,\e)$ (so that $\tau_2$ in (\ref{def:tau2}) is bounded by $\e^{\alpha/(d+1)}$ when $|z_1-z_2|\leq \e$ and $n$ is sufficiently large), we obtain that the first term in the decomposition (\ref{conv quadratic part}) of $\la (\bw^n_s)^2,f\ra$ is bounded by a constant (independent of $n,\e$) times
\begin{equation}
n^{-a} +\e^{\alpha/(d+1)}+ \e^{1/4}+\e^{\alpha/(2d+2)} + n^{-\beta(d-1)}\e^{(\alpha-d)/(2d+2)}.
\end{equation}
Letting $n$ tend to infinity in the above expression, we can write that the third term in the decomposition (\ref{conv quadratic part}) is bounded by a constant times
\begin{equation}
\e^{\alpha/(d+1)}+ \e^{1/4}+\e^{\alpha/(2d+2)} + \e^{(\alpha-1)/4}\mathbf{1}_{\{d=1\}}.
\end{equation}
Finally, the second term in the decomposition (\ref{conv quadratic part}) tends to $0$ by the assumption that $\bM^{n_k}$ converges to $M^\infty$. As in the fixed radius case, we can therefore conclude from \eqref{quadratic part} that
\begin{equation}\label{conv selection part 2}
\lim_{k\rightarrow \infty} \E\bigg[\bigg|\int_0^t \bigg\{\frac{u\sigma V_1}{\alpha}\la (\bw^{n_k}_s)^2,f\ra +\epsilon_n\bigg\}\rmd s - \frac{u\sigma V_1}{\alpha}\int_0^t \la (w^\infty_s)^2,f\ra \rmd s\bigg|\bigg]=0.
\end{equation}
(Note that this convergence is independent of the representatives of the different densities that we choose.)

Combining \eqref{limit motion}, \eqref{conv selection part} and \eqref{conv selection part 2}, we obtain \eqref{L1 limit2} and we can therefore conclude that ${\cal Z}^f$ is a martingale with respect to the natural filtration of $M^\infty$. As we already mentioned, when $d\geq 2$ this is sufficient to conclude that $M^\infty$ satisfies the equations stated in Theorem~\ref{th:stable}$(ii)$.

To identify the quadratic variation of ${\cal Z}^f$ and the bracket process between ${\cal Z}^f$ and ${\cal Z}^g$ when $d=1$, we proceed exactly as in the fixed radius case and therefore we do not provide all the details. Setting
\begin{equation}
W_t^n(f):= \la \bw_t^n,f\ra - \la \bw_0^n,f\ra - \int_0^t \cL^n \Psi_{\mathrm{Id},\varphi_f}(M^n_s)\rmd s,\qquad t\geq 0,
\end{equation}
we know from the paragraph \textbf{1)} on tightness that for every $n\geq 1$, $W^n(f)$ is a zero-mean martingale with predictable quadratic variation ${\cal Q}^n$ given in \eqref{expression Qn} (with $F=\mathrm{Id}$). As a consequence, for every $n\geq 1$, $0\leq t<t'$, $m\in \N$, $0\leq t_1<\cdots<t_m\leq t$ and $h_1,\ldots,h_m\in C_b(\cM_\lambda)$,
\begin{equation}\label{id qv 2}
\E\bigg[ \Big(\big(W^n_{t'}(f)\big)^2- \big(W^n_t(f)\big)^2 - {\cal Q}^n_{t'}+{\cal Q}^n_t\Big)\bigg(\prod_{i=1}^m h_i\big(\bM^n_{t_i}\big)\bigg)\bigg]=0.
\end{equation}
Recall from our calculations in \textbf{1)} that all summands in the expression for $W_t^n(f)$ are bounded uniformly in $n\geq 1$ and $t$ in a compact time interval. Furthermore, the same calculations as those we performed to obtain the limit of the selection part of $\cL^n\Psi_{\mathrm{Id},\varphi_f}(M^n_s)$ (see in particular \eqref{conv selection part} and \eqref{conv selection part 2}) show that
\begin{equation}\label{conv qv stable}
\lim_{k\rightarrow \infty}\E\bigg[\bigg|{\cal Q}^{n_k}_t - \frac{4u^2}{\alpha-1}\int_0^t \la w^\infty_s(1-w_s^\infty),f^2 \ra \rmd s\bigg|\bigg]=0,
\end{equation}
uniformly over compact intervals of time. Letting $n\rightarrow \infty$ in \eqref{id qv 2} along the converging subsequence, we arrive at
\begin{equation}\label{final qv stable}
\E\bigg[ \bigg(\big({\cal Z}^f_{t'}\big)^2- \big({\cal Z}^f_t\big)^2 - \frac{4u^2}{\alpha-1}\int_t^{t'}\la w^\infty_s(1-w_s^\infty),f^2 \ra \rmd s\bigg)\bigg(\prod_{i=1}^m h_i\big(M^\infty_{t_i}\big)\bigg)\bigg]=0.
\end{equation}
This allows us to identify the predictable quadratic variation of the martingale ${\cal Z}^f$ as
\begin{equation}
[{\cal Z}^f]_t= \frac{4u^2}{\alpha-1}\int_0^t\la w^\infty_s(1-w^\infty_s),f^2\ra \rmd s, \qquad t\geq 0.
\end{equation}
By the analogue of \eqref{small increments} (with $n^{-1/3}$ replaced by $n^{-\beta}$), every jump of $W^n(f)$ is bounded by $u_n\mathrm{Vol}(S_{\varphi_f})$ independently of the size of the radius of the event, where we recall that $u_n=un^{-\gamma}$. Consequently, ${\cal Z}^f$ has a.s. continuous trajectories. Since ${\cal Z}^f_0=0$, we can use the Dubins-Schwarz Theorem (or rather its extension since we do not know whether $[{\cal Z}^f]_\infty=+\infty$, see Remark~\ref{rk:DS}) to conclude that ${\cal Z}^f$ is a time-changed Brownian motion, solution to the stochastic differential equation
\begin{equation}
dW_t = \frac{2u}{\sqrt{\alpha-1}} \sqrt{\la w_t^\infty(1-w_t^\infty), f^2\ra}\, dB^f_t,
\end{equation}
where $B^f$ denotes standard Brownian motion. The bracket process between ${\cal Z}^f$ and ${\cal Z}^g$ is then obtained by the same kind of calculations, writing first the bracket process for a fixed $n$ and then identifying the limit by letting $n_k\rightarrow \infty$. We thus obtain that any limit of a subsequence of $(\bM^n)_{n\geq 1}$ satisfies the set of equations stated in Theorem~\ref{th:stable}$(i)$.

\medskip
\textbf{3) Uniqueness of the limit.}

\medskip
The argument is exactly the same as in the corresponding part of the proof of Theorem~\ref{th:fixed}. Indeed, by another modification of the results of Chapter~7 in \cite{LIA2009} (replacing Brownian motion by the symmetric $\alpha$-stable process $(\zeta_t)_{t\geq 0}$ generated by $\Da$ -- see Lemma~\ref{lemma generator alpha-stable}), we obtain that any solution to the limiting system of equations stated in Theorem~\ref{th:stable} is dual through the set of relations (\ref{dual formula}) to a system of particles following independent symmetric $\alpha$-stable processes (with the same law as $\zeta$), and branching independently at rate $u\sigma V_1/\alpha$ into two particles starting at the location of their parent. In one dimension, each pair of particles also coalesces at a rate $4u^2/(\alpha-1)$ times the local time at zero of their separation, independently of the other pairs. Since the set of all test functions of the form \eqref{test function again} is separating, we can again conclude that there is at most one solution to the system of equations of Theorem~\ref{th:stable}. Hence, this solution exists and the full sequence $(\bM^n)_{n\geq 0}$ converges to it in $D_{\cM_\lambda}[0,\infty)$.

\subsection{Proof of Theorem~\ref{th: conv duals alpha}}

Most of the proof is identical to that of Theorem~\ref{th: conv duals}.
That the only possible limit for $(\Xi_t^n)_{t\geq 0}$
is the system of branching (and in one dimension coalescing) symmetric
$\alpha$-stable processes described in the theorem, again
follows from an adaptation of Chapter~7 of \cite{LIA2009}, in which the only change is
that Brownian motion is replaced by the stable process generated by $\Da$
(see (\ref{generator stable motion})) and we have added natural selection/branching of particles. This gives us the analogue of Lemma~\ref{lem: conv dual} in the case of stable radii, whose proof is exactly the same as that of Lemma~\ref{lem: conv dual}.
\begin{lemma}\label{lem: conv dual stable}
The finite dimensional distributions of the system of scaled processes $\Xi^n$ converge as $n\rightarrow\infty$ to those of the system of branching and coalescing $\alpha$-stable motions $\Xi^\infty$,
described in the statement of Theorem~\ref{th: conv duals alpha}. In particular, the only possible
limit point for the sequence $(\Xi^n)_{n\geq 1}$ is $\Xi^\infty$.
\end{lemma}

Next, we have to show that the sequence $(\Xi^n)_{n\geq 1}$ is tight. Let again $\bfP_\psi$ denote the probability measure on $D_{\cM_p(\R^d)}[0,\infty)$ under which for each $n\geq 1$, the locations of the atoms of $\Xi_0^n$ have density $\psi$. As in the proof of Theorem~\ref{th: conv duals}, after showing that the compact containment condition holds if we replace $\R^d$ by its one-point compactification $\widehat{\R^d}$ and consider each $\Xi^n$ as taking its values in $\cM_p(\widehat{\R^d})$, we shall use Theorem~3.9.1 in \cite{EK1986} to deduce the tightness of $(\Xi^n)_{n\geq 1}$ in $D_{\cM_p(\widehat{\R^d})}[0,\infty)$ from the tightness of $(\Phi_{\exp,\ln f}(\Xi^n))_{n\geq 1}$ in $D_{[0,1]}[0,\infty)$ for every $f\in C^\infty(\widehat{\R^d})$ with values in $[0,1]$. More precisely, we show that for any such function $f$, every $T>0$, every sequence of stopping times $(\tau_n)_{n\geq 1}$ bounded by $T$ and every $\e>0$, there exists $\delta= \delta(f,T,\psi,\e)>0$ such that
\begin{equation}\label{tightness again}
\limsup_{n\rightarrow \infty} \, \bfP_{\psi}\bigg[\sup_{0\leq t \leq \delta} \bigg|\prod_{i=1}^{N^n_{\tau_n+t}} f\big(\xi_{\tau_n+t}^{n,i}\big) - \prod_{i=1}^{N^n_{\tau_n}} f\big(\xi_{\tau_n}^{n,i}\big)\bigg|>\e\bigg]\leq \e.
\end{equation}
Once these properties have been shown, we can use Lemma~\ref{lem: conv dual stable} and Corollary~3.9.3 in \cite{EK1986} to conclude that $(\Xi^n)_{n\geq 1}$ is tight in $D_{\cM_p(\R^d)}[0,\infty)$ and converges to $\Xi^\infty$.

Again, we proceed in four steps. First, by exactly the same arguments
as in the proof of Theorem~\ref{th: conv duals}, for every $T>0$ and every $\e>0$, there exists $K>0$ such that for every $n\in \N$ we have
\begin{equation}\label{tight 1}
\bfP_{\psi}[A_n]:=\bfP_{\psi}\bigg[\sup_{0\leq s\leq T+1} |\Xi^n_s|\leq K \bigg] \geq 1-\frac{\e}{4},
\end{equation}
which, in particular, grants us the compact containment condition since the set of all point measures on the compact space $\widehat{\R^d}$ with total mass less than $K$ is compact. Furthermore, there exists $\delta_1\in (0,1)$, independent of the subinterval of $[0,T]$ considered, such that
\begin{equation}\label{tight 2}
\bfP_{\psi}\big[\hbox{at least }1\hbox{ particle created in }(\tau_n,\tau_n+\delta_1]\, ;\, A_n\big] \leq \frac{\e}{4}.
\end{equation}

As before, the difficulty will be to control the coalescence (\emph{i.e.}, the events in which two or more particles are removed and replaced by one or two ``parental'' particles), but suppose for a moment that there is no
change in the number of particles in the interval $(\tau_n, \tau_n+\delta_2]$ and write
$I_n$ for the indexing set of the particles in $\Xi_{\tau_n}^n$. Then, exactly as before,
we can write
\begin{equation}\label{diff products 2}
\bigg|\prod_{i \in I_n} f\big(\xi_{\tau_n+t}^{n,i}\big) - \prod_{i\in I_n} f\big(\xi_{\tau_n}^{n,i}\big)\bigg|\leq C \|\nabla f\| \sum_{i\in I_n}\big|\xi_{\tau_n+t}^{n,i} - \xi_{\tau_n}^{n,i}\big|,
\end{equation}
and it suffices to consider the motion of a single particle to control the evolution of the whole set of particles. This is slightly more involved than in the fixed radius case.

Let $(Z^n_t)_{t\geq 0}$ be a L\'evy process, independent of $(\xi_t^n)_{t\geq 0}$ and with infinitesimal generator
\begin{equation}
D^n\phi(x) := u (1+s_n)\int_0^{n^{-\beta}} \frac{1}{r^{d+1+\alpha}} \int_{\R^d}\frac{V_r(x,y)}{V_r}\, \big(\phi(y)-\phi(x)\big)\, \rmd y\, \rmd r,
\end{equation}
for every $\phi\in C_0(\widehat{\R^d})$ and $x\in \widehat{\R^d}$. Then the process $(X_t)_{t\geq 0}$ defined by $X_t=\xi_t^n+Z_t^n$ has generator $(1+s_n)\Da$, where $\Da$ was shown in
Lemma~\ref{lemma generator alpha-stable} to be the generator of a symmetric stable process (indeed, observe that the
jump rates of $\xi^n$ and $Z^n$ depend only on the jump size $|y-x|$, hence the fact that the intensity measure of the jumps of $X$ is the sum of the intensity measures of $\xi^n$ and $Z^n$).
Using the strong Markov property and standard results on the growth of L\'evy
processes, see \emph{e.g.}~\cite{pruitt:1981}, we have for any $\eta, \delta>0$, and any stopping time $T_n$
\begin{equation}
\bfP_\psi\bigg[\sup_{t\in [0,\delta]}\, \big|X_{T_n+t} - X_{T_n}\big| >\eta\bigg]<C\frac{\delta}{\eta^\alpha}
\end{equation}
for a constant $C$ which is independent of $\eta,\delta$ and $T_n$.

Since
\begin{align}
&\bfP_\psi\bigg[\sup_{t\in [0,\delta]}\, \big|\xi_{T_n+t}^n - \xi_{T_n}^n\big| >\eta \bigg]\nonumber\\
&\leq \bfP_\psi\bigg[\sup_{t\in [0,\delta]}\, \big|X_{T_n+t} - X_{T_n}\big| >\eta \bigg]
+ \bfP_\psi\bigg[\sup_{t\in [0,\delta]}\, \big|Z_{T_n+t}^n - Z_{T_n}^n\big| >\eta \bigg], \label{decomposition}
\end{align}
it remains to show that
\begin{equation}
\bfP_\psi\bigg[\sup_{t\in [0,\delta]}\, \big|Z_{T_n+t}^n - Z_{T_n}^n\big| >\eta \bigg]\rightarrow 0,
\qquad\mbox{as }n\rightarrow \infty.
\end{equation}
Now, by construction,
the process $(Z_t^n)_{t\geq 0}$ has finite predictable quadratic variation, whose time derivative when $Z_t^n=x$ is
\begin{align}
(1+s_n) u & \int_0^{n^{-\beta}}\frac{1}{r^{d+1+\alpha}} \int_{\R^d} \frac{V_r(x,y)}{V_r}\, \big(f(y)-f(x)\big)^2 \, \rmd y\rmd r \nonumber\\
& = (1+s_n) u \int_0^{n^{-\beta}}\frac{1}{r^{d+1+\alpha}} \int_{\R^d} \frac{V_r(x,y)}{V_r}\, \big((y-x).\nabla f(x) + \mathcal{O}(|y-x|^2)\big)^2\, \rmd y\rmd r \nonumber\\
& \leq C \int_0^{n^{-\beta}} \frac{1}{r^{d+1+\alpha}} \int_{B(x,2r)}|y-x|^2\, \rmd y\rmd r = C'n^{-\beta(2-\alpha)},
\end{align}
where the Taylor expansion is justified since $V_r(x,y)=0$ if $|x-y|>2r$ and we are concentrating
on radii $r\leq n^{-\beta}$, and the first integral on the right hand side vanishes by
rotational symmetry. Hence, we can conclude that for any $\eta,\delta$,
\begin{equation}
\lim_{n\rightarrow \infty}\bfP_\psi\bigg[\sup_{t\in [0,\delta]}\, \big|Z_{T_n+t}^n - Z_{T_n}^n\big|
>\eta \bigg] =0.
\end{equation}
Coming back to (\ref{decomposition}), and taking $T_n=\tau_n$ and $\eta$ fixed, we can conclude that there
exists $\delta_3\in(0,\delta_2]$ such that for $n$ large enough,
\begin{equation}
\bfP_\psi\bigg[\sup_{t\in [0,\delta_3]}\, \big|\xi_{\tau_n+t}^n - \xi_{\tau_n}^n\big| >\eta \bigg] \leq \frac{\e}{4K}.
\end{equation}
Choosing $\eta = \e/(KC \|\nabla f\|)$ and recalling (\ref{diff products 2}), we obtain that for
all sufficiently large $n$,
\begin{equation}\label{tight 3}
\bfP_{\psi}\bigg[\sup_{t\in [0,\delta_3]}\bigg|\prod_{i \in I_n} f\big(\xi_{\tau_n+t}^{n,i}\big)
- \prod_{i\in I_n} f\big(\xi_{\tau_n}^{n,i}\big)\bigg|> \e\, ;\, A_n, B_{\delta_3}^c, C_{\delta_3}^c\bigg]\leq \frac{\e}{4},
\end{equation}
where as in the fixed radius case, $B_{\delta}^c$ is the event that there is no branching event in $(\tau^n,\tau^n+\delta]$ and $C_\delta^c$ is the event that there is no coalescence in $(\tau^n,\tau^n+\delta]$.

Finally, tightness will be proven if we can show that coalescence
events cannot accumulate. In particular, since we have controlled the total number of particles and
the probability of branching, we just need to control the probability that two particles
coalesce. The result will be based on the following lemma, in which we use again the interpretation of the replacement of a particle by its ``parent'' as a jump by this particle (or ancestral lineage - when there are two parents, we choose one of them uniformly at random).
\begin{lemma}\label{lem: limit coal}
Let $(\hat{\xi}^1_{n^\gamma t})_{t\geq 0}$ and $(\hat{\xi}^2_{n^\gamma t})_{t\geq 0}$ be two independent copies of the jump process obtained by following the (unscaled) position of one particle on the timescale $(n^\gamma t,\, t\geq 0)$,
and let $\zeta^n_t=\hat{\xi}^2_{n^\gamma t}-\hat{\xi}^1_{n^\gamma t}$
denote their difference.
Then, for every $t\geq 0$ we have:

$(i)$ When $d=1$, there exists $C(t)>0$ such that
\begin{equation}
\limsup_{n\rightarrow \infty} \, \bfE_\psi\bigg[\frac{1}{n^\gamma}\, \int_0^{n^{1-\gamma}t} \frac{1}{2^\alpha \vee |\zeta^n_s|^\alpha}\, \rmd s\bigg] \leq C(t).
\end{equation}
Furthermore, the function $t\mapsto C(t)$ can be chosen such that $C(t)\downarrow 0$ as $t\rightarrow 0$.

$(ii)$ When $d\geq 2$,
\begin{equation}
\lim_{n\rightarrow \infty} \bfE_\psi\bigg[\frac{1}{n^\gamma}\, \int_0^{n^{1-\gamma}t} \frac{1}{2^\alpha \vee |\zeta^n_s|^\alpha}\, \rmd s\bigg] =0.
\end{equation}
\end{lemma}
We defer the proof of Lemma~\ref{lem: limit coal} until after the end of the proof
of Theorem~\ref{th: conv duals alpha}.

Suppose that we start with a sample of two (non independent) particles at some (unscaled) separation $z_0\in \R^d$.
As before, we work on the timescale $n^\gamma$ so that a single particle jumps at
rate $\mathcal{O}(1)$ and we suppose the two particles $\xi^1$ and $\xi^2$ are currently at
locations $0$ and $z$ (in fact, only their separation matters). Then, the infinitesimal generator
$\Gamma$ of the difference walk $(\xi_{n^\gamma t}^2-\xi_{n^\gamma t}^1)_{t\geq 0}$ (until it reaches a cemetery state $\Delta$, say the point ``infinity'' in $(\widehat{\R^d})^2$, corresponding to the two walks having coalesced) is equal, for every given $\phi\in C_0((\widehat{\R^d})^2)$, to
\begin{align}
&\Gamma \phi(z)\nonumber\\
& = 2 u(1+s_n) \int_{\R^d} \bigg\{ \int_1^{\infty} \frac{1}{r^{d+1+\alpha}}\int_{\R^d} \ind_{\{0\notin B(x,r)\}}\ind_{\{z\in B(x,r)\}}\frac{ \ind_{\{y\in B(x,r)\}}}{V_r}\, \rmd x \rmd r \nonumber \\
& \qquad + (1-u_n) \int_1^{\infty} \frac{1}{r^{d+1+\alpha}}\int_{\R^d} \ind_{\{0\in B(x,r)\}}\ind_{\{z\in B(x,r)\}}\frac{ \ind_{\{y\in B(x,r)\}}}{V_r}\, \rmd x \rmd r  \bigg\}(\phi(y)-\phi(z))\rmd y \nonumber \\
 & \quad +u^2n^{-\gamma}(1+s_n) \int_{\R^d} \bigg\{ \int_1^{\infty} \frac{1}{r^{d+1+\alpha}}\int_{\R^d} \ind_{\{0\in B(x,r)\}}\ind_{\{z\in B(x,r)\}}\frac{\ind_{\{y\in B(x,r)\}}}{V_r}\, \rmd x\rmd r \bigg\}\nonumber\\
 & \qquad \qquad \qquad \qquad\qquad \qquad \qquad \qquad\qquad \qquad\qquad \qquad\qquad \qquad \times (\phi(\Delta)-\phi(z))\, \rmd y, \nonumber\\
& = 2u(1+s_n) \int_{\R^d} \bigg\{\int_1^\infty \frac{1}{r^{d+1+\alpha}}\, \frac{V_r(y,z)}{V_r}\, \rmd r\bigg\}(\phi(y)-\phi(z))\, \rmd y \nonumber \\
 & \quad - 2u^2n^{-\gamma}(1+s_n) \int_{\R^d}\bigg\{\int_1^\infty \frac{1}{r^{d+1+\alpha}}\frac{V_r(0,y,z)}{V_r}\, \rmd r \bigg\}(\phi(y)-\phi(z))\, \rmd y
\nonumber
\\ & \quad +u^2n^{-\gamma}(1+s_n) \int_{\R^d}\bigg\{\int_1^\infty \frac{1}{r^{d+1+\alpha}}\frac{V_r(0,y,z)}{V_r}\, \rmd r \bigg\}(\phi(\Delta)-\phi(z))\, \rmd y,\label{perturbation}
\end{align}
where $V_r(0,y,z)$ denotes the volume of the intersection $B(0,r)\cap B(y,r)\cap B(z,r)$.

From the first two terms on the r.h.s. of \eqref{perturbation}, we see that until coalescence
we can couple the difference walk (on the timescale $n^\gamma$) with the
difference $(\zeta_t^n)_{t\geq 0}$ between two independent random walks, each jumping according to the law
of a single walk but with each jump $z\mapsto y$ ``cancelled'' with probability
\begin{equation}
\Delta_n(z,y) = \frac{2u^2n^{-\gamma}(1+s_n) \int_1^\infty \frac{1}{r^{d+1+\alpha}}\frac{V_r(0,y,z)}{V_r}\, \rmd r}{2u(1+s_n)\int_1^\infty \frac{1}{r^{d+1+\alpha}}\, \frac{V_r(y,z)}{V_r}\, \rmd r}.
\end{equation}
(One can check that these two descriptions give rise to the same jump times and embedded chain.)
Each time we cancel a jump, with probability one half it was a coalescence in the original system (compare the second and third terms on the r.h.s. of \eqref{perturbation}),
but the key point is that if there are no cancelled jumps, then there was no coalescence.

It therefore suffices to show that we can find $\delta_2\in (0,\delta_1]$ such that, for sufficiently large $n$,
the probability that an event is cancelled in the interval $[0,\delta_2 n^{1-\gamma}]$ is
smaller than $\e/(4K(K-1))$.

Now, according to the expression on the right hand side of (\ref{perturbation}),
when the two particles lie at separation $z\in \R^d$, a cancelled event occurs at instantaneous rate
\begin{align}
2u^2n^{-\gamma}&(1+s_n) \int_{\R^d}\bigg\{\int_1^\infty \frac{1}{r^{d+1+\alpha}}\frac{V_r(0,y,z)}{V_r}\, \rmd r \bigg\} \rmd y \nonumber\\
& \leq 2u^2n^{-\gamma}(1+s_n) \int_{\R^d}\int_{1\vee \frac{|z|}{2}\vee \frac{|y|}{2}} \frac{1}{r^{d+1+\alpha}}\, \rmd r\rmd y = C_1n^{-\gamma} \big(2\vee |z|\big)^{-\alpha}.
\end{align}
Hence, (using the coupling with $(\zeta_t^n)_{t\geq 0}$), the probability of having no event
cancelled up to time $n^{1-\gamma}t$ (corresponding to time $nt$ in original units) is equal to
\begin{align}
&\bfE_\psi \bigg[\exp\bigg\{-\int_0^{n^{1-\gamma}t} 2u^2n^{-\gamma}(1+s_n) \int_{\R^d}\bigg\{\int_1^\infty \frac{1}{r^{d+1+\alpha}}\frac{V_r(0,y,\zeta_s^n)}{V_r\, }\, \rmd r \bigg\} \rmd y\rmd s\bigg\}\bigg] \\
& \geq \bfE_\psi\bigg[\exp\bigg\{-C_1n^{-\gamma}\int_0^{n^{1-\gamma}t} \big(2\vee |\zeta^n_s|\big)^{-\alpha}\, \rmd s \bigg\}\bigg]\geq 1- C_1\, \bfE_\psi\bigg[n^{-\gamma} \int_0^{n^{1-\gamma}t} \frac{\rmd s}{\big(2\vee |\zeta^n_s|\big)^{\alpha}}\bigg].\nonumber
\end{align}
But Lemma~\ref{lem: limit coal} shows that we can indeed find $\delta_2>0$ such that
\begin{equation}\label{events disappear}
\limsup_{n\rightarrow \infty}\, \bfE_\psi\bigg[n^{-\gamma} \int_0^{n^{1-\gamma}\delta_2} \frac{\rmd s}{\big(2\vee |\zeta^n_s|\big)^{\alpha}}\bigg] \leq \frac{\e}{2C_1K(K-1)}.
\end{equation}
Consequently,
\begin{equation}\label{proba of coal stable}
\bfP_\psi[\hbox{at least 1 coalescence in }(\tau_n,\tau_n+\delta_2];\, A_n,\, B_{\delta_2}^c]\leq \frac{K(K-1)}{2} \frac{\e}{2K(K-1)}=\frac{\e}{4},
\end{equation}
which was the last result we needed to complete the proof of tightness and therefore of Theorem~\ref{th: conv duals alpha}.

\begin{proof}[Proof of Lemma~\ref{lem: limit coal}]
As before, we shall exploit the fact that $(\zeta_t^n)_{t\geq 0}$ is ``nearly'' a symmetric
$\alpha$-stable process.
Indeed, the intensity at which $(\zeta_t^n)_{t\geq 0}$ jumps by some vector $y$ is independent of its current location and equal to
\begin{equation}
2(1+s_n) \bigg(\int_1^\infty \frac{1}{r^{d+1+\alpha}}\, \frac{V_r(0,y)}{V_r}\, \rmd r\bigg)\, \rmd y.
\end{equation}
Writing $(Z_t^n)_{t\geq 0}$ for a jump process, independent of $(\zeta_t^n)_{t\geq 0}$, starting at $0$ and with jump intensity
\begin{equation}
2(1+s_n) \bigg(\int_0^1 \frac{1}{r^{d+1+\alpha}}\, \frac{V_r(0,y)}{V_r}\, \rmd r\bigg)\, \rmd y,
\end{equation}
then the generator of the process $(X_t)_{t\geq 0}$, where $X_t = \zeta_t^n + Z^n_t$,
is precisely $2(1+s_n)$ times the
operator $\Da$ defined in (\ref{generator stable motion}), which we already checked
corresponds to a symmetric $\alpha$-stable process.
Once again, the idea is that the jumps of $(Z^n_t)_{t\geq 0}$ (which are bounded by $2$) do not
contribute much to the evolution of $(X_t)_{t\geq 0}$. More precisely, let us show that there exists
$C>0$ such that for every $n$ large enough and every $s\geq 1$,
\begin{equation}\label{control Z}
\bfP_\psi\bigg[\frac{|Z^n_s|}{\sqrt{s}} > (\log n)^2\bigg] \leq C e^{-(\log n)^2/d}.
\end{equation}
To this end, observe first that since the law of $Z^n_s$ is invariant under rotation, we can write that
\begin{equation}\label{one dim}
\bfP_\psi\bigg[\frac{|Z^n_s|}{\sqrt{s}} > (\log n)^2\bigg] \leq d\, \bfP_\psi\bigg[\frac{|Z_s^{n(1)}|}{\sqrt{s}} > \frac{(\log n)^2}{d}\bigg] = 2d\, \bfP_\psi\bigg[\frac{Z_s^{n(1)}}{\sqrt{s}} > \frac{(\log n)^2}{d}\bigg],
\end{equation}
where $Z_s^{n(1)}$ denotes the first coordinate of $Z^n_s$. Now, $(Z_s^{n(1)})_{s\geq 0}$
is again a symmetric L\'evy process with jumps bounded by $2$, and so Theorem 25.3 in \cite{SAT1999}
shows that for every $s,q\geq 0$, $\E[\exp(q Z_s^{n(1)})]<\infty$. In this case, it is known that
the characteristic exponent $\Psi^n$ of $(Z_s^{n(1)})_{s\geq 0}$, given here by a formula of the form
\begin{equation}
\Psi^n(q) = \int_{[-2,2]}\big(1-e^{iqx} + iqx\ind_{\{|x|< 1\}}\big) \, \mathfrak{m}^n(\rmd x),
\end{equation}
has an analytic extension to the half-plane with negative imaginary part, and we have
\begin{equation}
\bfE_\psi\Big[e^{qZ_s^{n(1)}}\Big] = e^{s\psi^n(q)}, \qquad \hbox{with}\quad \psi^n(q) = -\Psi^n(-iq).
\end{equation}
As a consequence, the Markov inequality gives us that
\begin{equation}\label{step 2}
\bfP_\psi\bigg[\frac{Z_s^{n(1)}}{\sqrt{s}} > \frac{(\log n)^2}{d}\bigg] \leq e^{-(\log n)^2/d + s\psi^n(1/\sqrt{s})}.
\end{equation}
Since the measure $\mathfrak{m}^n$ has support in $[-2,2]$, we can write that when $q$ is small
\begin{align}
\psi^n(q) & = -\int_{[-2,2]} \big(1- [1+qx +q^2x^2/2 + \mathcal{O}(q^3x^3)] + qx\ind_{\{|x|<1\}}\big) \, \mathfrak{m}^n(\rmd x) \nonumber\\
& = q \int_{[-2,2]} x\ind_{\{|x|\geq 1\}}\mathfrak{m}^n(\rmd x) + \frac{q^2}{2}\int_{[-2,2]}x^2\mathfrak{m}^n(\rmd x) + \mathcal{O}(q^3),
\end{align}
where the first term on the right is zero, by symmetry. Furthermore, $s_n\rightarrow 0$ and so $\mathfrak{m}^n$ converges to some finite $\mathfrak{m}$. Consequently, there exists a constant $C>0$ such that for every $s\geq 1$, $\psi^n(1/\sqrt{s})\leq C/s$. Together with (\ref{one dim}) and (\ref{step 2}), this gives us (\ref{control Z}).

It will be convenient to suppose that $\zeta_0=0$, but notice that there will be no loss of
generality in so-doing, since for $n$ sufficiently
large, $\zeta_0$ will be bounded by $(\log n)^2$ and so, for $s>1$, can be absorbed into our bound for $Z_s$.
Similarly, we can, and do,
replace $2^\alpha\wedge |\zeta_s^n|^\alpha$ by $1\wedge |\zeta_s^n|^\alpha$ in the
denominator of our integrand.

Based on these considerations, let us return to the integral of interest when $d\geq 2$.
Fixing $a\in (0,\gamma)$ and splitting the integral with respect to time into $\int_{[0,n^a]}+ \int_{[n^a, n^{1-\gamma}t]}$, we obtain
\begin{align}
&\bfE_\psi \bigg[\frac{1}{n^\gamma}\, \int_0^{n^{1-\gamma}t} \frac{1}{1 \vee |\zeta^n_s|^\alpha}\, \rmd s\bigg]  = \mathcal{O}(n^{a-\gamma}) + \frac{1}{n^\gamma}\int_{n^a}^{n^{1-\gamma}t}\bfE_\psi\bigg[\frac{1}{1\vee |X_s - Z^n_s|^\alpha}\bigg]\, \rmd s \nonumber\\
& \leq  C n^{a-\gamma} + n^{-\gamma}\int_{n^a}^{n^{1-\gamma}t}\bfP_\psi\bigg[\frac{|Z^n_s|}{\sqrt{s}}>(\log n)^2\bigg] \, \rmd s + n^{-\gamma}\int_{n^a}^{n^{1-\gamma}t}\bfE_\psi\bigg[\frac{\ind_{\{|Z^n_s|\leq (\log n)^2\sqrt{s}\}}}{1\vee |X_s -Z^n_s|^\alpha}\bigg]\, \rmd s \nonumber\\
& \leq C n^{a-\gamma} + Cn^{1-2\gamma}te^{-(\log n)^2/d} +  n^{-\gamma}\int_{n^a}^{n^{1-\gamma}t}\bfE_\psi\bigg[\frac{\ind_{\{|Z^n_s|\leq (\log n)^2\sqrt{s}\}}}{1\vee |X_s -Z^n_s|^\alpha}\bigg]\, \rmd s.  \label{estimate coal rate}
\end{align}
Since the first two terms on the right tend to $0$ as $n\rightarrow \infty$,
it now suffices to show that the last term remains bounded when $n$ is large.

By Lemma 5.3 in \cite{BW1998}, if $(p_s^\alpha)_{s\geq 0}$ denotes the transition density of
$(X_t)_{s\geq 0}$, we have, for every $s>0$ and $x\in \R^d$,
\begin{equation}\label{alpha semigroup}
p_s^\alpha(0,x) =: p_s^\alpha(x) = s^{-d/\alpha} p_1^\alpha(xs^{-1/\alpha})
\end{equation}
and there exists $C_{d,\alpha}>0$ (independent of $x$) such that
\begin{equation}\label{bound semigroup}
0\leq p_1^\alpha(x) \leq C_{d,\alpha}\big(1+|x|^{d+\alpha}\big)^{-1}.
\end{equation}
Hence, for any $s\geq n^a$ and any $z\in \R^d$ such that $|z|\leq (\log n)^2 \sqrt{s}$, we can write
\begin{align}
&\bfE_\psi \bigg[\frac{1}{1\vee |X_s - z|^\alpha}\bigg]  \leq s^{-d/\alpha}\int_{\R^d} \frac{1}{(1\vee |x-z|^\alpha)(1+|xs^{-1/\alpha}|^{d+\alpha})}\, \rmd x \\
& \leq s^{-d/\alpha} \int_{B(z,1)}\frac{1}{1+|xs^{-1/\alpha}|^{d+\alpha}}\, \rmd x + s^{-d/\alpha}\int_{B(z,1)^c}\frac{1}{|x-z|^\alpha (1+|xs^{-1/\alpha}|^{d+\alpha})}\, \rmd x \nonumber\\
& \leq Cs^{-d/\alpha} + C's^{-d/\alpha}\int_{B(0,s^{1/\alpha})\setminus B(z,1)} \frac{\rmd x}{|x-z|^\alpha} + C''s^{-d/\alpha}\int_{B(0,s^{1/\alpha})^c}\frac{\rmd x}{|x-z|^\alpha |xs^{-1/\alpha}|^{d+\alpha}}. \nonumber
\end{align}
But since $s \geq n^a$ and $|z|\leq (\log n)^2\sqrt{s}$, we have
\begin{equation}
|z|s^{-1/\alpha}\leq (\log n)^2 s^{\frac{1}{2}-\frac{1}{\alpha}}\leq (\log n)^2 n^{-a(2-\alpha)/(2\alpha)}\rightarrow 0,
\end{equation}
and so the second term on the right is bounded (after a change to polar coordinates) by
\begin{equation}
C' s^{-d/\alpha}\int_1^{s^{1/\alpha}} \rho^{d-1 -\alpha} d\rho = C' s^{-1},
\end{equation}
while the third term is bounded by
\begin{equation}
C'' s^{-d/\alpha} s^{1 + d/\alpha}\int_{s^{1/\alpha}}^\infty \rho^{d-1 -2\alpha -d}d\rho = C'' s^{-1}.
\end{equation}
Since all the constants depend on neither $z$ (in the range considered) nor $s$, we deduce
that the right hand side of (\ref{estimate coal rate}) is bounded by
\begin{equation}
C'n^{a-\gamma} + Cn^{1-2\gamma}t e^{-(\log n)^2/d} + C''n^{-\gamma}\big(n^{-a(d-\alpha)/d} +\log n + \log t\big) \rightarrow 0 \qquad \hbox{as }n\rightarrow \infty,
\end{equation}
which proves $(ii)$. \hfill

\medskip
The only point that differs when $d=1$ is that $1-d/\alpha >0$ and so
\begin{equation}
n^{-\gamma} \int_{n^a}^{n^{1-\gamma}t} s^{-1/\alpha}\rmd s \leq Cn^{-\gamma}n^{(1-\frac{1}{\alpha})(1-\gamma)} t^{1-\frac{1}{\alpha}}.
\end{equation}
An easy check confirms that $(1-\frac{1}{\alpha})(1-\gamma)-\gamma = 0$, and so $C(t)$ exists and is
proportional to $t^{1-\frac{1}{\alpha}}$. Since $\alpha>1$, we also have that
$C(t)\downarrow 0$ as $t\rightarrow 0$.
\end{proof}

\medskip
\noindent\textbf{Acknowledgment.} We are extremely grateful to the two reviewers (one in particular, who invested a lot of time in carefully reading and commenting different versions of the paper) and the Associate Editor for their valuable suggestions to significantly improve the presentation of the results and some of the arguments. AME was supported in part by EPSRC Grant EP/I01361X/1, AV was supported in part by the \emph{chaire Mod\'elisation Math\'ematique et Biodiversit\'e} of Veolia
Environnement-\'Ecole Polytechnique-Museum National d'Histoire Naturelle-Fondation X and FY was supported in part by EPSRC Grant EP/I028498/1.

\appendix

\section{Continuity estimates in the fixed radius case} \label{continuity fixed}

In this section, we state the continuity estimates for the scaled measures $M^n_T$ required in the proof of Theorem~\ref{th:fixed}. Because their proof is an adaptation of the (long and slightly more involved) proof of Proposition~\ref{prop: spatial continuity}$(ii)$, we do not give it here and instead refer to Appendix~\ref{continuity stable}. These estimates have the same flavour as the one dimensional estimates derived in \cite{MT1995} for the convergence of the local densities of $1$'s in the long range voter or contact process.

\begin{proposition}\label{prop: spatial continuity A}
Under the conditions of Theorem~\ref{th:fixed}, for every $T>0$ there exist $a,\lambda,v,C>0$ such that for every $n\geq 1$, $z_1,z_2\in \R^d$ such that $|z_1-z_2|<1$ and $\epsilon\in (0,1)$,
\begin{align}
  \E\bigg[\bigg|\frac 1 {V_{\epsilon}}& \int_{\R^d} w^n_T(x)
  (\ind_{\{|x-z_1|<\epsilon\}}-\ind_{\{|x-z_2|<\epsilon\}}) \rmd x\bigg|\bigg] \nonumber\\
  & \le Cn^{-a} + C\tau + C\big(|z_1-z_2|^{1/4} +\tau^{1/2}\big) e^{\lambda(|z_1|+\epsilon)} + Cn^{-(d-1)/6}\tau^{(2-d)/4},
\label{wish bound 2A}
\end{align}
where
\[
\tau=\tau(n,z_1,z_2) = n^{-v} \vee |z_1-z_2|^{2/(d+1)},
\]
and $\epsilon$ can depend on $n$ (as long as $\epsilon_n\leq 1$).

\end{proposition}

\section{Continuity estimates in the stable radius case}\label{continuity stable}

Our aim in this section is to obtain some continuity estimates for the measure $M^n_T$ (this time in the stable radius case), which are valid for fixed (large) $n$. Since in the stable radius case, we also need to compare the local densities of type-$1$ individuals over balls of radius $n^{-\beta}$ to the densities over balls of radius $\mathcal{O}(\log n)n^{-\beta}$, Proposition~\ref{prop: spatial continuity} below is more complete than Proposition~\ref{prop: spatial continuity A}. Lemma~\ref{lem: strong wish} will then follow as a corollary of item $(i)$.

\begin{proposition}\label{prop: spatial continuity}

Suppose the conditions of Theorem~\ref{th:stable} are satisfied. Fix $T>0$. Then,

\noindent $(i)$ There exist $a,C>0$ (dependent on $T$) such that for every $z\in \R^d$, $t\in [0,T]$, $n\geq 1$ and $n^{-\beta}\leq \epsilon_n <\epsilon_n'\leq 1$,
\begin{align}\nonumber
& \E\bigg[\bigg|\frac{1}{V_{\epsilon_n}}\int_{B(z,\epsilon_n)}w^n_t(y) \rmd y
- \frac{1}{V_{\epsilon_n'}}\int_{B(z,\epsilon_n')}w^n_t(y) \rmd y\bigg|\bigg] \\
\label{wish bound 1}
& \leq Cn^{-a} + C\tau_1+ C\epsilon_n'(\log n)^d \tau_1^{1-\frac{d+1}{\alpha}}
+(\log n)^{d/2}n^{\frac{\beta(\alpha-d)-\gamma}{2}}\Big[\epsilon_n^{'2}\,\tau_1^{1-\frac{2(d+1)}{\alpha}} \\
\nonumber
& \qquad \qquad \qquad \qquad \qquad \qquad \qquad \qquad \qquad \qquad \qquad \qquad + \epsilon_n'n^{-\frac{\beta(2-\alpha)d}{2(d+1)}}\tau_1^{1-\frac{d+1}{\alpha}}\Big]^{1/2},
\end{align}
where
\begin{equation}
\label{def:tau1}
\tau_1= \tau_1(n)=n^{-\beta(2-\alpha)/(2(d+1))}.
\end{equation}

\noindent $(ii)$ There exist $a,\lambda,C>0$ (dependent on $T$) such that for every $|z_1-z_2|<1$, $t\in [0,T]$, $n\geq 1$ and $\epsilon\in (0,1)$,
\begin{align}
  \E\bigg[\bigg|\frac{1}{V_{\epsilon}} &\int_{\R^d} w^n_t(x)
  (\ind_{\{|x-z_1|<\epsilon\}}-\ind_{\{|x-z_2|<\epsilon\}}) \rmd x\bigg|\bigg] \nonumber\\
  & \le Cn^{-a} + C\tau_2 + C\big(|z_1-z_2|^{1/4} +(\tau_2)^{1/2}\big) e^{\lambda(|z_1|+\epsilon)} + C\big(n^{-\beta(d-1)}\tau_2^{1-d/\alpha}\big)^{1/2},
\label{wish bound 2}
\end{align}
where
\begin{equation}
\label{def:tau2}
\tau_2=\tau_2(n,z_1,z_2) = n^{-\beta(2-\alpha)d/(4(d+1))} \vee |z_1-z_2|^{\alpha/(d+1)},
\end{equation}
and $\epsilon$ can depend on $n$ (as long as $\epsilon_n\leq 1$).
\end{proposition}
In particular, $(ii)$ implies uniform continuity of the limiting process of allele frequencies. That is:
\begin{corollary}
Suppose the conditions of Theorem~\ref{th:stable} are satisified and fix $T>0$. Then
for every $|z_1-z_2|<1$, $t\in[0,T]$ and $\epsilon\in (0,1)$,
\begin{align*}
\limsup_{n\rightarrow\infty}\ &
\E\bigg[\bigg|\frac 1 {V_{\epsilon}} \int_{\R^d} w^n_t(x)
  (\ind_{\{|x-z_1|<\epsilon\}}-\ind_{\{|x-z_2|<\epsilon\}}) \rmd x\bigg|\bigg]\\
&   \le C|z_1-z_2|^{(\alpha -1)/4}\ind_{\{d=1\}}+ C|z_1-z_2|^{\alpha/(d+1)}\ind_{\{d\geq 2\}} +
C\big(|z_1-z_2|^{1/4} \\
& \qquad \qquad \qquad \qquad\qquad \qquad\qquad \qquad+|z_1-z_2|^{\alpha/(2(d+1))}\big)
e^{\lambda(|z_1|+\epsilon)},
\end{align*}
where $C$ depends on $T$.
\end{corollary}

Before proving Proposition~\ref{prop: spatial continuity}, let us show how it
implies Lemma~\ref{lem: strong wish}.

\begin{proof}[Proof of Lemma~\ref{lem: strong wish}]
Set $\epsilon_n=n^{-\beta}$ and
$\epsilon_n'\in [n^{-\beta}, n^{-\beta}\log n]$ in $(i)$. Then
\[
\epsilon_n^{'2}\tau_1^{1-\frac{2(d+1)}{\alpha}}= (\log n)^2 n^{-2\beta}n^{\frac{\beta(2-\alpha)}{2(d+1)}- \frac{\beta(2-\alpha)}{\alpha}},
\]
and it is straightforward to check that the exponent of $n$ on the right hand side
is negative for any $\alpha\in (1,2)$.
Moreover,
\[
\epsilon_n' (\log n)^d \tau_1^{1-\frac{d+1}{\alpha}} \leq  (\log n)^a n^{-\beta(1-\frac{2-\alpha}{2}(\frac{1}{\alpha} - \frac{1}{d+1}))}
\]
for some $a>0$, and again one can check that the exponent of $n$ is negative in all dimensions.
Thus the right hand side of~(\ref{wish bound 1}) tends to zero and
the lemma follows.
\end{proof}

The rest of this section is devoted to the proof of Proposition~\ref{prop: spatial continuity}. Note that the different lemmas that appear in this proof will be shown later in Appendix \ref{sec:continuity_lemmas}.

\begin{proof}[Proof of Proposition~\ref{prop: spatial continuity}]
We define for $x\in\R^d$,
\[
  \sqcap_r(x) = \frac 1 {V_r} \ind_{\{|x|\le r\}},
\]
$\sqcap_r^{*k}$ to be the $k$-fold convolution of $\sqcap_r$ and
$\tilde w^n(x;r) = \frac 1 {V_r} \int_{B(x,r)} w^n(y) \rmd y$.
Recall the expression~(\ref{generator wn 2}) for the extended generator of $M^n$. For $\varphi \in \mathbb{L}^1(\R^d)$, we
follow our usual strategy of writing the
value of $\la w^n_T,\varphi\ra$ as a sum of drift and martingale terms (see the beginning of the proof of Theorem~\ref{th:stable}, where we can replace $\varphi\in C_c(\R^d)$ by $\varphi\in \mathbb{L}^1(\R^d)$ by a density argument): for any representative $w_t^n$ of the density of each $M_t^n$, we have
\begin{align}
  \la w^n_T,\varphi\ra =& \la w^n_0,\varphi\ra + {\cal M}^{n,\varphi}_T+u_n n^{1-\beta\alpha}
    \int_0^T \int_{\R^d} \int_{n^{-\beta}}^{\infty}
    \frac{1}{r^{d+1+\alpha}}\int_{B(x,r)^2}\frac{1}{V_r^2} \label{eq:decomp} \\
  & \qquad \times \Big\{
    w^n_t(y)(1+s_nw^n_t(z)) \la \ind_{B(x,r)}(1-w^n_t),\varphi\ra \nonumber \\
  & \qquad \qquad  -\big(1-w^n_t(y)+s_n(1-w^n_t(y)w^n_t(z))\big)
    \la \ind_{B(x,r)}w^n_t,\varphi\ra \Big\} \, \rmd y\rmd z\rmd r\rmd x\rmd t \nonumber \\
  =& \la w^n_0,\varphi\ra + {\cal M}^{n,\varphi}_T + u \int_0^T \int_{\R^d} \int_{n^{-\beta}}^{\infty}
    \frac{1}{r^{d+1+\alpha}}\int_{B(x,r)^2} \frac{1}{V_r^2} \big\{
    w^n_t(y)\la \ind_{B(x,r)},\varphi\ra   \nonumber\\
  & \  - \la \ind_{B(x,r)}w^n_t,\varphi\ra
    \nonumber + s_n (w^n_t(y)w^n_t(z)\la\ind_{B(x,r)},\varphi\ra
    - \la\ind_{B(x,r)}w^n_t,\varphi\ra) \big\} \, \rmd y\rmd z\rmd r\rmd x\rmd t \nonumber
\end{align}
(since $u_n n^{1-\beta\alpha}=u$), where $({\cal M}^{n,\varphi}_T)_{T\geq 0}$ is a mean zero martingale.
The first term in the integrand in~(\ref{eq:decomp}) is equal to:
\begin{align}
\int_{\R^d} &\int_{n^{-\beta}}^{\infty}
    \frac{1}{r^{d+1+\alpha}}\int_{B(x,r)} \frac{1}{V_r} \{
    w^n_t(y)\la \ind_{B(x,r)},\varphi\ra - \la \ind_{B(x,r)}w^n_t,\varphi\ra \}
    \rmd y\rmd r\rmd x  \nonumber \\
  =& \int_{\R^d} \int_{n^{-\beta}}^{\infty}
    \frac{1}{r^{d+1+\alpha}} \frac 1 {V_r} \int_{\R^d} \int_{\R^d}
    \ind_{\{|x-y|\le r\}} \ind_{\{|x-z|\le r\}}
    \{ w^n_t(y)\varphi(z) - w^n_t(z)\varphi(z) \} \rmd z\rmd y\rmd r\rmd x \nonumber \\
  =& \int_{n^{-\beta}}^{\infty} \int_{\R^d} \frac{V_r}{r^{d+1+\alpha}}
    \{ (\sqcap_r^{*2}*w^n_t)(z)\varphi(z) - w^n_t(z)\varphi(z) \} \rmd z\rmd r \nonumber \\
  =& \int_{\R^d} w^n_t(z) \int_{n^{-\beta}}^{\infty}
    \frac{V_r}{r^{d+1+\alpha}} \{(\sqcap_r^{*2}*\varphi)(z) - \varphi(z)\}\rmd r\rmd z.
\label{eq:decomp_stable1}
\end{align}
The second term in the integrand in~(\ref{eq:decomp}) is equal to
\begin{align*}
s_n \int_{\R^d}& \int_{n^{-\beta}}^{\infty} \frac{1}{r^{d+1+\alpha}}
    (\tilde w^n_t(x;r)^2 \la\ind_{B(x,r)},\varphi\ra
    - \la\ind_{B(x,r)}w^n_t,\varphi\ra) \rmd r\rmd x \nonumber \\
  =& s_n \int_{(\R^d)^2} \int_{n^{-\beta}}^{\infty} \frac{1}{r^{d+1+\alpha}}
    \ind_{\{|x-y|<r\}} (\tilde w^n_t(x;r)^2- w^n_t(y)) \varphi(y) \rmd y\rmd x\rmd r.
\end{align*}
Since $u_n^2 n^{1-\beta\alpha} = u^2 n^{-\gamma}
= u^2 n^{-(\alpha-1)/(2\alpha-1)}$,
the martingale term in~(\ref{eq:decomp}) has predictable quadratic variation
\begin{align*}
  [{\cal M}^{n,\varphi}]_T &= u^2 n^{-\gamma}
    \int_0^T \int_{\R^d} \int_{n^{-\beta}}^{\infty}
    \frac{1}{r^{d+1+\alpha}} \Big\{\tilde w^n_t(x;r)(1+s_n \tilde w^n_t(x;r))
    \la \ind_{B(x,r)}(1-w^n_t),\varphi\ra^2 \nonumber \\
  & \qquad \qquad + \big(1-\tilde w^n_t(x;r)+s_n(1-\tilde w^n_t(x;r)^2)\big)
    \la \ind_{B(x,r)}w^n_t,\varphi\ra^2 \Big\} \, \rmd r\rmd x\rmd t.
\end{align*}

It is convenient to replace this martingale problem by a mild version, obtained by
replacing $\varphi$ by the time dependent function $\zeta^n_t(x,z,\epsilon)$ chosen to
solve
\[ \partial_t \zeta^n_t(x;z,\epsilon)
  = \int_{n^{-\beta}}^{\infty} \frac{uV_r}{r^{d+1+\alpha}}
    \Big[(\sqcap_r^{*2}*\zeta^n_t(\cdot;z,\epsilon)) (x) - \zeta^n_t(x;z,\epsilon) \Big]\, \rmd r\]
with initial condition $\zeta^n_0(\cdot;z,\epsilon)$. That is
$\zeta^n_t(\cdot;z,\epsilon)$ is the density at time $t$ of the
$d$-dimensional L\'evy process, $(X^n_t)_{t\geq 0}$, with initial distribution
$\zeta^n_0(\cdot;z,\epsilon)$, zero drift, no Brownian component,
and L\'evy measure
\[ \nu^n(\rmd x) = \int_{n^{-\beta}}^\infty \frac{uV_r}{r^{d+1+\alpha}}
  \sqcap_r^{*2}(x) \rmd r \rmd x \]
for $x\in\R^d$ (in particular, $\zeta^n_t(x,z,\epsilon)\in \mathbb{L}^1(\R^d)$).
Here we assume that for any $n\in \N$, $z\in \R^d$ and $\epsilon>0$, $\zeta^n_0(\cdot;z,\epsilon) = \zeta^n_0(\cdot -z;0,\epsilon)$ and that the support of $\zeta^n_0(\cdot;0,\epsilon)$ is included in $B(0,\epsilon)$. Of course, the particular example we have in mind is $\zeta^n_0(\cdot;z,\epsilon)=\frac 1 {V_\epsilon} \ind_{\{|\cdot-z|<\epsilon\}}$.
The parameter $\epsilon$ can be taken to depend on $n$.
We observe that $\nu^n$ is radially symmetric.
Let
\begin{align*}
  a^n(x;r) &= r^{-d} \int_{\R^d} \ind_{\{|x-y|<r\}}
    (\tilde w^n(y;r)^2- w^n(x)) \rmd y \\
  b^n(x;r) &= \tilde w^n(x;r)(1+s_n \tilde w^n(x;r)) \\
  c^n(x;r) &= 1-\tilde w^n(x;r)+s_n(1-\tilde w^n(x;r)^2).
\end{align*}
Notice that $a^n$, $b^n$ and $c^n$ are all uniformly (in $n$, $x$ and $r$)
bounded between constants.  Suppose that we know the exponential decay of
$\zeta^n_{T-t}(\cdot;z,\epsilon)$ (which we prove in Lemma~\ref{lem:levy_exp3}),
then substituting in the martingale problem in the usual way, we obtain
\begin{align}
  \la w^n_T, \zeta^n_0(\cdot;z,\epsilon) \ra
  &= \la w^n_0, \zeta^n_T(\cdot;z,\epsilon) \ra  + {\cal M}^{n,\zeta_0^n(\cdot;z,\epsilon)}_T  \nonumber\\
 & \qquad + u s_n \int_0^T \int_{n^{-\beta}}^{\infty} \frac{1}{r^{1+\alpha}}
    \int_{\R^d} a^n_t(x;r) \zeta^n_{T-t}(x;z,\epsilon) \rmd x\rmd r\rmd t , \label{drift}\\
  \left[{\cal M}^{n,\zeta_0^n(\cdot;z,\epsilon)}\right]_T &= u^2 n^{-\gamma}
    \int_0^T \int_{n^{-\beta}}^{\infty} \frac{1}{r^{d+1+\alpha}} \int_{\R^d}
    \Big\{ b^n_t(x;r) \la \ind_{B(x,r)}(1-w^n_t),\zeta^n_{T-t}(\cdot;z,\epsilon)
    \ra^2 \nonumber\\
  & \qquad \qquad + c^n_t(x;r) \la
    \ind_{B(x,r)}w^n_t,\zeta^n_{T-t}(\cdot;z,\epsilon)
    \ra^2 \Big\} \rmd x \rmd r \rmd t \nonumber\\
  =& u^2 n^{-\gamma} \int_0^T \int_{n^{-\beta}}^{\infty}
    \frac{1}{r^{d+1+\alpha}} \int_{\R^d}\bigg\{
    b^n_t(x;r) \left(\int_{B(x,r)} (1-w^n_t(y)) \zeta^n_{T-t}(y;z,\epsilon) \rmd y
    \right)^2 \nonumber\\
  &\qquad\qquad
    + c^n_t(x;r) \left(\int_{B(x,r)} w^n_t(y) \zeta^n_{T-t}(y;z,\epsilon) \rmd y
    \right)^2\bigg\} \rmd x\rmd r\rmd t. \label{martingale}
\end{align}

In order to control the different terms appearing in~(\ref{drift}) and~(\ref{martingale}), we are
going to need to establish continuity estimates for $\zeta^n$.
In preparation for this, note that $(X^n_t)_{t\geq 0}$ is a continuous time random walk with jump rate
\[ A = \int_{n^{-\beta}}^\infty \frac{uV_r}{r^{d+1+\alpha}} \rmd r
  = V_1 n^{\alpha\beta}. \]
To describe the corresponding jump chain, let
$R_k$ be i.i.d. $\R^d$-valued random variables
distributed according to
$\frac{V_1} A r^{-(1+\alpha)} \ind_{\{r>n^{-\beta}\}} \rmd r$,
$Z_{1,k}$ and $Z_{2,k}$
be independent uniformly distributed random variables in $B(0,1)$, and
$Y_k=R_k (Z_{1,k}+Z_{2,k})$. Then we can write
\begin{equation}
  X^n_t = X^n_0+\sum_{k=1}^{K_t} Y_k,
\label{eq:xy_stable}
\end{equation}
where $K_t$ is a Poisson random variable with parameter $At$.
We define $f_Y$ as the density of $Y_1$, $f^{*k}_Y$ to be the $k$-fold
convolution of $f_Y$,
\begin{align}
  q^{n,\{k\}}_t(x) &= f^{*k}_Y(x) \P[K_t=k]
  = e^{-At} \frac{(At)^k}{k!} f^{*k}_Y(x) \nonumber\\
  q^n_t(x) &= \sum_{k=1}^\infty q^{n,\{k\}}_t(x). \label{def q}
\end{align}
Then,
\[ \zeta^n_t(x;z,\epsilon) = \zeta^n_0(x;z,\epsilon) e^{-At}
  + (\zeta^n_0(\cdot;z,\epsilon)*q^n_t(\cdot))(x). \]
Our estimates will involve splitting into two cases, according to whether the walk has taken greater
or fewer than $L$ steps in the interval $[0,t]$ and so it will be convenient to define
$q^{n,I}_t=\sum_{k\in I} q^{n,\{k\}}$
for $I\subset [1,\infty)$,
$\zeta^{n,\{k\}}_t(\cdot ; z,\epsilon)=\zeta^n_0(\cdot;z,\epsilon)*q^{n,\{k\}}_t(\cdot)$, and
$\zeta^{n,I}_t=\sum_{k\in I} \zeta_t^{n,\{k\}}$ for $I\subset [0,\infty)$.

Since the number of jumps made by the walk in $[0,t]$ has mean proportional to $n^{\alpha\beta}$,
with probability tending to one as $n\rightarrow\infty$ it will take at least $n^{c\alpha\beta}$
steps for any $c\in (0,1)$. We define
$c_1:= (\alpha-1)/(2\alpha) \in (0,1)$ and set
\[  L = n^{c_1\alpha\beta/2}.\]

In Section~\ref{sec:continuity_lemmas}, we shall prove a sequence of lemmas that control the behaviour
of the random walk. In particular, we establish the following. For every $t\geq 0$, let $q_t$ be the density function of value at time $t$ of the symmetric $\alpha$-stable process starting at 0 and with Laplace exponent
\[
\psi(\theta):=\int_{\R^d}\big(e^{i\theta\cdot x}-1\big) \nu(\rmd x),
\]
where
\[
\nu(\rmd x):=\int_0^\infty \frac{uV_r}{r^{d+1+\alpha}}\sqcap_r^{*2}(x)\rmd r\rmd x.
\]
(Note that this process is the one appearing in Lemma~\ref{lemma generator alpha-stable}.)
\begin{lemma}
Let $\|f\|_\lambda = \sup_x |f(x)| e^{\lambda |x|}$.
Let $c_2\in (0,\alpha)$ be a constant.
Recall $L=n^{c_1\alpha\beta/2}$, with $c_1=\frac{\alpha-1}{2\alpha}$.
For $x,y,z\in\R^d$ and $n$,
\begin{enumerate}
\item[(i)]
If $M\ge 2$ and $t\in [n^{-c_2\beta(2-\alpha)/(2(d+1))},T]$, then
\[ |q^{n,[M,\infty)}_t(x)-q_t(x)| \le C_{d,T} n^{-\beta(2-\alpha)d/(2(d+1))}
  + C_d n^{\beta d} (a^{M-1}+\P[K_t<M]) \]
  for some $a\in (0,1)$ independent of $M$ and $T$.
Furthermore,
\[ |q^{n,[L,\infty)}_t(x)-q_t(x)| \le C_{d,T} n^{-\beta(2-\alpha)d/(2(d+1))}.\]
\item[(ii)] If $t>0$, then $|q_t(x)-q_t(y)| \le C t^{-(d+1)/\alpha} |x-y|$.
\item[(iii)]
If $t\in [n^{-c_2\beta(2-\alpha)/(2(d+1))},T]$, then
\[ |q^{n,[L,\infty)}_t(x)-q^{n,[L,\infty)}_t(y)|
  \le C t^{-(d+1)/\alpha} |x-y| + C_{d,T} n^{-\beta(2-\alpha)d/(2(d+1))}.
\]
\item[(iv)] If $\lambda>0$, $t\le T$ and $|x|\ge 1$, then
  $q^{n,[1,\infty)}_t(x) \le C_{\lambda,T} e^{-\lambda (|x|-1)}$.
\item[(v)] If $\lambda>0$, $t\in [n^{-c_2\beta(2-\alpha)/(2(d+1))},T]$
and $|y-z|\le 1$, then
\begin{align*}\|\zeta^{n,[L,\infty)}_t(\cdot;y,\epsilon)
    -\zeta^{n,[L,\infty)}_t(\cdot;z,\epsilon)\|_\lambda
  \le  C_{\lambda,d,T} e^{\lambda \epsilon}& (t^{-(d+1)/(2\alpha)} |y-z|^{1/2}\\
      &\qquad+ n^{-\beta(2-\alpha)d/(4(d+1))}) e^{\lambda |z|},
\end{align*}
where $\epsilon$ can depend on $n$.
\end{enumerate}
\label{lem:levy_exp3}
\end{lemma}
Recall the definitions
\begin{align*}
\tau_1 & = n^{-\beta(2-\alpha)/(2(d+1))}, \\
  \tau_2 &= n^{-\beta(2-\alpha)d/(4(d+1))} \vee |z_1-z_2|^{\alpha/(d+1)},
\end{align*}
The quantity $\tau_1$ (\emph{resp.}, $\tau_2$) will be used in the bounds needed to prove
Proposition~\ref{prop: spatial continuity}$(i)$ (\emph{resp.}, $(ii)$).
Observe that for $t\ge\tau_2$ and $|z_1-z_2|<1$,
the estimate on the right hand side of Lemma~\ref{lem:levy_exp3}$(v)$ is
\[ \le C_{\lambda,d,T} (|z_1-z_2|^{1/2}+\tau_2) e^{\lambda \epsilon}e^{\lambda|z_1|}. \]

Since the organisations of the proofs are similar, we shall show
Proposition~\ref{prop: spatial continuity}$(i)$ and $(ii)$ in parallel. In both cases, we set
\[
\zeta^n_0(\cdot;z,\epsilon) := \frac{1}{V_\epsilon}\, \ind_{B(z,\epsilon)}
\]
(although most of the proof does not require a specific form for $\zeta_0^n$), and we estimate
\begin{align*}
(i) & \qquad \la w^n_T, \zeta^n_0(\cdot;z,\epsilon_n)-\zeta^n_0(\cdot;z,\epsilon_n')\ra,\\
(ii) & \qquad \la w^n_T, \zeta^n_0(\cdot;z_1,\epsilon)-\zeta^n_0(\cdot;z_2,\epsilon)\ra
\end{align*}
for the range of parameters stated in Proposition~\ref{prop: spatial continuity}, using~(\ref{drift})
and~(\ref{martingale}).

\subsection{Drift terms}
Let us split the different terms into the cases in which $K_t$, the number of jumps of $X^n$ by time $t$, is less than or larger than $L$. This first gives (using the fact that the function $a^n_t$ is bounded uniformly in $n,t,x,r$):
\begin{align}
\bigg|u s_n & \int_0^T\int_{n^{-\beta}}^{\infty}\frac{1}{r^{1+\alpha}}\int_{\R^d}a_t^n(x;r)\big(\zeta^{n,[0,L)}_{T-t}(x;z,\epsilon_n)-\zeta^{n,[0,L)}_{T-t}(x;z,\epsilon_n')\big)\rmd x\rmd r\rmd t \bigg| \nonumber \\
& \leq Cus_n \int_0^T \int_{n^{-\beta}}^{\infty}\frac{1}{r^{1+\alpha}}\int_{\R^d} \big(\zeta^{n,[0,L)}_{T-t}(x;z,\epsilon_n)+\zeta^{n,[0,L)}_{T-t}(x;z,\epsilon_n')\big)\rmd x\rmd r\rmd t \nonumber\\
& \leq Cus_n n^{\alpha\beta}\int_0^T \P[K_t<L]\, \rmd t \leq Cn^{-(1-c_1)\alpha\beta} \label{1}
\end{align}
by Lemma \ref{lem:Kt1} (which controls $\P[K_t<L]$)
and the fact that, by definition, $s_n n^{\alpha\beta}\equiv \sigma$.
The same estimate holds for $(ii)$ and the corresponding integral.

Next, let us split the remaining integral into an integral over large and small times. We can write
\begin{align*}
\zeta_{T-t}^{n,[L,\infty)}(x;z,\epsilon_n) &- \zeta_{T-t}^{n,[L,\infty)}(x;z,\epsilon_n') = \int_{\R^d}\big(\zeta_0^n(x';z,\epsilon_n)- \zeta_0^n(x';z,\epsilon_n')\big) q_{T-t}^{n,[L,\infty)}(x-x')\rmd x' \\
& = \int_{\R^d}\zeta_0^n(x';z,\epsilon_n)\big(q_{T-t}^{n,[L,\infty)}(x-x')-q_{T-t}^{n,[L,\infty)}(x-z)\big)\rmd x' \\
& \qquad - \int_{\R^d}\zeta_0^n(x';z,\epsilon_n')\big(q_{T-t}^{n,[L,\infty)}(x-x')-q_{T-t}^{n,[L,\infty)}(x-z)\big)\rmd x'.
\end{align*}
(The extra terms cancel since $\int_{\R^d}\zeta_0^n(x',z',\epsilon_n)\rmd x'=1$ for all choices of $\epsilon_n$.)
Since the second term above will be bounded in the same way as the first term, let us just consider
the first one. We have by Lemma~\ref{lem:levy_exp3}$(iii)$ and $(iv)$ (recalling also that the
support of $\zeta^n_0(\cdot;z,\epsilon)$ is contained in $B(z,\epsilon)$):
\begin{align}
C&us_n\int_0^{T-\tau_1}\int_{n^{-\beta}}^\infty \frac{1}{r^{1+\alpha}}\int_{\R^d}\int_{\R^d}\zeta_0^n(x';z,\epsilon_n)\big|q_{T-t}^{n,[L,\infty)}(x-x')-q_{T-t}^{n,[L,\infty)}(x-z)\big|\rmd x'\rmd x\rmd r\rmd t \nonumber \\
& \leq Cs_nn^{\alpha\beta}\int_{\tau_1}^T \int_{B(z,\log n)}\int_{\R^d}\zeta_0^n(x';z,\epsilon_n)\big(t^{-(d+1)/\alpha}|z-x'| \nonumber\\
& \qquad \qquad \qquad \qquad\qquad \qquad\qquad \qquad \qquad \qquad\qquad \qquad+ C_{d,T}n^{-\beta(2-\alpha)d/(2(d+1))}\big) \rmd x'\rmd x\rmd t \nonumber\\
& \qquad +  C's_nn^{\alpha\beta}\int_{\tau_1}^T \int_{B(z,\log n)^c}\int_{\R^d}\zeta_0^n(x';z,\epsilon_n)e^{-|x-z|} \rmd x'\rmd x\rmd t \nonumber\\
& \leq C\epsilon_n (\log n)^d \int_{\tau_1}^T t^{-(d+1)/\alpha}\rmd t + C_{d,T}T(\log n)^d n^{-\beta(2-\alpha)d/(2(d+1))}\big) + C'T\frac{(\log n)^{d-1}}{n}\nonumber \\
& \leq C \Big(\epsilon_n (\log n)^d \tau_1^{1-\frac{d+1}{\alpha}} + (\log n)^d n^{-\beta(2-\alpha)d/(2(d+1))} + \frac{(\log n)^{d-1}}{n}\Big).\label{2}
\end{align}

For $(ii)$, the corresponding calculation is different and uses Lemma \ref{lem:levy_exp3}$(v)$ with an arbitrary $\lambda>0$:
\begin{align}
\bigg| & u s_n \int_0^{T-\tau_2} \int_{n^{-\beta}}^{\infty} \frac{1}{r^{1+\alpha}} \int_{\R^d} a^n_t(x;r) \big[\zeta^{n,[L,\infty)}_{T-t}(x;z_1,\epsilon)-\zeta^{n,[L,\infty)}_{T-t}(x;z_2,\epsilon)\big] \rmd x\rmd r\rmd t \bigg| \nonumber\\
  &\le C u s_nn^{\alpha\beta} \sup_{t\in [0,T-\tau_2]} \big\|\zeta^{n,[L,\infty)}_{T-t}(x;z_1,\epsilon) -\zeta^{n,[L,\infty)}_{T-t}(x;z_2,\epsilon)\big\|_\lambda \int_0^{T-\tau_2} \int_{\R^d} e^{-\lambda |x|} \rmd x\rmd t  \nonumber\\
  &\le C_{\lambda,d,T} (|z_1-z_2|^{1/2}+\tau_2) e^{\lambda \epsilon} e^{\lambda|z_1|}.\label{3}
\end{align}

Finally, it remains to bound the integral corresponding to small $(T-t)$'s. For $(i)$, we obtain
\begin{align}
   \bigg|& u s_n \int_{T-\tau_1}^T \int_{n^{-\beta}}^{\infty}
    \frac{1}{r^{1+\alpha}} \int_{\R^d} a^n_t(x;r)
    \big[\zeta^{n,[L,\infty)}_{T-t}(x;z,\epsilon_n)-\zeta^{n,[L,\infty)}_{T-t}(x;z,\epsilon_n')\big] \rmd x\rmd r\rmd t \bigg|  \nonumber\\
  &\le C s_n \int_{T-\tau_1}^T \int_{n^{-\beta}}^{\infty}
    \frac{1}{r^{1+\alpha}} \int_{\R^d}
    \big(\zeta^{n,[L,\infty)}_{T-t}(x;z,\epsilon_n) + \zeta^{n,[L,\infty)}_{T-t}(x;z,\epsilon_n')\big)\rmd x\rmd r\rmd t \nonumber\\
  &\le C s_n n^{\alpha\beta} \tau_1 = C\tau_1. \label{4}
\end{align}
The same result obviously holds for $(ii)$, with $\tau_1$ replaced by $\tau_2$.

For the terms involving the initial condition $w_0^n$, similar arguments using Lemma~\ref{lem:levy_exp3}$(i)$ and $(v)$, and Lemma~\ref{lem:Kt1} lead to
\begin{align*}
\big|\la w_0^n,\zeta^n_T(\cdot;z,\epsilon_n) -\zeta^n_T(\cdot;z,\epsilon_n')\ra \big|\leq Ce^{-n^{c_1\alpha\beta/2}} + Cn^{-\beta(2-\alpha)d/(2(d+1))},
\end{align*}
and
\begin{align*}
\big|\la w_0^n,\zeta^n_T(\cdot;z_1,\epsilon) -\zeta^n_T(\cdot;z_2,\epsilon)\ra \big|\leq Ce^{-n^{c_1\alpha\beta/2}} +  Ce^{\lambda(|z_1|+\epsilon)}\big(\tau_2+|z_1-z_2|^{1/2}\big).
\end{align*}

\subsection{Martingale terms}
Now we turn to the martingale terms. As before, we first consider the case $K_t<L$.
We shall estimate the term involving $b^n$, but the same approach can also be applied to
the terms involving $c^n$. We have
\begin{align}
\bigg| u^2& n^{-\gamma} \int_0^T \int_{n^{-\beta}}^{\infty}
    \frac{1}{r^{d+1+\alpha}} \int_{\R^d} b^n_t(x;r)
      \la \ind_{B(x,r)}(1-w^n_t),\zeta^{n,[0,L)}_{T-t}(\cdot;z,\epsilon_n)
      \ra^2 \rmd z\rmd r\rmd t\bigg|  \nonumber \\
  &\le C n^{-\gamma} \int_0^T \int_{n^{-\beta}}^{\infty}
    \frac{1}{r^{d+1+\alpha}} \int \zeta^{n,[0,L)}_t(y;z,\epsilon_n)
    \int \ind_{\{|y-x|<r\}} \nonumber\\
    &  \qquad \qquad \qquad \qquad\qquad \qquad\times \int \ind_{\{|y'-x|<r\}} \zeta^{n,[0,L)}_t(y';z,\epsilon_n)
    \rmd y'\rmd x\rmd y\rmd r\rmd t \nonumber \\
  &\le C n^{-\gamma} \int_0^T \int_{n^{-\beta}}^{\infty}
    \frac{1}{r^{d+1+\alpha}} \int \zeta^{n,[0,L)}_t(y;z,\epsilon)
    \int \ind_{\{|x-y|<r\}} \rmd x\rmd y\rmd r\rmd t \nonumber \\
  &\le C n^{-\gamma} \int_0^T \int_{n^{-\beta}}^{\infty}
    \frac{1}{r^{1+\alpha}} \int \zeta^{n,[0,L)}_t(y;z,\epsilon) \rmd y\rmd r\rmd t \nonumber \\
  &\le C n^{-\gamma} n^{\alpha\beta} \int_0^T \P[K_t<L]\rmd t
  \le C n^\beta n^{-(1-(\alpha-1)/(2\alpha))\alpha\beta}
  = C n^{-(\alpha-1)\beta/2} \label{5}
\end{align}
by Lemma~\ref{lem:Kt1}. Of course, this inequality holds for $(i)$ and $(ii)$.

Now we turn to
\begin{align*}
  \Big|u^2 n^{-\gamma} \int_0^{T-\tau_1} \int_{n^{-\beta}}^{\infty}
    \frac{1}{r^{d+1+\alpha}} \int_{\R^d} b^n_t(x;r)
    \la \ind_{B(x,r)}(1-w^n_t),&\zeta^{n,[L,\infty)}_{T-t}(\cdot;z,\epsilon_n)\\
    & \qquad -\zeta^{n,[L,\infty)}_{T-t}(\cdot;z,\epsilon_n') \ra^2 \rmd x\rmd r\rmd t\Big|.
\end{align*}
Once again we write
\begin{align*}
\zeta^{n,[L,\infty)}_t&(y;z,\epsilon_n) -\zeta^{n,[L,\infty)}_t(y;z,\epsilon_n') \\
&= \int_{\R^d}\zeta_0^n(x';z,\epsilon_n)\big(q_t^{n,[L,\infty)}(y-x') -q_t^{n,[L,\infty)}(y-z) \big)\rmd x' \\
& \qquad - \int_{\R^d}\zeta_0^n(x';z,\epsilon_n')\big(q_t^{n,[L,\infty)}(y-x') -q_t^{n,[L,\infty)}(y-z) \big)\rmd x'.
\end{align*}
This gives us
\begin{align*}
\Big|& \int_{\R^d} b^n_{T-t}(x;r) \la \ind_{B(x,r)}(1-w^n_{T-t}),\zeta^{n,[L,\infty)}_t(\cdot;z,\epsilon_n) -\zeta^{n,[L,\infty)}_t(\cdot;z,\epsilon_n') \ra^2 \rmd x\Big|\\
& \leq C\int_{(\R^d)^3}\ind_{\{|x-y|\leq r\}}\ind_{\{|x-y'|\leq r\}}\bigg[\bigg(\int_{\R^d}\zeta_0^n(x';z,\epsilon_n)\big|q_t^{n,[L,\infty)}(y-x')\\
& \qquad \qquad \qquad \qquad\qquad \qquad\qquad \qquad\qquad \qquad\qquad \qquad \qquad \qquad-q_t^{n,[L,\infty)}(y-z) \big|\rmd x'\bigg) \\
& \qquad  \times \bigg(\int_{\R^d}\zeta_0^n(x';z,\epsilon_n)\big|q_t^{n,[L,\infty)}(y'-x') -q_t^{n,[L,\infty)}(y'-z) \big|\rmd x' \bigg) + S^n_t \bigg]\rmd y'\rmd y\rmd x,
\end{align*}
where $S_t^n$ is the sum of the remaining three terms comprising the squared integral on the first line.
Since all these terms behave in the same way, we shall only bound the first one. Writing as before
$V_r(y,y')(\leq C_dr^d)$ for the volume of $B(y,r)\cap B(y',r)$, and using Fubini's theorem, we can
replace the integral over $x$ by $V_r(y,y')$. Next, as in our estimates of the drift, we split the
integrals over $y,y'$ according to whether or not $y,y'\in B(z,\log n)$.
This gives us the following first bound, using Lemma~\ref{lem:levy_exp3}$(iii)$:
\begin{align*}
&\int_{B(z,\log n)^2} V_r(y,y') \bigg(\int_{\R^d}\zeta_0^n(x;z,\epsilon_n)\Big(t^{-\frac{d+1}{\alpha}}|z-x|+ C_{d,T}n^{-\frac{\beta(2-\alpha)d}{2(d+1)}}\Big)\rmd x\bigg)\\
& \qquad \qquad \qquad \times \bigg(\int_{\R^d}\zeta_0^n(x';z,\epsilon_n)\Big(t^{-\frac{d+1}{\alpha}}|z-x'|+ C_{d,T}n^{-\frac{\beta(2-\alpha)d}{2(d+1)}}\Big)\rmd x'\bigg)\rmd y'\rmd y\\
& \leq r^d\int_{B(z,\log n)^2}\ind_{\{|y-y'|\leq 2r\}}\Big(\epsilon_n t^{-\frac{d+1}{\alpha}} + C_{d,T}n^{-\frac{\beta(2-\alpha)d}{2(d+1)}}\Big)^2\rmd y'\rmd y \\
& \leq r^d (r\wedge \log n)^d (\log n)^d \Big(\epsilon_n t^{-\frac{d+1}{\alpha}} + C_{d,T}n^{-\frac{\beta(2-\alpha)d}{2(d+1)}}\Big)^2.
\end{align*}
Integrating over $t$ and $r$, we obtain
\begin{align}
& n^{-\gamma}\int_{\tau_1}^T\int_{n^{-\beta}}^\infty\frac{1}{r^{d+1+\alpha}}r^d (r\wedge \log n)^d (\log n)^d \Big(\epsilon_n t^{-\frac{d+1}{\alpha}} + C_{d,T}n^{-\frac{\beta(2-\alpha)d}{2(d+1)}}\Big)^2\rmd r\rmd t \nonumber \\
&\leq Cn^{-\gamma}(\log n)^d\bigg(\int_{n^{-\beta}}^{\log n}r^{d-1-\alpha}\rmd r +(\log n)^d\int_{\log n}^\infty r^{-1-\alpha}\rmd r \bigg)\nonumber\\
& \qquad \qquad \times \bigg(\epsilon_n^2\int_{\tau_1}^T t^{-\frac{2(d+1)}{\alpha}}\rmd t+2\epsilon_n n^{-\frac{\beta(2-\alpha)d}{2(d+1)}} \int_{\tau_1}^T t^{-\frac{d+1}{\alpha}}\rmd t + T n^{-\frac{\beta(2-\alpha)d}{(d+1)}} \bigg) \nonumber\\
& \leq Cn^{-\gamma}(\log n)^d \big( (\log n)^{d-\alpha}+n^{\beta(\alpha -d)}\big)\Big[\epsilon_n^2 \tau_1^{1-\frac{2(d+1)}{\alpha}} + 2\epsilon_n n^{-\frac{\beta(2-\alpha)d}{2(d+1)}}\tau_1^{1-\frac{d+1}{\alpha}} + n^{-\frac{\beta(2-\alpha)d}{(d+1)}}\Big]. \label{6}
\end{align}

Secondly, considering the case where $y\in B(z,\log n)$ and $y' \in B(z,\log n)^c$ and using Points $(iii)$ and $(iv)$ in Lemma~\ref{lem:levy_exp3}, the corresponding integral is bounded by
\begin{align*}
& \int_{B(z,\log n)}\int_{B(z,\log n)^c}V_r(y,y')\bigg(\int_{\R^d}\zeta_0^n(x';z,\epsilon_n)\Big(t^{-\frac{d+1}{\alpha}}|z-x|+C_{d,T}n^{-\frac{\beta(2-\alpha)d}{2(d+1)}}\Big)\rmd x\bigg)\\
& \qquad \qquad \times \bigg(\int_{\R^d}\zeta_0^n(x';z,\epsilon_n)e^{-|z-y'|}\rmd x'\bigg)\rmd y'\rmd y \\
& \leq C\Big(t^{-\frac{d+1}{\alpha}}\epsilon_n+C_{d,T}n^{-\frac{\beta(2-\alpha)d}{2(d+1)}}\Big) \int_{B(z,\log n)^c}\int_{B(z,\log n)\cap B(y',2r)}r^d e^{-|z-y'|}\rmd y\rmd y'\\
& \leq C \Big(t^{-\frac{d+1}{\alpha}}\epsilon_n +C_{d,T}n^{-\frac{\beta(2-\alpha)d}{2(d+1)}}\Big) r^d(r\wedge \log n)^d \frac{(\log n)^{d-1}}{n}.
\end{align*}
Integrating over $t$ and $r$ as well, we obtain
\begin{align}
n^{-\gamma}&\int_{\tau_1}^T \int_{n^{-\beta}}^\infty \frac{1}{r^{d+1+\alpha}}\Big(t^{-\frac{d+1}{\alpha}}\epsilon_n +C_{d,T}n^{-\frac{\beta(2-\alpha)d}{2(d+1)}}\Big) r^d(r\wedge \log n)^d \frac{(\log n)^{d-1}}{n} \, \rmd r\rmd t \nonumber\\
&\leq n^{-\gamma}\frac{(\log n)^{d-1}}{n}\Big[\epsilon_n \tau_1^{1-\frac{d+1}{\alpha}}+ C_{d,T}n^{-\frac{\beta(2-\alpha)d}{2(d+1)}}\Big]\Big[(\log n)^{d-\alpha} + n^{\beta(\alpha-d)}\Big].\label{7}
\end{align}
The case where $y\in B(z,\log n)^c$ and $y'\in B(z,\log n)$ is treated in the same way. Finally, if $y,y'\in B(z,\log n)^c$, Lemma~\ref{lem:levy_exp3}$(iv)$ gives us the bound
\begin{align*}
\int_{(B(z,\log n)^c)^2}&V_r(y,y') \bigg( \int_{\R^d}\zeta_0^n(x;z,\epsilon_n)e^{-|z-y|}\rmd x\bigg) \bigg( \int_{\R^d}\zeta_0^n(x';z,\epsilon_n)e^{-|z-y'|}\rmd x'\bigg)\rmd y'\rmd y\\
& \leq Cr^d \int_{B(z,\log n)^c}\int_{B(z,\log n)^c\cap B(y,2r)}e^{-|z-y|}e^{-|z-y'|}\rmd y'\rmd y \\
& \leq Cr^d\big(1\wedge r^d\big)\frac{(\log n)^{d-1}}{n}.
\end{align*}
Integrating over $t$ and $r$ gives the bound
\begin{equation}\label{8}
n^{-\gamma}\int_{\tau_1}^T \int_{n^{-\beta}}^\infty \frac{1}{r^{d+1+\alpha}}\, r^d\big(1\wedge r^d\big)\frac{(\log n)^{d-1}}{n}\, \rmd r\rmd t\leq CTn^{-\gamma}\frac{(\log n)^{d-1}}{n}\, (n^{\beta(\alpha-d)}+1).
\end{equation}

For the corresponding bound for $(ii)$, the argument is again much shorter thanks to Point~$(v)$ in Lemma~\ref{lem:levy_exp3}:
\begin{align*}
  \lefteqn{\Big|\int_{\R^d} \int_{\R^d} \int_{\R^d}
    \ind_{\{|y-x|<r\}} \ind_{\{|z-x|<r\}} b^n_t(x;r) (1-w^n_t(y)) (1-w^n_t(z))}\\
  &\qquad (\zeta^{n,[L,\infty)}_{T-t}(y;z_1,\epsilon)
      -\zeta^{n,[L,\infty)}_{T-t}(y;z_2,\epsilon))
    (\zeta^{n,[L,\infty)}_{T-t}(z;z_1,\epsilon)
      -\zeta^{n,[L,\infty)}_{T-t}(z;z_2,\epsilon)) \rmd z\rmd y\rmd x\Big|\\
  &\le \int_{\R^d} (\zeta^{n,[L,\infty)}_{T-t}(y;z_1,\epsilon)
      -\zeta^{n,[L,\infty)}_{T-t}(y;z_2,\epsilon))
    \int_{\R^d} \int_{\R^d} \ind_{\{|y-x|<r\}} \ind_{\{|z-x|<r\}} \\
  &\qquad |\zeta^{n,[L,\infty)}_{T-t}(z;z_1,\epsilon)
      -\zeta^{n,[L,\infty)}_{T-t}(z;z_2,\epsilon)| \rmd z\rmd x\rmd y \\
  &\le 2\sup_{t\in [0,T-\tau_2]} \|\zeta^{n,[L,\infty)}_{T-t}(\cdot;z_1,\epsilon)
      -\zeta^{n,[L,\infty)}_{T-t}(\cdot;z_2,\epsilon)\|_\lambda \\
  &\qquad \int_{\R^d}\zeta^{n,[L,\infty)}_{T-t}(y;z_1,\epsilon)
    \int_{\R^d} \int_{\R^d} \ind_{\{|y-x|<r\}} \ind_{\{|z-x|<r\}} e^{-\lambda|z|}
    \rmd z\rmd x\rmd y \\
  &\le C \sup_{t\in [0,T-\tau_2]}
    \|\zeta^{n,[L,\infty)}_{T-t}(\cdot;z_1,\epsilon)
    -\zeta^{n,[L,\infty)}_{T-t}(\cdot;z_2,\epsilon)\|_\lambda
    \int_{\R^d}\zeta^{n,[L,\infty)}_{T-t}(y;z_1,\epsilon)(r^{2d}\wedge r^d)\rmd y\\
  &\le C_{\lambda,d,T} (|z_1-z_2|^{1/2}+\tau_2) e^{\lambda(|z_1|+\epsilon)}
    (r^{2d}\wedge r^d),
\end{align*}
which yields
\begin{align}
 \Big|u^2 n^{-\gamma} &\int_0^{T-\tau_2} \int_{n^{-\beta}}^{\infty}
    \frac{1}{r^{d+1+\alpha}} \int_{\R^d} b^n_t(x;r)
    \la \ind_{B(x,r)}(1-w^n_t),\zeta^{n,[L,\infty)}_{T-t}(\cdot;z_1,\epsilon)
    \nonumber\\
    & \qquad \qquad \qquad \qquad\qquad \qquad\qquad \qquad\qquad -\zeta^{n,[L,\infty)}_{T-t}(\cdot;z_2,\epsilon) \ra^2 \rmd x\rmd r\rmd t\Big| \nonumber \\
  &\le C_{\lambda,d,T} n^{-\gamma} \int_0^{T-\tau_2}
    \int_{n^{-\beta}}^{\infty} \frac{1}{r^{d+1+\alpha}} (r^{2d}\wedge r^d)
    (|z_1-z_2|^{1/2}+\tau_2) e^{\lambda(|z_1|+\epsilon)} \rmd r \rmd t \nonumber \\
  &\le C_{\lambda,d,T} n^{-\gamma}  (|z_1-z_2|^{1/2}+\tau_2)e^{\lambda(|z_1|+\epsilon)}
    \left( \int_{n^{-\beta}}^1 r^{d-1-\alpha} \rmd r+\int_1^\infty r^{-1-\alpha} \rmd r
    \right) \nonumber \\
  &\le C_{\lambda,d,T} n^{-\gamma}  (|z_1-z_2|^{1/2}+\tau_2)e^{\lambda(|z_1|+\epsilon)}
    ( n^{(\alpha-d)\beta} + C) \nonumber \\
  &\le C_{\lambda,d,T} n^{-(d-1)\beta}(|z_1-z_2|^{1/2}+\tau_2) e^{\lambda(|z_1|+\epsilon)}
\label{eq:mart1}
\end{align}
since $n^{(\alpha-1)\beta}n^{-\gamma}=1$.

For $t\in (T-\tau_1,T)$, we apply Lemma~\ref{lem:ft_sup} to obtain
\begin{align*}
  \bigg| \int_{\R^d}& b^n_t(x;r) \la \ind_{B(x,r)}(1-w^n_t),
    \zeta^{n,[L,\infty)}_{T-t}(\cdot;z,\epsilon_n) \ra^2 \rmd x\bigg|\nonumber \\
  &\le C \int \zeta^{n,[L,\infty)}_{T-t}(y;z,\epsilon_n) \int \int
    \ind_{\{|y-x|<r\}} \ind_{\{|y'-x|<r\}} \zeta^{n,[L,\infty)}_{T-t}(y';z,\epsilon_n)
    \rmd y'\rmd x \rmd y \nonumber \\
  &\le C_d \int \zeta^{n,[L,\infty)}_{T-t}(y;z,\epsilon_n) \int
    \ind_{\{|y-x|<r\}} (1 \wedge (((T-t)^{-d/\alpha} + e^{-n^{c_5}}) r^d)) \rmd x \rmd y \nonumber \\
  &\le C_d (r^d \wedge (((T-t)^{-d/\alpha} + e^{-n^{c_5}}) r^{2d})), \nonumber
\end{align*}
which implies that
\begin{align}
& \Big|u^2 n^{-\gamma} \int_{T-\tau_1}^T \int_{n^{-\beta}}^{\infty} \frac{1}{r^{d+1+\alpha}} \int_{\R^d} b^n_t(x;r)
    (\la \ind_{B(x,r)}(1-w^n_t),\zeta^{n,[L,\infty)}_{T-t}(\cdot;z,\epsilon_n) \ra^2 \rmd x\rmd r\rmd t\Big| \nonumber \\
& \leq Cn^{-\gamma}\int_0^{\tau_1}\int_{n^{-\beta}}^\infty \Big[r^{-1-\alpha} \wedge \big(t^{-d/\alpha} + e^{-n^{c_5}}\big)r^{d-1-\alpha}\Big]\rmd r\rmd t \nonumber\\
& = Cn^{-\gamma}\int_{n^{-\beta}}^\infty \int_0^{r^\alpha \wedge \tau_1}r^{-1-\alpha}\rmd t\rmd r + Cn^{-\gamma}\int_{n^{-\beta}}^\infty \int_{r^\alpha \wedge \tau_1}^{\tau_1}\big(t^{-d/\alpha} + e^{-n^{c_5}}\big)r^{d-1-\alpha}\rmd t\rmd r \nonumber \\
& = Cn^{-\gamma}\int_{n^{-\beta}}^{\tau_1^{1/\alpha}} r^{-1}\rmd r + C\tau_1 n^{-\gamma}\int_{\tau_1^{1/\alpha}}^\infty r^{-1-\alpha}\rmd r + Cn^{-\gamma}\int_{n^{-\beta}}^{\tau_1^{1/\alpha}}\int_{r^\alpha}^{\tau_1}t^{-d/\alpha} r^{d-1-\alpha}\rmd t\rmd r \nonumber\\
&\leq Cn^{-\gamma}\log n + Cn^{-\gamma +\beta(\alpha-d)}\tau_1^{1-d/\alpha}. \label{9}
\end{align}
The same bound holds for $(ii)$, with $\tau_1$ replaced by $\tau_2$.

Combining (\ref{5}), (\ref{6}), (\ref{7}), (\ref{8}) and (\ref{9}) yields (recall that $\epsilon_n\leq \epsilon_n'$)
\begin{align*}
 \left[{\cal M}^{n,\zeta_0^n(\cdot;z,\epsilon_n) -\zeta_0^n(\cdot;z,\epsilon_n')}\right]_T &
  \le Cn^{-(\alpha-1)\beta/2} + n^{-\gamma + \beta(\alpha-d)}\tau_1^{1-d/\alpha} \\
  & \qquad + n^{-\gamma}(\log n)^d\big((\log n)^{d-\alpha}+ n^{\beta(\alpha-d)}\big)\Big[\epsilon_n^{'2}\,\tau_1^{1-\frac{2(d+1)}{\alpha}} \\
  & \qquad \qquad  + 2 \epsilon_n' n^{-\frac{\beta(2-\alpha)d}{2(d+1)}}\tau_1^{1-\frac{d+1}{\alpha}} + n^{-\frac{\beta(2-\alpha)d}{(d+1)}}\Big],
  \end{align*}
while combining (\ref{5}), (\ref{eq:mart1}) and (\ref{9}) gives us
\begin{align*}
 \left[{\cal M}^{n,\zeta_0^n(\cdot;z_1,\epsilon)-\zeta_0^n(\cdot;z_2,\epsilon)}\right]_T
  & \le C n^{-(\alpha-1)\beta/2} + C_{\lambda,d,T} (|z_1-z_2|^{1/2}+\tau_2)
    e^{\lambda(|z_1|+\epsilon)} \\
    & \qquad + C_d n^{-\gamma+\beta(\alpha-d)}\tau_2^{(\alpha-d)/\alpha} .
\end{align*}

Now, by the Burkholder-Davis-Gundy inequality (\cite{Bu73}),
\[
 \E\left[\sup_{t\le T} \left|{\cal M}^{n,\zeta_0^n(\cdot;z,\epsilon_n) -\zeta_0^n(\cdot;z,\epsilon_n')}_t\right| \right]
  \le \left[{\cal M}^{n,\zeta_0^n(\cdot;z,\epsilon_n)-\zeta_0^n(\cdot;z,\epsilon_n')}\right]_T^{1/2}.
\]
Combining this and the estimate for the drift term yields the desired result.
\end{proof}

\subsection{Lemmas}
\label{sec:continuity_lemmas}
We define for $\theta\in\R^d$,
\[ \tilde q^{n,\{k\}}_t(\theta)=\E\big[e^{i\theta\cdot (X^n_t-X^n_0)}
  \ind_{\{K_t=k\}}\big], \]
and correspondingly $\tilde q^{n,I}_t(\theta)$ for $I\subset [0,\infty)$, as
well as $\tilde q^n_t(\theta)=\tilde q^{n,[0,\infty)}_t(\theta)$.
Recall the representation of $X^n$ using random walks in~(\ref{eq:xy_stable}).
As $X^n$ has independent and stationary increments, the L\'evy-Khintchine Formula
(see \emph{e.g.} Theorems 2.7.10 and 2.8.1 of \cite{SAT1999})
implies that
\[ \tilde q^{n,[0,\infty)}_t(\theta) = \E\big[e^{i \theta \cdot (X^n_t-X^n_0)}\big]
  = e^{t \psi^n(\theta)}, \]
where
\begin{equation}
  \psi^n(\theta) = \int_{\R^d} (e^{i\theta \cdot x}-1) \nu^n(\rmd x).
\label{eq:psin}
\end{equation}
Similarly, we define the limiting L\'evy measure
\[ \nu(\rmd x) = \int_0^\infty \frac{uV_r}{r^{d+1+\alpha}} \sqcap_r^{*2}(x) \rmd r \rmd x, \]
as well as the corresponding function $\psi$,
\begin{equation}\label{def psi}
\psi(\theta) = \int_{\R^d} (e^{i\theta \cdot x}-1) \nu(\rmd x).
\end{equation}
We observe that for all $t>0$, $|e^{t\psi^n(\theta)}|\le 1$ and hence
$|e^{t\psi(\theta)}|\le 1$.

\begin{lemma}
For all $n$, we have:
\begin{enumerate}
\item[(i)] For all $\theta\in\R^d$, $|\psi^n(\theta)-\psi(\theta)|
  \le \frac{4^d} 3 n^{-\beta(2-\alpha)} |\theta|^2$.
\item[(ii)] For $|\theta|\le n^\beta$,
  $-\psi^n(\theta) \ge c |\theta|^\alpha$ for some positive constant $c=c_d$ independent of $n$.
  Hence $-\psi(\theta) \ge c |\theta|^\alpha$ for all $\theta$.
\end{enumerate}
\label{lem:levy_exp1}
\end{lemma}
\begin{proof}
Since $\nu$ is radially symmetric,~(\ref{eq:psin}) implies
\begin{align*}
  \psi^n(\theta) &= \frac 1 2 \int_{\R^d}
    (e^{i\theta \cdot x}-2+e^{-i\theta \cdot x}) \nu^n(\rmd x)
  = \frac 1 2 \int_{\R^d}
    (e^{i\theta \cdot x/2}-e^{-i\theta \cdot x/2})^2 \nu^n(\rmd x) \\
  &=  -2 \int_{\R^d} \sin^2(\theta \cdot x/2) \nu^n(\rmd x)
  = -2 \int_{\R^d} \sin^2(\theta \cdot x/2) \int_{n^{-\beta}}^\infty
    \frac{V_r}{r^{d+1+\alpha}} \sqcap_r^{*2}(x) \rmd r \rmd x.
\end{align*}
The calculations above can easily be repeated for $X$ and $\psi$, then
\begin{align*}
  \frac 1 2 |\psi^n(\theta)-\psi(\theta)|
  &= \left| \int_0^{n^{-\beta}} \frac{V_r}{r^{d+1+\alpha}} \int_{\R^d}
    \sin^2(\theta \cdot x/2) \int \sqcap_r(y) \sqcap_r(x-y) \rmd y \rmd x \rmd r \right| \\
  &=  \int_0^{n^{-\beta}} \frac{V_r}{r^{d+1+\alpha}} \int_{\R^d}
    \frac{\sin^2(\theta \cdot x/2)}{V_r^2} \int
    \ind_{\{|y|<r\}} \ind_{\{|x-y|<r\}} \rmd y \rmd x \rmd r \\
  &\le \int_0^{n^{-\beta}} \frac{1}{r^{d+1+\alpha}} \int_{\R^d}
    \sin^2(\theta \cdot x/2) \int \ind_{\{|x|<2r\}} \rmd x \rmd r \\
  &= \int_0^{n^{-\beta}} \frac{1}{r^{d+1+\alpha}} \int_{|x|<2r}
    \sin^2(\theta \cdot x/2) \rmd x \rmd r.
\end{align*}
Since $|\sin(x)|\le |x|$ for all $x$, we have
\begin{align*}
  \frac 1 2 |\psi^n(\theta)-\psi(\theta)|
  &\le \frac 1 4 \int_0^{n^{-\beta}} \frac{1}{r^{d+1+\alpha}} \int_{|x|<2r}
    \left(\sum_{i=1}^d \theta_i x_i\right)^2 \rmd x \rmd r \\
  &\le \frac 1 4 \int_0^{n^{-\beta}} \frac{1}{r^{d+1+\alpha}}
    \int_{-2r}^{2r} \ldots \int_{-2r}^{2r}
    \left(\sum_{i=1}^d \theta_i x_i\right)^2 \rmd x_1 \ldots \rmd x_d \rmd r \\
  &\le \frac 1 4 \sum_{i=1}^d \theta_i^2 \int_0^{n^{-\beta}}
    \frac{1}{r^{d+1+\alpha}} \int_{-2r}^{2r} \ldots \int_{-2r}^{2r} x_i^2
    \rmd x_1 \ldots \rmd x_d \rmd r.
\end{align*}
The $(d+1)$-dimensional integral above is the same for all $i$ (by symmetry),
and is equal to
\[
  \int_0^{n^{-\beta}} \frac{(4r)^{d-1}}{r^{d+1+\alpha}} \int_{-2r}^{2r} x_1^2
    \rmd x_1 \rmd r
  = \int_0^{n^{-\beta}} \frac{4^{d-1}}{r^{2+\alpha}} \frac 2 3 (2r)^3 \rmd r
  = \frac{4^{d+1}} 3 \int_0^{n^{-\beta}} r^{1-\alpha} \rmd r
  =\! \frac{4^{d+1}} 3 n^{-\beta(2-\alpha)}.
\]
Hence
\[ |\psi^n(\theta)-\psi(\theta)|
  \le \frac{4^d} 3 \, n^{-\beta(2-\alpha)} |\theta|^2, \]
as required by $(i)$.

For $(ii)$, we have
\begin{align*}
  -\frac 1 2 \psi^n(\theta) &=\int_{\R^d} \sin^2(\theta \cdot x/2)
    \int_{n^{-\beta}}^\infty \frac{V_r}{r^{d+1+\alpha}} \sqcap_r^{*2}(x)\rmd r\rmd x \\
  &= \int_{\R^d} \sin^2(\theta \cdot x/2)
    \int_{n^{-\beta}}^\infty \frac{1}{r^{d+1+\alpha}V_r} \int_{\R^d}
    \ind_{\{|y|<r\}} \ind_{\{|x-y|<r\}} \rmd y \rmd r \rmd x \\
  &\ge c_0 \int_{\R^d} \sin^2(\theta \cdot x/2)
    \int_{n^{-\beta}}^\infty \frac{1}{r^{d+1+\alpha}}
    \ind_{\{|x|<r\}} \rmd r \rmd x,
\end{align*}
since the intersection of the disc $\{y:|y|<r\}$ and $\{y:|y-x|<r\}$
has volume larger than $c_0 V_r$ for some positive constant $c_0$
(dependent on $d$) if $|x|<r$.  For $d=1$ and $\theta_1>0$, we have
\begin{align*}
  -\frac 1 2 \psi^n(\theta_1)
  &\ge 2 c_0 \int_{n^{-\beta}}^\infty \int_0^r \sin^2(\theta_1 x/2)
    \frac{1}{r^{2+\alpha}} \rmd x \rmd r \\
  &= 2 c_0 \int_0^\infty \int_{n^{-\beta}\vee x}^\infty
    \sin^2(\theta_1 x/2) \frac{1}{r^{2+\alpha}} \rmd r \rmd x \\
  &= \frac{2c_0}{1+\alpha} \int_0^\infty
    \sin^2(\theta_1 x/2) (n^{-\beta}\vee x)^{-(1+\alpha)} \rmd x \\
  &\ge \frac{2c_0}{1+\alpha} \int_{n^{-\beta}}^\infty
    \sin^2(\theta_1 x/2) x^{-(1+\alpha)} \rmd x \\
  &= \frac{2c_0}{1+\alpha}\, \theta_1^\alpha \int_{\theta_1 n^{-\beta}}^\infty
    \sin^2(y/2) y^{-(1+\alpha)} \rmd y.
\end{align*}
Since $\theta_1\le n^\beta$, the integral in the above is bounded below by
a constant. By symmetry, with thus obtain that for any $\theta$ such that $|\theta|\leq n^\beta$,
\begin{equation}
  -\psi^n(\theta) \ge c |\theta|^\alpha
\label{eq:psin2}
\end{equation}
for some $c>0$.
We can carry out a similar calculation for $d\ge 2$. Since $\psi^n$ is
radially symmetric, it suffices to consider $\theta=(\theta_1,0,\ldots,0)$ with $\theta_1>0$:
\begin{align*}
  -\frac 1 2 \psi^n(\theta)
  &\ge c_0 \int_{n^{-\beta}}^\infty \int_0^r \int_{|x|=\rho}
    \sin^2(\theta_1 x_1/2) \frac{1}{r^{d+1+\alpha}} \rmd x d\rho \rmd r \\
  &= c_0 \int_0^\infty \int_{n^{-\beta}\vee\rho}^\infty \int_{|x|=\rho}
    \sin^2(\theta_1 x_1/2) \frac{1}{r^{d+1+\alpha}} \rmd x \rmd r d\rho \\
  &= \frac{c_0}{d+\alpha} \int_0^\infty \int_{|x|=\rho}
    \sin^2(\theta_1 x_1/2) (n^{-\beta}\vee\rho)^{-(d+\alpha)} \rmd x d\rho \\
  &\ge \frac{c_0}{d+\alpha} \int_{n^{-\beta}}^\infty \int_{|x|=\rho}
    \sin^2(\theta_1 x_1/2) \rho^{-(d+\alpha)} \rmd x d\rho \\
  &= \frac{c_0}{d+\alpha} \int_{n^{-\beta}}^\infty \int_{|y|=1}
    \sin^2(\rho \theta_1 y_1/2) \rho^{-(1+\alpha)} \rmd y d\rho \\
  &= \frac{c_0}{d+\alpha}\, \theta_1^\alpha \int_{\theta_1 n^{-\beta}}^\infty
    \int_{|y|=1} \sin^2(r y_1/2) r^{-(1+\alpha)} \rmd y \rmd r.
\end{align*}
Since $\theta_1=|\theta|\le n^\beta$, the double integral in the above is bounded
below by a constant. Therefore we arrive at the same estimate as in~(\ref{eq:psin2})
and we have proved $(ii)$.
\end{proof}

\begin{lemma}
\begin{enumerate}
\item[(i)] Let $c_2\in (0,\alpha)$ be a constant. If
  $n^{-c_2\beta(2-\alpha)/(2(d+1))}\le t\le T$, then
\[ \int_{|\theta|\leq n^\beta} |(e^{t(\psi^n(\theta)-\psi(\theta))}-1)
    e^{t\psi(\theta)}| d\theta \le C_{d,T} n^{-\beta(2-\alpha)d/(2(d+1))}.
\]
\item[(ii)] Let $Z_r$ be a uniform random variable on $B(0,r)\subset \R^d$, then
\[ \E\big[e^{i\theta\cdot Z_r}\big]=\frac{2^{d/2}\Gamma(d/2+1)}{|r \theta|^{d/2}}
  J_{d/2}(|r \theta|), \]
where $J_{d/2}$ is the Bessel function of the first kind of order $d/2$.
\item[(iii)] If $M\ge 2$, then under the assumptions of $(i)$ there exist positive
$a$ (with $a<1$) and $C_d$, independent of $M$, such that for all $t>0$,
\[ \int_{|\theta|\ge n^\beta} |\tilde q^{n,[M,\infty)}_t(\theta)| d\theta
  \le C_d n^{\beta d} a^{M-1}. \]
\end{enumerate}
\label{lem:levy_exp2}
\end{lemma}
\begin{proof}
Let $\epsilon=n^{-\beta(2-\alpha)d/(d+1)}$.
For $|\theta|\le \sqrt{\epsilon n^{\beta(2-\alpha)}}
= n^{\beta(2-\alpha)/(2(d+1))}$,
Lemma~\ref{lem:levy_exp1}$(i)$ implies for $t\le T$ and sufficiently large $n$,
\[ |e^{t(\psi^n(\theta)-\psi(\theta))}-1| \le C t\big| \psi^n(\theta)-\psi(\theta)\big|
  \le C_{d,T} \epsilon. \]
Hence
\begin{align*}
&\int_{|\theta|\le n^{\beta}}  |(e^{t(\psi^n(\theta)-\psi(\theta)}-1)
    e^{t\psi(\theta)}| d\theta \\
  & \qquad \le \int_{|\theta|\le \sqrt{\epsilon n^{\beta(2-\alpha)}}}
    |(e^{t(\psi^n(\theta)-\psi(\theta)}-1) e^{t\psi(\theta)}| d\theta
  + \int_{\sqrt{\epsilon n^{\beta(2-\alpha)}}<|\theta|\le n^\beta}
    (e^{t\psi^n(\theta)}+e^{t\psi(\theta)}) d\theta \\
  &\qquad \le C_{d,T} (\epsilon n^{\beta(2-\alpha)})^{d/2} \epsilon
  + C_d \int_{\sqrt{\epsilon n^{\beta(2-\alpha)}}}^{n^\beta}
    r^{d-1} e^{-c t r^\alpha} \rmd r
\end{align*}
by Lemma~\ref{lem:levy_exp1}$(ii)$. The first term is equal to $C_{d,T}n^{-\frac{\beta(2-\alpha)d}{2(d+1)}}$. Since
\[
tr^\alpha \geq n^{(\alpha-c_2)\beta(2-\alpha)/(2(d+1))}
\]
in the integral, the second term is bounded by $Cn^{\beta d}e^{-cn^b}$ (with $b=(\alpha-c_2)\beta(2-\alpha)/(2(d+1))>0$). Both estimates combined give us $(i)$.

For $(ii)$, we use Theorem 4.15 of \cite{SW71},
which states that the Fourier
transform of the indicator function on the unit ball in $d$ dimensions is
\[
  \int_{\R^d} \ind_{[0,1]}(|x|) e^{i\theta\cdot x} \rmd x
  = \bigg|\frac\theta{2\pi}\bigg|^{-d/2} J_{d/2}(|\theta|).
\]
Hence, dividing by the volume of the unit ball in $d$ dimensions, which
is $\pi^{d/2}/\Gamma(d/2+1)$, yields
\begin{equation}
  \E\big[e^{i\theta\cdot Z_1}\big] = \frac{2^{d/2}\Gamma(d/2+1)}{|\theta|^{d/2}}
  J_{d/2}(|\theta|).
\label{eq:tilde_rho}
\end{equation}
Scaling $Z_1$ by a factor of $r$ gives us the desired result.

For $(iii)$, we recall from~(\ref{eq:xy_stable}) the representation of $X^n$
using random walks with step size $Y_k$.
Let $R$ be an $\R$-valued r.v. distributed
according to $\frac{V_1} A r^{-(1+\alpha)} \ind_{\{r>n^{-\beta}\}} \rmd r$,
$Z$ be a uniformly distributed random variable in $B(0,1)$
and $\tilde\rho(\theta)=\E[e^{i\theta\cdot Z}]$. Then $\tilde\rho$ is given
by~(\ref{eq:tilde_rho}) is real and
\[
  \tilde q^{n,[M,\infty)}_t(\theta)
  = \E_{K_t} \big[(\E_R[\tilde\rho(R\theta)^2])^{K_t} \ind_{\{K_t\ge M\}}\big],
\]
where $\E_{K_t}$ and $\E_R$ are expectations taken
with respect to $K_t$ and $R$, respectively.
Observe that
\[ \E_R[\tilde\rho(R\theta)^2] = n^{-\alpha\beta} \int_{n^{-\beta}}^\infty
  r^{-(1+\alpha)} \tilde\rho(r\theta)^2 \rmd r. \]

First, we show $|\tilde\rho(v)| = |\E[e^{i v\cdot Z}]|$ is bounded above
by a constant $a\in(0,1)$ for $|v|\ge 1$ uniformly.
Since $Z$ is radially symmetric about 0, we have
$\tilde\rho(v)=\E[\cos(v\cdot Z)] = \E[\cos(v_1 Z^{(1)})]$, where
$v_1$ and $Z^{(1)}$ denote the first coordinate of $v$ and $Z$, respectively.
It suffices to consider $v_1\ge 1$.
Let $\delta_1$ be a small positive constant.
If $|v_1 Z^{(1)} - n \pi| \ge \delta_1$ for all $n\in\Z$,
then $|\cos(v_1 Z^{(1)})|\le \cos \delta_1 <1$.
Let $I_n=((n\pi-\delta_1)/v_1, (n\pi+\delta_1)/v_1)$, then
\[
  \P\big[|v_1 Z^{(1)} - n \pi| < \delta_1 \mbox{ for some } n\in\Z\big]
  = \sum_{n=-\infty}^\infty \P\big[Z^{(1)} \in I_n\big].
\]
Since $-1\le Z^{(1)}\le 1$, the intervals $I_n$ for which
the probabilities on the right hand side above are non-empty and have total
length $\le 2\delta_1$. These intervals do not overlap.
The way to arrange non-overlapping intervals $J_n$ of total length
$2\delta_1$ so that the probability $\sum_n \P[Z^{(1)}\in J_n]$
is maximised is to take
$J_1=[-1,-1+\delta_1],J_2=[1-\delta_1,1]$ and $J_n=\emptyset$ otherwise.
Therefore
\[ \P\big[|v_1 Z^{(1)} - n \pi| < \delta_1 \mbox{ for some } n\in\Z\big]
  \le 2 \P\big[Z^{(1)}\ge [1-\delta_1,1]\big] \le 2 \delta_2 \]
for some $\delta_2\in (0,1/4)$ if we pick a sufficiently small $\delta_1$.
This implies
\begin{align*}
\E\big[\cos(v_1 Z^{(1)})\big]
&= \E\big[\cos(v_1 Z^{(1)}) \ind_{|v_1 Z^{(1)} - n \pi| \ge \delta_1}\big]
  + \E\big[\cos(v_1 Z^{(1)}) \ind_{|v_1 Z^{(1)} - n \pi| < \delta_1}\big] \\
&\le (\cos\delta_1) \P\big[|v_1 Z^{(1)} - n \pi| \ge \delta_1\big]
  + \P\big[|v_1 Z^{(1)} - n \pi| < \delta_1\big]
\le a
\end{align*}
for some $a\in (0,1)$. This estimate implies
\[ \E_R[\tilde\rho(R\theta)^2] \le a \]
for $|\theta|\geq n^\beta$.

Second, plugging in $\theta=n^\beta \xi$ yields
\begin{align*}
  \E_R[\tilde\rho(R n^\beta \xi)^2]
  &= n^{-\alpha\beta}
    \int_{n^{-\beta}}^\infty r^{-(1+\alpha)} \tilde\rho(r n^\beta \xi)^2 \rmd r
  = \int_1^\infty x^{-(1+\alpha)} \tilde\rho(x\xi)^2 \rmd x \\
  &= 2^d\Gamma(d/2+1)^2 \int_1^\infty x^{-(1+\alpha)}
    \frac{J_{d/2}(|x\xi|)^2}{|x\xi|^d} \rmd x
  \le C_d \int_1^\infty x^{-(1+\alpha)} |x\xi|^{-(d+1)} \rmd x \\
  &\le C_d |\xi|^{-(d+1)},
\end{align*}
where we use the fact $|J_\nu(z)|<C z^{-1/2}$ for $\nu>0$
(\cite{AS72}, p. 362, 9.1.61). The two estimates above imply
that there exist $a\in (0,1)$ and $C_d>0$ (both independent of $M$) such that for $|\xi|\ge n^{-\beta}$,
\[ \E_R[\tilde\rho(R n^\beta \xi)^2] \le a \wedge C_d |\xi|^{-(d+1)}. \]
We use this to estimate
\begin{align*}
  \int_{|\theta|\ge n^\beta}& |\tilde q^{n,[M,\infty)}_t(\theta)| d\theta\\
&  = n^{\beta d}\int_{|\xi|\ge 1} |\tilde q^{n,[M,\infty)}_t(n^\beta \xi)|d\xi
  \\
  &\le n^{\beta d}\int_{|\xi|\ge 1} (a \wedge C_d |\xi|^{-(d+1)})^M d\xi \\
  &\le n^{\beta d}\left(\int_{1\le|\xi|\le (C_d/a)^{1/(d+1)}} a^M d\xi
    + \int_{|\xi|>(C_d/a)^{1/(d+1)}} (C_d |\xi|^{-(d+1)})^M d\xi
    \right) \\
  &\le n^{\beta d}\left( C_d a^M (C_d/a)^{d/(d+1)}
    + \int_{(C_d/a)^{1/(d+1)}}^\infty (C_d r^{-(d+1)})^M r^{d-1} \rmd r
    \right).
\end{align*}
We take $\rho=r/C_d^{1/(d+1)}$ (hence $C_d r^{-(d+1)}=\rho^{-(d+1)}$) to obtain
\begin{align*}
  \int_{|\theta|\ge n^\beta} |\tilde q^{n,[M,\infty)}_t(\theta)| d\theta
  &\le C_d' n^{\beta d}\left(a^{M-d/(d+1)} + \int_{1/a^{1/(d+1)}}^\infty
    (\rho^{-(d+1)})^M (\rho C_d^{1/(d+1)})^{d-1} d\rho \right) \\
  &\le C_d' n^{\beta d} \left(a^{M-1}+ \int_{1/a^{1/(d+1)}}^\infty
    \rho^{-(d+1)(M-1)-2} d\rho\right) \\
  &\le C_d' n^{\beta d} \big(a^{M-1}+ (1/a^{1/(d+1)})^{-(d+1)(M-1)-1}\big) \\
  &\le C_d' n^{\beta d} a^{M-1}
\end{align*}
if $M\ge 2$. Hence we have established $(iii)$.
\end{proof}

\begin{lemma}
Let $c_3\in (0,1)$ be a constant. If $M=n^{c_3\alpha\beta/2}$
and $n^{-(1-c_3)\alpha\beta}\leq t\leq T$, then
$\P[K_t<M]\le C e^{-n^{c_3\alpha\beta/2}}$. Hence,
$\int_0^T \P[K_t<M] \rmd t \le C_T n^{-(1-c_3)\alpha\beta}$.
\label{lem:Kt1}
\end{lemma}
\begin{proof}
By a standard tail estimate for the $Poisson(V_1 n^{\alpha\beta} t)$
random variable $K_t$, since $M\leq V_1 n^{\alpha\beta}t$ we can write
\begin{align*}
  \P[K_t<M]
  &\le e^{-V_1 n^{\alpha\beta} t} \bigg(\frac{e V_1 n^{\alpha\beta} t} M\bigg)^M \\
  &= \exp(-V_1 n^{\alpha\beta} t + M (1+\log V_1+ \log (n^{\alpha\beta}t)
    - \log M)).
\end{align*}
The dominant term in the exponent above is
$V_1 n^{\alpha\beta} t$, which is $\ge V_1 n^{c_3\alpha\beta}$, hence
\[ \P[K_t<M] \le C e^{-n^{c_3\alpha\beta/2}}. \]
This establishes the estimate on $\P[K_t<M]$. The estimate on its integral
follows easily by splitting the integral over $[0,n^{-(1-c_3)\alpha\beta})$
and $[n^{-(1-c_3)\alpha\beta},T]$.
\end{proof}

Finally we turn to the proof of our key lemma.

\begin{proof}[Proof of Lemma~\ref{lem:levy_exp3}]
Recall from~(\ref{eq:xy_stable}) the representation of $X^n$
using random walks with step size $Y_k$: conditioned on $R_k$, which
has density $n^{-\alpha\beta} r^{-(1+\alpha)}\ind_{\{r>n^{-\beta}\}}\rmd r$,
$Y_k|R_k=r$ has density $\sqcap^{*2}_r(x)$. Recall also the definition of $q^n_t$ given in (\ref{def q}) and let $q$ be the density of the limiting $\alpha$-stable process with Laplace exponent $\psi$ defined in (\ref{def psi}).
We write
\begin{align*}
2\pi| q^{n,[M,\infty)}_t(x)&-q_t(x)|
 = \left| \int_{\R^d} (e^{t\psi^n(\theta)}
    -\E[e^{i\theta\cdot (X^n_t-X^n_0)} \ind_{\{K_t<M\}}] -e^{t\psi(\theta)})
    e^{-i\theta\cdot x} d\theta \right|  \\
  &\le \left|\int_{|\theta|<n^\beta} (e^{t(\psi^n(\theta)-\psi(\theta))}-1)
    e^{t\psi(\theta)} e^{-i\theta\cdot x} d\theta \right|
  + \left|\int_{|\theta|< n^\beta} \P[K_t<M] d\theta \right| \\
  & \qquad\qquad
  + \left|\int_{|\theta|\ge n^\beta} e^{t\psi(\theta)} d\theta \right|
  + \left|\int_{|\theta|\ge n^\beta} \tilde q^{n,[M,\infty)}_t(\theta) d\theta \right|.
\end{align*}
Lemma~\ref{lem:levy_exp2} implies that
the first and fourth terms are bounded above by
\[ C_{d,T} n^{-\beta(2-\alpha)d/(2(d+1))}, \ C_d n^{\beta d} a^{M-1}, \]
respectively, where we also use $t\ge n^{-c_2\beta(2-\alpha)/(2(d+1))}$.
The second term is bounded above by
\[ C_d \P[K_t<M] n^{\beta d}. \]
Lemma~\ref{lem:levy_exp1}$(ii)$ implies that the third term is bounded by
\begin{align*}
&\int_{|\theta|\ge n^\beta} e^{t\psi(\theta)} d\theta
  \le \int_{|\theta|\ge n^\beta} e^{-c t|\theta|^\alpha} d\theta
  \le C_d \int_{n^\beta}^\infty r^{d-1}
    \exp\big(-c n^{-c_2\beta(2-\alpha)/(2(d+1))} r\big) \rmd r \\
  &\leq  C_d \int_{n^\beta}^\infty
    \exp\big(-c n^{-c_2\beta(2-\alpha)/(2(d+1))} r + (d-1)\log r\big) \rmd r \\
  &\le C_d \int_{n^\beta}^\infty \exp\big(-c n^{-c_2\beta(2-\alpha)/(2(d+1))} r/2\big)
    \rmd r
  \le C_d \exp\Big(-\frac c 2 \, n^{\beta-c_2\beta(2-\alpha)/(2(d+1))}\Big).
\end{align*}
Combining the estimates for these four terms
yields the desired result in $(i)$ for the case $M\ge 2$.
Using Lemma~\ref{lem:Kt1}, the estimate
for $L=n^{c_1 \alpha\beta/2}$ follows easily (noting that $n^{-(1-c_1)\alpha\beta}$ is always smaller than $n^{-c_2\beta (2-\alpha)/[2(d+1)]}$ whenever $c_2<1$).

For $(ii)$, we observe that it was shown in
Lemma~\ref{lemma generator alpha-stable}
that the process $\eta_t$ with generator~(\ref{generator stable motion})
is a symmetric $\alpha$-stable process, hence
$\eta_t \stackrel{d}{=} t^{1/\alpha}\eta_1$. Let $f_{\eta_t}$ be the
density function of $\eta_t$. By Proposition 5.28.1 of \cite{SAT1999},
since $\int_{\R^d} |e^{t\psi(\theta)}| |\theta|^m d\theta<\infty$ for all $m>0$,
$f_{\eta_t}$ is $C^m$ for all $m>0$. In particular, this means that
the first derivative of $f_{\eta_t}$ is uniformly bounded, therefore
$f_{\eta_1}$ is uniformly continuous. This means that
\[ |f_{\eta_t}(x)-f_{\eta_t}(y)|
  = t^{-d/\alpha} |f_{\eta_1}(t^{-1/\alpha} x)-f_{\eta_1}(t^{-1/\alpha} y)|
  \le C t^{-(d+1)/\alpha} |x-y|. \]
Hence
\[
|q_t(x)-q_t(y)| \le C t^{-(d+1)/\alpha} |x-y|,
\]
as desired in $(ii)$. Part $(iii)$ follows easily from $(i)$ and $(ii)$.

Let $f_{X_k}$ denote the density of $X_k=\sum_{i=1}^k Y_i$. Since the density of $Y_1$ is radially symmetric and decreasing in $|x|$, the same properties hold for $f_{X_k}$.
Let $X_{k,1}$ denote the first coordinate of
$X_k$, then for $x_1\in [1,\infty)$ and $\lambda>0$,
\begin{align*}
f_{X_k}(x_1,0,\ldots,0)) & \le \P[X_{k,1} \ge x_1-1] \\
& \le e^{-\lambda (x_1-1)} \E\big[e^{(\lambda,0,\ldots,0)\cdot X_k}\big] = e^{-\lambda (x_1-1)} \E\big[e^{(\lambda,0,\ldots,0)\cdot Y_1}\big]^k.
\end{align*}
We would like to estimate $\E[e^{(\lambda,0,\ldots,0)\cdot Y_1}]$, for which
we calculate, using Lemma~\ref{lem:levy_exp2}$(ii)$,
\begin{align*}
  \E\big[e^{(\lambda,0,\ldots,0)\cdot Y_1}\big]-1
  &= n^{-\alpha\beta} \int_{n^{-\beta}}^\infty \frac 1 {r^{1+\alpha}} \left(
    \frac{2^d\Gamma(d/2+1)^2}{(r\lambda)^d}\, J_{d/2}(r\lambda)^2-1\right) \rmd r \\
  &= n^{-\alpha\beta} \lambda^\alpha \int_{\lambda n^{-\beta}}^\infty
    \frac 1 {\rho^{1+\alpha}} \left(
    \frac{2^d\Gamma(d/2+1)^2}{\rho^d}\, J_{d/2}(\rho)^2-1\right) d\rho.
\end{align*}
From \cite{AS72}, p.362, 9.1.69, Bessel functions are
related to generalised hypergeometric functions in the following way
\[ \Gamma\Big(\frac d 2+1\Big) J_{d/2}(x) (x/2)^{-d/2} = \,_0F_1\Big(\frac d 2+1;-x^2/4\Big)
  := 1+\sum_{n=1}^\infty \frac 1 {(\frac d 2+1)\ldots(\frac d 2+n)}
    \frac{(-x^2/4)^n}{n!}. \]
Hence
\[ \Big|\Big(\Gamma\Big(\frac d 2+1\Big) J_{d/2}(\rho) (\rho/2)^{-d/2}\Big)^2-1\Big| \le C_d \rho^2 \]
for $\rho\in [0,1]$. This implies
\begin{align*}
  \E\big[e^{(\lambda,0,\ldots,0)\cdot Y_1}\big]-1
  &\le n^{-\alpha\beta} \lambda^\alpha \left(
    \int_0^1 C_d \rho^{1-\alpha} d\rho
    + 2\int_1^\infty \rho^{-(1+\alpha)} d\rho \right) \\
  &\le C_\lambda n^{-\alpha\beta},
\end{align*}
where we also use $\big|\frac{2^{d/2}\Gamma(d/2+1)}{\rho^{d/2}} J_{d/2}(\rho)\big|
=|\E[e^{i\rho\cdot Z_1}]|\le 1$ in the first inequality.
Hence
\[ \E\big[e^{(\lambda,0,\ldots,0)\cdot Y_1}\big] \le 1+ C_\lambda n^{-\alpha\beta}
  \le e^{C_\lambda n^{-\alpha\beta}}, \]
which means
\[ f_{X_k}((x_1,0,\ldots,0))
  \le e^{-\lambda (x_1-1)} e^{C_\lambda n^{-\alpha\beta} k}. \]
Plugging the above into the random walk representation yields
\[
  q^{n,[1,\infty)}_t(x) \le \E_{K_t}\big[e^{-\lambda (x_1-1)}
    e^{C_\lambda n^{-\alpha\beta} K_t} \ind_{\{K_t\ge 1\}}\big]
  \le e^{-\lambda (x_1-1)} \exp\big(
    V_1 n^{\alpha\beta}t (e^{C_\lambda n^{-\alpha\beta}}-1)\big)
\]
since $K_t\sim Poisson(V_1 n^{\alpha\beta}t)$.
Since $n^{\alpha\beta} (e^{C_\lambda n^{-\alpha\beta}}-1))\to C_\lambda$
as $n\to\infty$, we have for $t\le T$ and $|x|\ge 1$,
\[
q^{n,[1,\infty)}_t(x) \le C_{\lambda,T} e^{-\lambda (|x|-1)},
\]
as desired in part $(iv)$.

For part $(v)$, we obtain,
\begin{align}
  \lefteqn{\zeta^{n,[L,\infty)}_t(x;y,\epsilon)
    -\zeta^{n,[L,\infty)}_t(x;z,\epsilon)} \nonumber \\
  &= \int (\zeta^{n}_0(x';y,\epsilon)
    -\zeta^{n}_0(x';z,\epsilon)) q^{n,[L,\infty)}_t(x-x') \rmd x' \nonumber \\
  &= \int (\zeta^{n}_0(x'-y;0,\epsilon)
    -\zeta^{n}_0(x'-z;0,\epsilon)) q^{n,[L,\infty)}_t(x-x') \rmd x'
    \nonumber \\
  &= \int \zeta^{n}_0(x';0,\epsilon)
    (q^{n,[L,\infty)}_t(x-y-x')-q^{n,[L,\infty)}_t(x-z-x')) \rmd x'.
\label{eq:qn_norm}
\end{align}
For $t\in [n^{-c_2\beta(2-\alpha)/(2(d+1))},T]$ and $|y-z|\le 1$, we have
\begin{align*}
  \lefteqn{\sup_x |q^{n,[L,\infty)}_t(y-x)-q^{n,[L,\infty)}_t(z-x)|
    e^{\lambda |x|} } \\
  &\le \sup_{x:|x-z|<2} \big|q^{n,[L,\infty)}_t(y-x)-q^{n,[L,\infty)}_t(z-x)\big|
    e^{\lambda |x|} \\
  &  \qquad \qquad  + \sup_{x:|x-z|\ge 2} \big|q^{n,[L,\infty)}_t(y-x)-q^{n,[L,\infty)}_t(z-x)\big|
    e^{\lambda |x|} \\
  &\le C_{\lambda,d,T} \big[(t^{-(d+1)/\alpha} |y-z|
    + n^{-\beta(2-\alpha)d/(2(d+1))}) e^{\lambda |z|} \\
  & \quad  + \sup_{x:|x-z|\geq 2} \min(t^{-(d+1)/\alpha} |y-z|
      + n^{-\beta(2-\alpha)d/(2(d+1))},
    e^{-2\lambda |x-y|}+e^{-2\lambda |x-z|}) e^{\lambda |x|}\big],
\end{align*}
where we use $(iii)$ for the first term, and $(iii)$ and $(iv)$ (applied with $2\lambda$) for the second.
Hence,
\begin{align*}
\sup_x&\, \big|q^{n,[L,\infty)}_t(y-x)-q^{n,[L,\infty)}_t(z-x)\big|
    e^{\lambda |x|}  \\
  &\le C_{\lambda,d,T} \big[(t^{-(d+1)/\alpha} |y-z|
    + n^{-\beta(2-\alpha)d/(2(d+1))}) e^{\lambda |z|} \\
  & \qquad + \sup_x(t^{-(d+1)/\alpha} |y-z|
      + n^{-\beta(2-\alpha)d/(2(d+1))})^{1/2}
    (e^{-2\lambda |x-y|}+e^{-2\lambda |x-z|})^{1/2} e^{\lambda |x|}\big] \\
  &\le C_{\lambda,d,T} \big[(t^{-(d+1)/\alpha} |y-z|
      + n^{-\beta(2-\alpha)d/(2(d+1))}
    + (t^{-(d+1)/\alpha} |y-z|\\
    &  \qquad \qquad \qquad \qquad \qquad \qquad \qquad \qquad \qquad \qquad \qquad
       + n^{-\beta(2-\alpha)d/(2(d+1))})^{1/2}\big] e^{\lambda |z|} \\
  &\le C_{\lambda,d,T} (t^{-(d+1)/(2\alpha)} |y-z|^{1/2}
      + n^{-\beta(2-\alpha)d/(4(d+1))}) e^{\lambda |z|}.
\end{align*}
Plugging this estimate into~(\ref{eq:qn_norm}) yields
\begin{align*}
\sup_x&\, \big|\zeta^{n,[L,\infty)}_t(x;y,\epsilon)
    -\zeta^{n,[L,\infty)}_t(x;z,\epsilon)\big| e^{\lambda |x|} \\
  &\le \sup_x \int \zeta^{n}_0(x';0,\epsilon)
    \big|q^{n,[L,\infty)}_t(x-x'-y)-q^{n,[L,\infty)}_t(x-x'-z)\big|
    e^{\lambda |x-x'|}e^{\lambda(|x|-|x-x'|)} \rmd x' \\
  &\le C_{\lambda,d,T} \big(t^{-(d+1)/(2\alpha)} |y-z|^{1/2}
      + n^{-\beta(2-\alpha)d/(4(d+1))}\big) e^{\lambda |z|}
    \int \zeta^{n}_0(x';0,\epsilon) e^{\lambda |x'|} \rmd x' \\
  &\le C_{\lambda,d,T} e^{\lambda \epsilon}\big(t^{-(d+1)/(2\alpha)} |y-z|^{1/2}
      + n^{-\beta(2-\alpha)d/(4(d+1))}\big) e^{\lambda |z|},
\end{align*}
as desired. Note that we used the assumption that
the support of $\zeta_0^n(\cdot;0,\epsilon)$ is contained in $B(0,\epsilon)$
to bound $e^{\lambda|x'|}$ by $e^{\lambda \epsilon}$.
Note also that this calculation holds even if $\epsilon=\epsilon_n$
depends on $n$.
\end{proof}

\begin{lemma}
There exists $c_5>0$ such that for all $t>0$,
\[ \sup_x \zeta^{n,[L,\infty)}_t(x;z,\epsilon)
  \le C_d (t^{-d/\alpha} + e^{-n^{c_5}}),
\]
where $\epsilon$ can depend on $n$.
\label{lem:ft_sup}
\end{lemma}
\begin{proof}
Let $\tilde\zeta^n_0(\theta)
  =\int_{\R^d} e^{i\theta\cdot x}\zeta_0(x;z,\epsilon) \rmd x$,
then $|\tilde\zeta^n_0(\theta)|\le 1$ regardless of $\epsilon$.
Let $\tilde\zeta^{n,[L,\infty)}_t(\theta)
  =\tilde q^{n,[L,\infty)}_t(\theta) \tilde\zeta^n_0(\theta)$,
where we recall that $\tilde q^{n,[L,\infty)}_t(\theta)
  =\E[e^{i\theta\cdot (X^n_t-X^n_0)}\ind_{\{K_t\ge L\}}]$.
Then
\begin{align*}
  \zeta^{n,[L,\infty)}_t&(x;z,\epsilon) \\
  &=  \frac 1 {2\pi} \int_{\R^d}
    \tilde\zeta^{n,[L,\infty)}_t(\theta) e^{-i\theta\cdot x} d\theta \\
  &\le \left| \frac 1 {2\pi} \int_{|\theta|<n^\beta}
    \tilde q^{n,[L,\infty)}_t(\theta) \tilde\zeta^n_0(\theta)
    e^{-i\theta\cdot x} d\theta \right|
  + \left| \frac 1 {2\pi} \int_{|\theta|\ge n^\beta}
    \tilde q^{n,[L,\infty)}_t(\theta) \tilde\zeta^n_0(\theta)
    e^{-i\theta\cdot x} d\theta \right| \\
  &\le \frac 1 {2\pi} \int_{|\theta|<n^\beta}
    |e^{t\psi^n(\theta)}-\tilde q^{n,[0,L)}_t(\theta)| d\theta
  + \frac 1 {2\pi} \int_{|\theta|\ge n^\beta}
    |\tilde q^{n,[L,\infty)}_t(\theta)| d\theta.
\end{align*}
Since $|\tilde q^{n,[0,L)}_t(\theta)|
= |\E[e^{i\theta\cdot (X^n_t-X^n_0)}\ind_{\{K_t<L\}}]| \le \P[K_t<L]$, we
apply Lemmas~\ref{lem:levy_exp1}$(ii)$, \ref{lem:Kt1}
and~\ref{lem:levy_exp2}$(iii)$ to each term above to obtain
\[
  \zeta^{n,[L,\infty)}_t(x;z,\epsilon)
  \le C_d \left(\int_{\R^d} e^{-c_4 t|\theta|^\alpha} d\theta
   + n^{\beta d} e^{-n^{(c_1/2) \alpha\beta}} + n^{\beta d} a^{L-1} \right)
\]
for some $c_4>0$ and $a\in (0,1)$.
Let $f(t)=\int_{\R^d} e^{-c_4 t|\theta|^\alpha} d\theta$, then
$f(t)=t^{-d/\alpha} f(1)$. Hence,
\[
  \zeta^{n,[L,\infty)}_t(x;z,\epsilon)
  \le C_d \left(t^{-d/\alpha} \int_{\R^d} e^{-c_4|\theta|^\alpha} d\theta
   + e^{-n^{c_5}} \right),
\]
for some $c_5>0$. This implies the desired result.
\end{proof}


\begin{thebibliography}{plain}
\bibitem{AS72} Abramowitz, M. and Stegun, I. A., eds.: Handbook of Mathematical Functions with Formulas, Graphs, and Mathematical Tables. \emph{New York: Dover Publications}, 1972.
\bibitem{ALD1978} Aldous, D.: Stopping times and tightness. \emph{Ann. Probab.} \textbf{6}, (1978), 335--340.
\bibitem{bah/pardoux:2013} Bah, B. and Pardoux, E.: Lambda-lookdown model with selection. \emph{Stoch. Process. Appl.} \textbf{125}, (2015), 1089--1126.
\bibitem{brunet/derrida:2001} Brunet, E. and Derrida, B.: Effect of microscopic noise on front propagation. \emph{J. Stat. Phys.} \textbf{103}, (2001), 269--282.
\bibitem{barton/depaulis/etheridge:2002} Barton, N. H., Depaulis, F. and Etheridge, A. M.: Neutral evolution in spatially continuous populations. \emph{Theor. Pop. Biol.} \textbf{61}, (2002), 31--48.
\bibitem{biswas/etheridge/klimek:2018}Biswas, N., Etheridge, A. M. and Klimek, A.: The spatial Lambda-Fleming-Viot process with fluctuating selection, \emph{Ar$\chi$v preprint} 1802.08188.
\bibitem{BEV2010} Barton, N. H., Etheridge, A. M. and V\'eber, A.: A new model for evolution in a spatial continuum. \emph{Electron. J. Probab.} \textbf{15}, (2010), 162--216.
\bibitem{BEV2013} Barton, N. H., Etheridge, A. M. and V\'eber, A.: Modelling evolution in a spatial continuum. \emph{JSTAT}, (2013), P01002.
\bibitem{BEV2012} Berestycki, N., Etheridge, A. M. and V\'eber, A.: Large-scale behaviour of the spatial $\Lambda$-Fleming-Viot process. \emph{Ann. Inst. H. Poincar\'e Probab. Statist.} \textbf{49}, (2013), 374--401.
\bibitem{BW1998} Biler, P. and Woyczynski, W. A.: Global and exploding solutions for nonlocal quadratic evolution problems. \emph{SIAM J. Appl. Math.} \textbf{59}, (1998), 845--869.
\bibitem{BIL1995} Billingsley, P.: Probability and Measure. \emph{Wiley}, New York, 1995.
\bibitem{Bu73} Burkholder, D. L.: Distribution function inequalities for martingales. \emph{Ann. Probab.} \textbf{1}, (1973), 19--42.
\bibitem{conlon/doering:2005} Conlon, J. G. and Doering, C. R.: On travelling waves for the stochastic Fisher-Kolmogorov-Petrovsky-Piscunov equation. \emph{J. Stat. Phys.} \textbf{120}, (2005), 421--477.
\bibitem{CR2013} Cabr{\'e}, X. and Roquejoffre, J.-M.: The influence of fractional diffusion in Fisher-KPP equations. \emph{Communications in Mathematical Physics} \textbf{320}, (2013), 679--722.
\bibitem{etheridge:2008} Etheridge, A. M.: Drift, draft and structure: some mathematical models of evolution.  \emph{Banach Center Publ.} \textbf{80}, (2008), 121--144.
\bibitem{ETH2011} Etheridge, A. M.: Some Mathematical Models from Population Genetics: \'Ecole d'\'et\'e de Probabilit\'es de Saint-Flour XXXIX-2009. \emph{Springer}, 2011.
\bibitem{EFP2016} Etheridge, A. M., Freeman, N. and Penington, S.: Branching Brownian motion, mean curvature flow and the motion of hybrid zones. \emph{Electron. J. Probab.} \textbf{22}, (2017), 103.
\bibitem{EFPS2016} Etheridge, A. M., Freeman, N., Penington, S. and Straulino, D.: Branching Brownian motion and selection in the spatial Lambda-Fleming-Viot process. \emph{Ann. Applied Probab.} \textbf{27}, (2017), 2605--2645.
\bibitem{EFS2017} Etheridge, A. M., Freeman, N. and Straulino, D.: Branching Brownian motion, the Brownian net and selection in the spatial $\Lambda$-Fleming-Viot process. \emph{Electron. J. Probab.} \textbf{22}, (2017), 39.
\bibitem{etheridge/kurtz:2014} Etheridge, A. M. and Kurtz, T. G.: Genealogical constructions of population models. \emph{Ann. Probab.} \textbf{47}, (2019), 1827--1910.
\bibitem{EK1986} Ethier, S. N. and Kurtz, T. G.: Markov processes: characterization and convergence. \emph{Wiley}, 1986.
\bibitem{EVA1997} Evans, S. N.: Coalescing Markov labelled partitions and a continuous sites genetics model with infinitely many types. \emph{Ann. Instit. H. Poincar\'e Probab. Statist.} \textbf{33}, (1997), 339--358.
\bibitem{felsenstein:1975} Felsenstein, J.: A pain in the torus: some difficulties with the model of isolation
by distance. \emph{Amer. Nat.} \textbf{109}, (1975), 359--368.
\bibitem{FIS1937} Fisher, R.: The wave of advance of advantageous genes. \emph{Annals of Eugenics}, \textbf{7}, (1937), 355--369.
\bibitem{FP2017} Forien, R. and Penington, S.: A central limit theorem for the spatial Lambda-Fleming-Viot process with selection. \emph{Electron. J. Probab.} \textbf{22}, (2017), 5.
\bibitem{foucart:2013} Foucart, C.: The impact of selection in the $\Lambda$-Wright-Fisher model. \emph{Electron. Commun. Probab.} \textbf{18}, (2013), 1--10. Erratum: \emph{Electron. Commun. Probab.} \textbf{15}, (2014), 1--3.
\bibitem{INW1968} Ikeda, N., Nagasawa, M. and Watanabe, S.: Branching Markov processes I. \emph{J. Math. Kyoto Univ.} \textbf{8}, (1968), 233--278.
\bibitem{KAL1976} Kallenberg, O.: Random measures. \emph{Academic Press}, London-New York-San Francisco, 1976.
\bibitem{KAL2002} Kallenberg, O.: Foundations of Modern Probability, 2nd edition. \emph{Springer}, New York, 2002.
\bibitem{KS2016} Khoshnevisan, D. and Schilling R. L.: From L\'evy-Type Processes to Parabolic SPDEs. \emph{Advanced Courses in Mathematics CRM Barcelona, Springer}, 2016.
\bibitem{kimura:1953} Kimura M.: ``Stepping stone'' model of population. \emph{Ann. Rept. Nat. Inst. Genetics Japan} \textbf{3}, (1953), 62--63.
\bibitem{KPP1937} Kolmogorov, A. N., Petrovsky, I. and Piscounov, N.: \'Etude de l'\'equation de la diffusion avec croissance de la quantit\'e de mati\`ere et son application \`a un probl\`eme biologique. \emph{Moscow Univ. Math. Bull.} \textbf{1}, (1937), 1--25.
\bibitem{krone/neuhauser:1997} Krone, S. M. and Neuhauser, C.: Ancestral processes with selection. \emph{Theor. Pop. Biol.} \textbf{51}, (1997), 210--237.
\bibitem{LL2010} Lawler, G. F. and Limic, V.: Random Walk: a modern introduction. \emph{Cambridge University Press}, 2010.
\bibitem{LIA2009} Liang, R. H.: Two continuum-sites stepping stone models in population genetics with delayed coalescence. \emph{PhD Thesis}, University of California, Berkeley, 2009.
\bibitem{MIL2015} Miller, L. R.: Evolution of highly fecund organisms. \emph{PhD Thesis}, University of Oxford, 2015.
\bibitem{mueller/mytnik/quastel:2011} Mueller, C., Mytnik, L. and Quastel, J.: Effect of noise on front propagation in reaction-diffusion equations of KPP type. \emph{Invent. Math.} \textbf{184}, (2011), 405--453.
\bibitem{mueller/sowers:1995} Mueller, C. and Sowers, R.: Random travelling waves for the KPP equation with noise. \emph{J. Functional Analysis} \textbf{128}, (1995), 439--498.
\bibitem{MT1995} Mueller, C. and Tribe, R.: Stochastic p.d.e.'s arising from the long range contact and long range voter processes. \emph{Probab. Theory Relat. Fields} \textbf{102}, (1995), 519--545.
\bibitem{neuhauser/krone:1997} Neuhauser, C. and Krone, S. M.: Genealogies of samples in models with selection. \emph{Genetics} \textbf{145}, (1997), 519--534.
\bibitem{PS1971} Port, S. C. and Stone, C. J.: Infinitely divisible processes and their potential theory. II. \emph{Ann. Instit. Fourier} \textbf{21}, (1971), 179--265.
\bibitem{pruitt:1981} Pruitt, W. E.: The growth of random walks and L\'evy processes. \emph{Ann. Probab.} \textbf{9}, (1981), 948--956.
\bibitem{REB1980} Rebolledo, R.: Sur l'existence de solutions \`a certains probl\`emes de semimartingales. \emph{C.R. Acad. Sci. Paris}, (1980), 290.
\bibitem{RY1999} Revuz, D. and Yor, M.: Continuous martingales and Brownian motion, 3rd edition. \emph{Springer}, 1999.
\bibitem{RR1966} Ridler-Rowe, C. J.: On first hitting times of some recurrent two-dimensional random walks. \emph{Z. Wahr. verw. Geb.} \textbf{5}, (1966), 187--201.
\bibitem{SAT1999} Sato, K.: L\'evy Processes and Infinitely Divisible Distributions. \emph{Cambridge University Press}, 1999.
\bibitem{shiga/shimizu:1980} Shiga, T. and Shimizu, A.: Infinite dimensional stochastic differential equations and their applications. \emph{J. Math. Kyoto Univ.} \textbf{20}, (1980), 395--416.
\bibitem{SW71} Stein, E. and Weiss, G.: Introduction to Fourier Analysis on Euclidean Spaces. \emph{Princeton University Press}, 1971.
\bibitem{VW2012} V\'eber, A. and Wakolbinger, A.: The spatial $\Lambda$-Fleming-Viot process: An event-based construction and a lookdown representation. \emph{Ann. Instit. H. Poincar\'e Probab. Statist.} \textbf{51}, (2015), 570--598.
\bibitem{YUL1925} Yule, G. U.: A mathematical theory of evolution, based on the conclusions of Dr. J.C. Willis. \emph{Philos. Trans. Roy. Soc. London Ser. B} \textbf{213}, (1924), 21--87.
\end{thebibliography}
\end{document}